\pdfoutput=1
\documentclass[reqno]{amsart}
\usepackage[margin=1.5in]{geometry}
\usepackage[utf8]{inputenc}
\usepackage[english]{babel}

\usepackage{csquotes}
\usepackage{graphicx}
\usepackage[utf8]{inputenc}
\usepackage[T1]{fontenc}
\usepackage[dvipsnames]{xcolor}
\usepackage{subfiles}

\usepackage[mathscr]{eucal}
\usepackage[shortlabels]{enumitem}
\usepackage[super]{nth}
\usepackage{leftindex}

\usepackage[backend=biber,
    style=alphabetic,
    sorting=nyt,
    giveninits=true,
    isbn=false,
    maxalphanames=4,
    date=year,
    maxbibnames=99,
    backref=false]{biblatex}
\renewbibmacro{in:}{%
  \ifentrytype{article}{}{\printtext{\bibstring{in}\intitlepunct}}} 
    
\DeclareRobustCommand{\SkipTocEntry}[5]{}

\usepackage{amsmath}
\usepackage{amssymb}

\usepackage{amsthm}
\usepackage{graphicx}

\usepackage{dsfont} 
\usepackage{mathtools}
\usepackage{xparse}
\usepackage{tikz-cd}
\tikzset{
    symbol/.style={%
        draw=none,
        every to/.append style={%
            edge node={node [sloped, allow upside down, auto=false]{$#1$}}}
    }
}
\usepackage{quiver}

\definecolor{leafgreen}{RGB}{48, 138, 3}

\usepackage{hyperref}
\hypersetup{
   pdfborder={0 0 0},
   colorlinks=true,
   allcolors=.,
   bookmarksdepth = 3
}

\usepackage[capitalize, nameinlink]{cleveref}
\crefformat{equation}{#2(#1)#3}
\Crefformat{equation}{#2(#1)#3}
\crefname{section}{Chapter}{Chapters}
\crefname{subsection}{Section}{Sections}
\crefname{subsubsection}{Section}{Sections}
\crefname{subappendix}{Section}{Sections}
\crefname{subsubappendix}{Section}{Sections}

\usepackage{microtype}

\theoremstyle{plain}
\newtheorem{theorem}{Theorem}[subsection]
\newtheorem{introtheorem}{Theorem}
\newtheorem{lemma}[theorem]{Lemma}
\newtheorem{proposition}[theorem]{Proposition}
\newtheorem{corollary}[theorem]{Corollary}

\newtheorem*{theorem*}{Theorem}

\theoremstyle{definition}
\newtheorem{definition}[theorem]{Definition}
\newtheorem{remark}[theorem]{Remark}
\newtheorem{example}[theorem]{Example}
\newtheorem{notation}[theorem]{Notation}
\newtheorem{construction}[theorem]{Construction}

\theoremstyle{remark}

\DeclareMathOperator{\End}{End}
\DeclareMathOperator{\coEnd}{coEnd}
\DeclareMathOperator{\Aut}{Aut}
\DeclareMathOperator{\id}{id}
\DeclareMathOperator*{\colim}{colim}

\DeclareMathOperator{\ev}{ev}

\DeclareMathOperator{\forget}{forget}
\DeclareMathOperator{\cofree}{cofree}
\DeclareMathOperator{\triv}{triv}
\DeclareMathOperator{\indec}{indec}
\DeclareMathOperator{\prim}{prim}

\DeclareMathOperator{\LMod}{LMod}
\DeclareMathOperator{\RMod}{RMod}
\DeclareMathOperator{\BMod}{BiMod}
\DeclareMathOperator{\MixMod}{MixMod}

\DeclareMathOperator{\LComod}{LComod}
\DeclareMathOperator{\RComod}{RComod}
\DeclareMathOperator{\CMon}{CMon}
\DeclareMathOperator{\Alg}{Alg}
\DeclareMathOperator{\coAlg}{coAlg}
\DeclareMathOperator{\CAlg}{CAlg}
\DeclareMathOperator{\coCAlg}{coCAlg}
\DeclareMathOperator{\Map}{Map}
\DeclareMathOperator{\Sym}{Sym}

\DeclareMathOperator{\fib}{fib}

\DeclareMathOperator{\tfib}{tfib}

\DeclareMathOperator{\cross}{cr}
\DeclareMathOperator{\segNerve}{N}
\DeclareMathOperator{\cosk}{cosk}

\DeclareMathOperator{\Bahr}{Bar}
\DeclareMathOperator{\Cobar}{Cobar}
\DeclareMathOperator{\Tot}{Tot}

\DeclareMathOperator{\const}{const}
\DeclareMathOperator{\pr}{pr}

\newcommand{\hide}[1]{}

\newlist{numberenum}{enumerate}{1}
\setlist[numberenum]{\upshape(\arabic*)}

\newcommand{\Fun}{\mathrm{Fun}}
\newcommand{\Lax}{\mathrm{Lax}}
\newcommand{\Laxunit}{\mathrm{Lax}^\mathrm{un}}
\newcommand{\FunL}{\mathrm{Fun}^{\mathrm{L}}}
\newcommand{\FunR}{\mathrm{Fun}^{\mathrm{R}}}

\newcommand{\Funtwo}{\mathrm{Fun}_2}

\newcommand{\FunLM}{\Fun_{\lten}}
\newcommand{\FunRM}{\Fun_{\rten}}
\newcommand{\FunBM}{\Fun_{\bten}}
\newcommand{\iFunLM}{\underline\Fun_{\lten}}
\newcommand{\iFunRM}{\underline\Fun_{\rten}}
\newcommand{\FunLL}[1]{\FunLM}
\newcommand{\FunRR}[1]{\FunRM}
\newcommand{\FunBB}[2]{\FunBM}
\newcommand{\iFunL}[1]{\iFunLM}
\newcommand{\iFunR}[1]{\iFunRM}
\newcommand{\Ortho}{\mathscr{O}\mathrm{rtho}}
\newcommand{\Funstr}[1]{\mathrm{Fun}^\mathrm{str}_{/ #1 \times #1}}
\newcommand{\Fin}{\mathrm{Fin}}

\newcommand{\LM}{\mathrm{LM}}
\newcommand{\RM}{\mathrm{RM}}
\newcommand{\BM}{\mathrm{BM}}
\newcommand{\cA}{\sA}
\newcommand{\cB}{\sB}
\newcommand{\cC}{\sC}
\newcommand{\cD}{\sD}
\newcommand{\cE}{\sE}
\newcommand{\cM}{\sM}

\newcommand{\cP}{\sP}
\newcommand{\sA}{\mathscr{A}}
\newcommand{\sB}{\mathscr{B}}
\newcommand{\sC}{\mathscr{C}}
\newcommand{\sD}{\mathscr{D}}
\newcommand{\sE}{\mathscr{E}}
\newcommand{\sF}{\mathscr{F}}

\newcommand{\sH}{\mathscr{H}}
\newcommand{\sK}{\mathscr{K}}
\newcommand{\sP}{\mathscr{P}}
\newcommand{\sM}{\mathscr{M}}
\newcommand{\sN}{\mathscr{N}}
\newcommand{\sO}{\mathscr{O}}

\newcommand{\sR}{\mathscr{R}}
\newcommand{\sU}{\mathscr{U}}
\newcommand{\sV}{\mathscr{V}}
\newcommand{\sW}{\mathscr{W}}
\newcommand{\sX}{\mathscr{X}}
\newcommand{\sY}{\mathscr{Y}}
\newcommand{\sZ}{\mathscr{Z}}

\newcommand{\Sp}{\mathrm{Sp}}

\newcommand{\diff}{\mathscr{D}\mathrm{iff}}
\newcommand{\qdiff}{\mathrm{q}\diff}
\newcommand{\fdiff}{\mathscr{D}\mathrm{iff}^{\mathrm{FSym}}}
\newcommand{\fdiffc}{\mathscr{D}\mathrm{iff}^{\mathrm{FSym}}_{\geq 1}}

\newcommand{\mlindiffc}{\mathscr{D}\mathrm{iff}^{\mathrm{mlin}}_{\geq 1}}
\newcommand{\mlindiffstc}{\mathscr{D}\mathrm{iff}^{\mathrm{mlin}}_{\mathrm{St}, \geq 1}}
\newcommand{\presl}{\mathscr{P}\mathrm{r}^{\mathrm{L}}}

\newcommand{\presr}{\mathscr{P}\mathrm{r}^{\mathrm{R}}}

\newcommand{\pressym}{\mathscr{P}\mathrm{r}^{\mathrm{Sym}}}

\newcommand{\pressymst}{\mathscr{P}\mathrm{r}^{\mathrm{Sym}}_{\mathrm{St}}}

\newcommand{\pressymc}{\mathscr{P}\mathrm{r}^{\mathrm{Sym}}_{\geq 1}}
\newcommand{\pressymstc}{\mathscr{P}\mathrm{r}^{\mathrm{Sym}}_{\mathrm{St}, \geq 1}}
\newcommand{\pressymdiffc}{\mathscr{P}\mathrm{r}^{\mathrm{Sym}}_{\mathrm{Diff},\geq 1}}
\newcommand{\Pfin}{\mathscr{P}^{\mathrm{fin}}}
\newcommand{\Spc}{\mathscr{S}}

\newcommand{\Cat}{\mathscr{C}\mathrm{at}}
\newcommand{\Cattwo}{\mathscr{C}\mathrm{at}_2}
\newcommand{\PreCattwo}{\mathscr{P}\mathrm{re}\mathscr{C}\mathrm{at}_2}
\newcommand{\Cattwofl}{\mathscr{C}\mathrm{at}^{\mathrm{fl}}_2}
\newcommand{\Catinfty}{\Cat}
\newcommand{\twocatLargeCat}{\widehat{\mathscr{C}\mathrm{at}}}
\newcommand{\CSS}{\mathscr{CS}\mathrm{eg}}

\newcommand{\twofSeg}{\PreCattwo}
\newcommand{\MonCat}{\mathscr{M}\mathrm{on}\mathscr{C}\mathrm{at}}
\newcommand{\RigSymMonCat}{\mathscr{R}\mathrm{ig}\mathscr{S}\mathrm{ym}\mathscr{M}\mathrm{on}}
\newcommand{\SymMonCat}{\mathscr{S}\mathrm{ym}\mathscr{M}\mathrm{on}}
\newcommand{\SymMonCattwo}{\mathscr{S}\mathrm{ym}\mathscr{M}\mathrm{on}_2}
\newcommand{\SymMonCattwofl}{\mathscr{S}\mathrm{ym}\mathscr{M}\mathrm{on}_2^\fl}
\newcommand{\Cocarttwo}{\mathscr{}\mathrm{co}\mathscr{C}\mathrm{art}_2}
\newcommand{\Cocart}{\mathscr{}\mathrm{co}\mathscr{C}\mathrm{art}}
\newcommand{\Cart}{\mathscr{C}\mathrm{art}}

\newcommand{\MMor}{\mathscr{M}\mathrm{or}}
\newcommand{\Hor}{\mathrm{Hor}}

\newcommand{\op}{\mathrm{op}}
\newcommand{\co}{\mathrm{co}}
\newcommand{\coop}{\mathrm{coop}}

\newcommand{\SymSeq}{\mathrm{SymSeq}}
\newcommand{\SymFunL}{\mathrm{SymFun}^\mathrm{L}}
\newcommand{\SymFun}{\mathrm{SymFun}}

\newcommand{\Ind}{\mathrm{Ind}}

\newcommand{\un}{\smallint}
\newcommand{\cartun}{\smallint^\vee}
\newcommand{\Tw}{\mathrm{Tw}}

\newcommand{\TwL}{\mathrm{Tw}^{\mathrm{L}}}
\newcommand{\Arlax}{\mathrm{Ar}^\mathrm{lax}}
\newcommand{\Ar}{\mathrm{Ar}}
\newcommand{\GrayTen}{\otimes_{\mathrm{lax}}}
\newcommand{\Exc}{\mathrm{Exc}}
\newcommand{\Homog}{\mathrm{Homog}}

\newcommand{\Pair}{\mathscr{P}\mathrm{air}}
\newcommand{\PairL}{\mathscr{P}\mathrm{air}^\mathrm{L}}
\newcommand{\mlin}{P_{\vec{1}}}

\newcommand{\inert}{\mathrm{inr}}

\newcommand{\Sph}{\mathbb{S}}
\newcommand{\unit}{\mathbf{1}}

\newcommand{\bbN}{\mathbb{N}}

\makeatletter
\newcommand{\oset}[3][0ex]{%
	\mathrel{\mathop{#3}\limits^{
			\vbox to#1{\kern-2\ex@
				\hbox{$\scriptstyle#2$}\vss}}}}
\makeatother
\newcommand{\eqarrow}{\oset[-.16ex]{\sim}{\longrightarrow}}

\newcommand{\laxfun}{\to}
\newcommand{\oplaxfun}{\laxfun}
\tikzset{
	rot90/.style={anchor=south, rotate=90, inner sep=.5mm}
} 

\newcommand{\angbr}[1]{{\langle #1 \rangle}}

\NewDocumentCommand\derprojlim{e{_}}{\mathchoice
{\varprojlim  \IfValueT{#1}{_{\mathclap{#1}}}{}^{\!1\!}\mathop{}}
{\varprojlim^1  \IfValueT{#1}{_{#1}}}
{\varprojlim^1  \IfValueT{#1}{_{#1}}}
{\varprojlim^1  \IfValueT{#1}{_{#1}}}}

\newcommand{\fl}{\mathrm{fl}}
\newcommand{\lten}{\mathrm{LM}}
\newcommand{\rten}{\mathrm{RM}}
\newcommand{\bten}{\mathrm{BM}}
\newcommand{\dual}{\mathrm{dual}}

\makeatletter
 \def\subsection{\@startsection{subsection}{1}%
 \z@{.7\linespacing\@plus\linespacing}{.5\linespacing}%
 {\normalfont\bfseries\centering}}
\makeatother

\renewcommand{\subset}{\subseteq}

\renewcommand{\epsilon}{\varepsilon}

\begin{document}

\title{On the chain rule in Goodwillie calculus}

\author{Max Blans}
\address{Mathematisch Instituut, Universiteit Utrecht, Utrecht, The Netherlands}
\email{m.a.blans@uu.nl}

\author{Thomas Blom}
\address{Max Planck Institute For Mathematics, Bonn, Germany}
\email{blom@mpim-bonn.mpg.de}

\begin{abstract}
    We prove a generalization of the Arone--Ching chain rule for Goodwillie derivatives by showing that for any pair of reduced finitary functors $F \colon \sD \to \sE$ and $G \colon \sC \to \sD$ between differentiable $\infty$-categories, there is an equivalence
    \[
    \partial_*(FG) \simeq \partial_*F \circ_{\partial_*{\id_\sD}} \partial_*G.
    \]
    This confirms a conjecture of Lurie.
    The proof of this theorem consists of two parts, which are of independent interest.
    We first show that the Goodwillie derivatives can be refined to a lax functor $\partial_* \colon \diff \to \pressymst$ from the $(\infty, 2)$-category of differentiable $\infty$-categories and reduced finitary functors to a certain $(\infty, 2)$-category of generalized symmetric sequences.
    Such a lax structure on the Goodwillie derivatives was long believed to exist, but has not been constructed prior to this work.
    We then finish the proof by studying the interaction of this lax functor with Koszul duality.
    In order to do so, we establish a new universal property of the bar-cobar adjunction.
    \end{abstract}

\maketitle

\tableofcontents

\section{Introduction}

Goodwillie calculus \cite{Goodwillie2003} attaches to any (reduced and finitary\footnote{We follow Goodwillie's convention of calling a functor that preserves filtered colimits \emph{finitary}.}) functor $F \colon \Spc_* \to \Spc_*$ from the category of pointed spaces to itself a sequence of \emph{derivatives}
\[
\partial_1 F, \qquad \partial_2 F, \qquad \partial_3 F, \qquad \partial_4 F, \qquad \ldots
\]
Here $\partial_n F$ is a spectrum with an action of the symmetric group $\Sigma_n$.
In \cite{chingThesis}, Ching shows that the symmetric sequence $\partial_*{\id_{\Spc_*}}$ can be given the structure of an operad.
This operad is Koszul dual to the cocommutative cooperad in spectra, and therefore can be thought of as a homotopical incarnation of the (desuspended) Lie operad.

Arone and Ching subsequently show in \cite{AroneChing2011} that for any (reduced finitary) functor $F \colon \Spc_* \to \Spc_*$, the derivatives $\partial_*F$ can be equipped with a left and right action of the operad $\partial_*{\id_{\Spc_*}}$.
They use this to formulate and prove a \emph{chain rule} for Goodwillie derivatives: if $F$ and $G$ are functors from $\Spc_*$ to itself, then there is an equivalence of symmetric sequences
\[
\partial_*(FG) \simeq \partial_*F \circ_{\partial_*{\id_{\Spc_*}}} \partial_*G.
\]
Here the right-hand side is a relative composition product, which can be computed as a two-sided bar construction.
They also prove a version of the chain rule for functors from spectra to spectra.
This takes the form $\partial_*(FG) \simeq \partial_*F \circ \partial_*G$, corresponding to the fact that $\partial_*{\id}_{\Sp}$ is the trivial operad.

It is natural to ask whether the chain rule also holds for functors between other categories than $\Spc_*$ and $\Sp$.
The right setting for this is laid out by Lurie in \cite[Chapter 6]{HA}, where he develops Goodwillie calculus in the $\infty$-categorical context.
In particular, he defines the derivatives of a functor $F \colon \sC \to \sD$, for $\sC$ and $\sD$ belonging to a large class of $\infty$-categories.
Lurie conjectures \cite[Conjecture 6.3.0.7]{HA} that the chain rule holds in this generality.
In this paper we prove that this is indeed the case:\footnote{In fact, Lurie only states his conjecture for functors between pointed compactly generated $\infty$-categories. We prove this conjecture for the strictly larger class of \emph{differentiable} $\infty$-categories. This class contains many other interesting examples, for example all presentable stable $\infty$-categories.}

\begin{introtheorem}
\label{thm: intro-chain-rule}
    Let $\sC$ be a pointed presentable $\infty$-category in which filtered colimits commute with finite limits (we call such $\infty$-categories differentiable).
    Then 
    \begin{enumerate}[label={\upshape(\arabic*)}]
        \item the derivatives $\partial_*{\id_\sC}$ admit a natural algebra structure for the composition product;
        \item for every reduced finitary functor $F \colon \sC \to \sD$ between differentiable $\infty$-categories, the derivatives $\partial_*F$ are naturally a $(\partial_*{\id_\sD},\partial_*{\id_{\sC}})$-bimodule; and
        \item for every composable pair of reduced finitary functors $\sC \xrightarrow{G} \sD \xrightarrow{F} \sE$ between differentiable $\infty$-categories, there is an equivalence
        \[
        \partial_*(FG) \simeq \partial_*F \circ_{\partial_*{\id}_\sD} \partial_*G
        \]
        of $(\partial_*{\id}_{\sE},\partial_*{\id_{\sC}})$-bimodules.
    \end{enumerate}
\end{introtheorem}

\begin{remark}
    In this result, the assumption that $\cC$, $\cD$ and $\cE$ are presentable can be replaced by the assumption that they admit finite limits and countable colimits. See \cref{thm:intro-quasi-diff} for a precise statement.
\end{remark}

\subsection{Main results}

We will now describe our general strategy for proving \cref{thm: intro-chain-rule} and the results that go into it.
Let us start by clarifying what we mean when we talk about the derivatives and the composition product in the context of functors between arbitrary differentiable categories (we will write category instead of $\infty$-category from now on).

As mentioned above, the derivatives of a functor $F \colon \Spc_* \to \Spc_*$ form a sequence $\{\partial_n F\}_{n \geq 1}$ of spectra with $\Sigma_n$-action.
For a general functor $F \colon \sC \to \sD$ between differentiable categories, this is no longer the case.
Instead, the $n$th derivative $\partial_n F$ is a symmetric functor of $n$ variables from $\Sp(\cC)$ to $\Sp(\cD)$ that preserves colimits in each variable separately.
The collection $\partial_* F = \{\partial_n F\}_{n \geq 0}$ of all these derivatives forms a \emph{functor symmetric sequence}, which we will call a \emph{symmetric sequence} for short:

\begin{definition}
	Let $\sA$ and $\sB$ be presentable categories.
	A \emph{symmetric sequence} $\sO$ from $\sA$ to $\sB$ is a sequence of functors
	\[
	\begin{tikzcd}
		\sO_n \colon \sA^{\times n}_{h\Sigma_n} \ar[r]& \sB, \qquad n \geq 0
	\end{tikzcd}
	\]
	such that for each $n$, the underlying functor of $\sO_n$ preserves colimits in each variable.
	We write $\SymFunL(\sA, \sB)$ for the category of symmetric sequences from $\sA$ to $\sB$.
\end{definition}

Observe that if $\sO$ is a symmetric sequence from $\Sp$ to $\Sp(\sD)$, then $\sO_n$ is completely determined by its evaluation at $n$ copies of the sphere spectrum $\sO_n(\Sph,\ldots,\Sph)$, which is an object of $\Sp(\cD)$ with a $\Sigma_n$-action.
It follows that a symmetric sequence from $\Sp$ to $\Sp(\sD)$ agrees with the usual notion of a symmetric sequence in $\Sp(\cD)$: a functor from $\Fin^\simeq$ to $\cD$.

In order to make sense of \cref{thm: intro-chain-rule}, we need to define a composition product for (functor) symmetric sequences.
The idea is to generalize the usual composition product for symmetric sequences.
For this, it will be convenient to use a coordinate-free indexing for (functor) symmetric sequences---using (non-empty) finite sets instead of natural numbers.
Given symmetric sequences $\sO$ from $\sD$ to $\sE$ and $\sP$ from $\sC$ to $\sD$, their \emph{composition product} will be given by the formula
\begin{equation}\label{eq: intro-composition-product}
(\sO \circ \sP)_I = \coprod_{E \in \mathrm{Part}(I)} \sO_E \circ \{\sP_J\}_{J \in E}.
\end{equation}
Here $I$ is a finite non-empty set and $\mathrm{Part}(I)$ denotes the set of partitions of $I$.
Observe the similarity with the formula for the composition product of (ordinary) symmetric sequences in a symmetric monoidal category.

Our first result is that this composition product is part of a monoidal structure, allowing us to make \cref{thm: intro-chain-rule} precise.

\begin{introtheorem}
\label{thm: intro-composition-product}
	Let $\sC$ be a presentable category.
	Then the composition product refines to a monoidal structure on the category $\SymFunL(\sC, \sC)$.
	Moreover, if $\sD$ is another presentable category, then
	$\SymFunL(\sC, \sD)$ is a bimodule category over $\SymFunL(\sC, \sC)$ and $\SymFunL(\sD, \sD)$.
\end{introtheorem}

\begin{remark}
    In the case where $\sC = \Sp$, a composition product monoidal structure on $\SymFunL(\sC, \sC)$ was already constructed in \cite{brantnerThesis}.
    This monoidal structure is equivalent to the one from \cref{thm: intro-composition-product}, see \cref{prop: equivalence-of-composition-products-spectra} below.
\end{remark}

To construct this monoidal structure, we use an alternative definition of symmetric sequences.
Let $\Sym$ denote the underlying functor of the free commutative algebra monad on $\presl$.
For any presentable category $\sC$, there is the formula
\[
\Sym(\sC) = \bigoplus_{n \geq 0} \sC^{\otimes n}_{h\Sigma_n},
\]
where $\otimes$ denotes the tensor product of presentable categories.
A symmetric sequence from $\sC$ to $\sD$ is then the same as a colimit preserving functor $\Sym(\sC) \to \sD$, and the category of such symmetric sequences is equivalent to $\FunL(\Sym(\sC), \sD)$.
Now it is a remarkable feature of the category $\presl$ that $\Sym(\sC)$ is not just the free commutative algebra, but also the \emph{cofree commutative coalgebra} on $\sC$.
This implies that there is an equivalence
\[
\FunL(\Sym(\sC), \sD) \simeq \FunL_{\coCAlg}(\Sym(\sC), \Sym(\sD)).
\]
If $\sC = \sD$, then the right-hand side has a monoidal structure given by functor composition and $\SymFunL(\sC,\sD)$ has a bimodule structure via pre- and postcomposition.
An explicit computation shows that the resulting monoidal product is given by formula \cref{eq: intro-composition-product}.

One of the main difficulties in generalizing the proof of the chain rule as given by Arone and Ching to arbitrary differentiable categories, is that they use model categorical arguments that are tailored to the categories of pointed spaces and spectra.
It is in particular not clear how to adapt their proofs to construct the algebra structure on $\partial_* {\id}_\sC$ and the bimodule structure on $\partial_*F$ in general.
However, they conjecture a more conceptual reason for the existence of these algebra and bimodule structures which does lend itself to generalization.
Write $\Fun^{\ast, \omega}(\sC, \sC)$ for the category of reduced finitary functors from $\sC$ to itself.
This category has a monoidal structure given by functor composition.
Suppose that the Goodwillie derivatives functor
\[
\begin{tikzcd}
    \partial_* \colon \Fun^{*, \omega}(\sC, \sC) \ar[r] & \SymFunL(\Sp(\sC), \Sp(\sC))
\end{tikzcd}
\]
admits a lax monoidal structure.
Then this would induce an algebra structure on $\partial_*{\id}_\sC$, since $\id_\sC$ is an algebra in $\Fun^{*, \omega}(\sC, \sC)$, being the unit of this monoidal category.
Moreover, since every functor $F \colon \sC \to \sC$ is trivially an $\id_\sC$-bimodule, $\partial_*F$ would automatically acquire the structure of a $\partial_*{\id}_\sC$-bimodule.
One of our main results is that such a lax monoidal structure on $\partial_*$ exists.

\begin{introtheorem}
\label{thm: intro-lax-monoidal-derivatives}
Let $\sC$ be a differentiable category.
The Goodwillie derivatives functor
\[
\begin{tikzcd}
    \partial_* \colon \Fun^{*, \omega}(\sC, \sC) \ar[r] & \SymFunL(\Sp(\sC), \Sp(\sC))
\end{tikzcd}
\]
admits a lax monoidal structure, which is strong monoidal if $\sC$ is a stable category.
\end{introtheorem}

Already in the case where $\sC = \Spc_*$, this result has interesting consequences.
For instance, it implies that if $T \colon \Spc_* \to \Spc_*$ is a reduced finitary monad on pointed spaces, then $\partial_*T$ has the structure of an operad in spectra.
This is something that has long been believed to be true, but has not been proved before.
Also observe that if $\sC$ is stable, the derivatives of a comonad will obtain the structure of a coalgebra since $\partial_*$ is strong monoidal in this case.

For the proof of \cref{thm: intro-chain-rule}, we must also consider functors $\sC \to \sD$ where $\sC$ is not necessarily equal to $\sD$.
The category $\Fun^{\ast, \omega}(\sC, \sD)$ has compatible left and right actions by the categories $\Fun^{\ast, \omega}(\sC, \sC)$ and $\Fun^{\ast, \omega}(\sD, \sD)$, and we would like these actions to be (laxly) intertwined by the functor $\partial_*$ with the actions of $\SymFunL(\sC, \sC)$ and $\SymFunL(\sD, \sD)$ on $\SymFunL(\sC, \sD)$.
For this, it is convenient to use the language of $(\infty, 2)$-categories (which we call $2$-categories throughout this paper).

Let $\sX$ and $\sY$ be $2$-categories.
Roughly speaking, a lax functor $F \colon \sX \to \sY$ is a functor that laxly preserves identities and composition of $1$-morphisms in the following sense: for each object $c \in \sX$, $F$ provides a $2$-morphism $\id_{F(c)} \to F(\id_c)$, and for each composable pair $g \colon c \to d$ and $f \colon d \to e$ of $1$-morphisms in $\sX$, $F$ provides a $2$-morphism
\[
\begin{tikzcd}
& F(d) \ar[dr, "F(f)"]& \\
F(c) \ar[ur, "F(g)"] \ar[rr, "F(f \circ g)"'name=a] & & F(e).
\arrow[from=1-2,to=a, Rightarrow, shorten <= 3pt, shorten >= 3pt]
\end{tikzcd}
\]
An important feature of such a functor is that for each object $x \in \sX$, it induces a lax monoidal functor $\sX(x, x) \to \sY(F(x), F(x))$ of endomorphism categories and for each pair $x, y \in \sX$, it induces a lax functor $\sX(x, y) \to \sY(F(x), F(y))$ of bimodule categories.

Let $\diff$ denote the $2$-category of differentiable categories and reduced finitary functors between them.
We construct a $2$-category $\pressymst$ with objects given by stable presentable categories and $1$-morphisms given by symmetric sequences.
Composition in this $2$-category is given by the composition product, as defined above.
We then prove the following theorem, which has \cref{thm: intro-lax-monoidal-derivatives} as an immediate corollary:

\begin{introtheorem}
\label{introtheorem: lax-2-functor}
    There exists a lax functor of $2$-categories
    \[
    \begin{tikzcd}
    \partial_* \colon \diff \ar[r] & \pressymst,
    \end{tikzcd}
    \]
    that sends a differentiable category $\sC$ to its stabilization $\Sp(\sC)$ and sends a reduced finitary functor $F \colon \sC \to \sD$ to the symmetric sequence $\partial_*F$.
    
    Moreover, the restriction of $\partial_*$ to the full subcategory $\diff_{\mathrm{St}} \subseteq \diff$ spanned by the stable differentiable categories is a \emph{strong} $2$-functor.
\end{introtheorem}

Parts (1) and (2) of \cref{thm: intro-chain-rule} as well as the chain rule for stable categories are formal consequences of this result.
It also provides a natural comparison map
\begin{equation}
\label{eq: intro-comparison-map}
\begin{tikzcd}
\partial_*F \circ_{\partial_*{\id}_\sD} \partial_*G \ar[r] & \partial_*(FG)
\end{tikzcd}
\end{equation}
for any pair of composable reduced finitary functors $\sC \xrightarrow{G} \sD \xrightarrow{F} \sE$ between differentiable categories.
It suffices to show this map is an equivalence to complete the proof of \cref{thm: intro-chain-rule}.
To do so, we use the theory of Koszul duality.

Let $\mathscr{M}$ be a monoidal category that admits both geometric realizations and totalizations.
Koszul duality as developed in \cite[\S 5.2]{HA} provides an adjunction
\[
\begin{tikzcd}[sep = large]
    \Alg^{\mathrm{aug}}(\mathscr{M}) \ar[r, shift left, "\Bahr"] & \coAlg^{\mathrm{aug}}(\mathscr{M}) \ar[l, shift left, "\Cobar"]
\end{tikzcd}
\]
between the categories of augmented algebras and coaugmented coalgebras in $\mathscr{M}$.
If $\unit \to Q$ is a coaugmented coalgebra, then $\Cobar(Q)$ is computed as the totalization of a cosimplicial object of the form
\[
\begin{tikzcd}
        \unit \ar[r, shift left] \ar[r, shift right] & Q \ar[l, shorten=0.4em] \ar[r, shift left=2] \ar[r] \ar[r, shift right=2] & Q \otimes Q \ar[l, shift left, shorten=0.4em] \ar[l, shift right, shorten=0.4em] \ar[r, shift left = 3] \ar[r, shift left]  \ar[r, shift right] \ar[r, shift right = 3] & Q \otimes Q \otimes Q \cdots, \ar[l, shift left = 2, shorten=0.4em] \ar[l, shorten=0.4em] \ar[l, shift right = 2, shorten=0.4em]
    \end{tikzcd}
\]
where the coface maps are given by the comultiplication and coaugmentation of $Q$, and the codegeneracy maps are given by the counit.
Similarly, $\Bahr(A)$ can be computed as the geometric realization of an analogous simplicial object.

One of the crucial insights of Arone and Ching is that there exists a connection between Goodwillie calculus and Koszul duality.
We write
\[
\begin{tikzcd}
\Sigma^\infty_\sC : \sC \ar[r, shift left] & \Sp(\sC) : \Omega^\infty_\sC \ar[l, shift left]
\end{tikzcd}
\]
for the stabilization adjunction of a differentiable category $\sC$.
Every adjunction gives rise to a \emph{Bousfield--Kan resolution}, which is a coaugmented cosimplicial object.
In the case of the stabilization adjunction this takes the form $\id_\sC \to \Omega^\infty_\sC (\Sigma^\infty_\sC \Omega^\infty_\sC)^\bullet \Sigma^\infty_\sC$.
Arone and Ching prove that this resolution induces an equivalence
\[
    \partial_* {\id}_\sC \simeq \Tot \left [
\begin{tikzcd}
        \partial_*(\Omega^\infty_\sC \Sigma^\infty_\sC) \ar[r, shift left] \ar[r, shift right] & \partial_*(\Omega^\infty_\sC (\Sigma^\infty_\sC\Omega^\infty_\sC) \Sigma^\infty_\sC) \ar[l, shorten=0.4em] \ar[r, shift left=2] \ar[r] \ar[r, shift right=2] & \partial_*(\Omega^\infty_\sC (\Sigma^\infty_\sC \Omega^\infty_\sC)^{\circ 2} \Sigma^\infty_\sC) \cdots \ar[l, shift left, shorten=0.4em] \ar[l, shift right, shorten=0.4em]
\end{tikzcd}
    \right ]
\]
upon taking derivatives.
By the stable chain rule and the fact that $\partial_*\Sigma^\infty_\sC \simeq \partial_*\Omega^\infty_\sC \simeq \unit$, it follows that the right-hand side is equivalent to $\Cobar\partial_*(\Sigma^\infty_\sC \Omega^\infty_\sC)$,
where $\partial_*(\Sigma^\infty_\sC \Omega^\infty_\sC)$ inherits a coalgebra structure from the comonad $\Sigma^\infty_\sC \Omega^\infty_\sC$.

As it is stated right now, this is only an equivalence of symmetric sequences.
We would like to lift it to an equivalence of algebras, where the algebra structure on $\partial_* {\id}_\sC$ comes from the lax monoidal structure on the derivatives functor.
To do this, we establish the following universal property of the functors $\Bahr$ and $\Cobar$.
This approach to Koszul duality was inspired by \cite{lurieThursdayKoszul}.

\begin{introtheorem}
\label{introtheorem:koszul-universal-property}
    Let $\mathscr{M}$ be a monoidal category that admits geometric realizations and totalizations.
    Suppose $A$ is an augmented algebra and $Q$ is a coaugmented coalgebra in $\mathscr{M}$.
    Then the mapping spaces $\Map_{\coAlg^{\mathrm{aug}}(\mathscr{M})}(\Bahr(A), Q)$ and $\Map_{\Alg^{\mathrm{aug}}(\mathscr{M})}(A, \Cobar(Q))$ are both naturally equivalent to the space
    \[
    \{\triv_Q(\unit)\} \times_{\LComod_Q(\mathscr{M})} \LComod_Q(\RMod_A(\mathscr{M})) \times_{\RMod_A(\mathscr{M})} \{ \triv_A(\unit) \}
    \]
    of compatibilities between the left $Q$-coaction and right $A$-action on the unit $\unit$.
\end{introtheorem}

We use this in the following way.
The functor $\Sigma^\infty_\sC$ has a canonical left coaction by $\Sigma^\infty_\sC \Omega^\infty_\sC$.
It trivially has a compatible right action by $\id_\sC$. 
Taking derivatives and again using that $\partial_*\Sigma^\infty_\sC \simeq \unit$, we obtain a compatible left coaction by $\partial_*(\Sigma^\infty_\sC \Omega^\infty_\sC)$ and right action by $\partial_*{\id}_\sC$ on the unit $\unit$.
By the universal property stated above, this yields a map
\[
\begin{tikzcd}
\partial_*{\id}_\sC \ar[r] & \Cobar \partial_*(\Sigma^\infty_\sC \Omega^\infty_\sC)
\end{tikzcd}
\]
of augmented algebras.
We then deduce that the underlying map comes from the Bousfield--Kan resolution, proving the following theorem:

\begin{introtheorem}
\label{thm: intro-koszul-goodwillie}
    Let $\sC$ be a differentiable category.
    Then there is an equivalence
    \[
    \partial_*{\id}_\sC \simeq \Cobar\partial_*(\Sigma^\infty_{\sC} \Omega^\infty_{\sC})
    \]
    of augmented algebras.
\end{introtheorem}

When all this is in place, the chain rule follows rather easily.
Indeed, by the naturality of the comparison map \cref{eq: intro-comparison-map}, we can reduce to the case $G = \Sigma^\infty_\sC$ and $F = \Omega^\infty_\sC$.
Since these functors are linear, we have
\[
\partial_*(\Sigma^\infty_\sC) \circ_{\partial_* {\id}_\sC} \partial_*(\Omega^\infty_\sC) \simeq \Bahr(\partial_*{\id}_\sC).
\]
This is equivalent to $\partial_*(\Sigma^\infty_\sC \Omega^\infty_\sC)$ by \cref{thm: intro-koszul-goodwillie} and the fact that $\Bahr$ and $\Cobar$ are inverse to each other in the category of (strongly positive) symmetric sequences.

Our chain rule admits a conceptual reformulation in terms of Morita categories.
To any monoidal category $\sV$ in which the monoidal structure preserves geometric realizations, one can assign its \emph{Morita $2$-category} $\MMor(\sV)$.
This is a $2$-category whose objects are algebras in $\sV$, whose 1-morphisms are bimodules, and where composition of 1-morphisms is given by the relative tensor product.
More generally, one can assign a Morita category to any suitable $2$-category, such as $\pressymst$.
Since the comparison map \cref{eq: intro-comparison-map} identifies the derivatives $\partial_*(FG)$ with the relative composition product of $\partial_* F$ and $\partial_* G$, the Goodwillie derivatives lift to a strong functor into the Morita $2$-category of $\pressymst$.

\begin{introtheorem}\label{thm: intro-Morita}
    Let $\sC$ be a differentiable category and write ${\BMod}_{({\partial_*{\id_\sC}},{\partial_*{\id_\sC}})}$ for the category of $\partial_*{\id_\sC}$-bimodules in $\SymFunL(\Sp(\sC), \Sp(\sC))$.
    Then the Goodwillie derivatives induce a strong monoidal functor
    \[\begin{tikzcd}
        \Fun^{*,\omega}(\sC,\sC) \ar[r] & {\BMod}_{({\partial_*{\id_\sC}},{\partial_*{\id_\sC}})},
    \end{tikzcd}\]
    where ${\BMod}_{({\partial_*{\id_\sC}},{\partial_*{\id_\sC}})}$ is endowed with the relative composition product.
    More generally, the Goodwillie derivatives refine to a strong functor
    \[\begin{tikzcd}
        \partial_* \colon \diff \ar[r] & \MMor(\pressymst)
    \end{tikzcd}\]
    of $2$-categories, where $\MMor(\pressymst)$ is the Morita $2$-category of $\pressymst$.
\end{introtheorem}

So far, we have assumed that our categories are presentable and that our functors preserve filtered colimits.
It turns out that these assumptions can be weakened.
Let us call a category $\cC$ \emph{quasi-differentiable} if it admits a zero object, finite limits, countable colimits, and if sequential colimits commute with finite limits in $\cC$.
A functor will be called \emph{quasi-finitary} if it preserves sequential colimits.
The 2-category of quasi-differentiable categories and reduced quasi-finitary functors between them will be denoted by $\qdiff$.

If $F \colon \cC \to \cD$ is a reduced quasi-finitary functor between quasi-differentiable categories, then its derivatives form a sequence
\[\partial_*F = \{\partial_n F \colon \Sp(\cC)^{\times n}_{h\Sigma_n} \to \Sp(\cD)\}_{n \geq 1}\]
of functors that preserve countable colimits in each variable.
There exists a 2-category $\Cat_{\mathrm{St}}^{\omega_1,\Sym}$ whose objects are stable categories that admit countable colimits and whose 1-morphisms are such sequences of functors.
Composition in this category is given by the same formula as the composition product \cref{eq: intro-composition-product}.
We prove the following chain rule for functors between quasi-differentiable categories. (In fact, this theorem is a relatively formal consequence of the chain rule for functors between differentiable categories.)

\begin{introtheorem}\label{thm:intro-quasi-diff}
    There exists a lax functor of 2-categories
    \[\partial_* \colon \qdiff \to \Cat_{\mathrm{St}}^{\omega_1,\Sym}\]
    that sends a quasi-differentiable category $\cC$ to its stabilization $\Sp(\cC)$ and sends a reduced quasi-finitary functor $F$ to its derivatives sequence $\partial_* F$.
    Moreover, for any pair of composable reduced quasi-finitary functors $\sC \xrightarrow{G} \sD \xrightarrow{F} \sE$, the induced comparison map
    \[\begin{tikzcd}\partial_* F \circ_{\partial_*{\id_\sD}} \partial_* G \ar[r] & \partial_*(FG) \end{tikzcd}\]
    is an equivalence.
    In particular, the Goodwillie derivatives refine to a strong functor
    \[\begin{tikzcd}
        \partial_* \colon \qdiff \ar[r] & \MMor(\Cat_{\mathrm{St}}^{\omega_1,\Sym})
    \end{tikzcd}\]
    into the Morita 2-category of $\Cat_{\mathrm{St}}^{\omega_1,\Sym}$.
\end{introtheorem}

\subsection{Applications} 

We will now sample a number of applications of our results.
All of these will appear in a sequel to this paper.

\subsubsection{Derivatives of the identity in pointed spaces}\

As we stated above, it is a theorem of Ching that the derivatives of the identity functor in pointed spaces $\partial_* {\id}_{\Spc_*}$ are equivalent to the desuspended spectral Lie operad.
We will reprove this in our setting by showing that $\partial_*(\Sigma^\infty \Omega^\infty)$ is the cocommutative cooperad in spectra.
This suffices by \cref{thm: intro-koszul-goodwillie} and the fact that the desuspended spectral Lie operad is Koszul dual to the cocommutative cooperad.

The cooperad structure on $\partial_*(\Sigma^\infty \Omega^\infty)$ was also computed by Arone and Ching in their paper on the chain rule.
They use a point set model for these derivatives that makes it easy to determine the cooperad structure.
It seems difficult to compare this model to ours, and we give an independent proof of this result.

\subsubsection{Derivatives of the identity in \texorpdfstring{$\sO$}{O}-algebras}\

Let $\sO$ be an operad in spectra with $\sO_0 = 0$ and $\sO_1 = \mathbb{S}$.
Let $\Alg_{\sO}$ be the category of $\sO$-algebras in spectra.
We will show that it is a straightforward consequence of our results that $\partial_*{\id}_{\Alg_{\sO}}$ is equivalent as an operad to $\sO$.
In a different setting, this was proved in \cite{ChingDayConvolution} and \cite{clark2021}.
Also cf.\ \cite{HarperHess} and \cite{PereiraThesis}.

\subsubsection{Reconstruction of Taylor towers}\

Suppose that $F \colon \sC \to \sD$ is a reduced finitary functor between differentiable categories. 
One can wonder whether it is possible to recover the functor $F$ from its derivatives (possibly with some extra structure).
Arone and Ching \cite{AroneChingClassification} give an affirmative answer to this question when the categories $\sC$ and $\sD$ are pointed spaces or spectra.
In this case, they show that the lift of the derivatives functor
\[
\begin{tikzcd}
\partial_* \colon \Fun^{\ast, \omega}(\sC, \sD) \ar[r] & {\BMod}_{({\partial_*{\id_\sD}},{\partial_*{\id_\sC}})}
\end{tikzcd}
\]
to the category of $(\partial_*{\id_\sD}, \partial_*{\id_\sC})$-bimodules in $\SymFunL(\Sp(\sC), \Sp(\sD))$ admits a right adjoint.
This adjunction induces a comonad $K$ on the category ${\BMod}_{({\partial_*{\id_\sD}},{\partial_*{\id_\sC}})}$ in the usual way, and for every $F \colon \sC \to \sD$, the derivatives $\partial_* F$ have a canonical $K$-comodule structure.
They then show that if the Goodwillie tower of $F$ converges, $F$ can be reconstructed from its derivatives together with the $(\partial_* {\id}_\sD,\partial_*{\id}_\sC)$-bimodule structure and the coaction of $K$.
We will show the same result holds for arbitrary differentiable categories $\sC$ and $\sD$.

\subsection{Related work}

In \cite[\S 6.3]{HA}, Lurie proves a dual version of the chain rule for functors between compactly generated categories. 
For any such category $\sC$, he constructs a stable $\infty$-operad $\Sp(\sC)^\otimes \to \Fin_*$ that encodes a coalgebra structure on the coderivatives $\partial^*{\id}_\sC$.
For any reduced finitary functor $F \colon \sC \to \sD$, he constructs a $\Delta^1$-family of stable $\infty$-operads that encodes left and right coactions of $\partial^*{\id}_\sC$ and $\partial^*{\id}_\sD$ on the coderivatives $\partial^*{\id}_\sD$.
The dual chain rule can then be formulated as a statement about $\Delta^2$-families of stable $\infty$-operads.
Lurie states that this theorem should be Koszul dual in an appropriate sense to the chain rule for ordinary derivatives (see \cite[Remark 6.3.0.17]{HA}), but a proof of this has not yet appeared.

Another approach to the chain rule making use of stable $\infty$-operads was given by Ching in \cite{ChingDayConvolution}.
Using a model for the derivatives in terms of Day convolution, he is able to construct a stable $\infty$-operad that captures the algebra structure on the derivatives of the identity, as well as a $\Delta^1$-family of stable $\infty$-operads corresponding to the bimodule structure on the derivatives of a functor.
He then formulates a version of the chain rule in this setting, which he conjectures to hold.

Comparing these results to the version of the chain rule proven in this paper would involve setting up a dictionary between (families of) certain kinds of stable $\infty$-operads and (co)algebras (and bi(co)modules over them) in symmetric sequences, which seems like a difficult task.

\subsection{Outline of the paper}

In \cref{sec: 2-categories}, we develop a number of methods to construct (lax) functors between $2$-categories.
This chapter is of a rather technical nature, and we advise the reader to skip it on a first reading and refer back to it when needed.

In \cref{sec:Lax-monoidal-derivatives}, we start by constructing the composition product of functor symmetric sequences, thus proving \cref{thm: intro-composition-product}. We then establish some of the basic properties of this monoidal structure that we need in the rest of the paper.
Next we construct the lax $2$-functor $\partial_* \colon \diff \to \pressymst$.
We end this chapter by proving the chain rule for stable categories, which completes the proofs of \cref{thm: intro-lax-monoidal-derivatives} and \cref{introtheorem: lax-2-functor}.

In \cref{sec: koszul-duality}, we establish the universal property of Koszul duality in terms of compatible actions and coactions on the unit (\cref{introtheorem:koszul-universal-property}).
We then use this to establish the identity $\partial_* {\id}_\sC \simeq \partial_*(\Sigma^\infty_\sC \Omega^\infty_\sC)$ (\cref{thm: intro-koszul-goodwillie}).
Finally, we prove the chain rule (\cref{thm: intro-chain-rule}).
This chapter is concluded with a reformulation of the chain rule in terms of Morita categories (\cref{thm: intro-Morita}) and a generalization of the chain rule to categories that are not necessarily presentable (\cref{thm:intro-quasi-diff}).

This paper has two appendices.
In \cref{appendix:Background-2-categories}, we develop all the $2$-category theory we need from scratch.
We have included it to make the paper more accessible to readers not familiar with $2$-categories, as well as to gather results that are hard to find or not available in the literature.

In \cref{appendix:orthofibrations}, we prove several technical results on (co)modules that are used in our treatment of Koszul duality.
In particular, we show that for any monoidal category $\sC$, there is a natural equivalence $\LComod_Q(\RMod_A(\sC)) \simeq \RMod_A(\LComod_Q(\sC))$.

\addtocontents{toc}{\SkipTocEntry}
\subsection*{Acknowledgments}

The authors would like to thank Fernando Abellán, Gregory Arone, Gijs Heuts and Joost Nuiten for numerous conversations on the contents of this work.
They thank Gijs Heuts for detailed comments on an earlier version of this paper.
They would also like to acknowledge a large intellectual debt to Gregory Arone and Michael Ching, whose original work on the chain rule was the main source of inspiration for this paper.
MB was supported by the European Research Council (ERC) through the grant “Chromatic homotopy theory of spaces”, grant no. 950048.
TB is grateful to the Copenhagen Centre for Geometry and Topology (DNRF151) and the Max Planck Institute for Mathematics in Bonn for their hospitality during the writing of this paper. He was supported by a grant from the Knut and Alice Wallenberg Foundation during part of the period in which this paper was written.

\addtocontents{toc}{\SkipTocEntry}
\subsection*{Conventions}

The term category will mean $(\infty,1)$-category throughout the paper.
Similarly, the term 2-category will mean $(\infty,2)$-category.
To emphasize the distinction between these notions, we will sometimes write 1-category instead of category.
This should not be confused with an ordinary $(1,1)$-category.
\newpage

\section{Functors between 2-categories}
\label{sec: 2-categories}

We start this chapter by describing a method to construct (op)lax functors between 2-categories, which we call ``twisting by adjunctions''. It is a direct generalization of the following fact proved in \cite[Lemma 3.10]{brantner2023pd}: given an adjunction
\[
\begin{tikzcd}[sep = large]
  \sC \ar[r, shift left, "F"] & \sD \ar[l, shift left, "G"]  
\end{tikzcd}
\]
between categories, the functor $F \circ - \circ G \colon \Fun(\sC,\sC) \to \Fun(\sD,\sD)$ can be made oplax monoidal with respect to functor composition. This method for constructing (op)lax functors of 2-categories plays a fundamental role in our construction of a lax structure on the Goodwillie derivatives in \cref{sec:Lax-monoidal-derivatives}.

In \cref{sec:Algebras-in-2-categories}, we study symmetric monoidal 2-categories and prove some basic properties. Notably, we show that in a symmetric monoidal 2-category $\sX$ where every object admits a dual, taking duals refines to a symmetric monoidal equivalence $\sX \simeq \sX^\op$ of 2-categories, and we show that the usual formula for free commutative algebras in a symmetric monoidal category \cite[\S 3.1.3]{HA} extends to 2-categories. We end with \cref{sec:The-2-category-PrL}, where we study one of the most important examples of a 2-category in this paper, namely the 2-category $\presl$ of presentable categories and colimit preserving functors.

The reader should feel free to skip this chapter and refer back to it as necessary; throughout the rest of the paper, the results of this chapter are mainly used as black boxes.

To accommodate readers who are not so familiar with the theory of 2-categories, we have included \cref{appendix:Background-2-categories}, which gives a mostly self-contained treatment of the 2-category theory needed to understand this paper. We will freely use the definitions and results from that appendix throughout this chapter. Note that we define 2-categories in \cref{definition:2-category} as complete Segal objects $\sX_\bullet$ in the category $\Cat$ with the property that $\sX_0$ is a space. In the literature, these often go by the name of \emph{complete twofold Segal spaces} and they are known to be equivalent to all the other common definitions of 2-categories \cite{GagnaHarpazea2022EquivalenceAllModels}. We will also have use for 2-categories that do not satisfy the completeness condition of a complete twofold Segal space. These will be called 2-precategories and are defined in \cref{definition:twofold-Segal-space}.

\subsection{Twisting by adjunctions}\label{sec:Twisting-by-adjunctions}

Suppose $\sX$ is a 2-(pre)category in the sense of \cref{definition:twofold-Segal-space} and that $f \colon x \to y$ is a 1-morphism in $\sX$. We say that $f$ admits a right adjoint if there exists a 1-morphism $g \colon y \to x$ together with unit and counit 2-morphisms $\eta \colon \unit_x \to gf$ and $\varepsilon \colon fg \to \unit_y$ satisfying the usual triangle identities. Given such an adjunction, we can define the ``twisting'' functors
\[
\begin{tikzcd}
f \circ - \circ g : \sX(x,x) \ar[r, shift left] & \sX(y,y) : g \circ - \circ f \ar[l, shift left],
\end{tikzcd}
\]
which are easily observed to be adjoint again via the unit and counit of $f \dashv g$. Moreover, for any $\phi,\psi \in \sX(x,x)$, the counit induces a natural map
\[
\begin{tikzcd}
f \circ \phi \circ \psi \circ g \ar[r] & (f \circ \phi \circ g) \circ (f \circ \psi \circ g),    
\end{tikzcd}
\]
suggesting that $f \circ - \circ g$ is oplax monoidal with respect to the monoidal structures on $\sX(x,x)$ and $\sX(y,y)$ given by composition (these monoidal structures are defined in \cref{def:endomorphisms-infty-category}).
The main result of this section asserts that this is indeed the case.

\begin{theorem}[Twisting by adjunctions]\label{theorem:Twisting-by-adjunctions}
    Let $\sX$ be a 2-precategory and let $f \dashv g$ be an adjunction in $\sX$. Then the functor
    \[
    \begin{tikzcd}
    f \circ - \circ g \colon \sX(x,x) \ar[r] & \sX(y,y)
    \end{tikzcd}
    \]
    can be refined to an oplax monoidal functor, where the oplax structure maps 
    \[f \circ \unit_x \circ g \to \unit_y \quad \text{and} \quad f \circ \phi \circ \psi \circ g \to (f \circ \phi \circ g) \circ (f \circ \psi \circ g) \]
    are given by the counit and unit of the adjunction $f \dashv g$, respectively.
\end{theorem}

\begin{remark}
    Observe that one can dualize \cref{theorem:Twisting-by-adjunctions} to show that
    \[\begin{tikzcd}
        g \circ - \circ f \colon \sX(y,y) \ar[r] & \sX(x,x)
    \end{tikzcd}\]
    refines to a lax monoidal functor, where the lax structure maps are given by
    \[\unit_x \to g \circ \unit_y \circ f \quad \text{and} \quad (g \circ \phi \circ f) \circ (g \circ \psi \circ f) \to g \circ \phi \circ \psi \circ f.\]
    Alternatively, since $g \circ - \circ f$ is right adjoint to $f \circ - \circ g$, one can apply \cite[Theorem 3.4.7]{linskensLaxMonoidalAdjunctions} to deduce that $g \circ f$ is lax monoidal.
\end{remark}

\begin{remark}
    In the case that $\sX = \Cat$, this result is also proved in \cite[Lemma 3.10]{brantner2023pd}.
\end{remark}

It turns out that our proof of \cref{theorem:Twisting-by-adjunctions} can be generalized in two directions: First, if we are given a collection of 1-morphisms $\{f_s \colon x_s \to y_s\}_{s \in S}$ in $\sX$ that admit right adjoints $\{g_s\}_{s \in S}$, then the functors
\[f_s \circ - \circ g_{s'} \colon \sX(x_s,x_{s'}) \rightarrow \sX(y_s, y_{s'}) \quad \text{and} \quad g_s \circ - \circ f_{s'} \colon \sX(y_s,y_{s'}) \to \sX(x_s,x_{s'})\]
can be organized into an oplax and a lax functor between appropriate 2-precategories, respectively. We refer the reader to \cref{definition:twofold-Segal-space} for the definition of a 2-precategory. Secondly, it turns out that the assumption that $f$ admits a right adjoint can be weakened to the assumption that certain left Kan extensions along $f$ exist. The precise statement is somewhat technical and can be found in \cref{theorem:Twisting-in-families} below.

\subsubsection{Lax arrow categories and comma categories}

Our proof of \cref{theorem:Twisting-by-adjunctions} will make use of the lax arrow category of a 2-category. It is constructed using the Gray tensor product $\GrayTen$ of 2-categories, for whose definition and basic properties we refer the reader to \cite{GagnaHarpazea2021GrayTensorProducts}.

\begin{definition}[Lax arrow categories]
    Let $\sX$ be a 2-category. Its \emph{lax arrow category} $\Arlax(\sX)$ is defined as the unique 2-category such that there is a natural equivalence
    \[\Map(\sY, \Arlax(\sX)) \simeq \Map(\sY \GrayTen [1], \sX).\]
\end{definition}

Since the Gray tensor product preserves colimits in each variable separately by \cite[Theorem 2.17]{GagnaHarpazea2021GrayTensorProducts}, we see that the lax arrow category $\Arlax(\sX)$ of any 2-category $\sX$ indeed exists. It follows from the universal property of $\Arlax(\sX)$ that its space of objects is given by $\Map(* \GrayTen [1], \sX) = \sX_1^\simeq$, while its 1-morphisms are given by maps $[1] \GrayTen [1] \to \sX$. In particular, given $f \colon x \to y$ and $g \colon z \to w$ in $\sX$, we see that the 1-morphisms in $\Arlax(\sX)$ from $f$ to $g$ are the lax squares in $\sX$ of the form
\[
\begin{tikzcd}
	x \ar[r] \ar[d,"f"'] & z \ar[d,"g"] \ar[shorten <=7pt, shorten >=7pt,dl,Rightarrow]\\
	y \ar[r] & w.
\end{tikzcd}
\]

Our proof of \cref{theorem:Twisting-by-adjunctions} hinges on the fact that the mapping categories of $\Arlax(\sX)$ can be explicitly computed as certain \emph{comma objects}.
Observe that a lax square
\[\begin{tikzcd}
	   \cE \ar[r] \ar[d] & \cC \ar[d] \ar[shorten <=7pt, shorten >=7pt,dl,Rightarrow]\\
	   \cD \ar[r,""] & \cB,
\end{tikzcd}\]
in $\Cat$ is equivalent to a commutative square of the form
\[\begin{tikzcd}
	   \cE \ar[r] \ar[d] & \cC \ar[d]\\
	   \Ar(\cB) \times_{\cB} \cD \ar[r,""] & \cB.
\end{tikzcd}\]
Here the pullback $\Ar(\cB) \times_{\cB} \cD$ is taken along the target projection $\Ar(\cB) \to \cB$.

\begin{definition}\label{definition:comma-object}
    Let a lax square
    \begin{equation}\label{equation:Lax-square-in-Cat}
    \begin{tikzcd}
	   \cE \ar[r] \ar[d] & \cC \ar[d] \ar[shorten <=7pt, shorten >=7pt,dl,Rightarrow]\\
	   \cD \ar[r,""] & \cB,
    \end{tikzcd}
    \end{equation}
    in $\Cat$ be given. We say that $\cE$ is a \emph{comma object of $\cC \to \cB \leftarrow \cD$} if the lax square \cref{equation:Lax-square-in-Cat} exhibits $\cE$ as the pullback \[\cC \times_{\cB} \Ar(\cB) \times_{\cB} \cD.\]
\end{definition}

In the work of Abellán--Martini \cite{AbellanMartini2024infty2TopoiDescent}, the following computation of the mapping categories of $\Arlax(\sX)$ is given.

\begin{proposition}[{\cite[Proposition 3.4.9]{AbellanMartini2024infty2TopoiDescent}}]\label{proposition:lax-arrow-category-and-comma-categories}
    Let $f \colon x \to y$ and $g \colon z \to w$ be 1-morphisms in $\sX$. Then there is a lax square
    \[\begin{tikzcd}
		{\Arlax(\sX)(f,g)} & {\sX(x,z)} \\
		{\sX(y,w)} & {\sX(x,w)}
		\arrow["{g_*}", from=1-2, to=2-2]
		\arrow["t"', from=1-1, to=2-1]
		\arrow[""{name=0, anchor=center, inner sep=0}, "s", from=1-1, to=1-2]
		\arrow[""{name=1, anchor=center, inner sep=0}, "{f^*}"', from=2-1, to=2-2]
		\arrow[shorten <=7pt, shorten >=7pt, Rightarrow, "\nu", from=1-2, to=2-1]
	\end{tikzcd}\]
    in $\Cat$ exhibiting $\Arlax(\sX)(f,g)$ as the comma category of $\sX(x,z) \xrightarrow{g_*} \sX(x,w) \xleftarrow{f^*} \sX(y,w)$, where $s$ and $t$ are the source and target projections. Moreover, given an object
    \[Q = \begin{tikzcd}
		x \ar[r,"h"] \ar[d,"f"'] & z \ar[d,"g"] \ar[shorten <=7pt, shorten >=7pt,dl,Rightarrow,"\alpha"]\\
		y \ar[r,"k"'] & w
	\end{tikzcd}\]
    of $\Arlax(\sX)(f,g)$, the component of the natural transformation $\nu$ at this $Q$ is given by the 2-morphism
    \[g_*s(Q) = gh \xRightarrow{\;\alpha\;} kf = f^*t(Q).\]
\end{proposition}

\subsubsection{Twisting by an adjunction}

Our proof of \cref{theorem:Twisting-by-adjunctions} uses the following remarkable fact about comma objects and adjunctions.

\begin{lemma}[Comma objects preserve adjunctions]\label{lemma:comma-preserves-adjunction}
    Let
    \[\begin{tikzcd}
	   \cE \ar[r,"F'"] \ar[d,"G'"'] & \cC \ar[d,"G"] \ar[shorten <=7pt, shorten >=7pt,dl,Rightarrow,"\nu"]\\
	   \cD \ar[r,"F"'] & \cB,
    \end{tikzcd}\]
    be a comma object in $\Cat$ and suppose that $F$ admits a left adjoint $L$. Then $F'$ admits a left adjoint $L'$. Moreover, this left adjoint $L'$ is fully faithful and the mate
    \[\begin{tikzcd} LG \ar[Rightarrow,r,"\eta"] & LGF'L' \ar[r,Rightarrow,"\nu"] & LFG'L' \ar[r,Rightarrow,"\varepsilon"] & G'L' \end{tikzcd}\]
    is an equivalence.
\end{lemma}

\begin{proof}
    The comma object $\cE$ is by \cref{definition:comma-object} equivalent to $\cC \times_\cB \Ar(\cB) \times_{\cB} \cD$. Under this equivalence, $F'$ corresponds to the projection onto $\cC$. Given $c$ in $\cC$, let $L'c$ denote the object $(c,\eta_{Gc} \colon Gc \to FLGc, LGc)$ of $\cE \simeq \cC \times_{\cB} \Ar(\cB) \times_{\cB} \cD$. A straightforward computation shows that for any $e$ in $\cE$, one has
    \[\Map_{\cE}(L'c, e) \simeq \Map_{\cC}(c,F'e) \times_{\Map_{\cB}(Gc,FG'e)} \Map_{\cD}(LGc,G'e).\]
    Since $\Map_{\cD}(LGc,G'e) \simeq \Map_{\cB}(Gc,FG'e)$ by assumption, the left-hand pullback is equivalent to $\Map_{\cC}(c,F'e)$. We therefore see that the objects $L'c$ assemble into a functor $L' \colon \cC \to \cE$ that is left adjoint to $F'$. It also follows that the unit $c \to F'L'c$ is simply the identity at $c$, hence $L'$ is fully faithful. For the final claim, we leave it to the reader to verify that componentwise, the mate
    \[\begin{tikzcd} LGc \ar[Rightarrow,r,"\eta"] & LGF'L'c \ar[r,Rightarrow,"\nu"] & LFG'L'c \ar[r,Rightarrow,"\varepsilon"] & G'L'c \end{tikzcd}\]
    can be identified with the composite
    \[\begin{tikzcd}
        LGc \ar[Rightarrow,r,"L\eta"] & LFLGc \ar[Rightarrow,r,"\varepsilon_L"] & LGc.
    \end{tikzcd}\]
    This map is invertible by the triangle identities.
\end{proof}

\begin{remark}
    One can also deduce \cref{lemma:comma-preserves-adjunction} directly from the 2-categorical universal property of comma objects (see for example Proposition 3.4.6 and Theorem 3.5.12 of \cite{RiehlVerity2022ElementsCategoryTheory}). Such a proof has the advantage of showing that \cref{lemma:comma-preserves-adjunction} is true in any 2-category and not just in $\Cat$.
\end{remark}

\begin{remark}
    \Cref{lemma:comma-preserves-adjunction} is remarkable since it is not true for ordinary pullbacks: the pullback of a right adjoint need not be a right adjoint. An example where this fails is the pullback square
    \[\begin{tikzcd}
        \varnothing \ar[r] \ar[d] \arrow["\lrcorner", very near start, phantom,dr] & \{0\} \ar[d] \\
        \{ 1 \} \ar[r] & {[1]}.
    \end{tikzcd}\]
\end{remark}

Combining \cref{lemma:comma-preserves-adjunction} and \cref{proposition:lax-arrow-category-and-comma-categories}, we can now prove \cref{theorem:Twisting-by-adjunctions}.

\begin{proof}[Proof of \cref{theorem:Twisting-by-adjunctions}]
    Let $\sX$ be a 2-precategory and $f \colon x \to y$ a 1-morphism that admits a right adjoint. By \cref{remark:Completion-2-precategories}, $\sX$ admits a fully faithful functor $\phi \colon \sX \to \sX'$ to a 2-category. It is clear that any functor of 2-precategories preserves adjunctions, hence we may assume without loss of generality that $\sX$ is a 2-category.
    
    By \cref{proposition:lax-arrow-category-and-comma-categories}, we see that
    \begin{equation}\label{diagram:lax-arrow-comma-category-2}
        \begin{tikzcd}
    		{\Arlax(\sX)(f,f)} & {\sX(x,x)} \\
    		{\sX(y,y)} & {\sX(x,y)}
    		\arrow["{f_*}", from=1-2, to=2-2]
    		\arrow["t"', from=1-1, to=2-1]
    		\arrow[""{name=0, anchor=center, inner sep=0}, "s", from=1-1, to=1-2]
    		\arrow[""{name=1, anchor=center, inner sep=0}, "{f^*}"', from=2-1, to=2-2]
    		\arrow[shorten <=7pt, shorten >=7pt, Rightarrow, "\nu", from=1-2, to=2-1]
    	\end{tikzcd}
     \end{equation}
    exhibits $\Arlax(\sX)(f,f)$ as a comma category. Note that since $f$ admits a right adjoint $g$, the functor $f^*$ admits a left adjoint $g^*$.
    By \cref{lemma:comma-preserves-adjunction}, we see that $s$ admits a left adjoint $l$.
    Moreover, we see that $f \circ - \circ g = g^*f_*$ is naturally equivalent to the functor $tl \colon \sX(x,x) \to \sX(y,y)$. 
    Since the source and target projections $\Arlax(\sX) \to \sX$ are functors of 2-categories, the functors $s$ and $t$ in diagram \cref{diagram:lax-arrow-comma-category-2} are strong monoidal. 
    By \cite[Theorem 3.4.7]{linskensLaxMonoidalAdjunctions}, it follows that $l$ and hence $tl = f \circ - \circ g \colon \sX(x,x) \to \sX(y,y)$ are oplax monoidal.

    Finally, to see that the oplax structure maps are indeed of the stated form, note that the left adjoint $l$ of $s$ sends $h \colon x \to x$ to the lax square
    \[\begin{tikzcd}
	   x \ar[r,"h"] \ar[d,"f"'] & x \ar[d,"f"] \ar[shorten <=7pt, shorten >=7pt,dl,Rightarrow,"\eta"]\\
	   y \ar[r,"fhg"'] & y.
	\end{tikzcd}\]
    Unwinding the proof of \cref{lemma:comma-preserves-adjunction} in the case of the comma object \cref{diagram:lax-arrow-comma-category-2}, we see that the unit of $l \dashv s$ is the identity and that the counit of this adjunction is given by
    \[\begin{tikzcd}
    	x & x & x & x \\
    	y & {y} & y & y,
    	\arrow[""{name=0, anchor=center, inner sep=0}, "h", from=1-1, to=1-2]
    	\arrow["f"', from=1-1, to=2-1]
    	\arrow["{\eta}", shorten <=7pt, shorten >=7pt, Rightarrow, from=1-2, to=2-1]
    	\arrow["f", from=1-2, to=2-2]
    	\arrow[""{name=1, anchor=center, inner sep=0}, "h", from=1-3, to=1-4]
    	\arrow["f"', from=1-3, to=2-3]
    	\arrow["\alpha", shorten <=7pt, shorten >=7pt, Rightarrow, from=1-4, to=2-3]
    	\arrow["f", from=1-4, to=2-4]
    	\arrow[""{name=2, anchor=center, inner sep=0}, "fhg"', from=2-1, to=2-2]
    	\arrow[""{name=3, anchor=center, inner sep=0}, "k"', from=2-3, to=2-4]
    	\arrow["{\id_h}"', shift left=0, bend left, shorten <=14pt, shorten >=14pt, Rightarrow, from=0, to=1]
    	\arrow["\beta" {yshift=0}, shift right=1, bend right, shorten <=14pt, shorten >=14pt, Rightarrow, from=2, to=3]
    \end{tikzcd}\]
    where $\beta$ denotes the composite $fhg \xRightarrow{\alpha} kfg \xRightarrow{\varepsilon} k$. Now let $\phi, \psi \colon x \to x$ be given. Using \cite[Proposition 3.2.7]{linskensLaxMonoidalAdjunctions}, we see that the oplax structure maps of $l$ are given by
    \[\begin{tikzcd}
    	x & x & x & x \\
    	y & {y} & y & y
    	\arrow[""{name=0, anchor=center, inner sep=0}, "\id_x", from=1-1, to=1-2]
    	\arrow["f"', from=1-1, to=2-1]
    	\arrow["{\eta}", shorten <=7pt, shorten >=7pt, Rightarrow, from=1-2, to=2-1]
    	\arrow["f", from=1-2, to=2-2]
    	\arrow[""{name=1, anchor=center, inner sep=0}, "\id_x", from=1-3, to=1-4]
    	\arrow["f"', from=1-3, to=2-3]
    	\arrow["\id", shorten <=7pt, shorten >=7pt, Rightarrow, from=1-4, to=2-3]
    	\arrow["f", from=1-4, to=2-4]
    	\arrow[""{name=2, anchor=center, inner sep=0}, "fg"', from=2-1, to=2-2]
    	\arrow[""{name=3, anchor=center, inner sep=0}, "\id_y"', from=2-3, to=2-4]
    	\arrow["\id"', shift left=0, bend left, shorten <=14pt, shorten >=14pt, Rightarrow, from=0, to=1]
    	\arrow["\varepsilon" {yshift=0}, shift right=1, bend right, shorten <=14pt, shorten >=14pt, Rightarrow, from=2, to=3]
    \end{tikzcd}\]
    and
    \[\begin{tikzcd}
    	x & x & x & x \\
    	y & {y} & y & y.
    	\arrow[""{name=0, anchor=center, inner sep=0}, "\phi \psi", from=1-1, to=1-2]
    	\arrow["f"', from=1-1, to=2-1]
    	\arrow["{\eta}", shorten <=7pt, shorten >=7pt, Rightarrow, from=1-2, to=2-1]
    	\arrow["f", from=1-2, to=2-2]
    	\arrow[""{name=1, anchor=center, inner sep=0}, "\phi \psi", from=1-3, to=1-4]
    	\arrow["f"', from=1-3, to=2-3]
    	\arrow["{\eta,\eta}", shorten <=7pt, shorten >=7pt, Rightarrow, from=1-4, to=2-3]
    	\arrow["f", from=1-4, to=2-4]
    	\arrow[""{name=2, anchor=center, inner sep=0}, "f \phi \psi g"', from=2-1, to=2-2]
    	\arrow[""{name=3, anchor=center, inner sep=0}, "f \phi gf \psi g"', from=2-3, to=2-4]
    	\arrow["{\id}"', shift left=0, bend left, shorten <=14pt, shorten >=14pt, Rightarrow, from=0, to=1]
    	\arrow["\beta" {yshift=0}, shift right=1, bend right, shorten <=14pt, shorten >=14pt, Rightarrow, from=2, to=3]
    \end{tikzcd}\]
    By definition of the counit of $l \dashv s$, the transformation $\beta$ is the composite
    \[\begin{tikzcd}
        f \phi \psi g \ar[r,Rightarrow,"{\eta,\eta}"] & f \phi gf \psi gfg \ar[r,Rightarrow,"\varepsilon"] & f \phi gf \psi g,
    \end{tikzcd}\]
    which by the triangle identities is just the unit $\eta \colon f \phi \psi g \Rightarrow f \phi g f \psi g$. Applying the target projection $t \colon \Arlax(\sX)(f,f) \to \sX(y,y)$, we obtain the desired description of the oplax structure maps of $tl = f \circ - \circ g$.
\end{proof}

In this proof, the only fact we used about the right adjoint $g$ is that it can be used to construct a left adjoint of $f^* \colon \sX(y,y) \to \sX(x,y)$. This suggests it should already go through under the assumption that certain left Kan extensions along $f$ exist.

\begin{definition}
    Let $\sX$ be a 2-precategory, let $f \colon x \to y$ be a 1-morphism and let $z$ be an object in $\sX$. We say that $f$ is \emph{left-extending with respect to $z$} if the functor $f^* \colon \sX(y,z) \to \sX(x,z)$ admits a left adjoint $f_!$.
\end{definition}

Diagrammatically, this means that for any $1$-morphism $h \colon x \to z$, its left Kan extension
\[\begin{tikzcd}
	x & z \\
	y
	\arrow["h", from=1-1, to=1-2]
	\arrow["f"', from=1-1, to=2-1]
	\arrow[""{name=0, anchor=center, inner sep=0},"f_!h"', from=2-1, to=1-2]
	\arrow[shorten >=2pt, Rightarrow, from=1-1, to=0]
\end{tikzcd}\]
along $f$ must exist.

\begin{example}
    Let $\sX$ be a 2-precategory and suppose that $f \colon x \to y$ is a 1-morphism that admits a right adjoint $g$. Then $f$ is left-extending with respect to any object $z$ of $\sX$, since $g^* \colon \sX(x,z) \to \sX(y,z)$ is left adjoint to $f^* \colon \sX(y,z) \to \sX(x,z)$.
\end{example}

By the exact same proof as \cref{theorem:Twisting-by-adjunctions}, we now obtain the following. Observe that $f_*$ denotes the postcomposition functor $\sX(x,x) \to \sX(x,y)$ and \textbf{not} the right Kan extension along $f$.

\begin{theorem}\label{theorem:Twisting-by-left-extending-functor}
    Let $\sX$ be a 2-precategory and suppose that $f \colon x \to y$ is a 1-morphism such that $f$ is left-extending with respect to $y$. Then the composite
    \[
    \begin{tikzcd}
    \sX(x,x) \ar[r, "f_*"] & \sX(x,y) \ar[r, "f_!"] & \sX(y,y) 
    \end{tikzcd}
    \]
    can be refined to an oplax monoidal functor.
\end{theorem}

\begin{remark}[The oplax structure maps]\label{remark:Unwinding-oplax-structure-maps-twisting}
    As in the case of \cref{theorem:Twisting-by-adjunctions}, one can unwind the proof of \cref{theorem:Twisting-by-left-extending-functor} to obtain an explicit description of the oplax structure maps of $f_!f_*$. For the unit, this is simply
    \[f_!f_*(\unit_x) \simeq f_!f^*(\unit_y) \xRightarrow{\varepsilon} \unit_y.\]
    The other oplax structure map is the map $f_!f_*(\phi \circ \psi) \to f_!f_*\phi \circ f_!f_*\psi$ obtained from the universal property of the left Kan extension $f_!(f \circ \phi \circ \psi)$ and the diagram
    \[\begin{tikzcd}
	x & x & x \\
	y & y & y.
	\arrow["\psi", from=1-1, to=1-2]
	\arrow["f"', from=1-1, to=2-1]
	\arrow["\phi", from=1-2, to=1-3]
	\arrow["\eta", shorten <=7pt, shorten >=7pt, Rightarrow, from=1-2, to=2-1]
	\arrow["f"{description}, from=1-2, to=2-2]
	\arrow["\eta", shorten <=7pt, shorten >=7pt, Rightarrow, from=1-3, to=2-2]
	\arrow["f", from=1-3, to=2-3]
	\arrow["{f_!(f \circ \psi)}"', from=2-1, to=2-2]
	\arrow["{f_!(f \circ \phi)}"', from=2-2, to=2-3]
\end{tikzcd}\]
    This map can also be constructed as the composite
    \begin{align*}
    f_!f_*(\phi \circ \psi) &\simeq f_!(f_*\phi \circ \psi) \xRightarrow{\eta} f_!(f^*f_!f_*\phi \circ \psi) \simeq f_!(f_!f_*\phi \circ f_*\psi) \\
    &\xRightarrow{\eta} f_!(f_!f_*\phi \circ f^*f_!f_*\psi) \simeq f_!f^*(f_!f_*\phi \circ f_!f_*\psi) \xRightarrow{\varepsilon} f_!f_*\phi \circ f_!f_*\psi.
    \end{align*}
\end{remark}

\begin{remark}
    Note that if $f$ admits a right adjoint $g$, then $g^*f_*$ and $f_*g^*$ are equivalent functors. This is not the case in the setting of \cref{theorem:Twisting-by-left-extending-functor}: if $f^* \colon \sX(x,x) \to \sX(y,x)$ admits a left adjoint $f_!$, then $f_!f^*$ and $f^*f_!$ need not be equivalent. \Cref{theorem:Twisting-by-left-extending-functor} only shows that the first of these is oplax monoidal.
\end{remark}

\begin{example}[Density comonad]
    In the situation of \cref{theorem:Twisting-by-left-extending-functor}, the oplax monoidal functor
    \[
    \begin{tikzcd}
        \sX(x,x) \ar[r, "f_!f_*"] & \sX(y,y) 
    \end{tikzcd}
    \]
    takes the unit $\unit_x$ to a coalgebra in $\sX(y,y)$. This coalgebra is classically known as the \emph{density comonad} of $f$.
\end{example}

\subsubsection{Twisting by families of adjunctions}

Suppose we are given a collection of 1-morphisms $\{f_a\colon x_a \to y_a\}_{a \in S}$ in $\sX$ that are pairwise left-extending with respect to each other's codomain; that is, for any $a,b \in S$, the 1-morphism $f_a$ is left extending with respect to $y_b$.
Then we obtain functors $(f_a)_! (f_b)_* \colon \sX(x_a,x_b) \to \sX(y_a,y_b)$, which by \cref{theorem:Twisting-by-left-extending-functor} are oplax monoidal whenever $a = b$. One may wonder whether these can be arranged into a single oplax functor between 2-precategories. We show that this is indeed the case in \cref{theorem:Twisting-in-families} below. While the precise result that we obtain is somewhat technical to state, it is a very powerful method for constructing (oplax) functors between 2-categories that will be used several times throughout this paper. In order to state the theorem we make use of \cref{construction:Modifying-space-of-objects}: here we describe, for a given 2-precategory $\sX$ and a map of spaces $\phi \colon S \to \sX_0$, that one can ``pull back'' $\sX$ along $\phi$ to obtain a 2-precategory $\phi^* \sX$ with space of objects $S$ and whose mapping categories agree with those of $\sX$.

\begin{theorem}[Twisting in families]\label{theorem:Twisting-in-families}
    Let $\sX$ be a 2-precategory and $\phi \colon S \to \sX_1^\simeq$ a map of spaces. Write $f_a \colon x_a \to y_a$ for $\phi(a)$. By composing $\phi$ with the source and target maps $s,t \colon \sX_1^\simeq \to \sX_0$, we obtain maps $\sigma, \tau \colon S \to \sX_0$ respectively. Suppose that for any pair $a,b \in S$, the 1-morphism $f_a$ is left-extending with respect to the target $y_b$ of $f_b$.
    Then there exists an oplax functor $\phi^\mathrm{tw} \colon \sigma^*\sX \oplaxfun \tau^*\sX$ such that for any $a,b \in S$, the functor $\sigma^*\sX(a,b) \to \tau^*\sX(a,b)$ on mapping categories is given by the composite
    \[
    \begin{tikzcd}
    \sX(x_a,x_b) \ar[r, "(f_b)_*"] & \sX(x_a,y_b) \ar[r, "(f_a)_!"] & \sX(y_a,y_b).
    \end{tikzcd}
    \]
\end{theorem}

\begin{remark}\label{remark:twisting-functor-into-X}
    Observe that we can compose $\phi^\mathrm{tw} \colon \sigma^* \sX \oplaxfun \tau^* \sX$ with the canonical map $\tau^* \sX \to \sX$ to obtain an oplax functor $\sigma^* \sX \oplaxfun \sX$.
\end{remark}

\begin{example}
    If we take $S = *$ in  \cref{theorem:Twisting-in-families}, then we recover \cref{theorem:Twisting-by-left-extending-functor}.
\end{example}

\begin{proof}[Proof of \cref{theorem:Twisting-in-families}]
    As in the proof of \cref{theorem:Twisting-by-adjunctions}, we may assume without loss of generality that $\sX$ is a 2-category. Since $\Arlax(\sX)_0 = \Map([1],\sX) \simeq \sX_1^\simeq$, the map $\phi \colon S \to \sX_1$ gives rise to a map $ \phi \colon S \to \Arlax(\sX)_0$. Pulling back $\Arlax(\sX)$ and the source-target projections $s,t \colon \Arlax(\sX) \to \sX$ along this map, we obtain a span
    \[\begin{tikzcd}
    	& {\phi^* \Arlax(\sX)} \\
    	{\tau^* \sX} && {\sigma^* \sX}
    	\arrow[from=1-2, to=2-1,"t"']
    	\arrow[from=1-2, to=2-3,"s"]
    \end{tikzcd}\]
    We claim that $s$ has a local left adjoint $l$ in the sense of \cref{ssec:local-adjunctions} and that the composite
    \[\begin{tikzcd}
        \sigma^* \sX \ar[r,"l"] & \phi^* \Arlax(\sX) \ar[r,"t"] & \tau^* \sX
    \end{tikzcd}\]
    has the desired properties of $\phi^\mathrm{tw}$. By (the dual of) \cref{prop: local-right-adjoint-is-lax}, the existence of $l$ follows if we can show that for any $a,b \in S$, the functor $s_{ab} \colon \Arlax(\sX)(f_a,f_b) \to \sX(x_a,x_b)$ admits a left adjoint $l_{ab}$. To see that this is the case, consider the comma object
    \[\begin{tikzcd}
		{\Arlax(\sX)(f_a,f_b)} & {\sX(x_a,x_b)} \\
		{\sX(y_a,y_b)} & {\sX(x_a,y_b)}
		\arrow["{(f_b)_*}", from=1-2, to=2-2]
		\arrow["t_{ab}"', from=1-1, to=2-1]
		\arrow[""{name=0, anchor=center, inner sep=0}, "s_{ab}", from=1-1, to=1-2]
		\arrow[""{name=1, anchor=center, inner sep=0}, "{f_a^*}"', from=2-1, to=2-2]
		\arrow[shorten <=7pt, shorten >=7pt, Rightarrow, from=1-2, to=2-1]
	\end{tikzcd}\]
    from \cref{proposition:lax-arrow-category-and-comma-categories}. It follows exactly as in the proof of \cref{theorem:Twisting-by-adjunctions} that the left adjoint $l_{ab}$ exists and that $t_{ab} \circ l_{ab}$ agrees with the composite
    \[
    \begin{tikzcd}
    \sX(x_a,x_b) \ar[r, "(f_b)_*"] & \sX(x_a,y_b) \ar[r, "(f_a)_!"] & \sX(y_a,y_b).
    \end{tikzcd} \qedhere
    \]
\end{proof}

\begin{remark}[The oplax structure maps]
    It follows from the proof of \cref{theorem:Twisting-in-families} and (the dual of) \cref{remark:lax-structure-maps-of-right-adjoint} that the oplax structure maps of $\phi^\mathrm{tw} \colon \sigma^* \sX \to \tau^* \sX$ can be described in exactly the same way as \cref{remark:Unwinding-oplax-structure-maps-twisting}.
\end{remark}

\subsection{Algebra in a 2-category}\label{sec:Algebras-in-2-categories}

Symmetric monoidal 2-categories are defined as functors $\sV \colon \Fin_* \to \Cattwo$ satisfying a Segal condition, analogous to the usual definition of a symmetric monoidal category. In this section, we will extend two well-known facts about symmetric monoidal categories to symmetric monoidal 2-categories: the fact that taking duals gives a symmetric monoidal equivalence $\sV \simeq \sV^\op$ if all objects of $\sV$ are dualizable, and the formula 
\[\Sym(x) = \coprod_{n \geq 0} x^{\otimes n}_{h \Sigma_n}\]
for the free commutative algebra on $x$. The results of this section will be used in \cref{sec:The-composition-product} to define the composition product of functor symmetric sequences and relate it to the usual composition product of symmetric sequences used to define operads.

Throughout this section, we will make use of a slightly different perspective on 2-categories: instead of considering complete Segal objects $\sX_\bullet$ in $\Cat$ such that $\sX_0$ is a space, we consider complete Segal objects in $\Cat$ for which the groupoid cores $\sX_\bullet^\simeq$ form a constant simplicial object in the category of spaces. We call these \emph{flipped 2-categories}. For details, the reader is referred to \cref{subsec:Flipped-2-categories}.

\subsubsection{Symmetric monoidal 2-categories}

Let $\Fin_*$ denote the category of finite pointed sets and $\angbr{n}$ the pointed set $\{ *,1,\ldots,n\}$. The usual definition of a symmetric monoidal category can be mimicked for 2-categories.

\begin{definition}[Symmetric monoidal 2-categories]
    A \emph{symmetric monoidal 2-category} is a functor $\sV = \sV_{\angbr{\bullet}} \colon \Fin_* \to \Cattwo$ such that $\sV_{\angbr{0}} \simeq *$ and the Segal maps
    \[
    \begin{tikzcd}
    \sV_{\angbr{n}} \ar[r] & \sV_{\angbr{1}} \times \ldots \times \sV_{\angbr{1}}
    \end{tikzcd}
    \]
    are equivalences. The \emph{underlying 2-category} of $\sV$ is defined as $\sV_{\angbr{1}}$. The full subcategory of $\Fun(\Fin_*,\Cattwo)$ spanned by the symmetric monoidal 2-categories is denoted by $\SymMonCattwo$.
\end{definition}

When no confusion can arise, we will often simply write $\sV$ for the underlying 2-category $\sV_{\angbr{1}}$ of $\sV$. Given objects $x$ and $y$ in (the underlying 2-category of) a symmetric monoidal 2-category $\sV$, their tensor product $x \otimes y$ is obtained by viewing $(x,y)$ as an object in $\sV_{\angbr{2}}$ and applying the functor $\sV_\angbr{2} \to \sV_\angbr{1} = \sV$ coming from the map $\angbr{2} \to \angbr{1}$ that sends $1$ and $2$ to $1$. The unit $\unit_\sV$ is the object of $\sV$ given by $* \simeq \sV_\angbr{0} \to \sV_\angbr{1} = \sV$.

We can now simply mimic the definition of the category of commutative algebras from \cite[Definition 2.1.3.1]{HA} in the 2-categorical setting. Let us write $\sV^\otimes \to \Fin_*$ for the cocartesian unstraightening in the sense of \cref{definition:cocartesian-fibration} of a symmetric monoidal 2-category $\sV$.

\begin{definition}[Commutative algebras in a symmetric monoidal 2-category]\label{definition:CAlg-2-category}
	Let $\sV$ be a symmetric monoidal 2-category. Then the 2-category $\CAlg(\sV)$ of commutative algebras in $\sV$ is defined as the full subcategory of 
	\[\Funtwo(\Fin_*,\sV^\otimes) \times_{\Funtwo(\Fin_*,\Fin_*)} \{\id_{\Fin_*}\}\]
	spanned by those sections $\Fin_* \to \sV^\otimes$ that take inerts to cocartesian lifts.
\end{definition}

In what follows, it will be convenient to take a slightly different perspective on symmetric monoidal 2-categories.
In \cref{subsec:Flipped-2-categories}, we observed that one can swap the simplicial coordinates in the definition of a 2-category, yielding the notion of a \emph{flipped 2-category} (\cref{definition:flipped-2-categories}).
This procedure of swapping simplicial coordinates yields an equivalence $(\cdot)^\fl \colon \Cattwo \eqarrow \Cattwofl \subset \Fun(\Delta^\op,\Cat)$ between the category $\Cattwo$ of 2-categories and the category $\Cattwofl$ of flipped 2-categories.
By furthermore swapping $\Fin_*$ and one copy of $\Delta^\op$ in the definition of a symmetric monoidal 2-category, one arrives at the following.

\begin{lemma}\label{lemma:Alternative-definition-symmetric-monoidal-2-category}
	Given a functor $\sW_\bullet \colon \Delta^\op \to \SymMonCat$, write $\sW_{\angbr{1},n}$ for the underlying category of $\sW_n$. The category $\SymMonCattwo$ is equivalent to the full subcategory of $\Fun(\Delta^\op, \SymMonCat)$ spanned by those functors $\sW_\bullet \colon \Delta^\op \to \SymMonCat$ for which $\sW_{\angbr{1}, \bullet}$ is a flipped 2-category.
\end{lemma}

\begin{proof}
	Given a symmetric monoidal 2-category $\sV \colon \Fin_* \to \Cattwo$, we can form the functor
	\[\begin{tikzcd}\sV^{\fl} \colon \Delta^\op \times \Fin_* \ar[r] & \Cat; \qquad ([n], \angbr{m}) \ar[r, maps to] & (\sV_{\angbr{m}})_n^\fl.\end{tikzcd}\]
	It follows from the equivalence $\Cattwo \simeq \Cattwofl$ of \cref{proposition:Equivalence-with-flipped-2-categories} that this yields an equivalence between $\SymMonCattwo$ and the full subcategory of $\Fun(\Delta^\op, \SymMonCat)$ described in the statement of this lemma.
\end{proof}

\begin{definition}
    A functor $\sW_\bullet \colon \Delta^\op \to \SymMonCat$ such that $\sW_{\angbr{1},\bullet}$ is a flipped 2-category will be called a \emph{flipped symmetric monoidal 2-category}. The full subcategory of $\Fun(\Delta^\op, \SymMonCat)$ spanned by the flipped symmetric monoidal 2-categories will be denoted by $\SymMonCattwofl$. Given a symmetric monoidal 2-category $\sV$, we write $\sV^\fl$ for the flipped symmetric monoidal 2-category corresponding to $\sV$ under the equivalence $\SymMonCattwo \simeq \SymMonCattwofl$ of \cref{lemma:Alternative-definition-symmetric-monoidal-2-category}.
\end{definition}

\begin{example}[Cartesian symmetric monoidal 2-categories]\label{example:Cartesian-2-monoidal-structure}
    Let $\sX$ be a 2-category that admits finite cartesian products. It follows by the same proof as \cref{proposition:Characterization-groupoidally-cocomplete-categories} that this is equivalent to the following:
    \begin{itemize}
        \item for any $n \geq 0$, the category $\sX^\fl_n$ admits finite products, and
        \item for any map $\alpha \colon [n] \to [m]$ in $\Delta$, the functor $\alpha^* \colon \sX_{[m]}^\fl \to \sX_{[n]}^\fl$ preserves finite products.
    \end{itemize}
    By \cite[Corollary 2.4.1.9]{HA}, the functor $\sX^\fl_\bullet \colon \Delta^\op \to \Cat$ lifts to a functor $\Delta^\op \to \SymMonCat$ that sends $[n]$ to the cartesian symmetric monoidal category $(\sX^\fl_n)^\times$. This is a flipped symmetric monoidal 2-category, so it follows from \cref{lemma:Alternative-characterization-P-algebras} that the cartesian product endows $\sX$ with the structure of a symmetric monoidal 2-category.
\end{example}

\subsubsection{Duality in symmetric monoidal 2-categories}

Let $\sV$ be a symmetric monoidal 2-category and let $x$ in $\sV$ be an object. We say that $x$ is \emph{dualizable} if there exists an object $x^\vee$ in $\sV$ together with 1-morphisms
\[\eta \colon \unit_\sV \to x^\vee \otimes x \quad \text{and} \quad \varepsilon \colon x \otimes x^\vee \to \unit_\sV\]
such that the composite maps
\begin{gather*}
    x \simeq \unit_\sV \otimes x \xrightarrow{\eta \otimes \id} x \otimes x^\vee \otimes x \xrightarrow{{\id}\otimes \varepsilon} x \otimes \unit_\sV \simeq x \\
    x^\vee \simeq x^\vee \otimes \unit_\sV \xrightarrow{{\id}\otimes \eta} x^\vee \otimes x \otimes x^\vee \xrightarrow{\varepsilon \otimes \id} \unit_{\sV} \otimes x^\vee \simeq x^\vee
\end{gather*}
are homotopic to the identity. Such a tuple $(x,x^\vee,\eta,\varepsilon)$ is called a \emph{duality datum}. Observe that a strong symmetric monoidal functor preserves any duality datum and hence dualizable objects. We let $\sV_\dual$ denote the full subcategory of $\sV$ spanned by dualizable objects. If every object of $\sV$ is dualizable, we say that $\sV$ is \emph{rigid}.

\begin{remark}\label{remark:dualizable-underlying-1-category}
    The definition of a duality datum only depends on the underlying symmetric monoidal 1-category of $\sV$. In particular, an object is dualizable in $\sV$ if and only if it is dualizable in $\sV^{\leq 1}$.
\end{remark}

The goal of this subsection is to show that, as for 1-categories, there is a symmetric monoidal equivalence $\sV_\dual \eqarrow \sV_\dual^\op$ that sends an object $x$ to its dual $x^\vee$. More precisely, we show the following.

\begin{proposition}\label{proposition:Rigid-2-category-self-dual}
    Let $\sV$ be a rigid symmetric monoidal 2-category. Then there is an equivalence of symmetric monoidal 2-categories $\sV \eqarrow \sV^\op$ given by sending an object $x$ to its dual $x^\vee$. This equivalence is natural in strong symmetric monoidal functors between rigid symmetric monoidal 2-categories $\sV$.
\end{proposition}

\begin{remark}
    Observe that in \cref{proposition:Rigid-2-category-self-dual} we obtain an equivalence $\sV \simeq \sV^\op$, and not an equivalence of $\sV$ with $\sV^\coop$. That is, taking duals only inverts the direction of the 1-morphisms and not of the 2-morphisms.
\end{remark}

For the proof of \cref{proposition:Rigid-2-category-self-dual}, we will use that the 1-categorical version of this statement, proved in \cite[Proposition 3.2.4]{Elliptic-1}, can be made natural.

\begin{proposition}\label{proposition:Rigid-1-category-self-dual-natural}
    Let $\sW$ be a rigid symmetric monoidal 1-category. Then there is a symmetric monoidal equivalence $\sW \eqarrow \sW^\op$ given by sending an object $x$ to its dual $x^\vee$. Moreover, this equivalence is natural in rigid symmetric monoidal categories $\sW$ and strong symmetric monoidal functors between them.
\end{proposition}

Before proving this proposition, let us deduce \cref{proposition:Rigid-2-category-self-dual} from it. The idea is to consider the flipped 2-category $\sV^\fl$ and apply \cref{proposition:Rigid-1-category-self-dual-natural} levelwise.

\begin{proof}[Proof of \cref{proposition:Rigid-2-category-self-dual}]
    Suppose $\sV$ is a rigid symmetric monoidal 2-category. By \cref{remark:dualizable-underlying-1-category}, $\sV^{\leq 1} = \sV^\fl_1$ is rigid as well. Since the degeneracy maps $\sV^\fl_0 \to \sV^\fl_n$ are strong symmetric monoidal and an equivalence on spaces of objects, we see that $\sV^\fl_n$ is rigid for any $n$.
    
    By \cref{lemma:Alternative-definition-symmetric-monoidal-2-category}, we may therefore identify rigid symmetric monoidal 2-categories with functors $\sW_\bullet \colon \Delta^\op \to \RigSymMonCat$ such that the underlying categories $\sW_{\angbr{1},\bullet}$ form a flipped 2-category. Here $\RigSymMonCat$ denotes the category of rigid symmetric monoidal categories and strong symmetric monoidal functors. The naturality of \cref{proposition:Rigid-1-category-self-dual-natural} means that there is a natural equivalence between the functor $(-)^\op \colon \RigSymMonCat \to \RigSymMonCat$ sending a rigid symmetric monoidal category to its opposite and the identity functor on $\RigSymMonCat$. 
    Whiskering this natural equivalence with the functor $\sW_\bullet \colon \Delta^\op \to \RigSymMonCat$ now yields the desired result, where we use that by \cref{remark:opposite-flipped-2-category}, the opposite of a flipped 2-category can be taken levelwise.
\end{proof}

For the proof of \cref{proposition:Rigid-1-category-self-dual-natural}, it will be useful to recall some basic facts on pairings of categories from \cite[Definition 5.2.1.5]{HA}, cf.\ \cite[Notation 3.35]{brantner2023pd}. Let us write $\Lambda^2_0$ for the span category $\bullet \leftarrow \bullet \rightarrow \bullet$.

\begin{definition}[Pairings of categories]
    A \emph{pairing} (of categories) consists of a triple $(\sC,\sD,\lambda \colon \sM \to \sC \times \sD)$ where $\sC$, $\sD$ and $\sM$ are categories and $\lambda$ is a right fibration. The category $\Pair$ is the full subcategory of $\Fun(\Lambda^2_0,\Cat)$ spanned by all pairings.
\end{definition}

Using straightening and the internal hom adjunction, any pairing $\lambda \colon \sM \to \sC \times \sD$ induces a functor $F_\lambda \colon \sC^\op \to \sP(\sD)$. We say that $\lambda$ is \emph{left representable} if $F_\lambda$ lands in $\sD$ and denote the corresponding functor by $\lambda^\mathrm{L} \colon \sC^\op \to \sD$. Observe that $F_\lambda(c) \in \sP(\sD)$ is representable if and only if $\sM \times_{\sC} \{c\}$ admits a terminal object; such an object will be called \emph{left universal}. Moreover, given a map of pairings
\[\begin{tikzcd}
    \cM \ar[r,"h"] \ar[d,"\lambda"] & \cM' \ar[d,"\mu"] \\
    \sC \times \sD \ar[r,"{f \times g}"] & \sC' \times \sD',
\end{tikzcd}\]
we obtain maps $F_\lambda(c)(d) \to F_{\mu}(f(c))(g(d))$. If $\lambda$ and $\mu$ are left representable, then these assemble into a natural transformation $g \circ \lambda^\mathrm{L} \Rightarrow \mu^\mathrm{L} \circ f$ of functors $\sC^\op \to \sD'$, which is invertible precisely if $h$ preserves left universal objects. This suggests that the category of left representable pairings and maps that preserve left universal objects is equivalent to $\Fun([1],\Cat)$, which is indeed the case by \cite[Proposition 2.2]{Li-Bland2015StackHigherInternal}.

\begin{definition}
    The category $\PairL$ is the (non-full) subcategory of $\Pair$ spanned by left representable pairings and maps of pairings that preserve left universal objects.
\end{definition}

\begin{proposition}[{\cite[Proposition 2.2]{Li-Bland2015StackHigherInternal}}]
    There is an equivalence $\PairL \simeq \Fun([1],\Cat)$ that sends a left representable pairing $\lambda \colon \sM \to \sC \times \sD$ to the functor $\lambda^\mathrm{L} \colon \sC^\op \to \sD$.
\end{proposition}

In particular, we obtain an equivalence between the category of commutative monoids in $\PairL$ and the category $\Fun([1],\SymMonCat)$. The proof of \cref{proposition:Rigid-1-category-self-dual-natural} uses that $\CMon(\PairL)$ admits a very convenient description. Observe that the cartesian product in $\PairL$ is obtained by restricting the cartesian product on $\Fun(\Lambda^2_0,\Cat)$. A bit of unpacking shows that $\CMon(\PairL)$ is the following (non-full) subcategory of $\CMon(\Fun(\Lambda^2_0,\Cat)) \simeq \Fun(\Lambda^2_0,\SymMonCat)$:
\begin{itemize}
    \item The objects of $\CMon(\PairL)$ are spans $\sW \leftarrow \sU \rightarrow \sW'$ of symmetric monoidal categories such that $\sU_{\angbr{1}} \to \sW_\angbr{1} \times \sW_\angbr{1}'$ is a left representable pairing with the property that
    \begin{itemize}
        \item the unit $\unit_\sU$ is left universal in $\sU_\angbr{1}$, and
        \item the collection of left universal objects in $\sU_\angbr{1}$ is closed under the tensor product.
    \end{itemize}
    \item The morphisms of $\CMon(\PairL)$ are those maps in $\Fun(\Lambda^2_0,\SymMonCat)$ for which the underlying morphism of pairings preserves left universal objects.
\end{itemize}

\begin{proof}[Proof of \cref{proposition:Rigid-1-category-self-dual-natural}]
    Let $\sW$ be a rigid symmetric monoidal category. As in the proof of \cite[Proposition 3.2.4]{Elliptic-1}, we can construct the pairing
    \[\begin{tikzcd}
        \lambda_\sW \colon (\sW \times \sW) \times_{\sW} \sW_{/\unit_\sW} \ar[r] & \sW \times \sW,
    \end{tikzcd}\] 
    where the pullback is taken along the tensor product $\otimes \colon \sW \times \sW \to \sW$.
    Observe that by \cite[Theorem 2.2.2.4]{HA}, there is a canonical symmetric monoidal structure on $\sW_{/\unit_\sW}$ such that $\sW_{/\unit_\sW} \to \sW$ is strong symmetric monoidal. This is easily seen to be natural in $\sW$, so we see that the pairing $\lambda_\sW$ is natural in $\sW$; more precisely, we have a functor
    \begin{align*}
        &\begin{tikzcd}[ampersand replacement=\&]
            \Lambda \colon \RigSymMonCat \ar[r] \& \Fun(\Lambda^2_0,\SymMonCat)
        \end{tikzcd}\\
        &\quad \begin{tikzcd}[ampersand replacement=\&]
            \sW \ar[r, maps to] \& (\lambda_\sW \colon (\sW \times \sW) \times_{\sW} \sW_{/\unit_\sW} \to \sW \times \sW)
        \end{tikzcd}
    \end{align*}
    where $\RigSymMonCat$ denotes the category of rigid symmetric monoidal categories.
    It follows as in the proof of \cite[Proposition 3.2.4]{Elliptic-1} that $\lambda_\sW \colon (\sW \times \sW) \times_{\sW} \sW_{/\unit_\sW} \to \sW \times \sW$ is a left representable pairing, where the left universal objects are those triples $(x,y, x \otimes y \to \unit_\sW)$ that exhibit $y$ as a dual of $x$. These are clearly closed under the tensor product, hence $\lambda_\sW$ defines an object of $\CMon(\PairL)$. We therefore obtain a strong symmetric monoidal functor $\lambda_{\sW}^\mathrm{L} \colon \sW^\op \to \sW$ that sends an object in $\sW$ to its dual. It follows as in the proof of \cite[Proposition 3.2.4]{Elliptic-1} that this is an equivalence.
    
    To see that this equivalence $\sW^\op \to \sW$ is natural in $\sW$, it suffices to show that the functor $\Lambda$ lands in the (non-full) subcategory $\CMon(\PairL)$, since $\Fun([1],\SymMonCat) \simeq \CMon(\PairL)$. We already verified that this is the case on objects, so it suffices to show that for any strong symmetric monoidal functor $f \colon \sW \to \sW'$, the induced functor
    \[
    \begin{tikzcd}
    (\sW \times \sW) \times_{\sW} \sW_{/\unit_{\sW}} \ar[r] & (\sW' \times \sW') \times_{\sW'} \sW'_{/\unit_{\sW'}}
    \end{tikzcd}
    \]
    preserves left universal objects. Since strong monoidal functors preserve duals, this is indeed the case. We conclude that $\Lambda$ defines a functor
    \[\begin{tikzcd} \RigSymMonCat \ar[r] & \CMon(\PairL) \simeq \Fun([1],\SymMonCat) \end{tikzcd}\]
    which is componentwise an equivalence $\sW^\op \eqarrow \sW$.
\end{proof}

\subsubsection{Free algebras in symmetric monoidal 2-categories}\label{subsec:Free-algebras-in-2-category}

We now show that, under suitable conditions, the explicit construction of free commutative algebras in a symmetric monoidal category also goes through in the setting of 2-categories. The strategy is to consider $\sV$ as a flipped 2-category and construct the free commutative algebra ``levelwise''. For this it will be useful to give an alternative description of the category $\CAlg(\sV)$.

\begin{construction}\label{construction:flipped-CAlg}
    Let $\sV^\fl$ be a flipped symmetric monoidal 2-category. Define $\CAlg(\sV^\fl_\bullet)$ as the composite
    \[\begin{tikzcd}[column sep = large]
        \Delta^\op \ar[r, "\sV^\fl"] & \SymMonCat \ar[r,"\CAlg(-)"] & \Cat.
    \end{tikzcd}\]
    Let $\CAlg(\sV^\fl)_n$ denote the full subcategory of $\CAlg(\sV^\fl_n)$ spanned by those objects that are in the image of the degeneracy map $\CAlg(\sV^\fl_0) \to \CAlg(\sV^\fl_n)$. These full subcategories form a simplicial object $\CAlg(\sV^\fl)$ in $\Cat$.
\end{construction}

\begin{lemma}\label{lemma:Alternative-characterization-P-algebras}
    The simplicial object $\CAlg(\sV^\fl)$ is a flipped 2-category. Moreover, it agrees with the flipped 2-category $\CAlg(\sV)^\fl$ corresponding to $\CAlg(\sV)$ under the equivalence $\Cattwo \simeq \Cattwofl$.
\end{lemma}

\begin{proof}
    Observe that by \cref{definition:CAlg-2-category} and the fact that pullbacks of (flipped) 2-categories are computed levelwise, the category $\CAlg(\sV)^\fl_n$ is equivalent to the pullback
    \[\Funtwo(\Fin_*,\sV^\otimes)^\fl_n \times_{\Funtwo(\Fin_*,\Fin_*)^\fl_n} \{\id_{\Fin_*} \}.\]
    Unwinding this expression using \cref{example:flipped-internal-hom}, it follows that this is precisely the category $\CAlg(\sV^\fl)_n$ from \cref{construction:flipped-CAlg}. We conclude that $\CAlg(\sV^\fl) \simeq \CAlg(\sV)^\fl$.
\end{proof}

Given a commutative algebra $\Fin_* \to \sV^\otimes$, we can precompose with $\{\angbr{1}\} \hookrightarrow \Fin_*$ to obtain an object in $\sV_\angbr{1}$. This yields a forgetful functor $U \colon \CAlg(\sV) \to \sV$.

\begin{proposition}[Free commutative algebras]\label{proposition:2-categorical-universal-property-of-sym}
    Let $\sV$ be symmetric monoidal 2-category whose underlying 2-category $\sV_\angbr{1}$ is groupoidally cocomplete in the sense of \cref{definition:Groupoidally-cocomplete-category} and suppose that the monoidal product $\otimes \colon \sV_\angbr{1} \times \sV_\angbr{1} \to \sV_\angbr{1}$ preserves groupoidal colimits in each variable separately. Then the forgetful functor $U \colon \CAlg(\sV) \to \sV$ admits a left adjoint $\Sym$. Moreover, for any $x$ in $\sV$, we have an equivalence
    \[U\Sym(x) \simeq \coprod_{n \geq 0} x^{\otimes n}_{h\Sigma_n}.\]
\end{proposition}

\begin{proof}
    Note that under the equivalence of \cref{lemma:Alternative-characterization-P-algebras}, the forgetful functor from $\CAlg(\sV)$ to $\sV$ corresponds to the levelwise forgetful functor $U_n \colon \CAlg(\sV^\fl)_n \to \sV^\fl_n$. Observe that the forgetful functor from the slightly larger category $\CAlg(\sV^\fl_n)$ to $\sV^\fl_n$ admits a left adjoint $F_n$ by \cite[Corollary 3.1.3.4]{HA} and \cref{proposition:Characterization-groupoidally-cocomplete-categories}. Since for any $\alpha \colon [m] \to [n]$, the functor $\alpha^* \colon \sV^\fl_n \to \sV^\fl_m$ is symmetric monoidal and preserves colimits indexed by groupoids, the functor $\alpha^* \colon \CAlg(\sV_n^{\fl}) \to \CAlg(\sV_m^{\fl})$ must preserve free commutative algebras. In particular, the diagram
    \[\begin{tikzcd}
		{\sV_n^\fl} & {\CAlg(\sV_n^{\fl})} \\
		{\sV_m^\fl} & {\CAlg(\sV_m^{\fl})}
		\arrow["{\alpha^*}"', from=1-1, to=2-1]
		\arrow[""{name=0, anchor=center, inner sep=0}, "{F_n}", from=1-1, to=1-2]
		\arrow["{\alpha^*}", from=1-2, to=2-2]
		\arrow[""{name=1, anchor=center, inner sep=0}, "{F_m}"', from=2-1, to=2-2]
		\arrow[shorten <=10pt, shorten >=10pt, Rightarrow, from=2-1, to=1-2]
	\end{tikzcd}\]
    commutes strictly for any $\alpha \colon [m] \to [n]$ and hence $F_n$ lands in the full subcategory $\CAlg(\sV^\fl)_n$ of $\CAlg(\sV^{\fl}_n)$. Moreover, by \cref{proposition:2-adjunctions-as-levelwise-adjunctions} it follows that these functors $F_n$ assemble into a left adjoint $\Sym$ of $U \colon \CAlg(\sV) \to \sV$.

    The explicit description of $U\Sym(x)$ now follows directly from the explicit description of the functors $F_n$ given in \cite[Example 3.1.3.14]{HA}. (Observe that loc.\ cit.\  assumes that the monoidal structure preserves all countable colimits in each of its variables. However, the proof also goes through if this only holds for colimits indexed by spaces.)
\end{proof}

\begin{remark}
    Given an $\infty$-operad $\cP^\otimes \to \Fin_*$ in the sense of \cite[Definition 2.1.1.10]{HA} and a symmetric monoidal 2-category $\sV$ satisfying the assumptions of \cref{proposition:2-categorical-universal-property-of-sym}, the same proof as above shows that free $\cP$-algebras in $\sV$ are computed in the underlying symmetric monoidal 1-category of $\sV$.
\end{remark}

\subsection{The 2-category \texorpdfstring{$\presl$}{Pr{\textasciicircum}L}}\label{sec:The-2-category-PrL}

We will now discuss one of the most important examples of a 2-category in this paper, namely the category $\presl$ of presentable categories and colimit preserving functors. A category $\sC$ is called \emph{presentable} if there exist a regular cardinal $\kappa$ and a small category $\sD$ that admits $\kappa$-small colimits such that $\sC \simeq \Ind_\kappa(\sD)$ \cite[Theorem 5.5.1.1]{HTT}. Recall the definition of a (locally) full subcategory of a 2-category from \cref{construction:1-full-subcategory}.

\begin{definition}[The 2-category $\presl$]
    Let $\twocatLargeCat$ denote the 2-category of (possibly large) categories. Then $\presl$ is defined as the locally full subcategory of $\twocatLargeCat$ spanned by the presentable categories and colimit preserving functors between them. Similarly, the 2-category $\presr$ is defined as the locally full subcategory of $\twocatLargeCat$ spanned by the presentable categories and accessible limit preserving functors.
\end{definition}

\begin{remark}
    The notation $\presl$ (resp.\ $\presr$) is conventionally used to denote the $1$-category of presentable categories and colimit preserving functors (resp.\ accessible limit preserving functors).
    We will commit a slight abuse of notation by using this notation for both the $1$-category and the $2$-category.
\end{remark}

Given presentable categories $\sC$ and $\sD$, we write $\FunL(\sC, \sD)$ for the category of colimit preserving functors.
This is the mapping category from $\sC$ to $\sD$ in $\presl$.
Similarly, we write $\FunR(\sC, \sD)$ for the mapping category in $\presr$.
By the adjoint functor theorem \cite[Corollary 5.5.2.9]{HTT}, a functor $F \colon \sC \to \sD$ lies in $\FunL(\sC, \sD)$ precisely if it is a left adjoint, and it lies in $\FunR(\sC, \sD)$ precisely if it is a right adjoint.
Recall that for any pair of presentable categories $\sC$ and $\sD$, Lurie constructs a tensor product $\sC \otimes \sD$ \cite[\S 4.8.1]{HA}.
This is characterized by the universal property
\[
\FunL(\sC \otimes \sD, \sE) \simeq \FunL(\sC, \FunL(\sD, \sE)).
\]
In other words, a colimit preserving functor $\sC \otimes \sD \to \sE$ is equivalent to a functor $\sC \times \sD \to \sE$ preserving colimits in both variables separately.
The tensor product  of presentable categories commutes with colimits in both variables.

\begin{proposition}\label{proposition:PrL-is-symmetric-monoidal-2-category}
    The tensor product of presentable categories makes $\presl$ into a symmetric monoidal 2-category.
\end{proposition}

\begin{proof}
    One can essentially mimic the construction of \cite[\S 4.8.1]{HA} and note that everything goes through 2-categorically. Throughout this proof, let $\twocatLargeCat$ denote the 2-category of (possibly large) categories and $\twocatLargeCat{}^\mathrm{L}$ the locally full subcategory of cocomplete categories and colimit-preserving functors. Note that the 2-category $\presl$ is a full subcategory of $\twocatLargeCat{}^\mathrm{L}$, so it suffices to show that $\twocatLargeCat{}^\mathrm{L}$ is a symmetric monoidal 2-category under the tensor product of cocomplete categories. The 2-category $\twocatLargeCat$ admits finite products, hence it is a symmetric monoidal 2-category under the cartesian product by \cref{example:Cartesian-2-monoidal-structure}. Let $\twocatLargeCat{}^\times \to \Fin_*$ denote the corresponding cocartesian fibration in the sense of \cref{definition:cocartesian-fibration}. Let $(\twocatLargeCat{}^\mathrm{L})^\otimes \subset \twocatLargeCat{}^\times$ be the locally full subcategory (see \cref{construction:1-full-subcategory}) defined as follows:
    \begin{enumerate}[(i)]
        \item an object $(\sC_1,\ldots,\sC_n)$ above $\langle n \rangle$ lies in $(\twocatLargeCat{}^\mathrm{L})^\otimes$ if and only if for every $i$, the category $\sC_i$ is cocomplete, and
        \item a 1-morphism $f \colon (\sC_1, \ldots, \sC_n) \to (\sD_1,\ldots,\sD_m)$ lying over $\alpha \colon \langle n \rangle \to \langle m \rangle$ belongs to $(\twocatLargeCat{}^\mathrm{L})^\otimes$ if and only if for every $1 \leq j \leq m$, the associated functor
        \[\prod_{\alpha(i) = j} \sC_i \to \sD_j\]
        preserves small colimits separately in each variable.
    \end{enumerate}
    To complete the proof, it suffices to show that $(\twocatLargeCat{}^\mathrm{L})^\otimes \to \Fin_*$ is a cocartesian fibration. This follows exactly as in the proof of Proposition 4.8.1.3 of \cite{HA}, using that the categories $\sP^\sK_{\sR}(\sC)$ of \cite[Proposition 5.3.6.2]{HTT} have a 2-categorical universal property and using \cref{lemma:When-locally-cocartesian-fibration-is-cocartesian} instead of \cite[Proposition 2.4.2.8]{HTT}.
\end{proof}

\begin{corollary}
The $2$-categories $\presl$ and $\presr$ both admit tensors and cotensors by small $1$-categories in the sense of \cref{definition:Tensors-and-cotensors}.
\end{corollary}
\begin{proof}
Let $I$ be a small $1$-category and let $\sC$ be a presentable category.
We write $\sP(I)$ for the presheaf category on $I$.
It is easily verified that both in $\presl$ and $\presr$, the presentable categories $\sC \otimes \sP(I)$ and $\Fun(I, \sC)$ provide respectively a tensoring and cotensoring of $\sC$ with $I$.
\end{proof}

\begin{corollary}
\label{cor: limits-and-colimits-presl-underlying}
    Let $I$ be a $1$-category. Then $I$-indexed limits and colimits in both $\presl$ and $\presr$ are computed in the underlying $1$-category.
\end{corollary}
\begin{proof}
    Since $\presl$ and $\presr$ admit all tensors and cotensors by small $1$-categories, this follows from \cref{corollary:cotensors-detect-colimits-underlying} and \cref{corollary:tensors-detect-limits-underlying}.
\end{proof}

Recall from \cite[Corollary 5.5.3.4]{HTT} that there is a canonical equivalence 
\begin{equation*}
    \presl \simeq (\presr)^{\mathrm{op}}
\end{equation*}
of $1$-categories that is the identity on objects and sends a functor to its right adjoint.
This equivalence is particularly useful for computing colimits (indexed by 1-categories) in $\presl$ and $\presr$.
Limits in both $\presl$ and $\presr$ are computed in $\Cat$ by Proposition 5.5.3.13 and Theorem 5.5.3.18 from \cite{HTT}.
In order to compute the colimit of a diagram in $\presl$, it follows from the equivalence above that we can do so by replacing all functors in the diagram by their right adjoints and then taking the limit in $\Cat$.
The same works for computing colimits in $\presr$, replacing right adjoints with their left adjoints.

This procedure for computing colimits in $\presl$ implies that certain limits and colimits in this category are equivalent.

\begin{proposition}
\label{prop: ambidexterity-prl}
    Let $X$ be a space and suppose that $F \colon X \to \presl$ is a functor.
    Then there is a canonical equivalence
    \[
    \colim_X F \simeq \lim_X F.
    \]
    The same holds for $\presr$.
\end{proposition}
\begin{proof}
    This is proved in Example 4.3.11 of \cite{HL-Ambidexterity} for the $1$-categories $\presl$ and $\presr$.
    The result therefore follows from \cref{cor: limits-and-colimits-presl-underlying}.
    We will sketch the argument for $\presl$ here; the one for $\presr$ goes exactly the same.

    Let $F \colon X \to \presl$ be a functor with $X$ a space. 
    We calculate the colimit of $F$ by passing to right adjoints (i.e.\ composing with the equivalence $\presl \simeq (\presr)^\op$ of 1-categories) and taking the limit of the resulting diagram $X^\op \to \presr$.
    However, since $X$ is a space, the functor $F$ lands in the groupoid core of $\presl$ and the passage to right adjoints is given by taking inverses.
    Composing with the equivalence $X \simeq X^\op$, we find that the resulting diagram $F \colon X \to \presr$ is equivalent to the diagram we started with, but now considered as taking values in $\presr$.
    Since limits in both $\presl$ and $\presr$ are computed in $\Cat$, we obtain the desired equivalence.
\end{proof}

\begin{example}
    The previous proposition implies that products and coproducts over arbitrary indexing sets agree in $\presl$, and are given by the product in $\Cat$.
    We will therefore use the symbol $\oplus$ to denote the product/coproduct in $\presl$.
    The proposition also implies that if $G$ is a discrete group and $\sC$ is a presentable category with a $G$-action, we have $\sC_{hG} \simeq \sC^{hG}$, where the colimit and limit are both taken in $\presl$.
\end{example}

We will use the following proposition in the next section to derive the formula for the composition product.

\begin{proposition}\label{prop: restriction-limit-right-adjoint-to-corestriction-colimit}
    Suppose that $p \colon X \to Y$ is a map of spaces and $F \colon Y \to \presl$ is a functor.
    Write
    \[
    \begin{tikzcd}
    p_! \colon \colim_X p^*F \ar[r] & \colim_Y F & \qquad p_* \colon \lim_Y F \ar[r] & \lim_X p^*F
    \end{tikzcd}
    \]
    for the colimit and limit comparison maps induced by $p$ in $\presl$.
    Under the equivalences $\colim_Y F \simeq \lim_Y F$ and $\colim_X p^*F \simeq \lim_X p^*F$ from the previous proposition, the functor $p_*$ is right adjoint to $p_!$.
\end{proposition}
\begin{proof}
    For any space $A$, we write $A^{\rhd}$ and $A^{\lhd}$ for the right and left cone on $A$ respectively. 
    Let $\Fun_{F}(X^{\rhd}, \presl)$ denote the category of functors $G \colon X^{\rhd} \to \presl$ together with an equivalence $G|_{X} \simeq p^*F$.
    We define $\Fun_{F}(X^{\lhd}, \presr)$ analogously.
    The proof of the previous proposition shows that we have an equivalence
    \[
        \Fun_{F}(X^{\rhd}, \presl)^\op \simeq \Fun_{F}(X^{\lhd}, \presr),
    \]
    given by combining $(X^{\rhd})^\op \simeq X^\lhd$ and $(\presl)^\op \simeq \presr$.
    This equivalence sends the object $X^{\rhd} \xrightarrow{p^{\rhd}} Y^\rhd \xrightarrow{\colim F} \presl$, corresponding to the comparison functor $p_!$, to the object $X^{\lhd} \xrightarrow{p^{\lhd}} Y^\lhd \xrightarrow{\lim F} \presr$, corresponding to the comparison functor $p_*$.
\end{proof}

\newpage

\section{Lax monoidal derivatives and the stable chain rule}\label{sec:Lax-monoidal-derivatives}

The goal of this chapter is to lift the Goodwillie derivatives functor to a lax $2$-functor
\[
\begin{tikzcd}
\partial_* \colon \diff \ar[r] & \pressymst
\end{tikzcd}
\]
from the $2$-category of differentiable categories and reduced finitary functors to the $2$-category of stable presentable categories and functor symmetric sequences.
We start by recalling the basic facts of Goodwillie calculus in \cref{sec: prelim-Goodwillie}.
In \cref{sec:The-composition-product}, we define the composition product of functor symmetric sequences and construct the $2$-category $\pressymst$.
We also prove some properties of the composition product that we will need later on.
\cref{sec: lax-structure-Goodwillie} is dedicated to constructing the lax $2$-functor $\partial_*$.
In \cref{sec: stable-chain-rule} we prove the stable chain rule, which says that the restriction of $\partial_*$ to the full subcategory on the \emph{stable} differentiable categories is a strong $2$-functor.

\subsection{Preliminaries on Goodwillie calculus}
\label{sec: prelim-Goodwillie}

We will now recall the basics of Goodwillie calculus.
Everything in this section -- except for the final lemma by Arone and Ching -- can already be found in some form in Goodwillie's original paper \cite{Goodwillie2003}.
The theory was developed in the $\infty$-categorical setting by Lurie in \cite[Chapter 6]{HA}, and for that reason we will mainly give references to this source.
Even for experts, who will find nothing new in this section, it might be helpful to skim its contents to familiarize oneself with our notation for various categories of (multilinear) symmetric functors.

\subsubsection{The Goodwillie tower}
To start off, we define a class of categories for which the theory of Goodwillie calculus is well-behaved.

\begin{definition}[Differentiable categories]
    We say a category $\sC$ is \emph{differentiable} if it is presentable, pointed, and for every filtered category $I$ the functor
    \[
    \begin{tikzcd}
    \colim \colon \Fun(I, \sC) \ar[r] & \sC
    \end{tikzcd}
    \]
    preserves finite limits.
\end{definition}

\begin{example}
    Any pointed compactly generated category is differentiable.
    More generally, any pointed compactly assembled category is differentiable (see \cite[\S 21.1.2]{SAG}).
\end{example}

\begin{example}
    Any stable presentable category is differentiable, since the colimit functor is exact in this case.
\end{example}

\begin{example}
    For any topos $\mathscr{X}$, the category of pointed objects $\mathscr{X}_{\ast /}$ is differentiable.
\end{example}

\begin{remark}
    Our definition of a differentiable category differs slightly from the one given by Lurie in \cite[Definition 6.1.1.6]{HA}, as he does not require a differentiable category to be presentable or pointed.
    Since we will only be dealing with pointed presentable categories, we have opted for this variant.
\end{remark}

Next, we introduce the kind of functors we will work with.

\begin{definition}[Reduced functors]
Let $\sC_1, \ldots, \sC_n$ and $\sD$ be presentable categories.
We call a functor $F \colon \prod_{1 \leq i \leq n} \sC_i \to \sD$ \emph{reduced} if $F(x_1, \ldots, x_n)$ is a final object whenever one of the $x_i$ is a final object.
\end{definition}

\begin{definition}[Finitary functors]
Let $\sC_1, \ldots, \sC_n$ and $\sD$ be presentable categories.
We call a functor $F \colon \prod_{1 \leq i \leq n} \sC_i \to \sD$ \emph{finitary} if it preserves filtered colimits in each variable separately.
\end{definition}

We will write $\Fun^*(\prod_{1 \leq i \leq k} \sC_i, \sD)$ for the full subcategory of $\Fun(\prod_{1 \leq i \leq k} \sC_i, \sD)$ on the reduced functors, $\Fun^{\omega}(\prod_{1 \leq i \leq k} \sC_i, \sD)$ for the full subcategory on the finitary functors, and $\Fun^{\ast, \omega}(\prod_{1 \leq i \leq k} \sC_i, \sD)$ for the full subcategory on the functors that are both reduced and finitary.

Throughout this paper, we will only apply Goodwillie calculus to functors that are reduced and finitary, even though the theory works in greater generality.
The reason for this restriction is that the chain rule only holds under these assumptions.
We will make it clear below exactly where these assumptions are needed.

\begin{definition}
    Let $\mathbf{P}(n)$ denote the poset of subsets of $\{1, \ldots, n\}$, ordered by inclusion.
    We call a functor $\mathcal{X} \colon \mathbf{P}(n) \to \sC$ an \emph{$n$-cube} in $\sC$.
    An $n$-cube is called \emph{cartesian} if it is a limit diagram.
    It is called \emph{strongly cocartesian} if it is left Kan extended from $\mathbf{P}_{\leq 1}(n)$, the subposet on the subsets of cardinality $\leq 1$.
\end{definition}

\begin{definition}[$n$-excisive functors]
    Let $\sC$ and $\sD$ be presentable categories and let $n \geq 0$.
    A functor $F \colon \sC \to \sD$ is \emph{$n$-excisive} if it sends strongly cocartesian $(n+1)$-cubes to cartesian $(n+1)$-cubes.
\end{definition}

We denote by $\Exc_n(\sC, \sD)$ the full subcategory of $\Fun(\sC, \sD)$ on the $n$-excisive functors.
We will use the notation $\Exc_n^{\ast, \omega}(\sC, \sD)$ for the full subcategory on the $n$-excisive functors that are additionally reduced and finitary.

The first fundamental theorem of Goodwillie calculus is that any reduced finitary functor between differentiable categories has a best approximation by an $n$-excisive functor.

\begin{theorem}[Goodwillie]
    Let $\sC$ and $\sD$ be differentiable categories.
    Then the inclusion $\Exc_n^{\ast, \omega}(\sC, \sD) \hookrightarrow \Fun^{\ast, \omega}(\sC, \sD)$ admits a left exact left adjoint
    \[
    P_n \colon \Fun^{\ast, \omega}(\sC, \sD) \to \Exc_n^{\ast, \omega}(\sC, \sD).
    \]
\end{theorem}
\begin{proof}
    This is proved in \cite[Theorem 6.1.1.10]{HA} without the reduced and finitary assumption.
    But since the construction of $P_n$ only involves finite limits and filtered colimits and $\sC$ and $\sD$ are differentiable, it is clear that $P_nF$ is reduced and finitary whenever $F$ is.
\end{proof}

For a  functor $F \colon \sC \to \sD$, we call $P_nF$ the \emph{$n$-excisive approximation} of $F$.
    
\begin{remark} \label{rem: multivariable-excisive-approx}
    There also is a version of the above theorem for functors in multiple variables.
    Let $\sC_1, \ldots, \sC_k$ and $\sD$ be differentiable categories, and let
    $\vec{n} = (n_1, \ldots, n_k)$ be a sequence of non-negative integers.
    Then we say a functor $F \colon \prod_{1 \leq i \leq k} \sC_i \to \sD$ is \emph{$\vec{n}$-excisive} if for every $1 \leq j \leq k$ and for every sequence of objects $\{x_i \in \sC_i\}_{i \neq j}$, the composite
    \[
    \begin{tikzcd}
        \sC_j \simeq \sC_j \times \prod_{i \neq j} \{x_i\} \ar[r, hook] & \sC_1 \times \cdots \times \sC_k \ar[r, "F"] & \sD 
    \end{tikzcd}
    \]
    is $n_j$-excisive.
    Let $\Exc_{\vec{n}}^{\ast, \omega}(\prod_{1 \leq i \leq k} \sC_i, \sD)$ be the full subcategory of $\Fun^{\ast, \omega}(\prod_{1 \leq i \leq k} \sC_i, \sD)$ spanned by the $\vec{n}$-excisive functors.
    Then by \cite[Proposition 6.1.3.6]{HA}, the inclusion of this subcategory has a left exact left adjoint
    \[
    \begin{tikzcd}
    P_{\vec{n}} \colon \Fun^{\ast, \omega}(\prod_{1 \leq i \leq k} \sC_i, \sD) \ar[r] & \Exc_{\vec{n}}^{\ast, \omega}(\prod_{1 \leq i \leq k} \sC_i, \sD).
    \end{tikzcd}
    \]
    This is given by composing the functors $P_{n_j}$ for $1 \leq j \leq k$, where we apply $P_{n_j}$ in the $j$th variable.
    We will write $\mlin$ for $P_{(1, \ldots, 1)}$, i.e.\ the functor applying $P_1$ in every variable.
\end{remark}

\begin{example}
    A functor $F \colon \sC \to \sD$ is $0$-excisive if and only if it is constant.
    The functor $P_0F$ is equivalent to the constant functor at $F(\ast)$, where $\ast$ is the zero object in $\sC$.
    If $F$ is reduced, this is just the constant functor at the zero object.
\end{example}

\begin{example}\label{example:explicit-description-P1}
A functor $F \colon \sC \to \sD$ is $1$-excisive if and only if it sends cocartesian squares to cartesian squares.
If $F$ is reduced, then $P_1F$ is the colimit of the sequential diagram
\[
\begin{tikzcd}
F \ar[r] & \Omega_\sD F \Sigma_\sC \ar[r] & \Omega_\sD^2 F \Sigma_\sC^2 \ar[r] & \Omega_\sD^3 F \Sigma_\sC^3 \ar[r] &\cdots
\end{tikzcd}
\]
In particular, $P_1\id_\sC \simeq \Omega^\infty_\sC \Sigma^\infty_\sC$.
We have a similar formula in the multivariable case.
Let $F \colon \prod_{1 \leq i \leq k} \sC_i \to \sD$ be a reduced multivariable functor.
Then
\[
\mlin F(x_1, \ldots, x_k) \simeq \colim_{n} \Omega_\sD^{kn} F(\Sigma_{\sC_1}^n x_1, \ldots, \Sigma_{\sC_k}^n x_k).
\]
\end{example}

\begin{definition}[Linear functors]
    A functor $F \colon \sC \to \sD$ between presentable categories is called \emph{linear} if it is both reduced and $1$-excisive.
    A multivariable functor $F \colon \prod_{1 \leq i \leq k} \sC_i \to \sD$ is called \emph{multilinear} if it is reduced and $(1, 1, \ldots, 1)$-excisive.
\end{definition}

It is not hard to prove that an $n$-excisive functor is also $(n+1)$-excisive.
So if $F \colon \sC \to \sD$ is a reduced finitary functor between differentiable categories, the universal property of the $(n+1)$-excisive approximation gives a canonical natural transformation $P_{n+1}F \to P_{n}F$.
Letting $n$ vary, we obtain a tower of functors under $F$:
\[
\begin{tikzcd}
    & & \vdots \ar[d] \\
    & & P_3 F \ar[d] \\
    & & P_2 F \ar[d] \\
F \ar[rr] \ar[urr] \ar[uurr]  & & P_1 F.
\end{tikzcd}
\]
This is known as the \emph{Goodwillie tower} of $F$.
We write $D_nF$ for the fiber of the map $P_nF \to P_{n-1}F$.

\begin{definition}[Homogeneous functors]
    Let $F \colon \sC \to \sD$ be a reduced finitary functor between differentiable categories.
    Then $F$ is called \emph{$n$-homogeneous} if it is $n$-excisive and $P_{n-1}F \simeq \ast$.
    We say $F$ is \emph{homogeneous} if it is $n$-homogeneous for some $n$.
    We write $\Homog_n^\omega(\sC, \sD)$ for the full subcategory of $\Fun^{\ast, \omega}(\sC, \sD)$ on the $n$-homogeneous functors.
\end{definition}

\begin{example}
    A functor is $1$-homogeneous precisely if it is linear.
\end{example}

\begin{example}
    As $P_n$ is left exact and $n$-excisive functors are closed under limits, we find that $D_n F$ is $n$-homogeneous.
\end{example}

\subsubsection{Symmetric functors and symmetric sequences}
One of the cornerstones of the theory of Goodwillie calculus is that an $n$-homogeneous functor is completely determined by its $n$th derivative, which is a symmetric multilinear functor in $n$ variables.
We will now introduce notation for the various kinds of symmetric functors we will use, and prove some basic results about them.

\begin{definition}[Symmetric functors]
    Let $\sC$ and $\sD$ be categories and let $n \geq 0$.
    We say a functor of the form $F \colon \sC^{\times n}_{h\Sigma_n} \to \sD$ is a \emph{symmetric functor (of arity $n$)}.
    Given such an $F$, we call the composite 
    \[
    \begin{tikzcd}
    \sC^{\times n} \ar[r] & \sC^{\times n}_{h\Sigma_n} \ar[r, "F"] & \sD
    \end{tikzcd}
    \]
    the underlying functor of $F$.
\end{definition}

 We write $\SymFun_n(\sC, \sD)$ for the category $\Fun(\sC^{\times n}_{h\Sigma_n}, \sD)$ of symmetric functors of arity $n$.
    For any property of multivariable functors, we will say a symmetric functor has this property if the underlying functor does.
So we will speak of reduced, multilinear, and finitary symmetric functors.

We will need sequences of symmetric functors of different arities.
Consider $\Cat$ as a symmetric monoidal category under the cartesian product. Given a category $\sC$, we let $\Sym^\times(\sC)$ denote the underlying category of the free commutative algebra on $\sC$ in $\Cat$.
There is the following formula:
\[
\Sym^\times(\sC) = \coprod_{n \geq 0} \sC^{\times n}_{h\Sigma_n}.
\]
Therefore, a functor $F \colon \Sym^\times(\sC) \to \sD$ to a category $\sD$ amounts to a sequence of symmetric functors $F_n \colon \sC^{\times n}_{h\Sigma_n} \to \sD$, one for each $n \geq 0$.
It will sometimes be convenient to use coordinate-free notation for such sequences, in the following sense.
We have an equivalence
\[
\Sym^\times(\sC) = \colim_{I \in \Fin^{\cong}} \sC^{\times I},
\]
where $\Fin^{\cong}$ denotes the category of finite sets and bijections.
Given $F \colon \Sym^\times(\sC) \to \sD$ and a finite set $I$, we write $F_I \colon \sC^{\times I} \to \sD$ for the $I$-component of $F$, which is an $\Aut(I)$-equivariant functor.

Let $\SymFun(\sC, \sD)$ denote the category $\Fun(\Sym^\times(\sC), \sD)$.
Suppose $\sC$ and $\sD$ are pointed and admit all colimits and finite limits. 
We say $F \colon \Sym^\times(\sC) \to \sD$ is reduced, finitary, or multilinear if each component $F_n$ is.
We write $\SymFun^{\ast, \omega}(\sC, \sD)$ and $\SymFun^{\mathrm{mlin}, \omega}(\sC, \sD)$ for the full subcategories on the reduced finitary and multilinear finitary symmetric functors, respectively.

There is a fully faithful embedding $\SymFun_n(\sC, \sD) \hookrightarrow \SymFun(\sC, \sD)$ with essential image those $F \colon \Sym^{\times}(\sC) \to \sD$ for which $F_m$ is the constant functor at the initial object whenever $m \neq n$.
We will say such an $F$ is concentrated in arity $n$.
The arity zero component of a functor $F \colon \Sym^{\times}(\sC) \to \sD$ is given by an object of $\sD$.
We say $F$ is \emph{positive} if $F_0$ is the initial object of $\sD$.
We write $\SymFun_{\geq 1}(\sC, \sD)$ for the full subcategory on the positive functors.

Any finitary symmetric functor has a best approximation by a reduced finitary symmetric functor.

\begin{proposition} \label{prop: reduction}
    Let $\sC$ and $\sD$ be differentiable categories and let $n \geq 0$.
    The inclusion functor $\SymFun_{n}^{\ast, \omega}(\sC, \sD) \hookrightarrow \SymFun_n^{\omega}(\sC, \sD)$ admits a right adjoint
    \[
    \begin{tikzcd}
    \mathrm{red} \colon \SymFun_n^{\omega}(\sC, \sD) \ar[r] & \SymFun_n^{\ast, \omega}(\sC, \sD).
    \end{tikzcd}
    \]
    We call $\mathrm{red}(F)$ the reduction of $F$.
\end{proposition}
\begin{proof}
    Without the finitary assumption, this is \cite[Proposition 6.1.4.3]{HA}.
    The functor $\mathrm{red}$ can be described as follows.
    Given $x_1, \ldots, x_n \in \sC$, we can form an $n$-cube in $\sD$ given by
    \[
    \mathbf{P}(n) \longrightarrow \sD, \quad S \longmapsto F(y_1, \ldots, y_n) \quad 
        \begin{cases}
        y_i = x_i \text{ if $i \not \in S$}, \\
        y_i = \ast \text{ else}.
        \end{cases}
    \]
    The object $\mathrm{red}(F)(x_1, \ldots, x_n)$ is given by the total fiber of this cube.
    Since $\sD$ is differentiable and the total fiber is obtained by taking finitely many finite limits, it follows that $\mathrm{red}(F)$ is finitary whenever $F$ is.
\end{proof}

\begin{remark}
    Taking the product of the right adjoints from the previous proposition for all $n \geq 0$, we obtain a functor
    \[
    \begin{tikzcd}
        \mathrm{red} \colon \SymFun^{\omega}(\sC) \ar[r] & \SymFun^{\ast, \omega}(\sC),
    \end{tikzcd}
    \]
    which is right adjoint to the inclusion $\SymFun^{\ast, \omega}(\sC, \sD) \hookrightarrow \SymFun^{\omega}(\sC, \sD)$.
\end{remark}

It is also possible to multilinearize a symmetric functor.

\begin{proposition} \label{prop: sym-multilinearization}
    Let $\sC$ and $\sD$ be differentiable categories and let $n \geq 0$.
    The inclusion functor $\SymFun_n^{\mathrm{mlin}, \omega}(\sC, \sD) \hookrightarrow \SymFun_n^{\ast, \omega}(\sC, \sD)$ admits a left exact left adjoint
    \[
    \begin{tikzcd}
    \mlin \colon \SymFun_n^{\ast, \omega}(\sC, \sD) \ar[r] & \SymFun_n^{\mathrm{mlin}, \omega}(\sC, \sD).
    \end{tikzcd}
    \]
    We call $\mlin(F)$ the multilinearization of $F$.
\end{proposition}
\begin{proof}
    By \cref{rem: multivariable-excisive-approx}, the inclusion $\Fun^{\mathrm{mlin}, \omega}(\sC^{\times n}, \sD) \hookrightarrow \Fun^{\ast, \omega}(\sC^{\times n}, \sD)$ admits a left exact left adjoint, given by $(1, \ldots, 1)$-excisive approximation:
    \[
    \begin{tikzcd}
    P_{\vec{1}} \colon \Fun^{\ast, \omega}(\sC^{\times n}, \sD) \ar[r] & \Fun^{\mathrm{mlin}, \omega}(\sC^{\times n}, \sD).
    \end{tikzcd}
    \]
    This adjoint is equivariant for the $\Sigma_n$-action on both sides, and we obtain the desired functor by taking fixed points.
\end{proof}

\begin{remark}
    We can take the product of the left adjoints from this proposition to obtain a left exact functor
    \[
    \begin{tikzcd}
        \mlin \colon \SymFun^{\ast, \omega}(\sC, \sD) \ar[r] & \SymFun^{\mathrm{mlin}, \omega}(\sC, \sD),
    \end{tikzcd}
    \]
    which is left adjoint to the inclusion $\SymFun^{\mathrm{mlin}, \omega}(\sC, \sD) \hookrightarrow \SymFun^{\ast, \omega}(\sC, \sD)$.
\end{remark}

We also consider sequences of symmetric functors that preserve all colimits in every variable separately.

\begin{definition}[Functor symmetric sequences]
    Let $\sC$ and $\sD$ be presentable categories.
    We call a functor $F \colon \Sym^\times(\sC) \to \sD$ a \emph{functor symmetric sequence} (or \emph{symmetric sequence} for short) if the underlying functor of each component $F_n$ preserves colimits in each variable separately.
\end{definition}

We write $\SymFunL(\sC, \sD)$ for the full subcategory of $\Fun(\Sym^\times(\sC), \sD)$ on the symmetric sequences.
If $\sC$ and $\sD$ are pointed, we let $\SymFunL_{\geq 1}(\sC, \sD)$ denote the full subcategory of $\SymFunL(\sC, \sD)$ on the positive symmetric sequences.

\begin{remark}
\label{rem: terminology-sym-seq}
    Our terminology is not standard: a symmetric sequence in a category $\sC$ is usually defined as a sequence of objects with $\Sigma_n$-action in $\sC$, one for each $n \geq 0$.
    In this paper, we will refer to these as \emph{ordinary} symmetric sequences to distinguish them from functor symmetric sequences.
    We write $\SymSeq(\sC) = \Fun(\Fin^{\cong}, \sC)$ for the category of ordinary symmetric sequences in $\sC$.
    Here $\Fin^{\cong}$ denotes the category of finite sets and bijections.
\end{remark}

\begin{example} \label{ex: funsymseq-is-symseq-spectra}
    A functor $F \colon \Sp^{\times n}_{h\Sigma_n} \to \Sp$ that preserves colimits in each variable separately is completely determined by $F(\Sph, \ldots, \Sph)$ considered as a spectrum with $\Sigma_n$-action.
    It follows that giving a functor symmetric sequence from spectra to itself is equivalent to giving an ordinary symmetric sequence in spectra.
    Consequently, we have an equivalence of categories $\SymFunL(\Sp, \Sp) \simeq \SymSeq(\Sp)$.
\end{example}

It will be useful to know that the category of symmetric sequences can also be defined by only making use of constructions in the category $\presl$.
For a presentable category $\sC$, we will write $\Sym(\sC)$ for the free commutative algebra on $\sC$ in $\presl$ with respect to the tensor product of presentable categories.
There is the formula
\[
\Sym(\sC) \simeq \bigoplus_{n \geq 0} \sC^{\otimes n}_{h\Sigma_n},
\]
where $\oplus$ denotes the coproduct in $\presl$.
Since $\Sym^{\times}(\sC)$ is the free symmetric monoidal category on $\sC$, there is a canonical symmetric monoidal functor $\mathrm{can} \colon \Sym^{\times}(\sC) \to \Sym(\sC)$ which restricts to the standard inclusion $\sC \hookrightarrow \Sym(\sC)$ on $\sC$.
It follows from the universal property of the tensor product of presentable categories that for any presentable category $\sD$, restriction along $\mathrm{can}$ induces an equivalence
\[
\begin{tikzcd}
\mathrm{can}^* \colon \FunL(\Sym(\sC), \sD) \ar[r, "\sim"] & \SymFunL(\sC, \sD). 
\end{tikzcd}
\]
In other words, a colimit preserving functor $\Sym(\sC) \to \sD$ is the same as a symmetric sequence from $\sC$ to $\sD$.

Note that a functor to a stable category is linear if and only if it is exact.
So if $\sD$ is stable, then $\SymFunL(\sC, \sD) = \SymFun^{\mathrm{mlin}, \omega}(\sC, \sD)$.

\begin{proposition} \label{prop: multilinear-is-stable}
Let $\sC$ and $\sD$ be presentable categories. 
The functor
\[
\begin{tikzcd}
\SymFunL(\Sp(\sC), \Sp(\sD)) \ar[r] & \SymFun^{\mathrm{mlin}, \omega}(\sC, \sD)
\end{tikzcd}
\]
given by sending a symmetric sequence $F \colon \Sym^\times(\Sp(\sC)) \to \Sp(\sD)$ to the composite
\[
\begin{tikzcd}[sep = large]
\Sym^\times(\sC) \ar[r, "\Sym^\times(\Sigma_\sC^\infty)"] & \Sym^\times(\Sp(\sC)) \ar[r, "F"] & \Sp(\sD) \ar[r, "\Omega_\sD^\infty"] & \sD 
\end{tikzcd}
\]
is an equivalence of categories.
\end{proposition}
\begin{proof}
    It follows from \cite[Corollary 6.2.3.22]{HA} that for each $n \geq 0$ we have an equivalence
    \[
    \Fun^{\mathrm{mlin}, \omega}(\Sp(\sC)^{\times n}, \Sp(\sD)) \xrightarrow{\sim} \Fun^{\mathrm{mlin}, \omega}(\sC^{\times n}, \sD), \quad F \mapsto \Omega_\sD^\infty \circ F \circ (\Sigma_\sC^\infty)^{\times n}. 
    \]
    This is clearly equivariant for the $\Sigma_n$-action on both sides, and the result now follows by taking homotopy fixed points.
\end{proof}

\subsubsection{Derivatives}

We now give the classification of homogeneous functors in terms of symmetric multilinear functors and the definition of the derivatives functor.

\begin{construction}[Symmetric cross-effects] \label{constr: cross-effects}
Let $\sC$ and $\sD$ be differentiable categories.
We equip $\sC$ with its cocartesian symmetric monoidal structure.
Since $\Sym^\times(\sC)$ is the free symmetric monoidal category on $\sC$, there is a unique symmetric monoidal functor $p \colon \Sym^\times(\sC) \to \sC$ which restricts to the identity on $\sC$.
The underlying functor of $p$ sends a sequence $(x_1, \ldots, x_n)$ to the coproduct $x_1 \sqcup \cdots \sqcup x_n$.
The \emph{symmetric cross-effects} functor is the composite
\[
\begin{tikzcd}
\cross \colon \Fun^{\omega}(\sC, \sD) \ar[r, "p^*"] & \SymFun^{\omega}(\sC, \sD) \ar[r, "\mathrm{red}"] & \SymFun^{\ast, \omega}(\sC, \sD).
\end{tikzcd}
\]
We write $\cross_n$ for the arity $n$ component of this functor.
\end{construction}

\begin{remark}
    If $F \colon \sC \to \sD$ is a reduced functor, then $\cross_0(F) \simeq \ast$.
    Therefore, the symmetric cross effects restrict to a functor $\cross \colon \Fun^{\ast, \omega}(\sC, \sD) \to \SymFun^{\ast, \omega}_{\geq 1}(\sC, \sD)$.
\end{remark}

As both $p^*$ and $\mathrm{red}$ are right adjoint functors, the symmetric cross-effects functor admits a left adjoint $p_!$ given by left Kan extension along $p$.

\begin{theorem}[Classification of homogeneous functors] \label{thm: classification-homogeneous}
    Let $\sC$ and $\sD$ be differentiable categories.
    Then the symmetric cross effects functor restricts to an equivalence
    \[
    \begin{tikzcd}
       \cross_n \colon \Homog_n^\omega(\sC, \sD) \ar[r, "\sim"] & \SymFun^{\mathrm{mlin}, \omega}_n(\sC, \sD).
    \end{tikzcd}
    \]
\end{theorem}
\begin{proof}
    See \cite[Theorem 6.1.4.7]{HA}.
\end{proof}

The inverse of $\cross_n$ can be described explicitly.
If $\sD$ is a stable differentiable category, then it is given by $p_!$ restricted to $\SymFun^{\mathrm{mlin}, \omega}_n(\sC, \sD)$.
By a cofinality argument given in the proof of loc.\ cit.\ one can see that this left Kan extension sends a symmetric multilinear functor $F \colon \sC^{\times n}_{h\Sigma_n} \to \sD$ to the functor given by the formula
\[
x \longmapsto F(x, \ldots, x)_{h \Sigma_n}.
\]
For a general symmetric multilinear functor $F \colon \Sym^\times(\sC) \to \sD$, it follows from \cref{prop: multilinear-is-stable} that we can write $F \simeq \Omega_\sD^\infty G$ for some symmetric multilinear functor $G \colon \Sym^\times(\sC) \to \Sp(\sD)$.
The corresponding $n$-homogeneous functor is then given by
\[
x \longmapsto \Omega_\sD^\infty G(x, \ldots, x)_{h\Sigma_n}.
\]

\begin{definition}[Derivatives]
    Let $F \colon \sC \to \sD$ be a reduced finitary functor between differentiable categories and let $n \geq 1$.
    By the preceding discussion, the functor $D_nF$ is of the form
    \[
    x \longmapsto \Omega_\sD^\infty \partial_nF(\Sigma_\sC^\infty x, \ldots, \Sigma_\sC^\infty x)_{h\Sigma_n}
    \]
    for some finitary symmetric multilinear functor $\partial_n F \colon \Sp(\sC)^{\times n}_{h \Sigma_n} \to \Sp(\sD)$.
    We call $\partial_nF$ the \emph{$n$th derivative} of $F$.
\end{definition}

It is possible to calculate the $n$th derivative of a functor without first determining its $n$-excisive approximation.

\begin{proposition}[{\cite[Remark 6.1.3.23]{HA}}]
    Let $F \colon \sC \to \sD$ be a reduced finitary functor between differentiable categories.
    Then the $n$th derivative of $F$ is naturally equivalent to the multilinearization of the $n$th cross-effect of $F$:
    \[
    \partial_nF \simeq \mlin \cross_n F.
    \]
\end{proposition}

It follows from this proposition that the derivatives functor can be defined in the following way.
This is the definition that we will be using throughout the rest of this paper.

\begin{definition}
    Let $\sC$ and $\sD$ be differentiable categories.
    The \emph{Goodwillie derivatives functor} is given by the composite
    \[
    \begin{tikzcd}[sep = large]
    \partial_* \colon \Fun^{\ast, \omega}(\sC, \sD) \ar[r, "\cross"] & \SymFun^{\ast, \omega}_{\geq 1}(\sC, \sD) \ar[r, "\mlin"] & \SymFunL_{\geq 1}(\Sp(\sC), \Sp(\sD)).
    \end{tikzcd}
    \]
    Here, we have used \cref{prop: multilinear-is-stable} to identify the target of the multilinearization functor with a category of symmetric sequences.
\end{definition}

\subsubsection{The Goodwillie tower of composite functors}
We finish this section by recalling a lemma due to Arone and Ching, which will play a crucial role in our proof of the chain rule.

\begin{lemma}[cf.\ {\cite[Proposition 1.3.1]{AroneChing2011}}] \label{lem: Arone-Ching-lemma}
    Suppose $G \colon \sC \to \sD$ and $F \colon \sD \to \sE$ are functors between differentiable categories and let $n \geq 1$.
    Then
    \begin{enumerate}[label={\upshape(\arabic*)}]
        \item if $G$ is reduced, the canonical map $P_n(FG) \to P_n((P_nF)G)$ is an equivalence;
        \item if $F$ is finitary, the canonical map $P_n(FG) \to P_n(F(P_nG))$ is an equivalence.
    \end{enumerate}
    In particular, if $G$ is reduced and $F$ is finitary there is an equivalence
    \[
    \begin{tikzcd}
    \partial_n(FG) \ar[r, "\sim"] & \partial_n((P_nF)(P_nG)).
    \end{tikzcd}
    \]
\end{lemma}
\begin{proof}
We start by briefly recalling the construction of $P_n$.
Suppose $\sA$ is a differentiable category.
For $x \in \sA$, let $C(x) \colon \mathbf{P}(n+1) \to \sA$ be the strongly cocartesian $(n+1)$-cube determined by the conditions that $C_{\emptyset}(x) = x$ and $C_{S}(x) = \ast$ whenever $|S| = 1$.
This construction is natural in $x$, and we write $C$ for the functor sending an object $x$ to the cube $C(x)$.
Given any other differentiable category $\sB$ and functor $K \colon \sA \to \sB$, we write $T_n K$ for the functor
\[
T_n K = \lim_{\emptyset \neq S \in \mathbf{P}(n+1)} K \circ C_S.
\]
There is an evident natural transformation $f_K \colon K \to T_n K$ and $P_n K$ is defined as the sequential colimit
\[
\begin{tikzcd} [sep = large]
    P_n K = \colim ( K \ar[r, "f_K"] & T_n K \ar[r, "f_{T_n K}"] & T_n (T_n K) \ar[r, "f_{T_n (T_n K)}"] & \cdots).
\end{tikzcd}
\]

In order to prove (1), we first construct a natural map $a \colon T_n(F) G \to T_n(FG)$.
To construct a map of $(n+1)$-cubes $C(G(x)) \to G(C(x))$ that is natural in $x$, it suffices to do so after restricting to $\mathbf{P}_{\leq 1}(n+1)$, since $C(G(x))$ is left Kan extended from this subcategory.
But since $G$ is reduced there is an essentially unique equivalence $C(G(x))_{\leq 1} \to G(C(x))_{\leq 1}$ on these restrictions that is the identity on the initial vertex.
Applying $F$ to the resulting natural transformation and taking limits over $\mathbf{P}(n+1) \setminus \{ \emptyset \}$, we obtain the map $a$.
By construction, the composite
\[
\begin{tikzcd}
    FG \ar[r] & T_n(F)G \ar[r, "a"] & T_n(FG)
\end{tikzcd}
\]
is the canonical map.
By iterating $a$ and taking the colimit, we obtain a natural map $P_n(F)G \to P_n(FG)$ compatible with the canonical maps from $FG$.
Since the target of this map is $n$-excisive, the universal property of $P_n$ yields a map $r \colon P_n(P_n(F)G) \to P_n(FG)$ such that the composite
\[
\begin{tikzcd}
    P_n(FG) \ar[r] & P_n(P_n(F)G) \ar[r, "r"] & P_n(FG)
\end{tikzcd}
\]
is the identity. 
Note that by construction, $r$ is natural in both $F$ and $G$.

Now consider the diagram
\[
\begin{tikzcd}
P_n(FG) \ar[r] \ar[d] & P_n(P_n(F)G) \ar[r, "r"] \ar[d] & P_n(FG) \ar[d] \\
P_n(P_n(F)G) \ar[r] & P_n(P_n(P_n(F))G) \ar[r, "r"] & P_n(P_n(F)G).
\end{tikzcd}
\]
The two horizontal arrows on the left and all vertical arrows are induced by the canonical map $F \to P_n(F)$.
Since $P_n(F) \to P_n(P_n(F))$ is an equivalence, the middle vertical map is one as well.
Therefore $P_n(FG) \to P_n(P_n(F)G)$ is the retract of an equivalence, which implies that it is an equivalence itself.
This completes the proof of (1).

For (2), observe that the assembly map
\[
\begin{tikzcd}
F(\lim_{\emptyset \neq S \in \mathbf{P}(n+1)} G(C_S)) \ar[r] & \lim_{\emptyset \neq S \in \mathbf{P}(n+1)} FG(C_S).
\end{tikzcd}
\]
provides a natural transformation $F (T_nG) \to T_n(FG)$ such that the composite 
\[
\begin{tikzcd}
FG \ar[r] & F (T_nG) \ar[r] & T_n(FG)
\end{tikzcd}
\]
is the canonical map.
From here on, the argument proceeds exactly as it did for (1).
The only point that needs some care is the equivalence $\colim_k F (T_n^kG) \simeq F (P_nG)$ that is used in the proof.
This holds since $F$ is finitary.
\end{proof}

\begin{remark}
The assumptions on $F$ and $G$ cannot be omitted.
For reducedness, this is easy to see:
if $G$ is a constant functor at an object not equivalent to the zero object, then it is clear that $P_n(FG)$ need not be equivalent to $P_n((P_nF)G)$.
Kuhn provided a counterexample to (2) in case $F$ is not finitary, see \cite[Example 7.4]{kuhnChromatic}.
It is for this reason that we only work with reduced finitary functors throughout this paper.
\end{remark}

We will need the following multivariate version of this lemma.

\begin{corollary}\label{cor: multivariate Arone-Ching}
      Let $\sC$, $\sD$, and $\sE$ be differentiable categories.
    Suppose we have reduced finitary functors $F \colon \sD^{\times k} \to \sE$ and $G_i \colon \sC^{\times n_i} \to \sD$ for $1 \leq i \leq k$.
    Let $P_{\vec{n}}{F}$ denote the functor obtained by taking the $n$th excisive approximation in each variable separately, and similarly for the $G_i$.
    Then for $n \geq 1$ both of the maps
        \[
        \begin{tikzcd}
        P_{\vec{n}}(F(G_1, \ldots, G_k)) \ar[r] & P_{\vec{n}}((P_{\vec{n}}{F})(G_1, \ldots, G_k))
        \end{tikzcd}
        \]
        and
        \[
        \begin{tikzcd}
        P_{\vec{n}}(F(G_1, \ldots, G_k)) \ar[r] & P_{\vec{n}}(F(P_{\vec{n}}{G_1}, \ldots, P_{\vec{n}}{G_k}))
        \end{tikzcd}
        \]
    are equivalences.
\end{corollary}
\begin{proof}
        This follows by iteratively applying \cref{lem: Arone-Ching-lemma}.
\end{proof}

\begin{example}\label{example:derivatives-Sigma-Loops}
    We saw in \cref{example:explicit-description-P1} that for a differentiable category $\sC$, the unit $\eta \colon \id_\sC \to \Omega^\infty_\sC \Sigma^\infty_\sC$ becomes an equivalence upon applying $P_1$.
    One can show that the counit $\varepsilon \colon \Sigma^\infty_\sC \Omega^\infty_\sC \to \id_{\Sp(\sC)}$ similarly induces an equivalence $P_1(\Sigma^\infty_\sC \Omega^\infty_\sC) \simeq \id_{\Sp(\sC)}$.
    Namely, by \cref{prop: multilinear-is-stable} it suffices to show that $P_1(\Sigma^\infty_\sC \Omega^\infty_\sC) \Sigma^\infty_\sC \to \Sigma^\infty_\sC$ is an equivalence.
    Since $\Sigma^\infty_\sC$ preserves colimits, this is equivalent to the map $P_1(\Sigma^\infty_\sC \Omega^\infty_\sC \Sigma^\infty_\sC) \to P_1(\Sigma^\infty_\sC)$.
    By the triangle identities, the composite
    \[
    \begin{tikzcd}
      P_1(\Sigma^\infty_\sC) \ar[r, "\eta"] &  P_1(\Sigma^\infty_\sC \Omega^\infty_\sC \Sigma^\infty_\sC) \ar[r, "\varepsilon"] & P_1(\Sigma^\infty_\sC)
    \end{tikzcd}
    \]
    is invertible.
    Since $\id_\sC \to \Omega^\infty_\sC \Sigma^\infty_\sC$ is an equivalence after applying $P_1$, it follows from \cref{lem: Arone-Ching-lemma} and two-out-of-three that $P_1(\varepsilon) \colon P_1(\Sigma^\infty_\sC \Omega^\infty_\sC) \to P_1(\id_{\Sp(\sC)}) = \id_{\Sp(\sC)}$ must be an equivalence.
\end{example}
    
\subsection{The composition product}\label{sec:The-composition-product}

For symmetric sequences $F \in \SymFunL(\sD, \sE)$ and $G \in \SymFunL(\sC, \sD)$, we can form their \emph{composition product} $F \circ G$, which is a symmetric sequence from $\sC$ to $\sE$.
If $F$ and $G$ are positive and $I$ is a finite set, this is given by the formula
\[
(F \circ G)_I = \coprod_{E \in \mathrm{Part}(I)} F_{E} \circ \{ G_J \}_{J \in E},
\]
where $\mathrm{Part}(I)$ denotes the set of partitions of $I$.
The $\circ$ on the right-hand side denotes composition of multivariable functors.
Note that the classical composition product of ordinary symmetric sequences is given by a very similar formula.

In this section, we construct a 2-category $\pressym$.
Its objects are presentable categories and the mapping categories $\pressym(\sC, \sD)$ are given by $\SymFunL(\sC, \sD)$.
Composition in this 2-category is given by the composition product.
In particular, this shows that the composition product defines a monoidal structure on $\SymFunL(\sC, \sC)$.
The full subcategory $\pressymst \subseteq \pressym$ spanned by the stable presentable categories will be the target of the lax derivatives functor that we construct in the next section.
After having defined $\pressym$, we prove some basic properties of the composition product.
We end this section by comparing the monoidal structure on $\SymFunL(\Sp, \Sp)$ with the composition product of symmetric sequences of spectra as constructed in \cite{brantnerThesis}.

\subsubsection{The 2-category \texorpdfstring{$\pressym$}{PrSym}}
We will construct the composition product of symmetric sequences by means of cocommutative coalgebras in the $2$-category $\presl$.
\begin{definition}
    Let $\sX$ be a symmetric monoidal $2$-category.
    A \emph{cocommutative coalgebra} in $\sX$ is defined as a commutative algebra in $\sX^\op$.
    The $2$-category of cocommutative coalgebras in $\sX$ is given by
    \[
    \coCAlg(\sX) \coloneqq \CAlg(\sX^\op)^\op,
    \]
    where the latter was defined in \cref{definition:CAlg-2-category}.
\end{definition}
Recall from \cref{proposition:PrL-is-symmetric-monoidal-2-category} that the tensor product of presentable categories equips $\presl$ with the structure of a symmetric monoidal $2$-category.
Our definition of the 2-category $\pressym$ is based on the following observation.

\begin{proposition} \label{prop: sym-is-cofree}
    Let $\sC$ be a presentable category.
    Then $\Sym(\sC)$ is the underlying category of the cofree commutative coalgebra on $\sC$ in the $2$-category $\presl$.
\end{proposition}
\begin{proof}
    By \cref{prop: ambidexterity-prl}, the tensor product of presentable categories commutes with limits indexed by groupoids in both variables.
    The result now follows by applying \cref{proposition:2-categorical-universal-property-of-sym} to $(\presl)^\op$.
\end{proof}

As a consequence of this proposition, we get the following equivalence for presentable categories $\sC$ and $\sD$:
\begin{equation} \label{eq: sym-seq-coalg}
\FunL(\Sym(\sC), \sD) \simeq \FunL_{\mathrm{coCAlg}}(\Sym(\sC), \Sym(\sD)).
\end{equation}
The right-hand side denotes the mapping category in the $2$-category $\mathrm{coCAlg}(\presl)$.
Composition in this $2$-category then gives a functor
\[
\begin{tikzcd}
\FunL(\Sym(\sD), \sE) \times \FunL(\Sym(\sC), \sD) \ar[r, "- \circ -"] & \FunL(\Sym(\sC), \sE),
\end{tikzcd}
\]
and we will prove below that on objects this is given by the formula for the composition product.
This motivates the following definition.

\begin{definition}
We define $\pressym$ to be the full subcategory of the 2-category $\mathrm{coCAlg}(\presl)$ spanned by the cofree commutative coalgebras $\Sym(\sC)$.
We write $\pressymst \subseteq \pressym$ for the full subcategory spanned by the objects $\Sym(\sC)$ where $\sC$ is a presentable stable category.
\end{definition}

\begin{remark}
We will think of objects of $\pressym$ as presentable categories $\sC$ rather than cocommutative coalgebras $\Sym(\sC)$.
This is justified by the fact that there is an equivalence of cocommutative coalgebras $\Sym(\sC) \simeq \Sym(\sD)$ if and only if $\sC$ is equivalent to $\sD$.
It follows from \cref{eq: sym-seq-coalg} that the mapping category $\pressym(\sC, \sD)$ is equivalent to $\SymFunL(\sC, \sD)$.
\end{remark}

\begin{remark}
    In the language of the theory of comonads, $\pressym$ is called the coKleisli category of the comonad $\Sym$.
\end{remark}

\begin{definition}
\label{def: composition product}
    Let $\sC$ be a presentable category.
    Then $\SymFunL(\sC, \sC)$ is the endomorphism category of the object $\sC$ in $\pressym$, and as such acquires a monoidal structure.
    We call this the \emph{composition product} monoidal structure.
\end{definition}

\subsubsection{Properties of the composition product}

We now derive the formula for the composition product of symmetric sequences and use it to prove some of its basic properties.
In order to do so, we first need a better understanding of the comultiplication map on cofree commutative coalgebras in $\presl$.
Throughout this section, we will only work with the $1$-categories $\presl$ and $\presr$, so that we can make use of the equivalence $\presl \simeq (\presr)^\op$ (2-categorically, we would have $\presl \simeq (\presr)^\coop$).

Suppose $\sC$ is a presentable category.
Then $\Sym(\sC)$ is the free commutative algebra on $\sC$ in $\presl$, so that there is a multiplication functor $\mu \colon \Sym(\Sym(\sC)) \to \Sym(\sC)$.
This functor can be understood in the following way.
We have an equivalence $\Sym(\sC) \simeq \colim_{I \in \Fin^{\cong}} \sC^{\otimes I}$.
It is easy to deduce that
\[
\Sym(\Sym(\sC)) \simeq \colim_{I \to J} \sC^{\otimes I},
\]
where the colimit is taken over $\Fun([1], \Fin)^{\cong}$, the groupoid core of the arrow category of finite sets.
Under this equivalence, the multiplication map $\mu$ is the functor induced on colimits by $\ev_0 \colon \Fun([1], \Fin)^{\cong} \to \Fin^{\cong}$.

By \cref{prop: sym-is-cofree}, $\Sym(\sC)$ is also the cofree commutative coalgebra on $\sC$ in $\presl$.
This is equivalent to the statement that $\Sym(\sC)$ is the free commutative algebra on $\sC$ in $\presr$, since $\presl \simeq (\presr)^\op$.
Hence, we obtain a comultiplication functor $\Delta \colon \Sym(\sC) \to \Sym(\Sym(\sC))$ in $\presl$ and a multiplication functor $\nu \colon \Sym(\Sym(\sC)) \to \Sym(\sC)$ in $\presr$.
The map $\nu$ is again induced on colimits by $\ev_0 \colon \Fun([1], \Fin)^{\cong} \to \Fin^{\cong}$, but since we are now taking colimits in $\presr$, $\nu$ need not be equivalent to $\mu$.

These various multiplication and comultiplication functors are related in the following way.

\begin{proposition} \label{prop: adjoints-of-comultiplication}
    Let $\sC$ be a presentable category.
    Then $\mu$ is left adjoint to $\Delta$, and $\Delta$ is left adjoint to $\nu$:
    \[
\begin{tikzcd}[sep = huge]
    \Sym(\Sym(\sC)) 
    \ar[rr, "\mu", bend left = 25]  
    \ar[rr, "\nu", bend right = 25]  & & 
    \Sym(\sC) 
    \ar[ll, "\Delta"']. 
\end{tikzcd}
\]
\end{proposition}
\begin{proof}
The cofree commutative coalgebra on $\sC \in \presl$ is by definition the free commutative algebra on $\sC$ in $(\presl)^\op$.
Since there is a symmetric monoidal equivalence $(\presl)^\op \simeq \presr$ given by sending a functor to its right adjoint, it follows immediately that $\Delta$ is left adjoint to $\nu$.

Since $\Fin^{\cong}$ is a groupoid, it follows from \cref{prop: ambidexterity-prl} that in both $\presl$ and $\presr$ we have
\[
\colim_{I \in \Fin^{\cong}} \sC^{\otimes I} \simeq \lim_{I \in \Fin^{\cong}} \sC^{\otimes I}.
\]
The same holds for $\Fun([1], \Fin)^{\cong}$.
Note that the limit on the right can be evaluated in $\Cat$, since limits in both $\presl$ and $\presr$ are computed in this way.
As we explained above, $\nu$ is induced on colimits in $\presr$ by $\ev_0 \colon \Fun([1], \Fin)^{\cong} \to \Fin^{\cong}$.
By \cref{prop: restriction-limit-right-adjoint-to-corestriction-colimit}, the left adjoint $\Delta$ of $\nu$ is therefore the induced map on limits in $\Cat$ by $\ev_0$.
The same reasoning applied to $\mu$ shows that the right adjoint of $\mu$ is also the functor induced on limits by $\ev_0$ in $\Cat$.
It follows that the right adjoint of $\mu$ is equivalent to $\Delta$.
\end{proof}

We now use this proposition to derive the formula for the composition product.

\begin{notation}
    Let $\{I_j\}_{j \in J}$ be a collection of finite sets indexed by a finite set $J$.
    Suppose $F_J \colon \sD^J \to \sE$ and $\{G_{I_j} \colon \sC^{I_j} \to \sD\}_{j \in J}$ are symmetric functors.
    We write
    \[
    \begin{tikzcd}
    F_J \circ \{G_{I_j}\}_{j \in J} \colon \sC^{\sqcup_{j \in J} I_j} \ar[r] & \sE
    \end{tikzcd}
    \]
    for the composite, which is an $\Aut(J) \ltimes \prod_{j \in J} \Aut(I_j)$-equivariant functor.
    We will often use this notation to denote the underlying functor of this composite as well.
\end{notation}

\begin{proposition}[Formula for the composition product]\label{prop:formula-composition-product}
Let $F \colon \Sym(\sD) \to \sE$ and $G \colon \Sym(\sC) \to \sD$ be symmetric sequences between presentable categories.
Then for any finite set $I$, the underlying functor of the $I$-component of the composition product of $F$ and $G$ is given by the formula
\[
(F \circ G)_I \simeq \colim_{\alpha \colon I \to J} F_J \circ \{G_{\alpha^{-1}(j)}\}_{j \in J},
\]
where the colimit is taken over the category $(\Fin_{I/})^{\cong}$ whose objects are maps of finite sets $I \to J$ and whose morphisms are bijections $J \xrightarrow{\sim} J'$ making the triangle commute.
\end{proposition}
\begin{proof}
Unwinding the definition of the composition product, we find that $F \circ G$ is given by the composite
\[
\begin{tikzcd}[sep = large]
    \Sym(\sC) \ar[r, "\Delta"] & \Sym(\Sym(\sC)) \ar[r, "\Sym(G)"] & \Sym(\sD) \ar[r, "F"] & \sE.
\end{tikzcd}
\]
Since $\Sym^{\times}(\sC)$ is the free commutative algebra on $\sC$ in $\Cat$, we obtain a multiplication map $m \colon \Sym^{\times}(\Sym^{\times}(\sC)) \to \Sym^{\times}(\sC)$.
We have a canonical symmetric monoidal functor $\mathrm{can} \colon \Sym^{\times}(\sC) \to \Sym(\sC)$, so that we obtain a commutative square
\[
\begin{tikzcd}[sep = large]
    \Sym^{\times}(\Sym^{\times}(\sC)) \ar[r, "\mathrm{can}^2"] \ar[d, "m"] & \Sym(\Sym(\sC)) \ar[d, "\mu"] \\
    \Sym^{\times}(\sC) \ar[r, "\mathrm{can}"] & \Sym(\sC).
\end{tikzcd}
\]
By the universal property of the tensor product of presentable categories, restriction along $\mathrm{can}$ induces a fully faithful embedding
\[
\begin{tikzcd}[sep = large]
   \mathrm{can}^* \colon \FunL(\Sym(\sC), \sE) \ar[r, hook] & \Fun(\Sym^{\times}(\sC), \sE)
\end{tikzcd}
\]
with essential image sequences of functors that preserve colimits in each variable separately.
We will also write $F$ and $G$ for the images of these symmetric sequences under this embedding.
Restriction along $\mathrm{can}^2$ also gives a fully faithful embedding of $\FunL(\Sym(\Sym(\sC)), \sE)$ in $\Fun(\Sym^{\times}(\Sym^{\times}(\sC)), \sE)$.
We therefore obtain the following commutative diagram
\[
\begin{tikzcd}[sep = large]
    \FunL(\Sym(\Sym(\sC)), \sE) 
    \ar[d, hook, "(\mathrm{can}^2)^*"] & 
    \FunL(\Sym(\sC), \sE) \ar[l, "\mu^*"'] \ar[d, hook, "\mathrm{can}^*"] \\
    \Fun(\Sym^{\times}(\Sym^{\times}(\sC)), \sE) 
    & \Fun(\Sym^{\times}(\sC), \sE) \ar[l, "m^*"'].
\end{tikzcd}
\]
Since $\mu$ is left adjoint to $\Delta$, which is itself a colimit preserving functor, $\Delta^*$ is left adjoint to $\mu^*$.
The functor $m^*$ has a left adjoint $m_!$ given by left Kan extension along $m$.
Therefore, we get an induced Beck--Chevalley transformation 
\[
\begin{tikzcd}
m_! \circ (\mathrm{can}^2)^* \ar[r] & \mathrm{can}^* \circ \Delta^*.
\end{tikzcd}
\]
To show this is an equivalence, it suffices to show that applying $m_!$ to a functor in the essential image of $(\mathrm{can}^2)^*$ yields a functor in the essential image of $\mathrm{can}^*$.
A functor $\Sym^\times(\Sym^\times(\sC)) \to \sE$ lies in the image of $(\mathrm{can}^2)^*$ precisely when all of its components $\sC^{\times n} \to \sE$ preserve colimits in each variable. The image of $\mathrm{can}^*$ is characterized by the same condition.
Since the value of $m_!$ on objects is computed as a colimit, and colimit preserving functors are stable under colimits, we find that the Beck--Chevalley transformation is an equivalence.

The upshot is that the composition product of $F \circ G$ can also be described as the left Kan extension along $m \colon \Sym^{\times}(\Sym^{\times}(\sC)) \to \Sym^{\times}(\sC)$ of the composite
\[
\begin{tikzcd}[sep = large]
\Sym^{\times}(\Sym^{\times}(\sC)) \ar[r, "\Sym^{\times}(G)"] & \Sym^{\times}(\sD) \ar[r, "F"] & \sE.
\end{tikzcd}
\]
Since $m$ is induced on colimits by the map of groupoids $\ev_0 \colon \Fun([1], \Fin)^{\cong} \to \Fin^{\cong}$, we can apply straightening to obtain a pullback square of categories
\[
\begin{tikzcd}[sep = large]
    \Sym^{\times}(\Sym^{\times}(\sC)) \ar[r, "m"] \ar[d, "p"] & \Sym^{\times}(\sC) \ar[d, "q"] \\
    \Fun([1], \Fin)^{\cong} \ar[r, "\ev_0"] & \Fin^{\cong}, 
\end{tikzcd}
\]
where the vertical arrows are cocartesian fibrations.
The map $\ev_0$ is a cocartesian fibration as well, since it is a functor between groupoids.
The functor $m$ sends $(\ev_0 \circ p)$-cocartesian edges to $q$-cocartesian edges.

We now use that a left Kan extension along $m$ can be computed fiberwise.
Choose $I \in \Fin^{\cong}$.
The fiber of $q$ over $I$ is equivalent to $\sC^I$.
It is easily computed that the fiber of ${\ev_0} \circ p$ over $I$ is given by $\colim_{I \to J} \sC^I$, where the indexing category is $(\Fin_{I/})^{\cong}$ and the diagram is constant at $\sC^I$.
Note that the restriction of the left Kan extension $m_!(F \circ \Sym^{\times}(G))$ to  the fiber of $q$ over $I$ is the underlying functor of $(F \circ G)_I$.
By \cite[Proposition 4.3.3.10]{HTT}, this restriction can be computed by first restricting $F \circ \Sym^{\times}(G)$ to the fiber of ${\ev_0} \circ p$ over $I$ and then left Kan extending along the map
\[
\begin{tikzcd}
    \colim_{I \to J} \sC^I \ar[r] & \sC^I
\end{tikzcd}
\]
induced by $m$ on fibers.
This yields the desired formula.
\end{proof}

\begin{remark} \label{rem: composition-connected-symmetric-sequences}
Observe that if $F$ and $G$ are positive, the colimit in the formula for $F \circ G$ can be indexed by just the surjections of finite sets.
The formula then simplifies to
\[
(F \circ G)_I \simeq \coprod_{E \in \mathrm{Part}(I)} F_{E} \circ \{ G_J \}_{J \in E},
\]
where $\mathrm{Part}(I)$ denotes the set of partitions of $I$.
\end{remark}

It follows from this formula that positive symmetric sequences are closed under the composition product.

\begin{definition}
    Let $\pressymc \subseteq \pressym$ denote the locally full subcategory on the positive symmetric sequences.
    Similarly, we write $\pressymstc$ for the corresponding locally full subcategory of $\pressymst$.
\end{definition}

The following proposition summarizes the exactness properties of the composition product.
Recall that a cosimplicial diagram in a category is called \emph{$k$-coskeletal} if it is right Kan extended from the full subcategory $\Delta_{\leq k} \subseteq \Delta$ spanned by the linearly ordered sets of cardinality $\leq k$.
We will use the term \emph{finite totalization} for the limit over a cosimplicial diagram that is $k$-coskeletal for some $k \geq 0$. 
A functor $F \colon \sC \to \sD$ is said to \emph{preserve finite totalizations} if for any $k$-coskeletal cosimplicial diagram $X^\bullet$ in $\sC$, the coassembly map $F(\Tot X^\bullet) \to \Tot F(X^\bullet)$ is an equivalence. 
Note that this does not imply that $F(X^\bullet)$ is coskeletal.

\begin{proposition} \label{prop: exactness-composition-product}
    Let $\sC$, $\sD$ and $\sE$ be presentable categories.
    The composition product
    \[
    \begin{tikzcd}
    \circ \colon \SymFunL_{\geq 1}(\sD, \sE) \times \SymFunL_{\geq 1}(\sC, \sD) \ar[r] & \SymFunL_{\geq 1}(\sC, \sE)
    \end{tikzcd}
    \]
    preserves colimits in the first variable and sifted colimits in the second variable.
    If $\sD$ and $\sE$ are additionally stable, it also preserves aritywise finite totalizations in both variables.
\end{proposition}  
\begin{proof}
    Given positive integers $n$ and $k_1, \ldots, k_n$, we can consider the functor
    \[
    \begin{tikzcd}
     M \colon \SymFunL_{\geq 1}(\sD, \sE) \times \SymFunL_{\geq 1}(\sC, \sD) \ar[r] & \Fun(\sC^{k_1 + \cdots + k_n}, \sE),
    \end{tikzcd}
    \]
    sending a pair of symmetric sequences $(F, G)$ to the composite $F_n(G_{k_1}, \ldots, G_{k_n})$.
    This clearly preserves all colimits in the first variable.
    It also preserves sifted colimits in the second variable as $F_n$ preserves colimits in each variable separately and for a sifted category $S$ the diagonal functor $S \to S^{\times n}$ is right cofinal.

    Recall that $\Delta_{\leq k}$ receives a left cofinal functor from a punctured $(k+1)$-cube, which is a finite diagram.
    From this it follows that the condition of being right Kan extended from $\Delta_{\leq k}$ only involves finite limits, so that finite totalizations are preserved by any exact functor.
    Therefore $M$ preserves aritywise finite totalizations in the first variable if $\sD$ and $\sE$ are stable.
    In this case $M$ also preserves aritywise finite totalizations in the second variable, since $F_k$ is exact in each variable separately and the diagonal functor $\Delta \to \Delta^{\times k}$ is left cofinal.

    If $F$ and $G$ are positive symmetric sequences, then by \cref{rem: composition-connected-symmetric-sequences} the underlying functor of ${F} \circ {G}$ is in each arity a finite coproduct of functors of the form $F_n(G_{k_1}, \ldots, G_{k_n})$.
    Since colimits (and in the stable case, finite limits) of symmetric sequences are computed on the level of underlying functors, this completes the proof.
\end{proof}

\subsubsection{Turning a symmetric sequence into a functor}\label{subsec:Lambda}

Given a symmetric sequence $F \colon \Sym(\sC) \to \sD$, we can extract a functor $\Lambda_F \colon \sC \to \sD$, given by the formula
\begin{equation}
\Lambda_F(x) = \coprod_{n \geq 0} F_n(x, \ldots, x)_{h\Sigma_n}.
\end{equation}
This construction is compatible with the composition product, in the sense that $\Lambda_{F \circ G} \simeq \Lambda_F \circ \Lambda_G.$
In fact, it even extends to a functor $\pressym \to \Cat$ of $2$-categories.

\begin{proposition} \label{prop: functorization}
    Let $\Lambda \colon \pressym \to \Cat$ be the $2$-functor corepresented by $\ast \in \pressym$.
    Then $\Lambda$ is the identity on objects and sends a symmetric sequence $F \colon \Sym(\sC) \to \sD$ to the functor $\Lambda_F \colon \sC \to \sD$.
\end{proposition}
\begin{proof}
Suppose $\sC$ is a presentable category. 
Since $\Sym(\ast) \simeq \Spc$, it follows that we have an equivalence 
\[
\SymFunL(\ast, \sC) \simeq \sC.
\]
This shows that $\Lambda$ is the identity on objects.
Now suppose that $\sC$ and $\sD$ are presentable categories.
By \cref{prop: corepresentable-functors-of-2-cats}, the functor $\SymFunL(\sC, \sD) \to \Fun(\sC, \sD)$ induced by $\Lambda$ on mapping categories is adjoint to the composition functor
\[
\begin{tikzcd}
\SymFunL(\sC, \sD) \times \SymFunL(\ast, \sC) \ar[r, "\circ"] & \SymFunL(\ast, \sD) \simeq \sD.
\end{tikzcd}
\]
Let $x \in \sC$ and consider it as a symmetric sequence $\Sym(\ast) \to \sC$ and let $F \colon \Sym(\sC) \to \sD$ be a symmetric sequence. 
Then the composition product $F \circ x$ is naturally equivalent to $\Lambda_{F}(x)$ as can be read off from the formula in \cref{prop:formula-composition-product}.
Hence $\Lambda$ sends $F$ to the functor $\Lambda_F$.
\end{proof}

\begin{remark} \label{rem: Lambda-F-reduced-finitary}
    For $F \in \SymFunL(\sC, \sD)$, the functor $\Lambda_F \colon \sC \to \sD$ is finitary. 
    If $F$ is additionally positive (i.e.\ $F_0$ is zero) then $\Lambda_F$ is also reduced.
\end{remark}

We will also write $\Lambda \colon \SymFunL(\sC, \sD) \to \Fun^\omega(\sC, \sD)$ for the functor induced by $\Lambda$ on mapping categories.
In the next proposition, we will give an alternative description of this functor, which will be necessary to relate it to the cross-effects functor later on.
Recall that we write $p \colon \Sym^\times(\sC) \to \sC$ for the functor sending a sequence $(x_1, \ldots, x_n)$ to the coproduct $x_1 \sqcup \cdots \sqcup x_n$.
More precisely, equipping $\sC$ with the cocartesian symmetric monoidal structure, $p$ is the symmetric monoidal functor determined by the requirement that the composite $\sC \hookrightarrow \Sym^\times(\sC) \xrightarrow{p} \sC$ is the identity.

\begin{proposition} \label{prop: functorization-is-kan-extension}
Let $\sC$ and $\sD$ be presentable categories.
The functor
\[
\begin{tikzcd}
\Lambda \colon \SymFunL(\sC, \sD) \ar[r] & \Fun^\omega(\sC, \sD)
\end{tikzcd}
\]
is given by left Kan extension along $p \colon \Sym^\times(\sC) \to \sC$.
\end{proposition}
\begin{proof}
    Let $F \in \SymFunL(\sC, \sD)$ and suppose that $x \in \sC$.
    We consider $x$ as a symmetric sequence from $\ast$ to $\sC$. 
    This is the unique colimit preserving functor
    \[
    \begin{tikzcd}
        \Sym(\ast) \simeq \Spc \ar[r, "x"] & \sC
    \end{tikzcd}
    \]
    that sends $\ast$ to $x$.
    The composition product $F \circ x$, which is equivalent to $\Lambda_F(x)$, is given by the composite
    \[
    \begin{tikzcd}[sep = huge]
     \Sym(\ast) \ar[r, "\Delta"] & \Sym(\Sym(\ast)) \ar[r, "\Sym(x)"] & \Sym(\sC) \ar[r, "F"] & \sD.
    \end{tikzcd}
    \]
    Here $\Delta \colon \Sym(\ast) \to \Sym(\Sym(\ast))$ is the comultiplication functor.
    Under the equivalences $\Sym(\ast) \simeq \Spc$ and $\Sym(\Sym(\ast)) \simeq \Fun(\Fin^{\cong}, \Spc)$ it is given by restriction along the functor $\Fin^{\cong} \to \ast$.
    Since the above description of $F \circ x$ is natural in $x$, it follows that $\Lambda \colon \FunL(\Sym(\sC), \sD) \to \Fun^\omega(\sC, \sD)$ is given by precomposition with
    \[
    \begin{tikzcd}
    \sC \simeq \FunL(\Spc, \sC) \ar[r, "\Sym"] & \FunL(\Sym(\Spc), \Sym(\sC)) \ar[r, "\Delta^*"] & \FunL(\Spc, \Sym(\sC)) \simeq \Sym(\sC).
    \end{tikzcd}
    \]
    This composite is easily seen to be equivalent to the functor $S \colon \sC \to \Sym(\sC)$ given by $x \mapsto \coprod_{n \geq 0} x^{\otimes n}_{h\Sigma_n}$.
    In other words, it is the composite of the inclusion $\sC \to \Sym(\sC)$ with
    $ \mathrm{forget}_{\CAlg} \circ \mathrm{free}_{\CAlg} \colon \Sym(\sC) \to \Sym(\sC)$, where $\mathrm{free}_{\CAlg}$ is the free commutative algebra functor on the symmetric monoidal category $\Sym(\sC)$ and $\mathrm{forget}_{\CAlg}$ is its right adjoint.
    Note that by \cite[Example 3.2.4.4 and Proposition 3.2.4.7]{HA}, $S$ is canonically a symmetric monoidal functor if we equip $\sC$ with the cocartesian symmetric monoidal structure.
    
    Write $\mathrm{can} \colon \Sym^\times(\sC) \to \Sym(\sC)$ for the unique symmetric monoidal functor for which the composite $\sC \hookrightarrow \Sym^\times(\sC) \to \Sym(\sC)$ is the standard inclusion $\sC \hookrightarrow \Sym(\sC)$.
    Consider the diagram 
    \[
    \begin{tikzcd}
        \Sym^\times(\sC) \ar[rr, "\mathrm{can}"] \ar[dr, "p"'] & & \Sym(\sC). \\
        & \sC \ar[ur, "S"'] &
    \end{tikzcd}
    \]
    To finish the proof, it suffices to show that $S$ is the left Kan extension of $\mathrm{can}$ along $p$.
    Indeed, since post-composing with a colimit preserving functor preserves Kan extensions, this would imply that we have a natural equivalence $p_!(\mathrm{can} \circ F) \simeq F \circ S$ natural in $F \in \SymFunL(\sC, \sD)$, which is what we wanted to prove.
    We will start by writing down a natural transformation $\alpha \colon p_! (\mathrm{can}) \to S$, which we do by providing a map $\mathrm{can} \to S \circ p$.
    For this, it suffices to give the map after restricting along the inclusion $\sC \to \Sym^\times(\sC)$, since all functors in the diagram are symmetric monoidal and $\Sym^\times(\sC)$ is the free symmetric monoidal category on $\sC$.
    Upon restricting to $\sC$, $\mathrm{can}$ becomes the inclusion $\sC \to \Sym(\sC)$ and $S \circ p$ becomes $S$, and there is an obvious choice of natural transformation $X \to S(x) = \coprod_{n \geq 0} x^{\otimes n}_{h\Sigma_n}$ given by the inclusion of the $n=1$ component.
    
    We will proceed to show $\alpha$ is an equivalence.
    Let $x \in \sC$.  
    We can compute $p_!(\mathrm{can})(X)$ as the colimit of the functor
    \[
    \begin{tikzcd}
    \Sym^\times(\sC)_{/x} \ar[r] & \Sym^\times(\sC) \ar[r, "\mathrm{can}"] & \Sym(\sC).
    \end{tikzcd}
    \]
    There is a right cofinal functor $\Fin^{\cong} \to \Sym^\times(\sC)_{/X}$ that sends the set $\{1, \ldots, n\}$ to the pair 
    \[
    ((x, \ldots, x),\quad x \sqcup \cdots \sqcup x \xrightarrow{\nabla} x),
    \] 
    where there are $n$ copies of $x$ and $\nabla$ denotes the fold map.
    It follows that the colimit of the above functor is given by $\coprod_{n \geq 0} x^{\otimes n}_{h\Sigma_n}$, which is equivalent to $S(x)$, and one can check that $\alpha$ indeed gives an equivalence.
\end{proof}

\subsubsection{Symmetric sequences of spectra}

In the previous section, we constructed for each presentable category $\sC$ a composition product monoidal structure on $\SymFunL(\sC, \sC)$.
As we explained in \cref{ex: funsymseq-is-symseq-spectra}, this category has a simple description in the case where $\sC$ is the category of spectra: there is an equivalence $\SymFunL(\Sp, \Sp) \simeq \SymSeq(\Sp)$, where the latter denotes the category $\Fun(\Fin^{\cong}, \Sp)$ of ordinary symmetric sequences in spectra.
Brantner already gave a construction in \cite[\S 4.1.2]{brantnerThesis} of a composition product monoidal structure on $\SymSeq(\Sp)$, making use of an idea going back to the classical category theory literature (cf.\ \cite{trimble}).
The goal of this section is to prove that these two monoidal structures agree.
\begin{proposition} \label{prop: equivalence-of-composition-products-spectra}
    The composition product monoidal structure on $\SymSeq(\Sp)$ constructed in the previous section is equivalent to the monoidal structure constructed by Brantner.
\end{proposition}

\begin{remark}
    Write $\Spc$ for the category of spaces. 
    We also have $\SymFunL(\Spc, \Spc) \simeq \SymSeq(\Spc)$, and Brantner's construction gives a composition product monoidal structure on the latter category as well. 
    Moreover, the above theorem holds if we replace $\Sp$ by $\Spc$.
    The proof is analogous.
    Brantner and Heuts prove in a forthcoming paper \cite{BrantnerHeuts} that the category $\Alg(\SymSeq(\Spc))$ is equivalent to the category of (singly-coloured) $\infty$-operads in the sense of Lurie \cite[Chapter 2]{HA}.
\end{remark}

\begin{remark}
    A different construction of a composition product monoidal structure on $\SymSeq(\Sp)$ was given by Haugseng in \cite{HaugsengSymmetric}.
    We expect this monoidal structure to also be equivalent to the one we consider in this section, but we don't have a proof.
\end{remark}

Let us start by recalling Brantner's construction.
The category $\SymSeq(\Sp)$ admits a Day convolution monoidal structure with respect to the smash product on $\Sp$.
As such, it is the free stable presentably symmetric monoidal category on one object.
(Presentably symmetric monoidal here means that the underlying category is presentable, and that the monoidal product preserves small colimits in both variables separately.)
So for any stable presentably symmetric monoidal category $\sC$, evaluation in the sphere spectrum induces an equivalence
\[
\begin{tikzcd}
    \Fun^{\mathrm{L}, \otimes}(\SymSeq(\Sp), \sC) \ar[r, "\sim"] & \sC.
\end{tikzcd}
\]
The left-hand side denotes the category of symmetric monoidal colimit preserving functors.
Plugging in $\sC = \SymSeq(\Sp)$, we obtain an equivalence 
\[
\End^{\mathrm{L}, \otimes}(\SymSeq(\Sp)) \simeq \SymSeq(\Sp).
\]
Functor composition endows the left-hand side with a monoidal structure.
Brantner then takes the opposite of this monoidal structure as his definition of the composition product on $\SymSeq(\Sp)$.
It is easily verified that for positive symmetric sequences $\sO$ and $\sP$, this monoidal product is given by the formula
\[
(\sO \circ \sP)_I = \bigoplus_{E \in \mathrm{Part}(I)} \sO_{E} \otimes \bigotimes_{J \in E} \sP_J,
\]
which justifies calling it the composition product.

Let $\presl_{\mathrm{St}}$ denote the full subcategory of the $2$-category $\presl$ on the stable presentable categories.
This $2$-category admits a symmetric monoidal structure given by the tensor product of presentable categories (\cref{proposition:PrL-is-symmetric-monoidal-2-category}).
The unit object is $\Sp$.
Our proof of \cref{prop: equivalence-of-composition-products-spectra} will rely on dualizability in this symmetric monoidal $2$-category.
As a preparation for the proof, we need two lemmas.

\begin{lemma} \label{lem: symseq-sp-is-sym-sp}
    There is an equivalence $\SymSeq(\Sp) \simeq \Sym(\Sp)$ of symmetric monoidal categories, where $\SymSeq(\Sp)$ is equipped with Day convolution.
\end{lemma}
\begin{proof}
    As we explained above, $\SymSeq(\Sp)$ is the free stable presentably symmetric monoidal category on one object.
    It is easily seen that $\Sym(\Sp)$ satisfies the same universal property.
    This proves the lemma.
\end{proof}

\begin{lemma} \label{lem: Sym-dualizable}
    Let $\sC \in \presl_{\mathrm{St}}$ be a dualizable category.
    Then $\Sym(\sC)$ is again dualizable.
\end{lemma} 
\begin{proof}
    By \cref{remark:dualizable-underlying-1-category} and \cite[Proposition D.7.3.1]{SAG}, a stable presentable category is dualizable if and only if it is the retract in $\presl_{\mathrm{St}}$ of a compactly generated category.
    But if $\sC$ is a retract of $\sD$, then $\Sym(\sC)$ is a retract of $\Sym(\sD)$.
    It therefore suffices to show that $\Sym(\sC)$ is compactly generated whenever $\sC$ is.

    Given small finite cocomplete categories $\sC_0$ and $\sD_0$, it follows from \cite[§4.8.1]{HA} that there is a small finite cocomplete category $\sC_0 \otimes_\omega \sD_0$ together with a functor $\sC_0 \times \sD_0 \to \sC_0 \otimes_\omega \sD_0$ that preserves colimits in both variables such that the restriction functor
    \[
    \begin{tikzcd}
    \Fun^{\mathrm{fin}}(\sC_0 \otimes_\omega \sD_0, \sE_0) \ar[r, hook] & \Fun(\sC_0 \times \sD_0, \sE_0)
    \end{tikzcd}
    \]
    is fully faithful with essential image those functors $\sC_0 \times \sD_0 \to \sE_0$ that preserve finite colimits in both variables.
    Here the left hand side denotes the category of functors that preserve finite colimits, and $\sE_0$ is a small finite cocomplete category.
    A comparison of universal properties shows that $\Ind(\sC_0) \otimes \Ind(\sD_0) \simeq \Ind(\sC_0 \otimes_\omega \sD_0)$,
    so that by induction a finite tensor product of compactly generated categories is compactly generated.
    It follows that if $\sC$ is compactly generated, $\Sym(\sC)$ is computed as the colimit in $\presl$ of a diagram of compactly generated categories, such that all the functors in the diagram are equivalences and hence preserve compact objects.
    It is a standard fact that the colimit of such a diagram is again compactly generated, see for instance \cite[Proposition 2.7]{mathewGalois}.
\end{proof}

\begin{proof}[Proof of \cref{prop: equivalence-of-composition-products-spectra}]
    Let $\presl_{\mathrm{dual}}$ denote the full subcategory of the $2$-category $\presl_{\mathrm{St}}$ on the dualizable categories.
    This $2$-category inherits a symmetric monoidal structure, and the functor
    \[
    \begin{tikzcd}
    \mathbf{D} \colon (\presl_{\mathrm{dual}})^{\op} \ar[r, "\sim"] & \presl_{\mathrm{dual}}
    \end{tikzcd}
    \]
    that sends an object to its dual is a strong symmetric monoidal equivalence by \cref{proposition:Rigid-2-category-self-dual}.
    By \cref{lem: Sym-dualizable}, the free-forgetful adjunction restricts to an adjunction
    \[
    \begin{tikzcd}[sep = huge]
    \presl_{\mathrm{dual}} \ar[r, shift left, "\Sym"] & \CAlg(\presl_{\mathrm{dual}}). \ar[l, shift left, "\forget"]
    \end{tikzcd}
    \]
    Since $\Sym$ is also the underlying functor of the cofree commutative coalgebra functor (\cref{prop: sym-is-cofree}), the same holds for the cofree-forgetful adjunction.
    We have the following commutative diagram of right adjoints between $2$-categories:
    \[
    \begin{tikzcd}[sep = large]
        \CAlg(\presl_{\mathrm{dual}}) \ar[r, "\mathbf{D}", "\sim"'] \ar[d, "\forget"'] & \CAlg((\presl_{\mathrm{dual}})^{\op}) \ar[d, "\forget"] \\
        \presl_{\mathrm{dual}} \ar[r, "\mathbf{D}", "\sim"'] & (\presl_{\mathrm{dual}})^{\op}.
    \end{tikzcd}
    \]
    Passing to left adjoints, we obtain an equivalence $\mathbf{D}(\Sym(\sC)) \simeq \Sym(\mathbf{D}(\sC))$ of cocommutative coalgebras in $\presl_{\mathrm{St}}$.
     In particular, $\mathbf{D}(\Sym(\Sp)) = \Sym(\Sp)$.
     The top functor in this diagram therefore induces an equivalence of monoidal endomorphism categories
     \[
     \begin{tikzcd}
     \End_{\CAlg(\presl)^{\op}}(\Sym(\Sp)) \ar[r, "\sim"] & \End_{\coCAlg(\presl)}(\Sym(\Sp)).
     \end{tikzcd}
     \]
     The left-hand side is Brantner's monoidal category of symmetric sequences by \cref{lem: symseq-sp-is-sym-sp}, whereas the right-hand side is our monoidal category of symmetric sequences.
\end{proof}

\subsection{A lax structure on the Goodwillie derivatives}
\label{sec: lax-structure-Goodwillie}

Differentiable categories and reduced finitary functors can be organized into a $2$-category.

\begin{definition}
    The $2$-category $\diff$ is the locally full subcategory of $\Cat$ with differentiable categories as objects and reduced finitary functors as morphisms.
\end{definition}

In this section we prove the first main result of this paper:

\begin{theorem} \label{thm: lax-derivatives}
    There exists a lax $2$-functor
    \[
    \begin{tikzcd}
    \partial_* \colon \diff \ar[r] & \pressymstc
    \end{tikzcd}
    \]
    which sends a differentiable category $\sC$ to its stabilization $\Sp(\sC)$, and on mapping categories is given by the Goodwillie derivatives functor.
\end{theorem}

This theorem in particular implies that for any differentiable category $\sC$ the derivatives functor
\[
\begin{tikzcd}
\partial_* \colon \Fun^{\ast, \omega}(\sC, \sC) \ar[r] & \SymFunL_{\geq 1}(\Sp(\sC), \Sp(\sC))
\end{tikzcd}
\]
can be refined to a lax monoidal functor.
To give the strategy for the proof of \cref{thm: lax-derivatives}, let us explain how this lax monoidal structure can be constructed.

The derivatives functor factors through $\SymFun_{\geq 1}^{\ast, \omega}(\sC, \sC)$ as symmetric cross-effects followed by multilinearization.
A first idea might be to show that both these functors admit lax monoidal structures.
This fails on the basic point that the category $\SymFun_{\geq 1}^{\ast, \omega}(\sC, \sC)$ is not monoidal in any relevant sense: even though it is possible to write down the formula for the composition product of two functors in this category, this product will not be associative since the functors involved do not commute with the relevant colimits.

We circumvent this problem by formally adjoining finite colimits to the category $\sC$, thus forcing the composition product on $\SymFun_{\geq 1}^{\ast, \omega}(\sC, \sC)$ to become associative.
More precisely, by embedding $\sC$ in an appropriate category of presheaves $\Pfin(\sC)$ we exhibit $\SymFun_{\geq 1}^{\ast, \omega}(\sC, \sC)$ as a reflective localization of the category $\SymFunL_{\geq 1}(\Pfin(\sC), \Pfin(\sC))$ of \emph{formal symmetric sequences}.
The factorization of the derivatives functor can then be extended to the following commutative diagram:
\[
\begin{tikzcd}[sep = large]
    \Fun^{\ast, \omega}(\sC, \sC) \ar[r, "\cross"] \ar[dr, "\widetilde{\cross}"'] & \SymFun^{\ast, \omega}_{\geq 1}(\sC, \sC) \ar[r, "\mlin"] \ar[d, hook, shift left, "\iota"] & \SymFunL_{\geq 1}(\sC, \sC) \\
    & \SymFunL_{\geq 1}(\Pfin(\sC), \Pfin(\sC)) \ar[u, shift left, "L"] \ar[ur, "\widetilde{\mlin}"']
\end{tikzcd}
\]
Here $L$ is left adjoint to the fully faithful embedding $\iota$.
The derivatives functor now factors through the monoidal category of formal symmetric sequences as the composite $\widetilde{\mlin} \circ \widetilde{\cross}$.
One can then proceed by showing $\widetilde{\cross}$ is lax monoidal and $\widetilde{\mlin}$ is strong monoidal on the image of $\widetilde{\cross}$.
The lax monoidality of $\widetilde{\cross}$ follows since it has a left adjoint which is oplax monoidal.
For $\widetilde{\mlin}$ one can use that it is a localization, so that it suffices to show it is compatible with the monoidal structure.
This compatibility can be deduced from the lemma of Arone and Ching which says that $P_1(FG) \simeq P_1(P_1(F) P_1(G))$.

We will now carry out a 2-categorical version of the argument sketched here.

\subsubsection{Formal symmetric sequences}

The definition of formal symmetric sequences relies on the following proposition, which allows us to formally adjoin non-empty finite colimits to a category.

\begin{proposition} \label{prop: pfin}
    Let $\sC$ be a pointed category that admits filtered colimits.
    Then there exists a category $\Pfin(\sC)$ together with a functor $j \colon \sC \to \Pfin(\sC)$ such that:
    \begin{enumerate}[label={\upshape(\arabic*)}]
        \item The category $\Pfin(\sC)$ is cocomplete and pointed,
        \item The functor $j$ is fully faithful, reduced, finitary, and preserves all limits that exist in $\cC$,
        \item For any cocomplete pointed category $\sD$, restriction along $j$ induces an equivalence
        \[
        \begin{tikzcd}
            j^* \colon \FunL(\Pfin(\sC), \sD) \ar[r, "\sim"] & \Fun^{\ast, \omega}(\sC, \sD),
        \end{tikzcd}
        \]
        with inverse given by left Kan extension along $j$.
    \end{enumerate}
    If $\sC$ is additionally accessible, then $\Pfin(\sC)$ is presentable.
\end{proposition}
\begin{proof}
In the notation of \cite[Proposition 5.3.6.2]{HTT}, we define $\Pfin(\sC)$ to be $\mathscr{P}_{\mathscr{R}}^{\mathscr{K}}(\sC)$, where $\mathscr{R}$ is the collection of all small filtered categories together with the empty category and $\mathscr{K}$ is the collection of all small categories.
It follows immediately from loc.\ cit.\ that $\Pfin(\sC)$ is cocomplete and comes equipped with a fully faithful functor $j \colon \sC \to \Pfin(\sC)$ that preserves the initial object and filtered colimits, and satisfies the universal property stated as (3) in the proposition.

In the proof of \cite[Proposition 5.3.6.2]{HTT}, the category $\Pfin(\sC)$ is constructed as a full subcategory of the presheaf category $\Fun(\sC^\op,\widehat{\Spc})$ and $j$ is given by the Yoneda embedding, which lands in $\Pfin(\sC)$.
Here $\widehat{\Spc}$ denotes the (very large) category of large spaces.
Since the Yoneda embedding $\sC \to \Fun(\sC^\op,\widehat{\Spc})$ preserves all limits that exist in $\sC$, so does $j \colon \sC \to \Pfin(\sC)$.
As $\sC$ is pointed and $j$ also preserves the initial object, $\Pfin(\sC)$ has to be pointed as well and $j$ is reduced. 

For the final claim, suppose that $\sC$ is $\kappa$-accessible so that we can write $\sC \simeq \Ind_{\kappa}(\sC_{\kappa})$, where $\sC_{\kappa}$ is the full subcategory of $\sC$ spanned by the $\kappa$-compact objects.
By comparing universal properties and using \cref{lemma:k-small-filtered-colimits} below, we obtain an equivalence $\Pfin(\sC) \simeq \Ind_{\kappa}(\mathscr{P}_{\mathscr{R}'}^{\mathscr{K}'}(\sC_{\kappa}))$,
where $\mathscr{R}'$ denotes the collection of all filtered $\kappa$-small categories together with the empty category and $\mathscr{K}'$ denotes the collection of all $\kappa$-small categories.
It follows from the construction of $\mathscr{P}_{\mathscr{R}'}^{\mathscr{K}'}(\sC_{\kappa})$ in loc.\ cit.\ that this category is small.
Hence, $\Pfin(\sC)$ is $\kappa$-accessible as well.
Since it has all small colimits by construction, it is presentable.
\end{proof}

We used the following lemma in the proof of the previous proposition.

\begin{lemma}\label{lemma:k-small-filtered-colimits}
    Let $F \colon \sC \to \sD$ be a functor between categories that admit filtered colimits and let $\kappa$ be a regular cardinal. Then $F$ preserves filtered colimits if and only if $F$ preserves $\kappa$-filtered colimits and filtered $\kappa$-small colimits. 
\end{lemma}

\begin{proof}
    This follows by the same argument as the final part of the proof of \cite[Proposition 5.4.2.9]{HTT}: there it is shown that any filtered colimit can be decomposed as a $\kappa$-filtered colimit of $\kappa$-small filtered colimits.
\end{proof}

\begin{remark}
\label{rmk: pfin-as-adjoint}
    Let $\twocatLargeCat$ be the 2-category of possibly large categories and $\sE \subset \twocatLargeCat$ the locally full subcategory with objects the pointed presentable categories and morphisms the reduced finitary functors.
    The construction $\sC \mapsto \Pfin(\sC)$ provides a left adjoint to the inclusion $\presl_\ast \hookrightarrow \sE$.
    Here $\presl_\ast$ denotes the full subcategory of $\presl$ on the pointed presentable categories.
\end{remark}

\begin{definition}[Formal symmetric sequences]
Let $\sC$ and $\sD$ be pointed presentable categories.
A \emph{formal symmetric sequence} from $\sC$ to $\sD$ is a symmetric sequence from $\Pfin(\sC)$ to $\Pfin(\sD)$.
\end{definition}

\begin{proposition} \label{prop: symseq-to-formalsymseq}
    Let $\sC$ and $\sD$ be pointed presentable categories.
    There is an equivalence
    \[
    \begin{tikzcd}
        \SymFun^{\ast, \omega}(\sC, \Pfin(\sD)) \ar[r, "\sim"] & \SymFunL(\Pfin(\sC), \Pfin(\sD)),
    \end{tikzcd}
    \]
    given by left Kan extension along $\Sym^\times(j) \colon \Sym^\times(\sC) \to \Sym^\times(\Pfin(\sC))$.
    This restricts to an equivalence between the full subcategories on objects concentrated in positive arity.
\end{proposition}
\begin{proof}
    Let $n \geq 0$.
    By part (3) of \cref{prop: pfin}, left Kan extension along $j^{\times n}$ gives a fully faithful embedding
    \[
    \begin{tikzcd}
        \Fun^{\ast, \omega}(\sC^{\times n}, \Pfin(\sD)) \ar[r] & \Fun(\Pfin(\sC)^{\times n}, \Pfin(\sD)),
    \end{tikzcd}
    \]
    with essential image the functors that preserve colimits in each variable separately.
    This embedding is equivariant for the obvious $\Sigma_n$-action on both sides and the result follows by passing to $\Sigma_n$-fixed points and taking the product over all $n \geq 0$ (or $n \geq 1$ for the second statement). 
\end{proof}

Recall from \cref{definition:twofold-Segal-space} that a $2$-precategory is a Segal object $\sX$ in $\Cat$ such that $\sX_0$ is a space. 
Given a $2$-category $\sX$ together with a map of spaces $\phi \colon S \to \sX_0$, we show in \cref{construction:Modifying-space-of-objects} how to construct a $2$-precategory $\phi^*\sX$ that has $S$ as its space of objects and mapping categories given by
\[
\phi^*\sX(x, y) = \sX(\phi(x), \phi(y)).
\]
We can organize all categories of formal symmetric sequences between differentiable categories in a $2$-precategory.

\begin{definition}\label{definition:pre-category-Diff-FSym}
Let $\mathscr{P}\mathrm{r}_0$ denote the space of presentable categories, i.e.\ the space of objects of both $\presl$ and $\pressym$.
Let $\diff_0$ be the space of objects of $\diff$.
Consider the morphism $\alpha \colon \diff_0 \to \mathscr{P}\mathrm{r}_0$ sending a differentiable category $\sC$ to $\Pfin(\sC)$.
We let $\fdiff$ be the $2$-precategory $\alpha^*\pressym$,
so that $\fdiff$ has differentiable categories as its objects and mapping categories given by
\[
\fdiff(\sC, \sD) = \SymFunL(\Pfin(\sC), \Pfin(\sD)).
\]
We write $\fdiffc$ for the locally full subcategory of $\fdiff$ spanned by the  formal symmetric sequences.
\end{definition}

Suppose $\sD$ is a differentiable category.
As the functor $j \colon \sD \hookrightarrow \Pfin(\sD)$ preserves limits and filtered colimits, the adjoint functor theorem tells us that $j$ admits a left adjoint $h$:
\begin{equation}\label{equation:Yoneda-and-left-adjoint}
\begin{tikzcd}[sep = huge]
    \sD \ar[r, shift right, hook, "j"'] & \Pfin(\sD) \ar[l, shift right, "h"'],
\end{tikzcd}
\end{equation}
so that $\sD$ is a reflective localization of $\Pfin(\sD)$. 
The functor $h$ sends a formal colimit in $\Pfin(\sD)$ to the actual colimit in $\sD$.
Postcomposing with this adjunction yields an adjunction between $\SymFun^{\ast, \omega}_{\geq 1}(\sC, \sD)$ and $\SymFun^{\ast, \omega}_{\geq 1}(\sC, \Pfin(\sD))$.
The latter is equivalent to $\SymFunL_{\geq 1}(\Pfin(\sC), \Pfin(\sD))$ by \cref{prop: symseq-to-formalsymseq}, so that we obtain an adjunction (left adjoint on top)
\begin{equation} \label{eq: adj-symfun-fsymfun}
    \begin{tikzcd}[sep = huge]
        \SymFun^{\ast, \omega}_{\geq 1}(\sC, \sD) \ar[r, shift right, hook, "\iota"'] & \SymFunL_{\geq 1}(\Pfin(\sC), \Pfin(\sD)). \ar[l, shift right, "L"']
    \end{tikzcd}
\end{equation}

We can use this adjunction to define the composition product of two functors in $\SymFun^{\ast, \omega}_{\geq 1}(\sC, \sD)$.
\begin{definition} \label{def: symmetric-composition-product}
    Let $F \in \SymFun^{\ast, \omega}_{\geq 1}(\sD, \sE)$ and $G \in \SymFun^{\ast, \omega}_{\geq 1}(\sC, \sD)$.
    We define their \emph{composition product} as
    \[
    F \circ G \coloneqq L(\iota F \circ \iota G),
    \]
    where the composition product on the right is taken in $\fdiff$.
\end{definition}

\begin{remark}
The composition product defined here is given by the same formula as the composition product in symmetric sequences.
However, it fails to be associative since the functors involved don't commute with the necessary colimits, as we already mentioned.
\end{remark}

\subsubsection{Cross-effects as a lax functor}

We can lift the cross-effects functor to the category of formal symmetric sequences by forming the composite
\[
\begin{tikzcd}
    \widetilde{\cross} \colon \Fun^{\ast, \omega}(\sC, \sD) \ar[r, "\cross"] & \SymFun^{\ast, \omega}_{\geq 1}(\sC, \sD) \ar[r, "\iota"] & \SymFunL_{\geq 1}(\Pfin(\sC), \Pfin(\sD)).
\end{tikzcd}
\]
The goal of this section is to construct a lax functor $\widetilde{\cross} \colon \diff \to \fdiffc$, which is given by $\widetilde{\cross}$ on mapping categories.

Let $\sC$ and $\sD$ be differentiable categories.
Recall that the cross-effects functor is given as the composite of two right adjoints:
\[
\begin{tikzcd}
    \cross \colon \Fun^{\ast, \omega}(\sC, \sD) \ar[r, "p^*"] & \Fun^{\omega}_{\geq 1}(\Sym^{\times}(\sC), \sD) \ar[r, "\mathrm{red}"] & \SymFun^{\ast, \omega}_{\geq 1}(\sC, \sD),
\end{tikzcd}
\]
where the functor $p \colon \Sym^{\times}(\sC) \to \sC$ sends a sequence of objects to their coproduct.
The left adjoint of $\cross$ is given by left Kan extension along $p$.
We will use the notation $q \colon \Sym^{\times}(\Pfin(\sC)) \to \Pfin(\sC)$ for the functor sending a sequence of objects to their coproduct in $\Pfin(\sC)$.

\begin{proposition}\label{prop: commuting-p-and-q}
Let $\sC$ and $\sD$ be differentiable categories. 
Then the square
\[
\begin{tikzcd}[sep = large]
    \Fun^{\ast, \omega}(\sC, \sD) & \SymFun^{\ast, \omega}_{\geq 1}(\sC, \sD) \ar[l, "p_!"'] \\
    \Fun^{\ast, \omega}(\Pfin(\sC), \Pfin(\sD)) \ar[u, "h \circ - \circ j"] & \SymFunL_{\geq 1}(\Pfin(\sC), \Pfin(\sD)) \ar[l, "q_!"] \ar[u, "L"'].
\end{tikzcd}
\]
commutes, where $h$ denotes the left adjoint to $j$ and $L$ is the left adjoint from \cref{eq: adj-symfun-fsymfun}.
In particular, the composition of $h\circ - \circ j$ with $q_!$ is left adjoint to the functor $\widetilde{\cross}$.
\end{proposition}
\begin{proof}
By \cref{prop: symseq-to-formalsymseq}, there is an equivalence 
\[
\Fun^{\ast, \omega}(\Sym^\times(\sC), \Pfin(\sD)) \simeq \SymFunL(\Pfin(\sC), \Pfin(\sD))
\]
given by left Kan extension along $\Sym^\times(j)$.
Under this equivalence, $L$ becomes postcomposition by $h$, and $q_!$ becomes left Kan extension along $q \circ \Sym^\times(j)$.
Note that the functor $j^* \colon \Fun^{\ast, \omega}(\Pfin(\sC), \Pfin(\sD)) \to \Fun^{\ast, \omega}(\sC, \Pfin(\sD))$ is equivalent to left Kan extension along $h$, as restriction along a right adjoint is the same as left Kan extension along its left adjoint (see \cite[Lemma 5.2.6.6]{HTT}).
Next observe that the composite
\[
\begin{tikzcd}[sep = large]
\Sym^\times(\sC) \ar[r, "\Sym^\times(j)"] & \Sym^\times(\Pfin(\sC)) \ar[r, "q"] & \Pfin(\sC) \ar[r, "h"] & \sC
\end{tikzcd}
\]
is equivalent to $p$, as it is symmetric monoidal for the cocartesian monoidal structure on $\sC$ and restricts to the identity along $\sC \hookrightarrow \Sym^\times(\sC)$.
Combining all of this, we have reduced the problem to showing that $p_! \circ h_*$ is equivalent to $h_* \circ p_!$ as functors 
\[
\begin{tikzcd}
\Fun^{\ast, \omega}(\Sym^\times(\sC), \Pfin(\sD)) \ar[r] & \Fun^{\ast, \omega}(\sC, \sD).
\end{tikzcd}
\]
This holds since $h$ preserves colimits.
\end{proof}

Suppose for a moment that $\sC = \sD$ in the statement of this proposition, so that all categories in the square except the top right one admit a composition monoidal structure.
By \cref{prop: functorization} and \cref{prop: functorization-is-kan-extension}, the functor $q_!$ admits a strong monoidal structure.
By \cref{theorem:Twisting-by-adjunctions}, the functor $h \circ - \circ j$ admits an oplax monoidal structure, with the oplax structure map given by the unit of the adjunction between $j$ and $h$:
\[
\begin{tikzcd}
    h \circ F \circ G \circ j \ar[r, "\eta"] & (h \circ F \circ j) \circ (h \circ G \circ j).
\end{tikzcd}
\]
Since the square in \cref{prop: commuting-p-and-q} commutes, we have therefore produced an oplax monoidal structure on the composite $p_! \circ L$.
This then gives a lax monoidal structure on $\widetilde{\cross}$, since this functor is right adjoint to $p_! \circ L$.

We will now carry out this argument on the level of $2$-categories to produce the lax $2$-functor $\widetilde{\cross} \colon \diff \to \fdiffc$.

\begin{lemma}
Let $\alpha \colon \diff_0 \to \diff_0$ be the map sending $\sC$ to $\Pfin(\sC)$.
Then there is an oplax $2$-functor $H \colon \alpha^*\diff \to \diff$ that is the identity on objects and is given by
\[
\begin{tikzcd}
h \circ - \circ j \colon \Fun^{\ast, \omega}(\Pfin(\sC), \Pfin(\sD)) \ar[r] & \Fun^{\ast, \omega}(\sC, \sD)
\end{tikzcd}
\]
on mapping categories.
\end{lemma}
\begin{proof}
Recall that $(\diff_1)^{\cong}$ denotes the space of arrows in $\diff$.
Let $\beta \colon \diff_0 \to (\diff_1)^{\cong}$ be the map that sends $\sC$ to the functor $h \colon \Pfin(\sC) \to \sC$.
This map comes from the unit of the adjunction described in \cref{rmk: pfin-as-adjoint}.
Observe that the composite $s \circ \beta \colon \diff_0 \to \diff_0$ of $\beta$ with the projection on the source is homotopic to the morphism $\alpha$.
The composite $t \circ \beta \colon \diff_0 \to \diff_0$ with the target projection is homotopic to the identity.
Since the functors $h$ are all left adjoints, we obtain by \cref{theorem:Twisting-in-families} an oplax functor from $(s \circ \beta)^*\diff \simeq \alpha^*\diff$ to $(t \circ \beta)^*\diff \simeq \diff$ satisfying the required properties.
\end{proof}

\begin{proposition}
\label{prop: lax-cross-effects}
    There exists a lax $2$-functor $\widetilde{\cross} \colon \diff \to \fdiffc$ that is the identity on objects and is given by the composite
    \[
    \begin{tikzcd}
        \widetilde{\cross} \colon \Fun^{\ast, \omega}(\sC, \sD) \ar[r, "\cross"] & \SymFun^{\ast, \omega}_{\geq 1}(\sC, \sD) \ar[r, "\iota"] & \SymFunL_{\geq 1}(\Pfin(\sC), \Pfin(\sD))
    \end{tikzcd}
    \]
    on mapping categories.
\end{proposition}
\begin{proof}
    Let $\pressymdiffc \subset \pressymc$ be the full subcategory spanned by the differentiable categories.
    Let $\Lambda \colon \pressymdiffc \to \diff$ be the restriction of the $2$-functor constructed in \cref{prop: functorization} (note that this functor indeed lands in $\diff$ by \cref{rem: Lambda-F-reduced-finitary}).
    Since this functor is the identity on objects, we can pull back along the map $\alpha \colon \diff_0 \to \diff_0$ that sends $\sC$ to $\Pfin(\sC)$ to obtain a functor
    \[
    \begin{tikzcd}
    \alpha^*\Lambda \colon \alpha^*\pressymdiffc = \fdiffc \ar[r] & \alpha^* \diff.
    \end{tikzcd}
    \]
    Composing $\alpha^* \Lambda$ with the oplax $2$-functor $H \colon \alpha^* \diff \to \diff$ constructed in the previous lemma, we obtain an oplax $2$-functor $\fdiffc \to \diff$ which is the identity on objects and on mapping categories is given by
    \[
    \begin{tikzcd}
    \SymFunL_{\geq 1}(\Pfin(\sC), \Pfin(\sD)) \ar[r, "\Lambda"] & \Fun^{\ast, \omega}(\Pfin(\sC), \Pfin(\sD)) \ar[r, "h \circ - \circ j"] & \Fun^{\ast, \omega}(\sC, \sD).
    \end{tikzcd}
    \]
    By \cref{prop: functorization-is-kan-extension}, the functor $\Lambda$ in this composite is given by left Kan extension along the functor $q \colon \Sym^{\times}(\Pfin(\sC)) \to \Pfin(\sC)$ sending a sequence of objects to their coproduct in $\Pfin(\sC)$.
    By \cref{prop: commuting-p-and-q}, the composite of $h \circ - \circ j$ with $q_!$ is left adjoint to $\widetilde{\cross}$.
    By \cref{prop: local-right-adjoint-is-lax}, we can pass to the local right adjoint of $H \circ \alpha^*\Lambda$ to obtain the desired lax $2$-functor $\widetilde{\cross} \colon \diff \to \fdiffc$.
\end{proof}

\begin{remark}
\label{rem: lax-comparison-map-cr}
    Let $\sC$, $\sD$, and $\sE$ be differentiable categories and suppose that
    $F \colon \sD \to \sE$ and $G \colon \sC \to \sD$ are reduced finitary functors.
    Recall from \cref{def: symmetric-composition-product} that the composition product $\cross(F) \circ \cross(G)$ is defined as $L(\iota \cross(F) \circ \iota \cross(G))$.
    We obtain a natural comparison map $\cross(F) \circ \cross(G) \to \cross(FG)$ by applying $L$ to the lax comparison map $\widetilde{\cross}(F) \circ \widetilde{\cross}(G) \to \widetilde{\cross}(FG)$ coming from the previous proposition.
    Unraveling the definitions, we see that this map is given by the composite
    \[
    \begin{tikzcd}
        \cross(F) \circ \cross(G) \ar[r, "\eta"] & \cross p_! (\cross(F) \circ \cross(G)) \ar[r, "\cross(\delta)"] & \cross(p_!\cross(F) \circ p_!\cross(G)) \ar[r, "\cross(\epsilon \circ \epsilon)"] & \cross(FG).
    \end{tikzcd}
    \]
    The natural transformations $\eta$ and $\epsilon$ are the unit and counit of the adjunction between $p_!$ and $\cross$.
    The map $\delta$ is the comparison map coming from the oplax structure on the local left adjoint of $\widetilde{\cross}$.
    More precisely, given $A \in \SymFun^{\ast, \omega}(\sD, \sE)$ and $B \in \SymFun^{\ast, \omega}(\sC, \sD)$, the map $\delta \colon p_!(A \circ B) \to p_!A \circ p_!B$ is defined as the composite
    \[
    \begin{tikzcd}
        p_!(A \circ B) \simeq h q_!(A \circ B) j \ar[r, "\sim"] & h q_!A \circ q_!B j \ar[r, "\alpha"] & (h q_!A j) \circ  (h  q_! B j) \simeq p_!A \circ p_!B,
    \end{tikzcd}
    \]
    where $\alpha$ is the unit map for the adjunction between $j$ and $h$.
    Unraveling the definitions, we find that $\delta$ is given by the colimit assembly map from
    \[
      \coprod_{n \geq 1} \coprod_{m_1, \ldots, m_n \geq 1} A_n(B_{m_1}(x, \ldots, x), \ldots, B_{m_n}(x, \ldots, x))_{h(\Sigma_n \ltimes \Sigma_{m_1} \times \cdots \times \Sigma_{m_n})}
    \]
    to 
    \[
    \coprod_{n \geq 1} A_n \left( \coprod_{m \geq 1} B_1(x, \ldots, x)_{h \Sigma_m}, \ldots, \coprod_{m \geq 1} B_m(x, \ldots, x)_{h \Sigma_m} \right)_{h \Sigma_n}.
    \]
    In particular, if $A \in \SymFunL(\sD, \sE)$, then $\delta$ is an equivalence.
\end{remark}

\subsubsection{Multilinearization is compatible with composition}

The second part of the factorization of $\partial_*$ is the composite
\[
\begin{tikzcd}
\SymFunL_{\geq 1}(\Pfin(\sC), \Pfin(\sD)) \ar[r, "L"] & \SymFun^{\ast, \omega}_{\geq 1}(\sC, \sD) \ar[r, "\mlin"] & \SymFun^{\mathrm{mlin}, \omega}(\sC, \sD),
\end{tikzcd}
\]
which is a reflective localization.
Our goal is to construct a functor of $2$-categories that is given by this localization on mapping categories.
Let us assume for a moment that $\sC = \sD$ to sketch the argument.
In this case, we could ask for the localization to be compatible with the composition product monoidal structure on the source.
This would provide a monoidal structure on the target for which localization is strong monoidal.
It turns out $\mlin \circ L$ isn't compatible with the monoidal structure on the entire category $\SymFunL_{\geq 1}(\Pfin(\sC), \Pfin(\sC))$, but it does become compatible once we pass to a suitable monoidal subcategory.
This subcategory contains the image of our lift $\widetilde{\cross}$ of the cross-effects functor, so that we can compose to obtain a lax structure on $\partial_*$.
Note that this produces a monoidal structure on $\SymFun^{\mathrm{mlin}, \omega}(\sC, \sD) \simeq \SymFunL(\Sp(\sC), \Sp(\sD))$ that is a priori different from the composition monoidal structure we have already constructed.
We will show these two monoidal structures agree in the next section.

We now carry out a 2-categorical version of the argument sketched above.

\begin{definition}
\label{def: image-iota-category}
    Let $\sF$ be the smallest locally full subcategory of $\fdiffc$ for which each of the mapping categories $\sF(\sC, \sD)$ contains the image of 
    \[
    \begin{tikzcd}
    \iota \colon \SymFun^{\ast, \omega}_{\geq 1}(\sC, \sD) \ar[r] & \SymFunL_{\geq 1}(\Pfin(\sC), \Pfin(\sD)).
    \end{tikzcd}
    \]
\end{definition}

\begin{remark} \label{rem: composition-image-iota}
    Observe that $\sF(\sC, \sD)$ is the full subcategory of $\SymFunL_{\geq 1}(\Pfin(\sC), \Pfin(\sD))$ spanned by formal symmetric sequences $F \colon \Sym(\Pfin(\sC)) \to \Pfin(\sD)$ for which there exists a positive integer $n \geq 1$ and functors $F_k \in \SymFun^{\ast, \omega}_{\geq 1}(\sC_k, \sC_{k+1})$ for $1 \leq k < n$ where the $\sC_i$ are differentiable, $\sC_1 = \sC$, and $\sC_n = \sD$, such that
    \[
    F = \iota(F_n) \circ \cdots \circ \iota(F_1).
    \]
    The compositions are taken in the $2$-precategory $\fdiffc$.
\end{remark}

\begin{proposition} \label{prop: composition-product-compatible-with-multilinearization}
There exists a strong functor of $2$-precategories
\[
\begin{tikzcd}
\widetilde{\mlin} \colon \sF \ar[r] & \mlindiffc, 
\end{tikzcd}
\]
where $\mlindiffc$ is a 2-category with objects differentiable categories and mapping categories given by $\SymFun^{\mathrm{mlin}, \omega}_{\geq 1}(\sC, \sD)$.
The functor $\widetilde{\mlin}$ is the identity on objects and is given by
\[
\begin{tikzcd}
    \sF(\sC, \sD) \ar[r, "L"] & \SymFun^{\ast, \omega}_{\geq 1}(\sC, \sD) \ar[r, "\mlin"] & \SymFun^{\mathrm{mlin}, \omega}_{\geq 1}(\sC, \sD)
\end{tikzcd}
\]
on mapping categories.
\end{proposition}
\begin{proof}
   We will prove that the functors $\mlin L \colon \sF(\sC, \sD) \to \SymFun^{\mathrm{mlin}, \omega}_{\geq 1}(\sC, \sD)$ for $\sC$ and $\sD$ ranging over all differentiable categories form a compatible family of localizations in the sense of \cref{def: compatible family of localizations}.
    By \cref{proposition:local-localizations}, this will imply that there is a $2$-precategory $\sF'$ with objects differentiable categories and mapping categories given by $\SymFun^{\mathrm{mlin}, \omega}_{\geq 1}(\sC, \sD)$, together with a strong $2$-functor $\sF \to \sF'$ that is given by $\mlin L$ on mapping categories.
    We are then done, since we can define $\diff^{\mathrm{mlin}}$ to be the completion (\cref{remark:Completion-2-precategories}) of $\sF'$ and $\widetilde{\mlin}$ to be the composite $\sF \to \sF' \to \diff^{\mathrm{mlin}}$.
    
    Suppose $\sC$, $\sD$, and $\sE$ are differentiable categories and $F \in \sF(\sD, \sE)$ and $G \in \sF(\sC, \sD)$. 
    Let $\eta \colon F \longrightarrow \iota \mlin L F$ be the unit for the localization $\mlin L$ at $F$. 
    We will also write $\eta$ for the unit at $G$.
    We need to show that $F \circ \eta$ and $\eta \circ G$ are sent to equivalences by $\mlin L$.
    Note that $\eta$ factors as
    \[
    \begin{tikzcd}
    F \ar[r, "\alpha"] & \iota LF \ar[r, "\beta"] & \iota \mlin L F,
    \end{tikzcd}
    \]
    where $\alpha$ and $\beta$ are the unit maps for the localizations $L$ and $\mlin$ respectively.
    The unit map $\eta \colon G \to \iota \mlin L G$ factors in a similar way.
    We claim it suffices to show two things to finish the proof: first, that $\beta \circ \iota LG$ and $\iota LF \circ \beta$ are sent to equivalences by $\mlin L$ and second, that for any $M \in \sF(\sD, \sE)$ and $N \in \sF(\sC, \sD)$, the maps $\alpha \circ N \colon M \circ N \to \iota L M \circ N$ and $M \circ \alpha \colon M \circ N \to M \circ \iota L N$ are sent to equivalences by $\mlin L$.
    Indeed, by applying $\mlin L$ to the commutative diagram
    \[
    \begin{tikzcd}[sep = large]
        F \circ G \ar[r, "\alpha \circ G"] \ar[dr, "\eta \circ G"'] & \iota L F \circ G \ar[d, "\beta \circ G"] \ar[r, "\iota L F \circ \alpha"] & \iota L F \circ \iota L G \ar[d, "\beta \circ \iota L G"] \\
        & \iota \mlin L F \circ G \ar[r, "\iota \mlin L F \circ \alpha"] & \iota \mlin L F \circ \iota L G
    \end{tikzcd}
    \]
    it then follows by 2-out-of-3 that $\eta \circ G$ is sent to an equivalence, and similarly for $F \circ \eta$.

    We start by analyzing $\beta \circ \iota LG$ and $\iota LF \circ \beta$.
    These are sent by $\mlin L$ to the canonical maps
    \[
    \mlin(LF \circ LG) \to \mlin(\mlin LF \circ LG), \quad \mathrm{and} \quad \mlin(LF \circ LG) \to \mlin(LF \circ \mlin LG)
    \]
    respectively, both of which are equivalences by \cref{cor: multivariate Arone-Ching}.

    Before we consider $\alpha$, we make some remarks about the adjunction between $\iota$ and $L$.
    Recall that this adjunction is induced by post-composition with the adjunction
    \[
    \begin{tikzcd}[sep = huge]
        \sE \ar[r, shift right, hook, "j"'] & \Pfin(\sE) \ar[l, shift right, "h"'].
    \end{tikzcd}
    \]
    Post-composition with this adjunction also leads to an adjunction on the level of underlying functors:
    \[
  \begin{tikzcd}[sep = huge]
      \Fun^{\ast, \omega}(\coprod_{k \geq 1} \sD^{\times k}, \sE) \ar[r, shift right, hook, "\iota"'] & \Fun^{\ast, \omega}(\coprod_{k \geq 1} \sD^{\times k}, \Pfin(\sD)). \ar[l, shift right, "L"']
  \end{tikzcd}
    \]
    This adjunction is of course compatible with the one for symmetric functors, in the sense that if $F$ is a symmetric functor, the underlying functor of $\iota F$ will be $\iota$ applied to the underlying functor of $F$, and similarly for $L$.
    For the remainder of this proof, we will only work on the level of underlying functors. 
    This suffices, since equivalences are detected on underlying functors.
    So from now on, we will write $F$ and $G$ for the underlying functors of these symmetric sequences.
    
    With this out of the way, we return to the proof.
    Let $M \in \sF(\sD, \sE)$ and $N \in \sF(\sC, \sD)$.
    By \cref{rem: composition-image-iota} and the formula for the composition product, we find that $M$ and $N$ are in each arity finite coproducts of functors in the image of $\iota$, so that we can write $M_n = \coprod_i \iota M_n^i$.
    Then $LM_n \simeq \coprod_i M^i_n$, and the map $\alpha$ is given by the canonical comparison map
    \[
    \begin{tikzcd}
    \alpha \colon \coprod_i \iota M_n^i \ar[r] & \iota( \coprod_i M_n^i ).
    \end{tikzcd}
    \]
    Since composition of functors preserves colimits in the first variable and $L \circ \iota \simeq \id$, it follows from the formula for the composition product that $L(\alpha \circ N)$ is an equivalence.
    In particular, $\alpha \circ N$ is sent to an equivalence by $\mlin L$.
    
    Just like we did for $M$, we can write the arity $n$ component of $N$ as a finite coproduct $N_n = \coprod_{j} \iota N_n^j$.
    In arity $n$, the composition product $M \circ N$ is therefore a coproduct of functors of the form $\iota M_k^i(N_{n_1}^{j_1}, \ldots, N_{n_k}^{j_k})$ with $n_1 + \cdots + n_k = n$.
    The map $L(M \circ \alpha)$ is the coproduct over $i$ of the assembly maps
    \[
    \begin{tikzcd}
    \coprod_{j_1, \ldots, j_k} M_k^i(N_{n_1}^{j_1}, \ldots, N_{n_k}^{j_k}) \ar[r] & M_k^i(\coprod_{j_1} N_{n_1}^{j_1}, \ldots, \coprod_{j_k} N_{n_k}^{j_k}).
    \end{tikzcd}
    \]
   In order to show $M \circ \alpha$ is sent to an equivalence by $\mlin L$, it therefore suffices to prove the following statement:  given integers $k, n_1, \ldots, n_k \geq 1$ and $A \in \Fun^{\ast, \omega}(\sD^{\times k}, \sE)$, the composite functor
    \[
    \begin{tikzcd}
    \prod_{i = 1}^k \Fun^{\ast, \omega}(\sC^{\times n_i}, \sD) \ar[r, "A \circ -"] & \Fun^{\ast, \omega}(\sC^{\times n_1 + \cdots + n_k}, \sE) \ar[r, "\mlin"] & \Fun^{\mathrm{mlin}, \omega}(\sC^{\times n_1 + \cdots + n_k}, \sE)   
    \end{tikzcd}
    \]
    preserves finite coproducts in each variable.
    In fact, our proof will show it preserves all small colimits in each variable.
    Since $A$ preserves initial objects and filtered colimits in each variable and $\mlin$ preserves all colimits, it suffices to prove the functor preserves pushouts in each variable.
    Because we can work one variable at a time, we can reduce to the case $k=1$.
    So let $A \in \Fun^{\ast, \omega}(\sD, \sE)$ and let
    \[
    \begin{tikzcd}
        B_0 \ar[r] \ar[d] & B_1 \ar[d] \\
        B_2 \ar[r] & B_3
    \end{tikzcd}
    \]
    be a pushout square in $\Fun^{\ast, \omega}(\sC^{\times n}, \sD)$.
    Composing with $P_1 A$ yields a pullback square of functors since $P_1 A$ is excisive.
    Applying $\mlin$ to the result gives
    \[
    \begin{tikzcd}
        \mlin(P_1(A)B_0) \ar[r] \ar[d] & \mlin(P_1(A)B_1) \ar[d] \\
        \mlin(P_1(A)B_2) \ar[r] & \mlin(P_1(A)B_3),
    \end{tikzcd}
    \]
    which is a cartesian square of multilinear functors by left exactness of $\mlin$.
    But such a square is cocartesian as well, as the category $\Fun^{\mathrm{mlin}, \omega}(\sC^{\times n}, \sE)$ is stable.
    Since we have equivalences $\mlin(P_1(A)B_i) \simeq \mlin(AB_i)$ by \cref{cor: multivariate Arone-Ching}, this completes the proof.
\end{proof}

Recall that we defined the composition product of two functors in $\SymFun^{\ast, \omega}_{\geq 1}(\sC, \sD)$ in \cref{def: symmetric-composition-product}.
The preceding proposition then has the following consequence:

\begin{corollary} \label{cor: multilinearization-symmetric-composition}
    Let $F \in \SymFun^{\ast, \omega}_{\geq 1}(\sD, \sE)$ and $G \in \SymFun^{\ast, \omega}_{\geq 1}(\sC, \sD)$.
    Then there is a natural equivalence
    \[
    \mlin(F \circ G) \simeq \mlin(F) \circ \mlin(G),
    \]
    where the composition product on the right is taken in $\mlindiffc$.
\end{corollary}

\subsubsection{Equivalence of the two composition products}

In \cref{prop: composition-product-compatible-with-multilinearization}, we constructed the 2-category $\mlindiffc$ by localizing the mapping categories of $\fdiff_{\geq 1}$.
For stable categories $\sC$ and $\sD$, this localized mapping category is simply the category $\SymFunL(\sC,\sD)$ of symmetric sequences from $\sC$ to $\sD$.
In particular, the composition of 1-morphisms in $\mlindiffc$ gives a ``localized'' composition product for symmetric sequences between stable categories.
The goal of this section is to show that this localized composition product agrees with the composition product defined in \cref{def: composition product}.
More generally, we will show the following.

\begin{proposition}\label{proposition:equivalence-two-composition-products}
    There exists an equivalence of 2-categories
    \[\begin{tikzcd}
        \pressymstc \ar[r,"\sim"] & \mlindiffc
    \end{tikzcd}\]
    which sends a stable presentable category to itself and a which is the identity on mapping categories. In particular, for two symmetric sequences between stable categories, their composition product agrees with their localized composition product.
\end{proposition}

Let us start with the following observation.

\begin{lemma}\label{lemma:mlindiffstable-equals-mlindiff}
    Let $\mlindiffstc$ denote the full subcategory of $\mlindiffc$ spanned by the stable presentable categories. Then the inclusion
    \[\begin{tikzcd}
        \mlindiffstc \ar[r,hook] & \mlindiffc
    \end{tikzcd}\]
    is an equivalence of 2-categories.
\end{lemma}

\begin{proof}
    It suffices by \cref{lemma:fully-faithful-ess-surjective-is-equivalence} to show that this inclusion is essentially surjective.
    To this end, let $\sC$ be any differentiable category.
    In the stabilization adjunction
    \[
    \begin{tikzcd}[sep = large]
        \Sigma^\infty_\sC : \sC \ar[r, shift left, ""] & \Sp(\sC) : \Omega^\infty_\sC \ar[l, shift left, ""]
    \end{tikzcd}
    \]
    both functors are linear and preserve filtered colimits.
    In particular, we can view $\Sigma^\infty_\sC$ and $\Omega^\infty_\sC$ as 1-morphisms in $\mlindiffc$ that are concentrated in arity 1.
    We claim that these 1-morphisms are inverse to each other in $\mlindiffc$.
    To see this, it suffices to show that the unit and counit maps
    \[\begin{tikzcd}\eta \colon \id_\sC \ar[r] & \Omega^\infty \Sigma^\infty\end{tikzcd} \quad \text{and} \quad \begin{tikzcd} \varepsilon \colon \Sigma^\infty \Omega^\infty \ar[r] & \id_{\Sp(\sC)} \end{tikzcd}\]
    become equivalences after applying $P_1$. This was established in \cref{example:explicit-description-P1,example:derivatives-Sigma-Loops}.
    Since $\Sp(\sC)$ lies in the image of $\mlindiffstc \hookrightarrow \mlindiffc$, it follows that this inclusion is essentially surjective.
\end{proof}

\begin{construction}[The 2-precategory $\sF'_{\mathrm{St}}$]\label{construction:the-category-FSt}
    Note that in \cref{prop: composition-product-compatible-with-multilinearization}, we constructed a 2-precategory $\sF'$ by localizing the mapping categories of $\fdiff_{\geq 1}$ and then defined $\mlindiffc$ as the completion of $\sF'$.
    If we restrict to the full sub-2-precategory $\fdiff_{\mathrm{St},\geq 1}$ spanned by the stable presentable categories, then the localizations
    \[\begin{tikzcd}
        P_1L \colon \SymFunL_{\geq 1}(\Pfin(\sC), \Pfin(\sD)) \ar[r] & \SymFun^{\mathrm{mlin}, \omega}_{\geq 1}(\sC, \sD)
    \end{tikzcd}\]
    still form a weakly compatible family of localizations in the sense of \cref{remark:weakly-compatible-family-of-localizations}.
    We therefore obtain a 2-precategory $\sF'_{\mathrm{St}}$, which is essentially by construction the full sub-2-precategory of $\sF'$ spanned by the stable presentable categories.
    In particular, we see that the completion of $\sF'_{\mathrm{St}}$ is $\mlindiffstc$.
\end{construction}

The strategy for proving \cref{proposition:equivalence-two-composition-products} will be to construct a lax functor $\pressymstc \laxfun \fdiff_{\mathrm{St}, \geq 1}$ and use \cref{proposition:Recognizing-localizations} to identify $\pressymstc$ with the localization $\sF'_\mathrm{St}$.
For the following lemma, recall the adjunction
\[\begin{tikzcd}[sep = huge]
    \Pfin(\sD) \ar[r, shift left, "h"] & \sD \ar[l, shift left, hook', "j"]
\end{tikzcd}\]
from \cref{equation:Yoneda-and-left-adjoint}. We can view $h$ as a 1-morphism in $\pressymc$ concentrated in arity 1.

\begin{lemma}\label{lemma:comparison-composition-products-laxfun}
    There exists a lax functor $\phi \colon \pressymstc \laxfun \fdiff_{\mathrm{St}, \geq 1}$ that sends a stable presentable category to itself and which on mapping categories is the inclusion
    \[\begin{tikzcd}
        \SymFunL_{\geq 1}(\sC,\sD) \ar[r,hook,"j_*"] & \SymFun^{*,\omega}_{\geq 1}(\sC,\Pfin(\sD)) \simeq \SymFunL_{\geq 1}(\Pfin(\sC), \Pfin(\sD)).
    \end{tikzcd}\]
\end{lemma}

\begin{proof}
    We will construct the lax functor $\phi$ by ``twisting with $h \colon \Pfin(\sD) \to \sD$'' inside $\pressymc$ in the sense of \cref{theorem:Twisting-in-families}.
    In fact we will use a dual version of \cref{theorem:Twisting-in-families}: we will apply it to the category $(\pressymc)^\coop$ (the notation $\coop$ is defined in \cref{construction:opposites-of-2-categories}).
    Observe that for any pair of differentiable categories, the functor
    \[\begin{tikzcd}
    	\SymFunL_{\geq 1}(\Pfin(\sC),\Pfin(\sD)) \ar[r,"h_*"] & \SymFunL_{\geq 1}(\Pfin(\sC),\sD)
    \end{tikzcd}\]
    may be identified with
    \[\begin{tikzcd}
    	\SymFun_{\geq 1}^{*,\omega}(\sC,\Pfin(\sD)) \ar[r,"h_*"] & \SymFun^{*,\omega}_{\geq 1}(\sC,\sD)
    \end{tikzcd}\]
    This functor has a right adjoint $h^!$, given by postcomposition with $j$.
    This shows that in $(\pressymc)^\coop$, the 1-morphism $h \colon \Pfin(\sD) \to \sD$ is left extending with respect to $\Pfin(\sC)$. 
    Now consider the map $\psi \colon (\pressymstc)^\simeq \to \Ar(\pressymc)^\simeq$ that sends a stable presentable category $\sD$ to the 1-morphism $h \colon \Pfin(\sD) \to \sD$.
    Let $\sigma$ and $\tau$ denote the composite of $\psi$ with the source and target projections $s,t \colon \Ar(\pressymc)^\simeq \to (\pressymc)^\simeq$ respectively.
    Then $\tau$ is simply the inclusion $(\pressymstc)^\simeq \hookrightarrow (\pressymc)^\simeq$.
    On the other hand, $\sigma$ agrees with the map $\alpha$ used to define $\fdiff$ in \cref{definition:pre-category-Diff-FSym}, restricted to $(\pressymstc)^\simeq$.
    In particular, $\tau^*\pressymc = \pressymstc$, while $\sigma^*\pressymc = \fdiff_{\mathrm{St},\geq 1}$ is the 2-precategory of  formal symmetric sequences between stable presentable categories.
    By applying \cref{theorem:Twisting-in-families} to $(\pressymc)^\coop$, we obtain a lax functor $\pressymstc \to \fdiff_{\mathrm{St},\geq 1}$ which is given on mapping spaces by
    \[\begin{tikzcd}
    	\SymFunL_{\geq 1}(\sC,\sD) \ar[r,"{h^*}"] & \SymFunL_{\geq 1}(\Pfin(\sC),\sD) \ar[r,"{h^!}"] & \SymFunL_{\geq 1}(\Pfin(\sC),\Pfin(\sD)).
    \end{tikzcd}\]
    This is the desired lax functor.
\end{proof}

We now apply \cref{proposition:Recognizing-localizations} to prove \cref{proposition:equivalence-two-composition-products}.

\begin{proof}[Proof of \cref{proposition:equivalence-two-composition-products}]
    In order to apply \cref{proposition:Recognizing-localizations} to $\phi \colon \pressymstc \laxfun \fdiff_{\mathrm{St},\geq 1}$, we first need to verify that the composite
    \[\begin{tikzcd}
    \SymFunL_{\geq 1}(\sC,\sD) \ar[r, "h^*"] & \SymFunL_{\geq 1}(\Pfin(\sC),\sD) \ar[r,"h^!"] & \SymFunL_{\geq 1}(\Pfin(\sC),\Pfin(\sD))
    \end{tikzcd}\]
    admits a left adjoint, where $h_*$ and $h^!$ are in the proof of \cref{lemma:comparison-composition-products-laxfun}.
    Since $h^!$ is defined as the right adjoint of the postcomposition functor $h_*$, it suffices to show that $h^*$ has a left adjoint.
    This follows since we may identify $h^*$ with the inclusion $\SymFunL_{\geq 1}(\sC,\sD) \hookrightarrow \SymFun^{*,\omega}_{\geq 1}(\sC,\sD)$, which has a left adjoint given by multilinearization.
    In what follows, write $h_!$ for this left adjoint.
    Observe that both right adjoints $h^!$ and $h^*$ are fully faithful.

    We need to show that for a sequence
    \[\begin{tikzcd}
        \sC_0 \ar[r,"F_1"] & \sC_1 \ar[r,"F_2"] & \cdots \ar[r,"F_n"] & \sC_n
    \end{tikzcd}\]
    of composable 1-morphisms in $\pressymstc$, the lax structure map
    \[\begin{tikzcd}
        h^!h^*(F_n) \circ \cdots \circ h^!h^*(F_1) \ar[r] & h^!h^*(F_n \circ \cdots \circ F_1)
    \end{tikzcd}\]
    is an equivalence after applying $h_!h_*$.
    By (the dual of) \cref{remark:Unwinding-oplax-structure-maps-twisting}, this lax structure map is given by the composite
    \begin{align*}
        h^!h^* F_n& \circ \cdots \circ h^!h^* F_1 \xRightarrow{\;\eta\;} h^!h_*(h^!h^* F_n \circ \cdots \circ h^!h^* F_1) \simeq h_!(h_*h^!h^* F_n \circ \cdots \circ h^!h^* F_1) \\
        &\xRightarrow{\;\varepsilon\;} h^!(h^* F_n \circ h^!h^* F_{n-1} \circ \cdots \circ h^!h^* F_1) \simeq h^!(F_n \circ h_*h^!h^* F_{n-1} \circ \cdots \circ h^!h^* F_1) \\
        &\xRightarrow{\;\varepsilon\;} h^!(F_n \circ h^* F_{n-1} \circ \cdots \circ h^!h^* F_1) \simeq \cdots \xRightarrow{\;\varepsilon\;} \cdots\\
        &\xRightarrow{\;\varepsilon\;} h^!(F_n \circ \cdots \circ h^* F_1) \simeq h^!h^*(F_n \circ \cdots \circ F_1)
    \end{align*}
    Since $h^!$ is fully faithful, we see that the counit $\varepsilon$ is invertible. In particular, up to equivalence we may identify the lax structure map with the unit
    \[h^!h^* F_n \circ \cdots \circ h^!h^* F_1 \xRightarrow{\;\eta\;} h^!h_*(h^!h^* F_n \circ \cdots \circ h^!h^* F_1).\]
    Since $h_*$ is a reflective localization, we see that this unit becomes an equivalence after applying $h_!h_*$.
    The same argument shows that the lax structure maps
    \[\id_{\Pfin(\sC)} \xRightarrow{\;\eta\;} h^!h_*(\id_{\Pfin(\sC)}) \simeq h^!h_*(\id_{\sC})\]
    become equivalences after applying $h_!h_*$. We therefore conclude that the hypotheses of \cref{proposition:Recognizing-localizations} are satisfied and we obtain a weakly compatible family of localizations $L'$ such that $\pressymstc = L'\fdiff_{\mathrm{St},\geq 1}$.
    
    It is clear that this family of weakly compatible localizations agrees with the family of localizations $P_1L$ from \cref{prop: composition-product-compatible-with-multilinearization} on the full sub-2-precategory $\fdiff_{\mathrm{St},\geq 1}$ of $\fdiff_{\geq 1}$.
    We therefore see that the 2-precategory $\sF'_{\mathrm{St}}$ from \cref{construction:the-category-FSt} is equivalent to $\pressymstc$.
    Note that this implies that $\sF'_{\mathrm{St}}$ is complete, hence $\sF'_{\mathrm{St}} \simeq \mlindiffstc$.
    We already saw that the inclusion $\mlindiffstc \hookrightarrow \mlindiffc$ is an equivalence, so we conclude that $\pressymstc \simeq \mlindiffc$.
    It follows from our construction of this equivalence that it is indeed of the form stated in \cref{proposition:equivalence-two-composition-products}.
\end{proof}

\subsubsection{Conclusion of the proof}

We now combine all the results from this section to obtain the desired lax derivatives functor.

\begin{proof}[Proof of \cref{thm: lax-derivatives}]
    Observe that the lax functor $\widetilde{\cross} \colon \diff \to \fdiffc$ from \cref{prop: lax-cross-effects} lands in the $2$-precategory $\sF$ defined in \cref{def: image-iota-category}. By composing with the strong functor $\widetilde{\mlin} \colon \sF \to \mlindiffc$ from \cref{prop: composition-product-compatible-with-multilinearization}, and using the equivalence $\mlindiffc \simeq \pressymstc$ from \cref{proposition:equivalence-two-composition-products}, we obtain a lax functor
    \[
    \begin{tikzcd}
        \partial_* \colon \diff \ar[r] & \pressymstc.
    \end{tikzcd}
    \]
    It is clear from the description of $\widetilde{\cross}$ and $\widetilde{\mlin}$ that this functor satisfies the required properties.
\end{proof}

\begin{remark}
\label{rem: lax-comparison-map-derivatives}
    Let $\sC$, $\sD$, and $\sE$ be differentiable categories and suppose that $F \colon \sD \to \sE$ and $G \colon \sC \to \sD$ are reduced finitary functors.
    The lax functor $\partial_*$ provides a natural comparison map $\partial_*(F) \circ \partial_*(G) \to \partial_*(FG)$.
    It follows from the construction of $\partial_*$ that this map is obtained by applying $\mlin$ to the composite
    \[
    \begin{tikzcd}
        \cross(F) \circ \cross(G) \ar[r, "\eta"] & \cross p_! (\cross(F) \circ \cross(G)) \ar[r, "\delta"] & \cross(p_!\cross(F) \circ p_!\cross(G)) \ar[r, "\cross(\epsilon \circ \epsilon)"] & \cross(FG),
    \end{tikzcd}
    \]
    considered in \cref{rem: lax-comparison-map-cr}, and precomposing with the equivalence $\mlin(\cross(F) \circ \cross(G)) \simeq \mlin\cross(F) \circ \mlin\cross(G)$ from \cref{cor: multilinearization-symmetric-composition}.
\end{remark}

\cref{thm: lax-derivatives} has the following important corollary:

\begin{corollary}
\label{cor: lax-monoidal-derivatives-and-monad}
    Let $\sC$ be a differentiable category.
    Then the derivatives functor
    \[
    \begin{tikzcd}
    \partial_* \colon \Fun^{\ast, \omega}(\sC, \sC) \ar[r] & \SymFunL(\Sp(\sC), \Sp(\sC))
    \end{tikzcd}
    \]
    can be refined to a lax monoidal functor.
    In particular, for any reduced finitary monad $T \colon \sC \to \sC$, the derivatives $\partial_*T$ obtain the structure of an algebra in $\SymFunL(\Sp(\sC), \Sp(\sC))$.
\end{corollary}

\subsection{The stable chain rule}
\label{sec: stable-chain-rule}

We will now prove the following theorem.

\begin{theorem}[Stable chain rule] \label{thm: stable chain rule}
   Let $\sC \xrightarrow{G} \sD \xrightarrow{F} \sE$ be a pair of reduced finitary functors between differentiable categories.
   Suppose that $\sD$ is stable.
   Then the comparison map
   \[
   \begin{tikzcd}
   \partial_*(F) \circ \partial_*(G) \ar[r] & \partial_*(FG)
   \end{tikzcd}
   \]
   of the lax functor $\partial_*$ from \cref{thm: lax-derivatives} is an equivalence. 
\end{theorem}

This theorem has the following corollaries.

\begin{corollary}
   The restriction
    \[
    \begin{tikzcd}
    \partial_* \colon \diff_{\mathrm{St}} \ar[r] & \pressymstc
    \end{tikzcd}
    \]
    of the lax derivatives functor to the subcategory $\diff_{\mathrm{St}} \subset \diff$ spanned by the stable presentable categories is a strong $2$-functor.
\end{corollary}
\begin{proof}
    By \cref{lemma:Criterion-when-lax-functors-are-strict}, this follows from \cref{thm: stable chain rule} together with the fact that for a stable presentable category $\sD$ we have $\partial_*{\id_\sD} \simeq \unit$. The latter holds since $\id_\sD$ is a linear functor.
\end{proof}

We saw in \cref{cor: lax-monoidal-derivatives-and-monad} that the derivatives of a reduced finitary monad obtain the structure of an algebra in symmetric sequences.
In the stable case, the analogous statement for comonads also holds.

\begin{corollary}
    Let $\sC$ be a stable presentable category and let $Q \colon \sC \to \sC$ be a reduced finitary comonad.
    Then $\partial_*Q$ canonically acquires the structure of a coalgebra in $\SymFunL(\sC,\sC)$.
\end{corollary}
\begin{proof}
    By restricting the $2$-functor from the previous corollary to the appropriate mapping categories, we obtain a strong monoidal functor $\partial_* \colon \Fun^{\ast, \omega}(\sC, \sC) \to \SymFunL(\sC, \sC)$.
    Since strong monoidal functors preserve coalgebras, the result follows.
\end{proof}

\begin{example}
    Let $\sC$ be a differentiable category.
    The stabilization adjunction 
    \[
    \begin{tikzcd}
    \Sigma^\infty_\sC \colon \sC \ar[r, shift left] & \Sp(\sC) : \Omega^\infty_\sC \ar[l, shift left] 
    \end{tikzcd}
    \]
    induces a comonad $\Sigma^\infty_\sC \Omega^\infty_\sC$ on $\Sp(\sC)$.
    It is easy to see that this comonad is reduced.
    The functor $\Sigma^\infty_\sC$ is finitary since it is a left adjoint.
    To see that $\Omega^\infty_\sC$ is finitary as well, recall that $\Sp(\sC)$ is defined as the category $\Exc_1^{\ast}(\Spc^{\mathrm{fin}}_{\ast}, \sC)$ of reduced $1$-excisive functors from the category of finite pointed spaces to $\sC$.
    Since $\sC$ is differentiable, this category is closed in $\Fun(\Spc^{\mathrm{fin}}_{\ast}, \sC)$ under filtered colimits.
    It follows that $\Omega^\infty_\sC$ is finitary, since it is given by evaluation in $S^0$.
    Hence $\partial_*(\Sigma^\infty_\sC \Omega^\infty_\sC)$ obtains the structure of a coalgebra in $\SymFunL(\Sp(\sC), \Sp(\sC))$.
\end{example}

We need two lemmas for the proof of \cref{thm: stable chain rule}.
Recall that for differentiable categories $\sC$ and $\sD$, we have an adjunction 
\[
\begin{tikzcd}[sep = large]
\SymFun^{\ast, \omega}_{\geq 1}(\sC, \sD) \ar[r, shift left, "p_!"] & \Fun^{\ast, \omega}(\sC, \sD). \ar[l, shift left, "\cross"]
\end{tikzcd}
\]
The left adjoint $p_!$ is given by the formula $p_!F(x) = \coprod_{n \geq 1} F(x, \ldots, x)_{h\Sigma_n}$.

\begin{lemma}
\label{lem: cross-effects unit equiv}
    Let $\sC$ and $\sD$ be differentiable categories and suppose that $\sD$ is stable. 
    Let $F \in \SymFunL(\sC, \sD)$ be a symmetric sequence concentrated in finitely many arities.
    Then the unit $\eta \colon F \to \cross( p_! F)$ becomes an equivalence after multilinearization.
\end{lemma}
\begin{proof}
    The functor $p_!$ commutes with coproducts, since it is a left adjoint.
    Because $\sD$ is stable, $\cross \colon \Fun^{\ast, \omega}(\sC, \sD) \to \SymFun^{\ast, \omega}_{\geq 1}(\sC, \sD)$ commutes with finite coproducts.
    We can therefore reduce to the case where $F$ is concentrated in a single arity $n$.
    Then the functor $p_!F$ is $n$-homogeneous, and the map $\eta$ is an equivalence in degree $n$ by the equivalence between homogeneous and symmetric multilinear functors.
    For $k \neq n$, the multilinearization of $\cross_k(p_!F)$ is $\partial_k(p_!F)$, which vanishes.
\end{proof}

\begin{lemma}
\label{lem: lax-comparison-map-stabilization-functors}
    Let $\sC$ and $\sD$ be differentiable categories and let $F \colon \sC \to \Sp(\sD)$ be a reduced finitary functor.
    Then the lax comparison map
    \[
    \begin{tikzcd}
        \partial_*(\Omega^\infty_\sD) \circ \partial_*(F) \ar[r] & \partial_*(\Omega^\infty_\sD F)
    \end{tikzcd}
    \]
    coming from \cref{thm: lax-derivatives} is an equivalence.
\end{lemma}
\begin{proof}
    Observe that since $\Omega^\infty_\sD$ is a linear functor, $\cross(\Omega^\infty_\sD)$ is concentrated in arity $1$, where it is given by $\Omega^\infty_\sD$.
    Therefore, the composition product $\cross(\Omega^\infty_\sD) \circ \cross(F)$ (\cref{def: symmetric-composition-product}) is given in arity $n$ by $\Omega^\infty_\sD \cross_n(F)$. This is equivalent to $\cross_n(\Omega^\infty_\sD F)$, since cross effects are computed as a total fiber and hence commute with any limit preserving functor.
    So we have an equivalence $\cross(\Omega^\infty_\sD) \circ \cross(F) \simeq \cross(\Omega^\infty F)$.
    Consider the following diagram
    \[
    \begin{tikzcd}
        \cross(\Omega^\infty_\sD) \circ \cross(F) \ar[r, "\eta"] \ar[d, "\sim"] & \cross p_!(\cross(\Omega^\infty_\sD) \circ \cross(F)) \ar[r, "\cross(\delta)"] \ar[d, "\sim"] &\cross(p_!\cross(\Omega^\infty_\sD) \circ p_!\cross(F)) \ar[d, "\cross(\epsilon \circ \epsilon)"] \\
        \cross(\Omega^\infty_\sD F) \ar[r, "\eta"] & \cross p_! \cross(\Omega^\infty_\sD F) \ar[r, "\cross(\epsilon)"] & \cross(\Omega^\infty_\sD F) &
    \end{tikzcd}
    \]
    The maps $\eta$ and $\epsilon$ are the unit and counit of the adjunction between $\cross$ and $p_!$.
    The map $\delta$ was defined in \cref{rem: lax-comparison-map-cr}.
    By \cref{rem: lax-comparison-map-derivatives}, applying $\mlin$ to the composite starting in the top left corner and going right twice and down once yields the comparison map $\partial_*(\Omega^\infty_\sD) \circ \partial_*(F) \to \partial_*(\Omega^\infty_\sD F)$.
    It therefore suffices to show this composite is an equivalence.
    Note that the bottom horizontal composite is an equivalence by the triangle identities, so we are done if we can show the diagram commutes.
    The left square commutes by naturality of $\eta$.
    Observe that $p_! \cross(\Omega^\infty_\sD) = \Omega^\infty_\sD$.
    By \cref{rem: lax-comparison-map-cr}, the map $\delta$
    is given by the colimit assembly map
    \[
    \begin{tikzcd}
    \coprod_{n \geq 1} (\Omega^\infty \cross_nF(x, \ldots, x))_{h\Sigma_n} \ar[r] & \Omega^\infty\coprod_{n \geq 1} \cross_nF(
    x, \ldots, x)_{h\Sigma_n}.
    \end{tikzcd}
    \]
    For any functor $G$ and $n \geq 1$, the counit $\epsilon \colon p_! \cross_n G \to G$ is given by the composite
    \[
    \begin{tikzcd}
    \cross_nG(x, \ldots, x)_{h\Sigma_n} \ar[r] & G(x \sqcup \ldots \sqcup x)_{h\Sigma_n} \ar[r, "\nabla"] & G(x),
    \end{tikzcd}
    \]
    where the first map is the canonical map coming from the definition of the cross effects as a total fiber and the second map is induced by the fold map.
    From this it is easily seen that the right-hand square commutes as well, completing the proof.
\end{proof}

\begin{proof}[Proof of Theorem \ref{thm: stable chain rule}]
    Let $G \colon \sC \to \sD$ and $F \colon \sD \to \sE$ be a pair of reduced finitary functors between differentiable categories and suppose $\sD$ is stable.
    We need to show that the lax comparison map
    \[
    \mu \colon \partial_*F \circ \partial_*G \longrightarrow \partial_*(F \circ G)
    \]
    is an equivalence.
    We start by making a series of reductions.
    To show $\mu$ is an equivalence, it suffices to show it is an equivalence in each arity $n$.
    By naturality, we get a commutative square
    \[
    \begin{tikzcd}
        (\partial_*F \circ \partial_*G)_n \ar[r, "\mu"] \ar[d] & \partial_n(F \circ G) \ar[d] \\
        (\partial_*P_nF \circ \partial_*P_nG)_n \ar[r, "\mu"] & \partial_n(P_nF \circ P_nG).
    \end{tikzcd}
    \]
    The left vertical arrow is an equivalence because the arity $n$ term of the composition product of two symmetric sequences whose arity $0$ terms vanish depends only on the terms of arity $\leq n$ of both sequences.
    The right vertical arrow is an equivalence by Lemma \ref{lem: Arone-Ching-lemma}.
    So we can assume that $F$ and $G$ are $n$-excisive for some $n$.
    
    Observe that the functors $\partial_*(-) \circ \partial_*G$ and $\partial_*(- \circ G)$ both preserve finite limits.
    Again using the naturality of $\mu$, it follows by induction on the fiber sequences
    \[
    \begin{tikzcd}
    D_kF \ar[r] & P_k F \ar[r] & P_{k-1}F
    \end{tikzcd}
    \]
    that we can reduce to the case where $F$ is $k$-homogeneous for some $k$.
    Then $F \colon \sD \to \sE$ factors over $\Omega^\infty_\sD \colon \Sp(\sE) \to \sE$, so that by \cref{lem: lax-comparison-map-stabilization-functors} it moreover suffices to prove the theorem in the case where $\sE$ is stable.
    
    By \cref{prop: exactness-composition-product}, the composition product commutes with finite totalizations in the second variable.
    It follows from \cite[Proposition 3.37]{brantner2023deformation} that an $n$-excisive functor between stable presentable categories preserves finite totalizations.
    As the derivatives functor preserves finite limits, we find that both sides of the comparison map $\mu$ commute with finite totalizations in the $G$-variable.
    As $G$ is $n$-excisive, it follows from the Dold-Kan Theorem for stable categories \cite[\S 1.2.3]{HA} applied to the finite tower
    \[
    \begin{tikzcd}
    G = P_nG \ar[r] & P_{n-1}G \ar[r] & \cdots \ar[r] & P_1G
    \end{tikzcd}
    \]
    that $G$ is a finite totalization of finite coproducts of homogeneous functors.
    So we can reduce to the case where $G$ is such a coproduct. 
    Note that by the classification of homogeneous functors it now suffices to prove the theorem in the special case where $F = p_!A$ and $G = p_!B$ for $A \in \SymFunL(\sD, \sE)$ and $B \in \SymFunL(\sC, \sD)$ and both $A$ and $B$ are concentrated in finitely many arities. 
    We'll proceed to do so.

    Consider the following commutative diagram:
    \[
    \begin{tikzcd}
        A \circ B \ar[r, "\eta"] \ar[d, "\eta \circ \eta"] & \mathrm{cr}(p_!(A \circ B)) \ar[r, "\sim"] \ar[d, "\mathrm{cr}(p_!(\eta \circ \eta))"] & \mathrm{cr}(p_!A \circ p_!B) \ar[d, "\mathrm{cr}(p_! \eta \circ p_! \eta)"] \\
        \mathrm{cr}(p_!A) \circ \mathrm{cr}(p_! B) \ar[r, "\eta"] \ar[drr, "\nu"'] & \mathrm{cr}(p_!(\mathrm{cr}(p_!A) \circ \mathrm{cr}(p_!B))) \ar[r, "\sim"] & \mathrm{cr}(p_!\mathrm{cr}(p_!A) \circ p_!\mathrm{cr}(p_!B)) \ar[d, "\mathrm{cr}(\epsilon \circ \epsilon)"] \\ & & \mathrm{cr}(p_! A \circ p_! B)      
    \end{tikzcd}
    \]
    The object $\mathrm{cr}(p_!A) \circ \mathrm{cr}(p_! B)$ appearing in the bottom left corner is given by the (non-associative) composition product defined in \cref{def: symmetric-composition-product}.
    Observe that by \cref{rem: lax-comparison-map-derivatives}, the morphism $\nu$ becomes
    \[
    \begin{tikzcd}
    \mu \colon \partial_*p_!A \circ \partial_*p_!B \ar[r] & \partial_*(p_! A \circ p_! B)
    \end{tikzcd}
    \]
    after multilinearization, which is the morphism we want to show is an equivalence.
    The two horizontal arrows on the right are equivalences by \cref{rem: lax-comparison-map-cr} (in the notation of that remark, they are both induced by the comparison map $\delta$).
    The right vertical composite is homotopic to the identity by the triangle identities for the adjunction between $p_!$ and $\mathrm{cr}$.
    Applying $\mlin$ to the top left horizontal map $\eta$ in the diagram, we obtain an equivalence by \cref{lem: cross-effects unit equiv}. 
    Since $\mlin(\eta \circ \eta) \simeq \mlin(\eta) \circ \mlin(\eta)$ by \cref{cor: multilinearization-symmetric-composition}, this map is an equivalence as well.
    We are now done by 2-out-of-3 for equivalences.
\end{proof}
\newpage

\section{Koszul duality and the chain rule}
\label{sec: koszul-duality}

In this chapter we finally prove the chain rule.
As we explained in the introduction, the theory of Koszul duality (or bar-cobar duality) plays an important role in the proof.

In \cref{sec: bar-cobar-duality}, we construct the bar-cobar adjunction by showing that $\Bahr$ and $\Cobar$ satisfy a universal property in terms of compatible right actions and left coactions on the unit object.
This treatment of Koszul duality owes a large debt to \cite{lurieThursdayKoszul}.
In \cref{sec: koszul-symm-seq}, we show that bar-cobar duality gives an equivalence between augmented algebras and coaugmented coalgebras in strongly positive symmetric sequences.
In \cref{sec: koszul-dual-partial-id}, we prove that for every differentiable category $\sC$, there is an equivalence $\partial_*\id_\sC \simeq \Cobar \partial_*(\Sigma^\infty_\sC \Omega^\infty_\sC)$ of augmented algebras.
After all this preparatory work, we prove the chain rule in \cref{sec: chain-rule}.

\subsection{The bar-cobar adjunction}\label{sec: bar-cobar-duality}

\subsubsection{Endomorphism objects}

\begin{definition}[Endomorphism objects]
    Let $\sC$ be a monoidal category and let $\sM$ be a category that is right-tensored over $\sC$.
    Suppose $x \in \sM$.
    An \emph{endomorphism object} of $x$ is an object $e \in \sC$ together with a map $a \colon x \otimes e \to x$ such that for any $d \in \sC$, the following composite is an equivalence
    \[
    \begin{tikzcd}
    \Map_{\sC}(d, e) \ar[r, "x \otimes -"] & \Map_{\sM}(x \otimes d, x \otimes e) \ar[r, "a_*"] & \Map_{\sM}(x \otimes d, x).
    \end{tikzcd}
    \]
    If an endomorphism object exists, it is unique up to contractible choice.
    We will write $\End(x)$ for $e$.
\end{definition}

\begin{remark}
    This definition concerns endomorphism objects classifying right actions.
    Of course one can also define endomorphism objects for left actions, and all results from this section hold in that case as well.
\end{remark}

In practice, it is often easy to verify the condition from this definition.
The key fact is that an endomorphism object canonically obtains the structure of an algebra, which satisfies a universal property in the category $\Alg(\sC)$.

\begin{proposition}[{\cite[Corollary 4.7.1.41]{HA}}] \label{prop: endomorphism-object-as-algebra}
    Let $\sM$ be right-tensored over $\sC$.
    Suppose that $x \in \sM$ admits an endomorphism object $\End(x) \in \sC$.
    Then $\End(x)$ can be equipped with the structure of an algebra in $\sC$, such that its action on $x$ can be refined to a right module structure.
    Moreover, for every algebra $A \in \Alg(\sC)$, the canonical map
    \[
    \begin{tikzcd}
    \Map_{\Alg(\sC)}(A, \End(x)) \ar[r] & \RMod_A(\sM) \times_\sM \{x\}
    \end{tikzcd}
    \]
    is an equivalence.
\end{proposition}

\begin{remark}
    One can show that the two universal properties of endomorphism objects are compatible in the following way.
    Suppose that $E$ is an algebra in $\sC$ and $x$ is equipped with the structure of a right $E$-module.
    Then the pair $(E, x)$ satisfies the universal property of the previous proposition if and only if the action map $x \otimes E \to x$ coming from the module structure exhibits $E$ as the endomorphism object of $x$. 
    See \cite[Lemma 3.32]{brantner2023pd} for a proof.
\end{remark}

\begin{definition}[Coendomorphism objects]
    Let $\sM$ be left-tensored over $\sC$.
    Then a \emph{coendomorphism object} of $x \in \sM$ is an endomorphism object of $x \in \sM^{\mathrm{op}}$, where we regard $\sM^{\op}$ as left-tensored over $\sC^{\mathrm{op}}$.
    If such a coendomorphism object exists we denote it by $\coEnd(x)$.
    By the above proposition, any coendomorphism object acquires the structure of a coalgebra in $\sC$, and $x$ acquires the structure of a left comodule such that for any $Q \in \coAlg(\sC)$, the map
    \[
    \begin{tikzcd}
    \Map_{\coAlg(\sC)}(\coEnd(x), Q) \ar[r] & \LComod_{Q}(\sM) \times_{\sM} \{x\}
    \end{tikzcd}
    \]
    is an equivalence.
\end{definition}

Next, we will give some examples of endomorphism and coendomorphism objects that will play a role in the rest of this chapter.

\begin{example}[Comonad of an adjunction] \label{ex: comonad-adjunction}
    Suppose that
    \[
    \begin{tikzcd}
        \sC \ar[r, shift left=1, "F"] & \sD \ar[l, shift left = 1, "G"]
    \end{tikzcd}
    \]
    is a pair of adjoint functors.
    The category $\Fun(\sC, \sD)$ is left-tensored over $\Fun(\sD, \sD)$, where the monoidal structure and tensoring are given by functor composition.
    It is an easy verification that the unit
    \[
    \begin{tikzcd}
    F \ar[r, "\eta"] & F(GF) \simeq (FG)F
    \end{tikzcd}
    \]
    exhibits $FG$ as a coendomorphism object of $F$ \cite[Lemma 4.7.3.1]{HA}.
    By \cref{prop: endomorphism-object-as-algebra}, this endows $FG$ with the structure of a comonad on $\sC$. Dually, the counit of the adjunction exhibits $GF$ as the endomorphism object of $F$, giving $GF$ the structure of a monad on $\sC$.
\end{example}

\begin{definition}
    Let $\sC$ be a monoidal category with monoidal unit $\unit$.
    An \emph{augmented algebra} in $\sC$ is an algebra $A \in \Alg(\sC)$ together with a map of algebras $A \to \unit$.
    A \emph{coaugmented coalgebra} in $\sC$ is a coalgebra $Q \in \coAlg(\sC)$ together with a map of coalgebras $\unit \to Q$.
    We write $\Alg^{\mathrm{aug}}(\sC)$ for the category of augmented algebras and $\coAlg^{\mathrm{aug}}(\sC)$ for the category of coaugmented coalgebras in $\sC$.
\end{definition}

\begin{example}[Endomorphisms of the unit]
    Let $\sC$ be a monoidal category with unit $\unit$, and regard $\sC$ as right-tensored over itself.
    Then $\End(\unit) \simeq \unit$.
    Indeed, the endomorphism object of $\unit$ is an object $e$ with the property that giving a map $d \to e$ is equivalent to giving a map $\unit \otimes d \simeq d \to \unit$. 
    It is clear that $\unit$ itself satisfies this property.

    For any $A \in \Alg(\sC)$, it then follows from \cref{prop: endomorphism-object-as-algebra} that we get an equivalence
    \[
    \Map_{\Alg(\sC)}(A, \unit) \simeq \Map_{\Alg(\sC)}(A, \End(\unit)) \simeq \RMod_A(\sC) \times_\sC \{ \unit \}.
    \]
    Therefore, endowing $A$ with an augmentation is equivalent to equipping $\unit$ with a right $A$-module structure.
    The same proof shows that giving an augmentation of $A$ is equivalent to giving a \emph{left} $A$-module structure on $\unit$.
\end{example}

The following proposition describes the naturality of endomorphism objects.

\begin{proposition} \label{prop: naturality-endomorphism-objects}
    Let $f \colon \sC \to \sD$ be a lax monoidal functor and let $\sM \in \RMod_{\sC}(\Cat)$ and $\sN \in \RMod_{\sD}(\Cat)$.
    Let $g \colon \sM \to \sN$ be a lax $\sC$-linear functor, where we regard $\sN$ as a $\sC$-module via restriction of scalars along $f$.
    Suppose that $x \in \sM$ is an object such that both $x$ and $g(x)$ admit endomorphism objects.
    Then there is a unique algebra morphism $\phi \colon f(\End(x)) \to \End(g(x))$, such that for any $A \in \Alg(\sC)$, the map
    \[
    \begin{tikzcd}
    \RMod_A(\sM) \times_\sM \{ x \} \ar[r, "g_*"] & \RMod_{f(A)}(\sN) \times_\sN \{ g(x) \}
    \end{tikzcd}
    \]
    corresponds under the equivalence of \cref{prop: endomorphism-object-as-algebra} to the map
    \[
    \begin{tikzcd}
    \Map_{\Alg(\sC)}(A, \End(x)) \ar[r] & \Map_{\Alg(\sD)}(f(A), \End(g(x)))
    \end{tikzcd}
    \]
    induced by $\phi$.
\end{proposition}
\begin{proof}
The object $g(x)$ obtains the structure of a right $f(\End(x))$-module, since $g$ is a lax $\sC$-linear functor.
This module structure is classified by an algebra morphism $\phi \colon f(\End(x)) \to \End(g(x))$, and it is easily seen that this morphism has the required property.
\end{proof}

\subsubsection{Indecomposables, primitives, and Koszul duality}

\begin{definition}
    Let $\sC$ be a monoidal category and let $\epsilon \colon A \to \unit$ be an augmented algebra in $\sC$.
    Then restriction of scalars along $\epsilon$ induces a functor
    \[
    \triv_A \colon \sC \simeq \RMod_\unit(\sC) \to \RMod_A(\sC),
    \]
    the \emph{trivial $A$-module functor}.
\end{definition}

Note that since restriction of scalars is a left $\sC$-linear functor by \cref{lem: Restriction-of-scalars-C-linear} (where we equip $\RMod_A(\sC)$ with the standard left tensoring over $\sC$), we have a natural equivalence $\triv_A(x) \simeq x \otimes \triv_A(\unit)$.

Under mild assumptions on the category $\sC$, the functor $\triv_A$ admits a left adjoint.

\begin{proposition}
    Let $\sC$ be a monoidal category and let $A \in \Alg^{\mathrm{aug}}(\sC)$. 
    Suppose that $\sC$ admits geometric realizations of simplicial diagrams.
    Then the functor $\triv_A$ admits a left adjoint.
\end{proposition}
\begin{proof}
    Let $\sE \subset \RMod_A(\sC)$ be the full subcategory spanned by the $A$-modules $M$ for which the functor
    \[
    \Map_{\RMod_A}(M, \triv_A(-))
    \]
    is corepresentable.
    We need to show that $\sE = \RMod_A(\sC)$.
    Note that the free $A$-modules $x \otimes A$ all lie in $\sE$, since we have
    \[
    \Map_{\RMod_A}(x \otimes A, \triv_A(-)) \simeq \Map_{\sC}(x, \forget_A \triv_A(-)) \simeq \Map_{\sC}(x, -),
    \]
    that is, the functor is corepresented by $x$ in this case.
    Given an arbitrary $A$-module $M$, let $B_{\bullet}(M, A, A) = M \otimes A^{\otimes \bullet + 1}$ be the standard simplicial resolution of $M$ consisting of free $A$-modules, so that we have $M \simeq |B_{\bullet}(M, A, A)|$.
    Note that the functor $\triv_A$ admits a partially defined left adjoint $F \colon \sE \to \sC$.
    As $B_{\bullet}(M, A, A)$ lies entirely in $\sE$, we can apply $F$ to it to obtain a simplicial object in $\sC$.
    As $\sC$ admits geometric realizations, the colimit of $F(B_{\bullet}(M, A, A))$ exists and it is clear that this is the corepresenting object we are after.
\end{proof}

When the left adjoint to $\triv_A$ exists, we call it the \emph{indecomposables} functor:
\[
\begin{tikzcd}[sep = large]
    \RMod_A(\sC) \ar[r, shift left, "\indec_A"] & \sC \ar[l, shift left, "\triv_A"].
\end{tikzcd}
\]
Recall that if $A \in \Alg(\sC)$, $M \in \RMod_A(\sC)$, and $N \in \LMod_A(\sC)$, then we can form the \emph{bar construction} \cite[Construction 4.4.2.7]{HA}, which is a simplicial object
\[
B_\bullet(M, A, N) =
\begin{tikzcd}
          \cdots M \otimes A \otimes A \otimes N \ar[r, shift left=2] \ar[r] \ar[r, shift right=2] & M \otimes A \otimes N \ar[l, shift left, shorten=0.4em] \ar[l, shift right, shorten=0.4em] \ar[r, shift left] \ar[r, shift right] & M \otimes N. \ar[l, shorten=0.4em]
    \end{tikzcd}
\]
The face maps are given by the multiplication of $A$ and the right and left actions of $A$ on $M$ and $N$.
The degeneracies are given by the unit of $A$.
We will write $B(M, A, N)$ or $\Bahr(M, A, N)$ for its geometric realization, and we will sometimes refer to this as the bar construction as well.
Note that $B(M, A, N)$ is equivalent to the relative tensor product $M \otimes_A N$.
However, we warn the reader that this notation might be misleading: if the tensor product on $\sC$ does not commute with geometric realizations, then the relative tensor product $- \otimes_A -$ is not associative.

If $A$ and $\sC$ are as in the previous proposition, then it follows from the proof that for any $A$-module $M$, we have the formula
\[
\indec_A(M) \simeq B(M, A, \unit),
\]
where $\unit$ is a left $A$-module via the augmentation.
We write $\Bahr(A)$ (or sometimes $BA$ for brevity) for the bar construction $B(\unit, A, \unit)$, which is equivalent to $\indec_A(\triv_A(\unit))$.
The following proposition implies that $\Bahr(A)$ canonically admits the structure of a coalgebra in $\sC$.

\begin{proposition}
\label{prop: bar-is-coendomorphism-object}
    Let $\sC$ be a monoidal category that admits geometric realizations of simplicial objects and let $A \in \Alg^{\mathrm{aug}}(\sC)$.
    Then $\Bahr(A)$ is the coendomorphism object of $\triv_A(\unit)$, where we use the standard left-tensoring of $\RMod_A(\sC)$ over $\sC$.
\end{proposition}
\begin{proof}
    This follows immediately from the adjunction between $\triv_A$ and $\indec_A$, since for any $x \in \sC$ we have
    \begin{align*}
    \Map_{\sC}(\indec_A(\triv_A(\unit)), x) &\simeq \Map_{\RMod_A}(\triv_A(\unit), \triv_A(x)) \\ 
                                            &\simeq \Map_{\RMod_A}(\triv_A(\unit), x \otimes \triv_A(\unit)) \qedhere
    \end{align*}
\end{proof}

\begin{remark} \label{rem: coaugmentation-of-Bar(A)}
    As the forgetful functor $\RMod_A(\sC) \to \sC$ is compatible with the left-tensoring by $\sC$, it follows from \cref{prop: naturality-endomorphism-objects} that we obtain a canonical coalgebra morphism $\unit \simeq \coEnd(\unit) \to \Bahr(\unit)$. 
    This equips $\Bahr(A)$ with a coaugmentation.
\end{remark}

We call the coaugmented coalgebra $\Bahr(A)$ the \emph{Koszul dual} of $A$.

The preceding discussion has a counterpart for coaugmented coalgebras.
Suppose that $\eta \colon \unit \to Q$ is a coaugmented coalgebra.
Then $\eta$ induces the \emph{trivial $Q$-comodule} functor, which admits a right adjoint if $\sC$ admits totalizations of cosimplicial objects. We call this right adjoint the \emph{primitives} functor:
\[
\begin{tikzcd}[sep = large]
    \sC \ar[r, shift left, "\triv_Q"] & \LComod_Q(\sC) \ar[l, shift left, "\prim_Q"].
\end{tikzcd}
\]
For a left $Q$-comodule $N$, the primitives functor is given by the \emph{cobar construction} $C(\unit, Q, N)$ (we will also use the notation $\Cobar(\unit, Q, N)$).
In general, if $M$ is a right $Q$-comodule and $N$ is a left $Q$-comodule, then the cobar construction is the cosimplicial object
\[
C^\bullet(M, Q, N) =
\begin{tikzcd}
        M \otimes N \ar[r, shift left] \ar[r, shift right] & M \otimes Q \otimes N \ar[l, shorten=0.4em] \ar[r, shift left=2] \ar[r] \ar[r, shift right=2] & M \otimes Q \otimes Q \otimes N \cdots. \ar[l, shift left, shorten=0.4em] \ar[l, shift right, shorten=0.4em]
    \end{tikzcd}
\]
We write $\Cobar(Q)$ for the cobar construction $C(\unit, Q, \unit)$.
This is the endomorphism object of $\triv_Q(\unit) \in \LComod_Q(\sC)$, and hence canonically acquires the structure of an augmented algebra in $\sC$.
We call it the \emph{Koszul dual} of $Q$.

The main theorem of Koszul duality is that $\Bahr$ and $\Cobar$ define an adjunction between the categories of augmented algebras and coaugmented coalgebras in $\sC$.

\begin{theorem} \label{thm: koszul-left-module-right-comodule-structure-on-unit}
Let $\sC$ be a monoidal category admitting geometric realizations of simplicial objects and totalizations of cosimplicial objects.
Suppose that $A$ is an augmented algebra and $Q$ is a coaugmented coalgebra in $\sC$.
Then the spaces $\Map_{\coAlg^{\mathrm{aug}}(\sC)}(\Bahr(A), Q)$ and $\Map_{\Alg^{\mathrm{aug}}(\sC)}(A, \Cobar(Q))$ are both naturally equivalent to
\[
\{\triv_Q(\unit)\} \times_{\LComod_Q(\sC)} \LComod_Q(\RMod_A(\sC)) \times_{\RMod_A(\sC)} \{ \triv_A(\unit) \}.
\]
\end{theorem}
\begin{proof}
Write $\eta \colon \unit \to Q$ for the coaugmentation of $Q$. 
Observe that we have a cartesian square
\[
\begin{tikzcd}[sep = large]
\Map_{\coAlg^{\mathrm{aug}}(\sC)}(\Bahr(A), Q) \ar[r] \ar[d] & \Map_{\coAlg(\sC)}(\Bahr(A), Q) \ar[d] \\
\ast \ar[r, "\eta"] & \Map_{\coAlg(\sC)}(\unit, Q).
\end{tikzcd}
\]
The right-hand map is restriction along the coaugmentation of $\Bahr(A)$.
By \cref{prop: naturality-endomorphism-objects} and \cref{rem: coaugmentation-of-Bar(A)}, this map is equivalent to the canonical map
\[
\LComod_Q(\RMod_A(\sC)) \times_{\RMod_A(\sC)} \{ \triv_A(\unit) \} \to \LComod_Q(\sC) \times_{\sC} \{ \unit \},
\]
and $\eta$ corresponds to the left $Q$-comodule $\triv_Q(\unit)$.
Hence $\Map_{\coAlg^{\mathrm{aug}}(\sC)}(\Bahr(A), Q)$ is of the desired form.
A completely analogous argument shows $\Map_{\Alg^{\mathrm{aug}}(\sC)}(A, \Cobar(Q))$ is equivalent to
\[
\{\triv_Q(\unit)\} \times_{\LComod_Q(\sC)} \RMod_A(\LComod_Q(\sC)) \times_{\RMod_A(\sC)} \{ \triv_A(\unit) \}.
\]
We are now done since there is a natural equivalence
\[
\RMod_A(\LComod_Q(\sC)) \simeq \LComod_Q(\RMod_A(\sC)),
\]
by \cref{theorem:bicomod-interchange}.
\end{proof}

\begin{remark}
    This theorem says that giving a map of coaugmented coalgebras $\Bahr(A) \to Q$ (or equivalently, a map of augmented algebras $A \to \Cobar(Q)$), is equivalent to equipping $\unit$ with a compatible right $A$-module and left $Q$-comodule structure, such that forgetting the right $A$-action gives the coaugmentation of $Q$, and forgetting the left $Q$-coaction gives the augmentation of $A$.
\end{remark}

\begin{corollary}
\label{cor: bar-cobar-adjunction}
    Let $\sC$ be a monoidal category admitting geometric realizations of simplicial objects and totalizations of cosimplicial objects.
    Then the constructions $\Bahr$ and $\Cobar$ give rise to a pair of adjoint functors (left adjoint on top):
    \[
    \begin{tikzcd}[sep = large]
        \Alg^{\mathrm{aug}}(\sC) \ar[r, shift left, "\Bahr"] & \coAlg^{\mathrm{aug}}(\sC) \ar[l, shift left, "\Cobar"].
    \end{tikzcd}
    \]
\end{corollary}
\begin{remark}
    Such an adjunction was also constructed in \cite[Section 5.2]{HA} and \cite[Section 3.4]{brantner2023pd}. 
    It is a priori not clear that the adjunction considered in these sources is equivalent to the one from this corollary. 
    We will prove this in the next section.
\end{remark}
\begin{proof}
    Consider the functor $\Alg^{\mathrm{aug}}(\sC)^{\op} \times \coAlg^{\mathrm{aug}}(\sC) \to \Spc$ that sends a pair $(A, Q)$ to the space
    \[
    \{\triv_Q(\unit)\} \times_{\LComod_Q(\sC)} \LComod_Q(\RMod_A(\sC)) \times_{\RMod_A(\sC)} \{ \triv_A(\unit) \}.
    \]
    Adjointing over one factor to the right, we obtain a functor
    \[
    \begin{tikzcd}
    \coAlg^{\mathrm{aug}}(\sC) \ar[r] & \Fun(\Alg^{\mathrm{aug}}(\sC)^{\op}, \Spc),
    \end{tikzcd}
    \]
    and it follows from \cref{thm: koszul-left-module-right-comodule-structure-on-unit} that this functor lands in the image of the Yoneda embedding.
    Composing with the inverse of this embedding, we obtain the functor
    \[
    \begin{tikzcd}
    \Cobar \colon \coAlg^{\mathrm{aug}}(\sC) \ar[r] & \Alg^{\mathrm{aug}}(\sC).
    \end{tikzcd}
    \]
    It follows from the same theorem that this functor admits a left adjoint, given on objects by $A \mapsto \Bahr(A)$.
\end{proof}

The bar-cobar adjunction has the following naturality property.

\begin{proposition}
    Let $\sC$ and $\sD$ be monoidal categories that admit geometric realizations of simplicial objects.
    Suppose that $F \colon \sC \to \sD$ is a strong monoidal functor that preserves geometric realizations.
    Then the following square commutes
    \[
    \begin{tikzcd}[sep = large]
        \Alg^{\mathrm{aug}}(\sC) \ar[r, "\Bahr"] \ar[d, "F"] & \coAlg^{\mathrm{aug}}(\sC) \ar[d, "F"] \\
        \Alg^{\mathrm{aug}}(\sD) \ar[r, "\Bahr"] & \coAlg^{\mathrm{aug}}(\sD).
    \end{tikzcd}
    \]
    The analogous statement for $\Cobar$ is true if $F$ preserves totalizations.
\end{proposition}
\begin{proof}
    Since $F$ is a strong monoidal functor, it induces a natural map from
    \[
    \{\triv_Q(\unit_\sC)\} \times_{\LComod_Q(\sC)} \LComod_Q(\RMod_A(\sC)) \times_{\RMod_A(\sC)} \{\triv_A(\unit_\sC)\}
    \]
    to
    \[
    \{\triv_{F(Q)}(\unit_\sD)\} \times_{\LComod_{F(Q)}(\sD)} \LComod_{F(Q)}(\RMod_{F(A)}(\sD)) \times_{\RMod_{F(A)}(\sD)} \{\triv_{F(A)}(\unit_\sD)\},
    \]
    where we consider both as functors $\coAlg^{\mathrm{aug}}(\sC) \times \Alg^{\mathrm{aug}}(\sC)^{\op} \to \Spc$.
    Adjointing over one factor, we obtain a natural transformation as follows:
    \[
    \begin{tikzcd}
        \Alg^{\mathrm{aug}}(\sC)^{\op} \ar[r, ""{name=0}] \ar[d, "F"] & \sP(\coAlg^{\mathrm{aug}}(\sC)^{\op}) \\
        \Alg^{\mathrm{aug}}(\sD)^{\op} \ar[r, ""{name=1}] & \sP(\coAlg^{\mathrm{aug}}(\sD)^{\op}) \ar[u, "F^*"'].
        \arrow[shorten <=5pt, shorten >=5pt, Rightarrow, from=1, to=0]
    \end{tikzcd}
    \]
    Since both horizontal functors land in the full subcategories spanned by the representables by \cref{thm: koszul-left-module-right-comodule-structure-on-unit}, the Yoneda lemma yields a natural transformation $\alpha \colon \Bahr \circ F \to F \circ \Bahr$.
    We will show $\alpha$ is an equivalence.

    Let $A \in \Alg^{\mathrm{aug}}(\sC)$.
    By construction, the map $\alpha \colon \Bahr(F(A)) \to F(\Bahr(A))$ comes from the naturality of coendomorphism objects (\cref{prop: naturality-endomorphism-objects}).
    Inspecting the proof of \cref{prop: bar-is-coendomorphism-object}, we see that this map can be described as the Beck--Chevalley transformation $\indec_{F(A)} \circ F \to F \circ \indec_A$ of the right square in the diagram
    \[
    \begin{tikzcd}[sep = huge]
       \sC \ar[d, "F"] & \RMod_A(\sC) \ar[l, "\mathrm{forget}"'] \ar[d, "F"] & \sC \ar[l, "\triv_A"'] \ar[d, "F"] \\
       \sD & \RMod_{F(A)}(\sD) \ar[l, "\mathrm{forget}"'] & \sD \ar[l, "\triv_{F(A)}"']
    \end{tikzcd}
    \]
    evaluated at the object $\triv_A(\unit_\sC)$.
    We will show that this Beck--Chevalley transformation is an equivalence in general.
    Since both $\indec_{F(A)} \circ F$ and $F \circ \indec_A$ commute with geometric realizations, it suffices to prove this for free $A$-modules.
    
    Observe that the Beck--Chevalley transformation of the outer square in the diagram above is an equivalence, since both horizontal composites are equivalent to the identity functor.
    This transformation factors as the composite of the Beck--Chevalley transformations of the two inner squares in the diagram:
    \[
    \begin{tikzcd}
        \indec_{F(A)}(F(A) \otimes F(X)) \ar[r] & \indec_{F(A)}(F(A \otimes X)) \ar[r] & F(\indec_A(A \otimes X)).
    \end{tikzcd}
    \]
    The left-hand map is an equivalence since $F$ is strong monoidal.
    Hence the right-hand map is an equivalence by 2-out-of-3, completing the proof.
\end{proof}

\begin{remark}
\label{rem: koszul-unit-initial-final}
    For a monoidal category $\sC$ with unit $\unit$, the double slice category $\sC_{\unit // \unit}$ inherits a monoidal structure for which the forgetful functor $\sC_{\unit // \unit} \to \sC$ is strong monoidal \cite[Remark 2.2.2.5]{HA}. 
    It also preserves geometric realizations and totalizations.
    Moreover, it induces an equivalence on categories of augmented algebras 
    \[
    \begin{tikzcd}
    \Alg^{\mathrm{aug}}(\sC_{\unit // \unit}) \ar[r, "\sim"] & \Alg^{\mathrm{aug}}(\sC),
    \end{tikzcd}
    \]
    and similarly for coaugmented coalgebras.
    It therefore follows from the previous proposition that the bar-cobar adjunction for $\sC$ is equivalent to the one for $\sC_{\unit // \unit}$.
    As a consequence, we can always assume the unit object is both initial and final when considering this adjunction, since this is the case for $\sC_{\unit // \unit}$.
\end{remark}

\subsubsection{A comparison with Lurie's bar-cobar adjunction}

In \cite[\S 5.2]{HA}, Lurie gives a different construction of the bar-cobar adjunction.
He shows that $\Map_{\coAlg^{\mathrm{aug}}}(\Bahr(A), Q)$ and $\Map_{\Alg^{\mathrm{aug}}(\sC)}(A, \Cobar(Q))$ are both naturally equivalent to the space of lifts of the pair $(A,Q)$ to an augmented algebra object in the twisted arrow category $\Tw(\sC)$.
It turns out this agrees with the construction we gave above.
Write $\Bahr' \colon \Alg^{\mathrm{aug}}(\sC) \to \coAlg^{\mathrm{aug}}(\sC)$ for the left adjoint coming from Lurie's construction.

\begin{proposition}
\label{prop: comparison-bar-cobar-adjunctions}
    Assume that the tensor product on $\sC$ preserves geometric realizations.
    Then the bar-cobar adjunction constructed in terms of the twisted arrow category is equivalent to the bar-cobar adjunction constructed by means of (co)endomorphism objects.
    In particular, there is a natural equivalence $\Bahr' \simeq \Bahr$.
\end{proposition}

\begin{remark}
    We expect that this proposition still holds without the assumption that the tensor product commutes with geometric realizations.
    However, this weaker statement suffices for our applications.
\end{remark}

\begin{proof}[Proof of \cref{prop: comparison-bar-cobar-adjunctions}]
    By \cref{rem: koszul-unit-initial-final}, we can assume that the unit object of $\sC$ is both initial and final.
    Let $A \in \Alg^{\mathrm{aug}}(\sC)$. In \cite[Proposition 3.34]{brantner2023pd} it is proved that $\Bahr'(A)$ is the coendomorphism object of $\triv_A(\unit) \in \RMod_A(\sC)$, so that we have an equivalence $\Bahr'(A) \simeq \Bahr(A)$ of coaugmented coalgebras.
    To prove the proposition, it will suffice to make this equivalence natural in $A$.
    For this, we will use the language of orthofibrations from \cref{appendix:orthofibrations}.

    Consider the following commutative diagram
    \[
    \begin{tikzcd}[sep = large]
        & \BMod(\sC) \ar[d] \ar[r, "\Bahr'"] & \LComod(\RMod(\sC)) \ar[d] \\
        \Alg(\sC) \ar[ur, "\theta"] \ar[r, "\Delta"] & \Alg(\sC) \times \Alg(\sC) \ar[r, "\Bahr' \times \id"] & \coAlg(\sC) \times \Alg(\sC).
    \end{tikzcd}
    \]
    Here $\theta$ is the functor that sends an algebra $A$ to itself considered as an $(A, A)$-bimodule and $\Delta$ is the diagonal functor.
    Note that $\Alg^{\mathrm{aug}}(\sC) \simeq \Alg(\sC)$ and similarly for coalgebras since the unit is both initial and final in $\sC$, so that the source and target of the lower right horizontal functor are indeed correct.
    The top horizontal functor $\Bahr'$ is constructed in \cite[Theorem 3.26]{brantner2023pd}.
    There they also prove that the right square in the diagram commutes.
    It is for this construction that we need that the tensor product on $\sC$ commutes with geometric realizations.
    This functor $\Bahr'$ sends an $(A, B)$-bimodule $M$ to a $(\Bahr'(A), B)$-mixed module with underlying object $\Bahr(\unit, A, M)$ (as in \cref{appendix:orthofibrations}, we use the term ``mixed module'' for an object with a compatible left coaction and right action).

    Observe that since $\Bahr(\unit, A, A) \simeq \unit$, the composite $\Bahr' \circ \theta$ lands in the full subcategory $\LComod(\RMod(\sC))^{\unit} \subseteq \LComod(\RMod(\sC))$ on those mixed modules $M$ for which the underlying object in $\sC$ is $\unit$.
    This can equivalently be described as the full subcategory on the objects $(Q, M, A)$ that project to $(Q, \triv_Q(\unit))$ and $(A, \triv_A(\unit))$ along the projections
    \[
    \begin{tikzcd}
    \LComod(\sC) & \LComod(\RMod(\sC)) \ar[l] \ar[r] & \RMod(\sC),
    \end{tikzcd}
    \]
    since $\unit$ is both initial and final in $\sC$.
    From this it is easy to see that the projection $\LComod(\RMod(\sC))^{\unit} \to \coAlg(\sC) \times \Alg(\sC)$
    is the orthofibration (\cref{definition:orthofibration}) corresponding to the functor $H \colon \coAlg(\sC) \times \Alg(\sC)^\op \to \Spc$ sending $(Q, A)$ to
    \[
        \{\triv_Q(\unit)\} \times_{\LComod_Q(\sC)} \LComod_Q(\RMod_A(\sC)) \times_{\RMod_A(\sC)} \{\triv_A(\unit)\}.
    \]
    
    Hence, it follows from \cref{cor: spaces-ortho-section-end} that the diagram above provides an element of the end 
    \[
    \int_{A \in \Alg(\sC)} H(A, \Bahr'(A)).
    \]
    By \cref{prop: parametrized-yoneda}, this element corresponds to a map
    \[
    \begin{tikzcd}
    \alpha \colon \Map_{\coAlg(\sC)}(\Bahr'(A), Q) \ar[r] & H(A, Q).
    \end{tikzcd}
    \]
    that is natural in $A$ and $Q$.
    It follows from \cite[Proposition 3.34]{brantner2023pd} that the composite $\Bahr' \circ \theta$ sends $A$ to the left $\Bahr'(A)$-comodule $K(A)$ with underlying object $\triv_A(\unit) \in \RMod_A(\sC)$ that witnesses the fact that $\Bahr'(A)$ is the coendomorphism object of $\triv_A(\unit)$.
    Unraveling the definitions, we see that the natural transformation $\alpha$ sends a coalgebra morphism $f \colon \Bahr'(A) \to Q$ to the pushforward of $K(A)$ along $f$.
    Hence $\alpha$ is an equivalence by the universal property of the coendomorphism object.

    By construction of the functor $\Bahr$ (\cref{cor: bar-cobar-adjunction}), we have a natural equivalence $H(A, Q) \simeq \Map_{\coAlg(\sC)}(\Bahr(A), Q)$.
    Therefore, $\alpha$ provides an equivalence
    \[
        \Map_{\coAlg(\sC)}(\Bahr'(A), Q) \simeq \Map_{\coAlg(\sC)}(\Bahr(A), Q),
    \]
    natural in $A$ and $Q$.
    We are now done by an application of the Yoneda lemma.
\end{proof}

A useful fact that already made an appearance in the proof of the previous proposition is that the bar-cobar adjunction can be lifted to an adjunction on the level of (co)modules.
Let 
\[
\LMod^{\mathrm{aug}}(\sC) \coloneqq \LMod(\sC) \times_{\Alg(\sC)} \Alg^{\mathrm{aug}}(\sC),
\]
where we take the pullback along the forgetful functor $\Alg^{\mathrm{aug}}(\sC) \to \Alg(\sC)$.
Write $p \colon \LMod^{\mathrm{aug}}(\sC) \to \Alg^{\mathrm{aug}}(\sC)$ for the projection, which is a cartesian fibration.
Similarly, we obtain a cocartesian fibration $q \colon \LComod^{\mathrm{aug}}(\sC) \to \coAlg^{\mathrm{aug}}(\sC)$.

\begin{proposition} \label{prop: koszul-duality-modules}
    Suppose that $\sC$ is a monoidal category such that the tensor product preserves geometric realizations.
    Then there is a diagram
    \[
    \begin{tikzcd}[sep = huge]
    \LMod^{\mathrm{aug}}(\sC) \ar[r, shift left, "\Bahr"] \ar[d, "p"] & \LComod^{\mathrm{aug}}(\sC) \ar[l, shift left, "\Cobar"] \ar[d, "q"] \\
    \Alg^{\mathrm{aug}}(\sC) \ar[r, shift left, "\Bahr"] & \coAlg^{\mathrm{aug}}(\sC) \ar[l, shift left, "\Cobar"],
    \end{tikzcd}
    \]
    where both rows are adjunctions, such that if we restrict to the square containing just the left adjoints or just the right adjoints, the diagram commutes.
    If $A$ is an augmented algebra, then $\Bahr$ sends a left $A$-module $M$ to a left $\Bahr(A)$-comodule with underlying object $\indec_A(M)$.
    Similarly, if $Q$ is a coaugmented coalgebra, then $\Cobar$ sends a left $Q$-comodule $N$ to a left $\Cobar(Q)$-module with underlying object $\prim_Q(N)$.

    Moreover, $\Bahr \colon \LMod^{\mathrm{aug}}(\sC) \to \LComod^{\mathrm{aug}}(\sC)$ is a right $\sC$-linear functor and sends $p$-cocartesian edges to $q$-cocartesian edges.
\end{proposition}
\begin{proof}
    This proposition holds for the bar-cobar adjunction defined in terms of twisted arrow categories by \cite[Theorem 3.26]{brantner2023pd}.
    The result therefore follows from \cref{prop: comparison-bar-cobar-adjunctions}.
\end{proof}

If $M$ is a left $A$-module, this proposition gives a canonical $\Bahr(A)$-comodule structure on $\indec_A(M)$.
We call this comodule the \emph{Koszul dual} of $M$.

\subsection{Koszul duality for symmetric sequences}
\label{sec: koszul-symm-seq}

Let $\sC$ be a stable presentable category.
We now apply the theory of Koszul duality from the previous section to the category $\SymFunL(\sC, \sC)$ with the composition product (see \cref{def: composition product}).
We restrict our attention to the following subcategory.

\begin{definition}
    A symmetric sequence $F \colon \Sym(\sC) \to \sC$ is \emph{strongly positive} if $F_0 = 0$ and $F_1 = \id_\sC$.
    Let $\SymFunL_+(\sC, \sC)$ denote the full subcategory of $\SymFunL(\sC, \sC)$ on the strongly positive symmetric sequences.
    Since this subcategory is closed under the composition product and contains the unit, it inherits a monoidal structure.
\end{definition}

\begin{remark}
A symmetric sequence that satisfies the analogue of the strongly positive condition is often called reduced in the literature.
We have departed from this convention since it clashes with the terminology from Goodwillie calculus, where a functor $F \colon \sC \to \sC$ is called reduced if $P_0 F = *$.
\end{remark}

We write $\unit \in \SymFunL_+(\sC, \sC)$ for the unit symmetric sequence, which is $\id_\sC$ in arity $1$ and vanishes in all other arities.

\begin{remark}
For any algebra $A \in \Alg(\SymFunL_+(\sC, \sC))$, the unit map $\unit \to A$ is automatically given by the inclusion of the arity $1$ component.
Moreover, $A$ admits a unique augmentation, given by the projection $A \to \unit$ onto the arity $1$ component.
Similar remarks apply to the counit and coaugmentation of any coalgebra in $\SymFunL_+(\sC, \sC)$.
\end{remark}

The goal of this section is to prove the following theorem.

\begin{theorem} \label{thm: sseq-bar-cobar-equiv}
The bar-cobar adjunction gives an adjoint equivalence
\[
\begin{tikzcd}[sep = large]
\Alg^{\mathrm{aug}}(\SymFunL_+(\sC, \sC)) \ar[r, shift left, "\Bahr", "\sim"'] &  \coAlg^{\mathrm{aug}}(\SymFunL_+(\sC, \sC)). \ar[l, shift left=2, "\Cobar"]
\end{tikzcd}
\]
\end{theorem}

The argument we will give essentially goes back to \cite[Theorem 2.18]{ChingBarCobar}, where it is used to prove the analogous result for ordinary symmetric sequences.
We learned this version of the proof from Heuts \cite[Theorem 3.4]{heuts2024koszulduality}.

We start by establishing a finiteness property of cobar constructions in the category $\SymFunL_+(\sC)$, which we will use in the proof. 
Recall that a cosimplicial object is called \emph{$n$-coskeletal} if it is right Kan extended from the subcategory $\Delta_{\leq n} \subset \Delta$.

\begin{proposition} \label{prop: cobar-symseq-finite}
    Let $Q \in \coAlg^{\mathrm{aug}}(\SymFunL_+(\sC, \sC))$ and let $M$ be a right and $N$ a left $Q$-comodule in $\SymFunL(\sC, \sC)$ such that $M_0 = N_0 = 0$.
    Then the arity $n$ component of the cobar construction $C^\bullet(M, Q, N)$ is $(n-1)$-coskeletal.
    In particular, the cobar construction is arity-wise a finite totalization.
\end{proposition}
\begin{proof}
    Write $C^\bullet$ for $C^\bullet(M, Q, N)$.
    It suffices to show the fiber of $\Tot^n C^\bullet \to \Tot^{n-1} C^\bullet$ vanishes through arity $n$ (see for example \cite[Remark 2.8]{barwick2022ktheorypolynomialfunctors}).
    There is a standard formula for this fiber given by $\Omega^n \tfib \mathcal{X}_n$, where $\mathcal{X}_n$ is an $n$-cube sending a subset $S \subset \underline{n} = \{1, \ldots, n\}$ to  
    \[
    M \circ Q^{\circ \underline{n} \setminus S} \circ N.
    \]
    The maps in the cube are given by applying the counit $\epsilon \colon Q \to \unit$ to the various factors.
    For instance, in the case $n=2$ we get the square
    \[
    \begin{tikzcd}[sep = large]
       M \circ Q^{\circ 2} \circ N
        \ar[r, "\epsilon \circ \id"] \ar[d, "\id \circ \epsilon"] & 
        M \circ Q \circ N 
        \ar[d, "\epsilon"] \\
        M \circ Q \circ N 
        \ar[r, "\epsilon"] & 
        M \circ N
    \end{tikzcd}
    \]
    It suffices to show this total fiber vanishes through arity $n$.
    Let $\mathcal{X}'_n$ be the $n$-cube sending $S \subset \langle n \rangle$ to
    \[
    M \circ Q^{\circ \underline{n} \setminus S},
    \]
    and having the same maps as before, so that $\mathcal{X}_n = \mathcal{X}'_n \circ N$.
    Since $N_0 = 0$ and the composition product is exact in the first variable, it suffices to show the total fiber of $\mathcal{X}'_n$ vanishes through arity $n$.
    We prove this by induction on $n$.
    For $n=0$ there is nothing to prove, since all symmetric sequences involved vanish in arity $0$.
    Suppose the statement holds for a given $n \geq 0$.
    By the inductive formula for the total fiber of a cube and the fact that the composition product is exact in the first variable, we obtain
    \[
    \begin{tikzcd}
    \tfib \mathcal{X}'_{n+1} \simeq \fib (\tfib \mathcal{X}'_n \circ Q \ar[r, "\epsilon"] & \tfib \mathcal{X}'_n)
    \end{tikzcd}
    \]
    By hypothesis, $\tfib \mathcal{X}'_n$ vanishes through arity $n$.
    Since all symmetric sequences involved vanish in arity 0, it immediately follows from the formula for the composition product and the fact that $\epsilon$ is an isomorphism in arity one that $\tfib \mathcal{X}'_{n+1}$ vanishes through arity $n+1$.
\end{proof} 

It is possible to truncate a symmetric sequence at a given arity, and this operation will also play a role in the proof.
For $n \geq 1$, let $\SymFunL_{+, \leq n}(\sC, \sC)$ denote the full subcategory of $\SymFunL_+(\sC, \sC)$ spanned by symmetric sequences $F$ for which $F_k = 0$ whenever $k > n$.
The obvious truncation functor $\tau_{\leq n} \colon \SymFunL_+(\sC, \sC) \to \SymFunL_{+, \leq n}(\sC, \sC)$ is both left and right adjoint to the inclusion, making it both a reflective localization and colocalization.
It is easily seen to be compatible with the composition product in the sense of \cite[Definition 2.2.1.6]{HA}.
This equips $\SymFunL_{+, \leq n}(\sC, \sC)$ with a monoidal structure for which $\tau_{\leq n}$ is strong monoidal and the inclusion $\SymFunL_{+, \leq n}(\sC, \sC) \to \SymFunL_+(\sC, \sC)$ admits both a lax monoidal and oplax monoidal structure.
It follows that for any $A \in \Alg^{\mathrm{aug}}(\SymFunL_+(\sC, \sC))$ we obtain an algebra structure on $\tau_{\leq n} A$ together with an algebra morphism $A \to \tau_{\leq n}A$ which is an equivalence in arities $\leq n$.
Similarly, for any augmented coalgebra $Q$ we obtain a coalgebra structure on $\tau_{\leq n}Q$ together with a coalgebra map $\tau_{\leq n}Q \to Q$.

\begin{proof}[Proof of \cref{thm: sseq-bar-cobar-equiv}]
    We will show that both the unit and counit of the bar-cobar adjunction are equivalences.

    Let $A \in \Alg^{\mathrm{aug}}(\SymFunL_+(\sC, \sC))$ and write $\eta \colon A \to \Cobar(\Bahr(A))$ for the unit.
    By \cref{prop: koszul-duality-modules} we have an adjunction
    \[
    \begin{tikzcd}[sep = huge]
    \RMod_A(\SymFunL(\sC, \sC)) \ar[r, shift left, "\indec_A"] 
    & \RComod_{BA}(\SymFunL(\sC, \sC)) \ar[l, shift left, "\prim_{BA}"],
    \end{tikzcd}
    \]
    where $\indec_A(M) = \Bahr(M, A, \unit)$ and $\prim_{BA}(N) = \Cobar(N, BA, \unit)$.
    Consider the unit of this adjunction at the free right $A$-module $A$:
    \[
    \begin{tikzcd}
        A \ar[r] & \Cobar(\Bahr(A, A, \unit), \Bahr(A), \unit)) \simeq \Cobar(\Bahr(A)).
    \end{tikzcd}
    \]
    It is straightforward to verify that the underlying map of this unit is equivalent to the underlying map of $\eta$.
    So it is enough to show that this unit map in right $A$-modules is an equivalence.

    Observe that $\lim_n \Bahr(\tau_{\leq n} A, A, \unit) \simeq \Bahr(A, A, \unit)$, since in every arity this limit eventually becomes constant.
    It follows that the unit $A \to \prim_{BA}\indec_A(A)$ is the inverse limit of the unit maps at $\tau_{\leq n} A$.
    Therefore, it suffices to show that the unit map $\tau_{\leq n} A \to \prim_{BA} \indec_A (\tau_{\leq n} A)$ is an equivalence for all $n \geq 1$.
    We will prove this by induction on $n$.

    Since the composition product commutes with all limits and colimits in the left variable, the categories $\RMod_A(\SymFunL(\sC, \sC))$ and $\RComod_{BA}(\SymFunL(\sC, \sC))$ admit all limits and colimits, and these are preserved and reflected by the forgetful functors \cite[Corollaries 4.2.3.3 and 4.2.3.5]{HA}.
    In particular, both categories are stable, and the functors $\indec_A$ and $\prim_{BA}$ are exact.
    
    Consider the following diagram of right $A$-modules
    \[
    \begin{tikzcd}[sep = large]
        A_n \ar[r] \ar[d] & \Cobar(\Bahr(A_n, A, \unit), \Bahr(A), \unit)) \ar[d] \\
        \tau_{\leq n} A \ar[r] \ar[d] & \Cobar(\Bahr(\tau_{\leq n}A, A, \unit), \Bahr(A), \unit)) \ar[d] \\
        \tau_{\leq n-1} A \ar[r] & \Cobar(\Bahr(\tau_{\leq n - 1}A, A, \unit), \Bahr(A), \unit)).
    \end{tikzcd}
    \]
    Both columns are fiber sequences and all horizontal maps are components of the unit of the adjunction between $\indec_A$ and $\prim_{BA}$.

    The bottom horizontal arrow is an equivalence by induction, so in order to show the middle arrow is an equivalence, it suffices to show the top one is.
    The fact that the right $A$-module $A_n$ is concentrated in arity $n$ implies that it is a trivial module for degree reasons; we give a proof of this statement in \cref{lem: degree-trivial-module} below.
    Since the composition product preserves geometric realizations in the right variable by \cref{prop: exactness-composition-product}, triviality of $A_n$ implies that $\Bahr(A_n, A, \unit) \simeq A_n \circ \Bahr(A)$.
    The unit map at $A_n$ is now easily seen to come from the split augmented cosimplicial diagram $A_n \to A_n \circ C_\bullet(\Bahr(A), \Bahr(A), \unit).$
    Hence it is an equivalence.

    The proof that for any augmented coalgebra $Q$ the counit $\Bahr(\Cobar(Q)) \to Q$ is an equivalence goes exactly the same.
    Only the final step needs some extra care: to argue that $\Cobar(Q_n, Q, \unit) \simeq Q_n \circ \Cobar(Q)$, we use that $\Cobar(Q)$ is arity-wise a finite totalization by \cref{prop: cobar-symseq-finite} and that the composition product preserves these kinds of totalizations in the right variable by \cref{prop: exactness-composition-product}. 
\end{proof}

We used the following lemma in the proof.

\begin{lemma}
    \label{lem: degree-trivial-module}
    Let $A \in \Alg^{\mathrm{aug}}(\SymFunL_+(\sC, \sC))$ and let $M \in \RMod_A(\SymFunL(\sC, \sC))$ be a module concentrated in a single arity $n \geq 1$.
    Then $M$ is a trivial module.
\end{lemma}
\begin{proof}
    Let $\SymFunL_{\leq n}(\sC, \sC) \subset \SymFunL(\sC, \sC)$ be the full subcategory on the $n$-truncated symmetric sequences that are not necessarily equal to $\id_\sC$ in arity $1$.
    Just as was the case for $\SymFunL_{+, \leq n}(\sC, \sC)$, this category inherits a monoidal structure from $\SymFunL(\sC, \sC)$.
    The module $M$ lies in the image of the fully faithful inclusion 
    \[
    \begin{tikzcd}
    \RMod_{\tau_{\leq n} A}(\SymFunL_{\leq n}(\sC, \sC)) \ar[r, hook] & \RMod_A(\SymFunL(\sC, \sC)),
    \end{tikzcd}
    \] 
    and this inclusion preserves trivial modules, so it suffices to show $M$ is trivial in $\RMod_{\tau_{\leq n} A}(\SymFunL_{\leq n}(\sC, \sC))$.
    Since $M$ is concentrated in arity $n$, the formula for the composition product in $\SymFunL_{\leq n}(\sC, \sC)$ implies that the unit $\unit \to \tau_{\leq n} A$ induces an equivalence $M \xrightarrow{\sim} M \circ \tau_{\leq n}A$.
    This morphism therefore exhibits $M$ as a free $\tau_{\leq n}A$-module, with generators given by $M$ as well. 
    The universal property of free modules now allows us to construct a map $M \to \triv(M)$ of right $\tau_{\leq n}A$-modules that is the identity on underlying objects, and hence an equivalence.
 \end{proof}

\subsection{The Koszul dual of \texorpdfstring{$\partial_*{\id}$}{∂(id)}}
\label{sec: koszul-dual-partial-id}

We will now explain how the theory of Koszul duality ties in with Goodwillie calculus.
Let $\sC$ be a differentiable category.
Then $\id_\sC$ is a monad on $\sC$, and $\Sigma^\infty_\sC \Omega^\infty_\sC$ is a reduced finitary comonad on $\Sp(\sC)$.
By lax monoidality of derivatives and the stable chain rule, $\partial_*{\id_\sC}$ is an algebra and $\partial_*(\Sigma^\infty_\sC \Omega^\infty_\sC)$ is a coalgebra in $\SymFunL(\Sp(\sC), \Sp(\sC))$.
The main theorem of this section is:
\begin{theorem} \label{thm: koszul-dual-of-id}
    There is an equivalence of algebras
    \[
    \partial_*{\id_\sC} \simeq \Cobar \partial_*(\Sigma^\infty_\sC \Omega^\infty_\sC).
    \]
    In other words, $\partial_*{\id_\sC}$ is the Koszul dual algebra of $\partial_*(\Sigma_\sC^\infty \Omega_\sC^\infty)$.
\end{theorem}

In order to construct a map of algebras $\partial_*{\id_\sC} \to \Cobar\partial_*(\Sigma^\infty_\sC \Omega^\infty_\sC)$, it suffices by \cref{thm: koszul-left-module-right-comodule-structure-on-unit} to provide the unit symmetric sequence $\unit$ with a compatible left $\partial_*(\Sigma^\infty_\sC \Omega^\infty_\sC)$-comodule and right $\partial_*(\id_\sC)$-module structure.
In fact, we can already construct this structure on the level of functors.

\begin{construction}\label{construction:Koszul-duality-map}
The functor $\Sigma^\infty_\sC$ canonically has the structure of a left $\Sigma^\infty_\sC \Omega^\infty_\sC$-comodule.
As $\id_\sC$ is the unit in $\Fun(\sC, \sC)$, there is an equivalence
\[
\LComod_{\Sigma^\infty_\sC \Omega^\infty_\sC}(\Fun(\sC, \Sp(\sC))) \simeq \RMod_{\id_\sC}(\LComod_{\Sigma^\infty_\sC \Omega^\infty_\sC}(\Fun(\sC, \Sp(\sC)))).
\]
We consider $\Sigma^\infty_\sC$ as an object in the category on the right.
It follows from the lax monoidality of derivatives and the stable chain rule combined with \cref{lem: strong-lax-induced-map-mix-mod} that
\[
\partial_*\Sigma^\infty_\sC \in \RMod_{\partial_*(\id_\sC)}(\LComod_{\partial_*(\Sigma^\infty_\sC \Omega^\infty_\sC)}(\SymSeq(\Sp(\sC)))).
\]
Since $\Sigma^\infty_\sC$ is a linear functor, we have $\partial_*\Sigma^\infty_\sC \simeq \unit$. Therefore, this provides a compatible left $\partial_*(\Sigma^\infty_\sC \Omega^\infty_\sC)$-comodule and right $\partial_*(\id_\sC)$-module structure on $\unit$.
We write 
\[
\begin{tikzcd}
\alpha \colon \partial_*{\id_\sC} \ar[r] & \Cobar\partial_*(\Sigma^\infty_\sC \Omega^\infty_\sC)
\end{tikzcd}
\]
for the induced algebra morphism.
\end{construction}

Our goal is to show $\alpha$ is an equivalence.
Recall that for a differentiable category $\sC$, the \emph{Bousfield--Kan resolution} of $\id_\sC$ is a coaugmented cosimplicial object
\[
\begin{tikzcd}
    \id_\sC \ar[r] & \Omega_\sC^\infty \Sigma_\sC^\infty \ar[r, shift left] \ar[r, shift right] & \Omega_\sC^\infty(\Sigma_\sC^\infty \Omega_\sC^\infty)\Sigma_\sC^\infty \ar[l, shorten=0.4em] \ar[r, shift left=2] \ar[r] \ar[r, shift right=2] & \Omega_\sC^\infty(\Sigma_\sC^\infty \Omega_\sC^\infty)^2\Sigma_\sC^\infty \cdots. \ar[l, shift left, shorten=0.4em] \ar[l, shift right, shorten=0.4em]
\end{tikzcd}
\]
The coface and codegeneracy maps in this diagram use the unit and the counit of the adjunction between $\Omega_\sD^\infty$ and $\Sigma_\sD^\infty$.
This resolution has the property that applying $\Sigma_\sC^\infty$ to it yields the canonical split cobar resolution of $\Sigma_\sC^\infty$ by cofree $\Sigma_\sC^\infty \Omega_\sC^\infty$-comodules.
It should be mentioned that the construction of the Bousfield--Kan resolution is completely formal and can be carried out for any adjunction.

We will show below that the underlying map of $\alpha$ is given by taking derivatives of the Bousfield--Kan resolution and applying the stable chain rule.
The proof that $\alpha$ is an equivalence then amounts to the following crucial observation due to Arone and Ching.
\begin{proposition} \label{prop: Arone-Ching-Bousfield-Kan}
    Let $\sC \xrightarrow{F} \sD \xrightarrow{G} \sE$ be a pair of reduced finitary functors between differentiable categories.
    Then the map
    \[
    \begin{tikzcd}
        \partial_*(GF) \ar[r] & \Tot \partial_*(G \Omega_\sD^\infty (\Sigma_\sD^\infty \Omega_\sD^\infty)^{\bullet} \Sigma_\sD^\infty F)
    \end{tikzcd}
    \]
    obtained by applying $\partial_*$ to the Bousfield--Kan resolution is an equivalence.
\end{proposition}
\begin{proof}
This was originally proved as \cite[Theorem 4.1.1]{AroneChing2011}.
A proof in the $\infty$-categorical setting was given as \cite[Proposition B.4]{Heuts-Thesis}.
We will sketch the argument.
Working one arity at a time, we can reduce by \cref{lem: Arone-Ching-lemma} to the case where $F$ is $n$-excisive.
Inductively going up the Goodwillie tower of $F$ and using the fact that the functor $F \mapsto \partial_*(FG)$ is left exact, we can further reduce to the case where $F$ is homogeneous.
Then $F \simeq H \Sigma_\sD^\infty$ for some functor $H$.
The cosimplicial object now has the form
\[
\partial_*(H\Sigma_\sD^\infty \Omega_\sD^\infty (\Sigma_\sD^\infty \Omega_\sD^\infty)^{\bullet} \Sigma_\sD^\infty G),
\]
and the counit map $\epsilon \colon \Sigma_\sD^\infty \Omega_\sD^\infty \to \id$ in the leftmost factor provides an extra codegeneracy.
This completes the proof.
\end{proof}

What remains is to show that the underlying map of $\alpha$ indeed comes from the Bousfield--Kan resolution.
We start by recalling how this resolution is constructed.
The adjunction between $\Sigma_\sC^\infty$ and $\Omega_\sC^\infty$ factors in the standard way as a composite of adjunctions (left adjoints on top):
\[
\begin{tikzcd}[sep = large]
\sC \ar[r, shift left, "\Sigma_\sC^\infty"] & 
\LComod_{\Sigma_\sC^\infty \Omega_\sC^\infty}(\Sp(\sC)) \ar[r, shift left, "\forget"] \ar[l, shift left, "R"] &
\Sp(\sC) \ar[l, shift left, "\cofree"].
\end{tikzcd}
\]
We commit a slight abuse of notation by denoting the leftmost left adjoint in this diagram by $\Sigma_\sC^\infty$.
It is the functor taking an object $x \in \sC$ to $\Sigma_\sC^\infty x$ with the canonical $\Sigma_\sC^\infty \Omega_\sC^\infty$-coalgebra structure.
For an arbitrary $M \in \LComod_{\Sigma_\sC^\infty \Omega_\sC^\infty}(\Sp(\sC))$, we can write $M$ as the totalization of its canonical cobar resolution by cofree comodules: $M \simeq \Tot C^\bullet(\Sigma_\sC^\infty \Omega_\sC^\infty, \Sigma_\sC^\infty \Omega_\sC^\infty, M).$
Since $R \circ \cofree \simeq \Omega_\sC^\infty$ and $R$ preserves limits, we obtain the formula
\[
RM \simeq \Tot C^\bullet(\Omega_\sC^\infty, \Sigma_\sC^\infty\Omega_\sC^\infty, M). 
\]
The unit of the adjunction between $\Sigma_\sC^\infty$ and $R$ then provides the coaugmentation
\[
\begin{tikzcd}
\id_\sC \ar[r] & C^{\bullet}(\Omega_\sC^\infty, \Sigma_\sC^\infty\Omega_\sC^\infty, \Sigma_\sC^\infty),
\end{tikzcd}
\]
and this is precisely the Bousfield--Kan resolution of $\id_\sC$.

\begin{lemma} \label{lem: coassembly-lemma}
    The underlying map of $\alpha \colon \partial_*{\id_\sC} \to \Cobar \partial_*(\Sigma_\sC^\infty \Omega_\sC^\infty)$, which was defined in \cref{construction:Koszul-duality-map}, is given by the composite
    \[
    \begin{tikzcd}
    \partial_*{\id_\sC} \ar[r] & \Tot \partial_*(\Omega_\sD^\infty (\Sigma_\sD^\infty \Omega_\sD^\infty)^{\bullet} \Sigma_\sD^\infty) \ar[r, "\simeq"] & \Tot C^\bullet(\partial_*\Omega_\sD^\infty, \partial_*(\Sigma_\sD^\infty \Omega_\sD^\infty), \partial_*\Sigma_\sD^\infty),
    \end{tikzcd}
    \]
    where the first map is given by applying $\partial_*$ to the Bousfield--Kan resolution and the second map comes from the stable chain rule.
\end{lemma}
\begin{proof}
The totalization of the Bousfield--Kan resolution is the underlying functor of the monad $R\Sigma_\sC^\infty$. Therefore, by \cref{ex: comonad-adjunction}, it is the endomorphism object of $\Sigma_\sC^\infty$, where the latter is regarded as an object in $\LComod_{\Sigma_\sC^\infty \Omega_\sC^\infty}(\Fun(\sC, \Sp(\sC)))$:
\[
\End(\Sigma_\sC^\infty) = \Tot C^\bullet(\Omega_\sC^\infty, \Sigma_\sC^\infty\Omega_\sC^\infty, \Sigma^\infty_\sC).
\]
By \cref{prop: naturality-endomorphism-objects}, there is a canonical algebra morphism  
\[
\begin{tikzcd}
\partial_*C(\Omega_\sC^\infty, \Sigma_\sC^\infty \Omega_\sC^\infty, \Sigma_\sC^\infty) = \partial_*\End(\Sigma_\sC^\infty) \ar[r, "\beta"] & \End(\partial_*\Sigma_\sC^\infty) = \Cobar \partial_*(\Sigma_\sC^\infty \Omega_\sC^\infty).
\end{tikzcd}
\]
The map $\alpha$ is then given by composing $\beta$ with the map $\partial_*{\id_\sC} \to \partial_*\End(\Sigma_\sC^\infty)$ coming from the Bousfield--Kan resolution.
It now suffices to show that the underlying map of $\beta$ is the coassembly map for the totalization of the Bousfield--Kan resolution.
Note that the underlying map of $\beta$ classifies the action
\[
\begin{tikzcd}
\partial_*(\Sigma^\infty_\sC) \circ \partial_*(\End(\Sigma_\sC^\infty)) \ar[r, "\mathrm{lax}"] & \partial_*(\Sigma_\sC^\infty \circ \End(\Sigma_\sC^\infty)) \ar[r] & \partial_*\Sigma_\sC^\infty,
\end{tikzcd}
\]
where the second map is given by applying derivatives to the action of $\End(\Sigma_\sC^\infty)$ on $\Sigma_\sC^\infty$.
Unraveling the definitions, we see it suffices to show that this action is homotopic to the composite
\[
\begin{tikzcd}
    \partial_*(\Sigma^\infty_\sC) \circ \partial_*(\End(\Sigma_\sC^\infty)) \ar[r, "\id \circ \gamma"] & \partial_*(\Sigma^\infty_\sC) \circ \End(\partial_*\Sigma_\sC^\infty) \ar[r] & \partial_*\Sigma^\infty_\sC,
\end{tikzcd}
\]
where $\gamma$ is the coassembly map and the second map is the canonical action.
This can readily be deduced from the following three observations
\begin{enumerate}
    \item The action of $\End(\Sigma_\sC^\infty)$ on $\Sigma_\sC^\infty$ is given by the composite
    \[
    \begin{tikzcd}
    \Sigma_\sC^\infty \circ \End(\Sigma_\sC^\infty) \ar[r] & \Sigma_\sC^\infty \Omega_\sC^\infty \Sigma_\sC^\infty \ar[r, "\epsilon"] & \Sigma_\sC^\infty,
    \end{tikzcd}
    \]
    where the first map is induced by the map of cosimplicial diagrams 
    \[
    \begin{tikzcd}
    C^\bullet(\Omega_\sC^\infty, \Sigma_\sC^\infty \Omega_\sC^\infty, \Sigma_\sC^\infty) \ar[r, "\tau"] & C^\bullet(\Omega_\sC^\infty, \id_{\Sp(\sC)}, \Sigma_\sC^\infty),
    \end{tikzcd}
    \] 
    obtained by applying the counit $\epsilon \colon \Sigma_\sC^\infty \Omega_\sC^\infty \to \id_{\Sp(\sC)}$ in the middle factor.
    \item The action of $\End(\partial_*\Sigma_\sC^\infty)$ on $\partial_*\Sigma_\sC^\infty$ is induced by the map $\partial_*(\tau)$.
    \item For an arbitrary functor $F \colon \sC \to \sC$, the composite 
    \[
    \begin{tikzcd}
    \partial_*(\Sigma_\sC^\infty) \circ \partial_*(\Omega_\sC^\infty) \circ \partial_*(F) \ar[r, "\mathrm{lax}"] & \partial_*(\Sigma_\sC^\infty \Omega_\sC^\infty F) \ar[r, "\epsilon"] & \partial_* F
    \end{tikzcd}
    \]
    is an equivalence. \qedhere
\end{enumerate}
\end{proof}

We have now proved the main theorem of this section.

\begin{proof}[Proof of \cref{thm: koszul-dual-of-id}]
    Combine \cref{prop: Arone-Ching-Bousfield-Kan} and \cref{lem: coassembly-lemma}.
\end{proof}

\subsection{The chain rule}
\label{sec: chain-rule}

In this section we prove the chain rule, which expresses the derivatives $\partial_*(FG)$ of the composite of two functors in terms of the derivatives of $\partial_*(F)$ and $\partial_*(G)$.

\begin{theorem}[Chain rule] \label{thm: unstable-chain-rule}
    Let $\sC \xrightarrow{G} \sD \xrightarrow{F} \sE$ be a pair of reduced finitary functors between differentiable categories. Then there is an equivalence
    \[
    \partial_*F \circ_{\partial_*{\id_\sD}} \partial_*G \simeq \partial_*(FG).
    \]
    of $(\partial_*{\id_\sE}, \partial_*{\id_\sC})$-bimodules.
\end{theorem}

Again, we start by constructing a comparison map.

\begin{construction}\label{construction:The-lax-comparison-map}
Let $\sC \xrightarrow{G} \sD \xrightarrow{F} \sE$ be a pair of reduced finitary functors between differentiable categories.
The composite $FG$ is the colimit of the bar construction $B_\bullet(F, \id_\sD, G)$, which is a constant simplicial diagram.
Taking derivatives and using the lax monoidal structure, we obtain a map of simplicial diagrams.
On colimits, this induces
\[
\begin{tikzcd}
\phi \colon \partial_*F \circ_{\partial_*{\id_\sD}} \partial_*G \ar[r] & \partial_*(FG),
\end{tikzcd}
\]
which is a map of $(\partial_*{\id_\sE}, \partial_*{\id_\sC})$-bimodules.
\end{construction}

We will prove \cref{thm: unstable-chain-rule} by showing $\phi$ is an equivalence.
The proof relies on strong finiteness properties of the derivatives of the Bousfield--Kan resolution, which we will first establish.

\begin{proposition}[cf.\ {\cite[Remark 4.1.3]{AroneChing2011}}] \label{prop: bousfield-kan-finite}
    Let $\sC \xrightarrow{G} \sD \xrightarrow{F} \sE$ be a pair of composable reduced finitary functors between differentiable categories.
    Then the cosimplicial object
    \[
    \begin{tikzcd}
        \partial_n(F\Omega_\sD^\infty \Sigma_\sD^\infty G) \ar[r, shift left] \ar[r, shift right] & \partial_n(F(\Omega_\sD^\infty \Sigma_\sD^\infty)^2G) \ar[l, shorten=0.4em] \ar[r, shift left=2] \ar[r] \ar[r, shift right=2] & \partial_n(F(\Omega_\sD^\infty \Sigma_\sD^\infty)^3 G) \cdots \ar[l, shift left, shorten=0.4em] \ar[l, shift right, shorten=0.4em]
    \end{tikzcd}
    \]
    obtained by applying $\partial_n$ to the Bousfield--Kan resolution is $(n-1)$-coskeletal.
\end{proposition}
\begin{proof}
    Applying $\partial_*$ to the Bousfield--Kan resolution and using the stable chain rule, we obtain the cobar construction 
    \[
    C^\bullet(\partial_*(F\Omega_{\sD}^\infty), \partial_*(\Sigma_{\sD}^\infty \Omega_{\sD}^\infty), \partial_*(\Sigma_{\sD}^\infty G)).
    \]
    By \cref{example:derivatives-Sigma-Loops}, we see that $\partial_*(\Sigma^\infty_\sC \Omega^\infty_\sC)$ lies in $\SymFunL_+(\Sp(\sC), \Sp(\sC))$, so the result follows from \cref{prop: cobar-symseq-finite}.
\end{proof}

We can now prove the chain rule.

\begin{proof}[Proof of \cref{thm: unstable-chain-rule}]
    Resolving both $F$ and $G$ by their Bousfield--Kan resolutions, it follows from \cref{prop: Arone-Ching-Bousfield-Kan} and the stable chain rule that $\partial_*F \circ_{\partial_*{\id_\sD}} \partial_*G$ is equivalent to
\[
        \Tot \left [ \partial_*(F\Omega_\sD^\infty (\Sigma_\sD^\infty \Omega_\sD^\infty)^{\bullet}) \circ \partial_*\Sigma_\sD^\infty \circ_{\partial_*{\id_\sD}} \partial_*\Omega_\sD^\infty \circ \partial_*((\Sigma_\sD^\infty \Omega_\sD^\infty)^{\bullet} \Sigma_\sD^\infty G) \right ].
\]
Here we use that the Bousfield--Kan resolution is a finite totalization in every arity (\cref{prop: bousfield-kan-finite}), that the composition product commutes with finite totalizations in both variables (\cref{prop: exactness-composition-product}), and that the diagonal $\Delta \to \Delta \times \Delta$ is left cofinal.
By naturality of the map $\phi \colon \partial_*F \circ_{\partial_*{\id_\sD}} \partial_*G \to \partial_*(FG)$, we can therefore reduce to the case $F = \Sigma_\sD^\infty$ and $G = \Omega_\sD^\infty$.
But in that case we have
\[
\partial_*\Sigma_\sD^\infty \circ_{\partial_*{\id_\sD}} \partial_*\Omega_\sD^\infty \simeq \Bahr(\partial_*{\id_\sD}) \simeq \partial_*(\Sigma_\sD^\infty \Omega_\sD^\infty).
\]
by \cref{thm: sseq-bar-cobar-equiv} and \cref{thm: koszul-dual-of-id}.
One can see that this equivalence is given by the map $\phi$ by a diagram chase similar to the one in the proof of \cref{lem: coassembly-lemma}.
\end{proof}

\begin{corollary}
    Let $\sC$ and $\sD$ be differentiable categories.
    Then the following diagram commutes
    \[
    \begin{tikzcd}[sep = large]
        \Fun^{\ast, \omega}(\sC, \sD) \ar[r, "\partial_*"] \ar[d, "\Sigma_\sD^\infty \circ -"] & \LMod_{\partial_*{\id_\sD}}(\SymFunL(\Sp(\sC), \Sp(\sD))) \ar[d, "\indec_{\partial_*{\id_\sD}}"] \\
        \LComod_{\Sigma_\sD^\infty \Omega_\sD^\infty}(\Fun^{\ast, \omega}(\sC, \Sp(\sD))) \ar[r, "\partial_*"] & \LComod_{\partial_*(\Sigma_\sD^\infty \Omega_\sD^\infty)}(\SymFunL(\Sp(\sC), \Sp(\sD))),
    \end{tikzcd}
    \]
    where the right vertical arrow is the Koszul duality functor from \cref{prop: koszul-duality-modules}.
\end{corollary}
\begin{proof}
    By the chain rule, we have a natural equivalence of left $\partial_*(\Sigma_\sD^\infty\Omega_\sD^\infty)$-comodules
    \[
    \partial_*(\Sigma_\sD^\infty F) \simeq \partial_*\Sigma_\sD^\infty \circ_{\partial_*{\id_\sD}} \partial_*F.
    \]
    But the right-hand side is precisely $\indec_{\partial_*{\id_\sD}}(\partial_*F)$.
\end{proof}

\begin{remark}
    This corollary has a version for right modules and comodules, where one precomposes with $\Omega_\sD^\infty$ instead of postcomposing with $\Sigma_\sD^\infty$. 
\end{remark}

\subsubsection{The chain rule as a strong 2-functor}

The chain rule proved in \cref{thm: unstable-chain-rule} suggests that the lax functor $\partial_* \colon \diff \to \pressymstc$ can be refined to a strong functor if we replace the composition product of symmetric sequences by the relative composition product of bimodules.
The purpose of this subsection is to show that this can indeed be done.
The idea is to replace $\pressymstc$ by its \emph{Morita category}: this is a 2-precategory whose objects are algebras in the monoidal categories $\SymFunL_{\geq 1}(\sC,\sC)$, where $\sC$ ranges over all presentable categories, while the 1-morphisms are bimodules in $\SymFunL_{\geq 1}(\sC,\sD)$ over these algebras.
The composition of 1-morphisms is then given by the relative composition product.

We start by stating a general existence result and universal property of such a Morita category.
Given 2-precategories $\sX$ and $\sY$, let us write $\Lax(\sX,\sY)$ for the category of lax functors from $\sX$ to $\sY$, defined as a full subcategory of $\Fun_{/\Delta^\op}(\un \sX, \un \sY)$ (see \cref{definition:lax-functors}).
If $F \colon \sX \to \sY$ is a lax functor such that for any object $x$ in $\sX$, the map $\unit_{F(x)} \to F(\unit_x)$ from \cref{map:unit-2-cell-of-lax-functor} is an equivalence, then we call $F$ a \emph{unital lax functor}.
The full subcategory of $\Lax(\sX,\sY)$ spanned by the unital lax functors will be denoted $\Laxunit(\sX,\sY)$.

\begin{proposition}\label{proposition:Universal-property-Morita}
    Let $\sX$ be a 2-precategory and suppose that
    \begin{enumerate}[(1)]
        \item for any pair of objects $x$ and $y$ in $\sX$, the mapping category $\sX(x,y)$ admits geometric realizations, and
        \item for any three objects $x$, $y$ and $z$ in $\sX$, the composition map
        \[\sX(y,z) \times \sX(x,y) \to \sX(x,z)\]
        preserves geometric realizations in each variable separately.
    \end{enumerate}
    Then there exists a 2-precategory $\MMor(\sX)$, the \emph{Morita 2-precategory of $\sX$}, with the universal property that for any 2-precategory $\sY$, there is a natural equivalence
    \[\Lax(\sY,\sX) \simeq \Laxunit(\sY,\MMor(\sX)).\]
\end{proposition}

\begin{proof}
    This essentially follows from \cite[Theorem B]{Blom2024StraighteningEveryFunctor}. That statement is proved in the more general context of double categories (see \cref{example:Double-categories}), so let us explain how to deduce the version for 2-precategories from it. The inclusion of 2-precategories into double categories admits a right adjoint $\Hor$ called the \emph{horizontal fragment} (see for example \cite[\S 2.2]{Ruit2023FormalCategoryTheory}). It is not hard to see that for any 2-precategory $\sY$ and any double category $\mathbb{D}$, the counit $\Hor(\mathbb{D}) \to \mathbb{D}$ induces equivalences
    \[\Lax(\sY,\Hor(\mathbb{D})) \simeq \Lax(\sY,\mathbb{D}) \quad \text{and} \quad \Laxunit(\sY,\Hor(\mathbb{D})) \simeq \Laxunit(\sY,\mathbb{D}).\]
    In particular, if we define $\MMor(\sX)$ as $\Hor(\mathbb{M}\mathrm{or}(\sX))$, where $\mathbb{M}\mathrm{or}$ denotes the Morita double category from \cite[Theorem B]{Blom2024StraighteningEveryFunctor}, then we obtain a natural equivalence
    \[\Lax(\sY,\sX) \simeq \Laxunit(\sY,\MMor(\sX)). \qedhere\]
\end{proof}

The universal property of $\MMor(\sX)$ tells us that the space of objects of $\MMor(\sX)$ is given by the space of lax functors from $*$ to $\sX$. By \cref{example:Lax-functors-out-of-point}, such a lax functor consists of the data of an object $x$ in $\sX$ together with an algebra in the monoidal category $\sX(x,x)$. Similarly, a 1-morphism in $\MMor(\sX)$ is a lax functor from $[1]$ to $\sX$; that is, a functor $(\Delta_{/[1]})^\op \to \un \sX$ over $\Delta^\op$ that preserves cocartesian lifts of inerts. By \cref{proposition:approximation-to-BM-operad,proposition:bimodule-mapping-categories}, this is the same data as two objects $x$ and $y$ in $\sX$, algebras $A$ in $\sX(x,x)$ and $B$ in $\sX(y,y)$ and a $(B,A)$-bimodule $M$ in $\sX(x,y)$. It follows by \cite[Definitions A.11 \& A.16]{Blom2024StraighteningEveryFunctor} that the composition of 1-morphisms
\[\begin{tikzcd}
    A \ar[r,"M"] & B \ar[r,"N"] & C
\end{tikzcd}\]
in $\MMor(\sX)$ is given by the relative composition product
\[N \circ_B M = \Bahr(N,B,M).\]

The chain rule (\cref{thm: unstable-chain-rule}) can be reformulated as the statement that the Goodwillie derivatives assemble into a \emph{strong} functor into $\MMor(\pressymstc)$ (as opposed to a unital lax one):

\begin{theorem}\label{theorem:Derivatives-as-strong-functor}
    The lax functor $\partial_* \colon \diff \to \pressymstc$ refines to a strong functor
    \[\begin{tikzcd}
        \partial_* \colon \diff \ar[r] & \MMor(\pressymstc)
    \end{tikzcd}\]
    into the Morita 2-precategory of $\pressymstc$. On objects, this functor sends a differentiable category $\sC$ to the algebra $\partial_*{\id_\sC}$ in $\SymFunL_{\geq 1}(\sC,\sC)$, and on 1-morphisms it is given by sending $F \colon \sC \to \sD$ to the $(\partial_*{\id_\sD},\partial_*{\id_\sC})$-bimodule $\partial_* F$.
\end{theorem}

\begin{proof}
    Applying \cref{proposition:Universal-property-Morita} to $\pressymstc$, we obtain a unital lax functor $\diff \to \MMor(\pressymstc)$ having the stated behaviour on objects and 1-morphisms.
    To see that it is a strong functor, it suffices by \cref{lemma:Criterion-when-lax-functors-are-strict} to show that for any pair of composable functors $\sC \xrightarrow{G} \sD \xrightarrow{F} \sE$, the lax structure map
    \[\partial_* F \circ_{\partial_*{\id_{\sD}}} \partial_* G \to \partial_*(FG)\]
    is an equivalence.
    An inspection of the proof of \cite[Theorem B]{Blom2024StraighteningEveryFunctor} shows that this is precisely the map constructed in \cref{construction:The-lax-comparison-map}, which we show is an equivalence in \cref{thm: unstable-chain-rule}.
\end{proof}

\begin{remark}
    In \cref{theorem:Derivatives-as-strong-functor}, we obtain a strong functor $\diff \to \MMor(\pressymstc)$ into a 2-precategory. By composing with the universal functor from $\MMor(\pressymstc)$ into its completion (see \cref{remark:Completion-2-precategories}), we obtain a strong functor into an actual 2-category.
\end{remark}

By restricting this result to the endofunctors of a differentiable category $\sC$, we obtain the following.

\begin{corollary}
    Let $\sC$ be a differentiable category and write ${\BMod}_{({\partial_*{\id_\sC}},{\partial_*{\id_\sC}})}$ for the category of $\partial_*{\id_\sC}$-bimodules in $\SymFunL(\Sp(\sC), \Sp(\sC))$.
    Then the Goodwillie derivatives refine to a strong monoidal functor
    \[\begin{tikzcd}
        \Fun^{*,\omega}(\sC,\sC) \ar[r] & {\BMod}_{({\partial_*{\id_\sC}},{\partial_*{\id_\sC}})},
    \end{tikzcd}\]
    where ${\BMod}_{({\partial_*{\id_\sC}},{\partial_*{\id_\sC}})}$ is endowed with the relative composition product.
    More generally, given two differentiable categories $\sC$ and $\sD$, the Goodwillie derivatives define a strong functor
    \[\begin{tikzcd}[column sep = normal]
    	{\Fun^{*,\omega}(\sD,\sD)} & {\Fun^{*,\omega}(\sC,\sD)} & {\Fun^{*,\omega}(\sC,\sC)} \\
    	{{\BMod}_{({\partial_*{\id_\sD}},{\partial_*{\id_\sD}})}} & {{\BMod}_{({\partial_*{\id_\sD}},{\partial_*{\id_\sC}})}} & {{\BMod}_{({\partial_*{\id_\sC}},{\partial_*{\id_\sC}})}}
    	\arrow["\curvearrowright"{description, style={font=\large}}, draw=none, from=1-1, to=1-2]
    	\arrow[from=1-1, to=2-1]
    	\arrow[from=1-2, to=2-2]
    	\arrow["\curvearrowleft"{description, style={font=\large}}, draw=none, from=1-3, to=1-2]
    	\arrow[from=1-3, to=2-3]
    	\arrow["\curvearrowright"{description, style={font=\large}}, draw=none, from=2-1, to=2-2]
    	\arrow["\curvearrowleft"{description, style={font=\large}}, draw=none, from=2-3, to=2-2]
\end{tikzcd}\]
    of bitensored categories.
\end{corollary}

\begin{proof}
    This follows directly from \cref{def:endomorphisms-infty-category,proposition:bimodule-mapping-categories}.
\end{proof}

\subsubsection{A chain rule for non-presentable categories}

The aim of this section is to weaken the assumptions under which the chain rule can be proved.
More precisely, we will show that there is already a chain rule for functors between categories $\cC$ such that
\begin{itemize}
    \item $\sC$ is pointed;
    \item $\sC$ admits countable colimits;
    \item $\sC$ admits finite limits;
    \item sequential colimits commute with finite limits in $\sC$.
\end{itemize}
We will call a category satisfying these assumptions \emph{quasi-differentiable}.
Moreover, we will call a functor that preserves sequential colimits \emph{quasi-finitary}.

\begin{remark}\label{remark:sequential-vs-countable-fil}
    Any countable filtered category $I$ admits a right cofinal functor from $\bbN$.
    To see this, note that by \cite[Proposition 5.3.1.18]{HTT} there is a cofinal functor $P \to I$ where $P$ is a directed poset.
    An inspection of the proof shows that $P$ is countable if $I$ is.
    Choosing a bijection $\phi \colon \bbN \cong P$ of sets, we recursively construct an increasing function $f \colon \bbN \to P$ as follows: set $f(0) = \phi(0)$ and let $f(n)$ be any element of $P$ such that $f(n) > \phi(n)$ and $f(n) > f(n-1)$.
    It is then easily verified that $f$ is cofinal.
    
    It therefore follows that a functor is quasi-finitary if and only if it preserves all countable filtered colimits.
    Similarly, sequential colimits commute with finite limits in $\sC$ precisely if countable filtered colimits commute with finite limits in $\sC$.
\end{remark}

We write $\Fun^{\ast, \mathrm{seq}}(\sC, \sD)$ for the category of reduced, quasi-finitary functors from $\sC$ to $\sD$.
These are the kind of functors to which we will generalize the chain rule.

\begin{definition}\label{def:quasi-derivatives}
Let $\sC$ and $\sD$ be quasi-differentiable categories. The \emph{$n$th derivative functor} is defined as
\[
\begin{tikzcd}
    \partial_n \colon \Fun^{\ast, \mathrm{seq}}(\sC, \sD) \ar[r, "\mlin \circ \cross_n"] & \SymFun^{\ast, \mathrm{seq}, 1{\text -}\mathrm{exc}}_{n}(\sC, \sD) \ar[r, "\sim"] & \SymFun_n^{\omega_1{\text -}\mathrm{sm}}(\Sp(\sC), \Sp(\sD)).
\end{tikzcd}
\]
The category on the right is the category of symmetric functors $\Sp(\sC)^{\times n} \to \Sp(\sD)$ that preserve countable colimits in each variable separately.
The equivalence on the right follows from \cite[Corollary 6.2.3.22]{HA} together with the fact that a reduced, 1-excisive, sequential colimit preserving functor between stable categories is the same as a countable colimit preserving functor.
\end{definition}

Note that by enlarging our universe, we may assume without loss of generality that all our categories are small; see \cref{rem: enlarging-the-universe} below.
We will now prove the chain rule for reduced quasi-finitary functors between small quasi-differentiable categories.
For any small quasi-differentiable category $\cC$, the category $\Ind_{\omega_1}(\cC)$ is differentiable.
This will allow us to deduce a chain rule for quasi-differentiable categories from the one we already proved for differentiable categories.

\begin{proposition}\label{prop: ind-omega1-differentiable}
    Let $\cC$ be a small category that admits countable colimits and let $\cD$ be a cocomplete category. Then
    \begin{enumerate}[\upshape{(}1\upshape{)}]
        \item\label{item2:omega1} the Yoneda embedding $y \colon \cC \hookrightarrow \Ind_{\omega_1}(\cC)$ preserves countable colimits and all limits that exist in $\cC$;
        \item\label{item4:omega1} a functor $\cC \to \cD$ is (reduced) quasi-finitary if and only if its left Kan extension $\Ind_{\omega_1}(\cC) \to \cD$ is (reduced) finitary;
        \item\label{item5:omega1} a functor $\cC \to \cD$ preserves countable colimits if and only if its left Kan extension $\Ind_{\omega_1}(\cC) \to \cD$ preserves all small colimits;
        \item\label{item1:omega1} the category $\Ind_{\omega_1}(\cC)$ is presentable; and
        \item\label{item3:omega1} if $\cC$ is quasi-differentiable, then $\Ind_{\omega_1}(\cC)$ is differentiable.
    \end{enumerate}
\end{proposition}

\begin{proof}
    Item \ref{item2:omega1} follows from Propositions 5.1.3.2 and 5.3.5.14 of \cite{HTT}, item \ref{item5:omega1} follows from \cite[Example 5.3.6.8]{HTT} and item \ref{item1:omega1} is proved in \cite[Theorem 5.5.1.1]{HTT}.
    
    The ``if'' direction of \ref{item4:omega1} follows from item \ref{item2:omega1}.
    For the ``only if'' direction, let $f \colon \cC \to \cD$ be quasi-finitary and note that the square
    \[\begin{tikzcd}
    	{\Fun(\bbN,\cC)} & \cC \\
    	{\Fun(\bbN,\cD)} & \cD
    	\arrow["\colim", from=1-1, to=1-2]
    	\arrow["{f_*}"', from=1-1, to=2-1]
    	\arrow["f", from=1-2, to=2-2]
    	\arrow["\colim"', from=2-1, to=2-2]
    \end{tikzcd}\]
    commutes by assumption.
    Write $\overline{f}$ for the left Kan extension of $f$ along $\cC \to \Ind_{\omega_1}(\cC)$.
    Then $\overline{f}$ preserves $\omega_1$-filtered colimits.
    By applying $\Ind_{\omega_1}$ to the top part of the previous square and using \cref{lem: commute-fun-ind-countable} below, we see that the square
    \[\begin{tikzcd}
    	{\Fun(\bbN,\Ind_{\omega_1}(\cC))} & \Ind_{\omega_1}(\cC) \\
    	{\Fun(\bbN,\cD)} & \cD
    	\arrow["\colim", from=1-1, to=1-2]
    	\arrow["{\overline{f}_*}"', from=1-1, to=2-1]
    	\arrow["\overline{f}", from=1-2, to=2-2]
    	\arrow["\colim"', from=2-1, to=2-2]
    \end{tikzcd}\]
    commutes.
    It follows that $\overline{f}$ preserves sequential colimits, hence all filtered colimits by \cref{remark:sequential-vs-countable-fil,lemma:k-small-filtered-colimits}.
    This proves \ref{item4:omega1}.

    For item \ref{item3:omega1}, note that since the Yoneda embedding $\cC \to \Ind_{\omega_1}(\cC)$ preserves both initial and terminal objects, it follows that $\Ind_{\omega_1}(\cC)$ is pointed if $\cC$ is.
    Hence to prove item \ref{item3:omega1}, it suffices to show that for any finite category $K$, the limit functor
    \[{\lim} \colon \Ind_{\omega_1}(\Fun(K,\cC)) \simeq \Fun(K,\Ind_{\omega_1}(\cC)) \to \Ind_{\omega_1}(\cC)\]
    preserves filtered colimits.
    Since it preserves $\omega_1$-filtered colimits, by \cite[Lemma 5.3.5.8]{HTT} it is left Kan extended from $\Fun(K,\cC)$.
    The composite ${\lim} \colon \Fun(K,\cC) \to \cC \hookrightarrow \Ind_{\omega_1}(\cC)$ is quasi-finitary since $\cC$ is quasi-differentiable and the Yoneda embedding preserves countable colimits.
    It therefore follows from \ref{item4:omega1} that
    ${\lim} \colon \Fun(K,\Ind_{\omega_1}(\cC)) \to \Ind_{\omega_1}(\cC)$
    preserves filtered colimits.
\end{proof}

We used the following lemma in the proof of the previous proposition.

\begin{lemma} \label{lem: commute-fun-ind-countable}
    Let $\sC$ be a small category that admits countable colimits.
    Suppose that $I$ is a countable category.
    Then the Yoneda embedding induces an equivalence
    \[
    \begin{tikzcd}
    \Ind_{\omega_1}(\Fun(I, \sC)) \ar[r, "\sim"] & \Fun(I, \Ind_{\omega_1}(\sC)).
    \end{tikzcd}
    \]
\end{lemma}

\begin{proof}
   The Yoneda embedding gives a fully faithful functor $\Fun(I,\sC) \hookrightarrow \Fun(I,\Ind_{\omega_1}(\cC))$ that preserves countable colimits by \ref{item2:omega1} of \cref{prop: ind-omega1-differentiable}. By the Makkai--Pitts criterion \cite[Proposition 1.7]{Harr2023DerivedCategoryLocally}, it suffices to show that the image of $\Fun(I,\cC)$ in $\Fun(I,\Ind_{\omega_1}(\cC))$ consists of $\omega_1$-compact objects $X$ with the property that the functors $\Map(X,-)$ jointly detect equivalences in $\Fun(I,\Ind_{\omega_1}(\cC))$.
   
  To see that these objects are $\omega_1$-compact, observe that since $I$ is countable, its twisted arrow category $\TwL(I)$ is also countable. Let $X$ be an object of $\Fun(I,\cC)$ and $Y_\bullet \colon J \to \Fun(I,\Ind_{\omega_1}(\cC))$ an $\omega_1$-filtered diagram. Then
   \[
   \Map(X,\colim_j Y_j) \simeq \lim_{(i \to i') \in \TwL(I)} \colim_j \Map(X(i), Y_j(i')) \simeq \colim_j \Map(X,Y_j)
   \]
   by \cite[Prop.\ 5.1]{GepnerHaugsengea2017LaxColimitsFree}, where we use that $\omega_1$-filtered colimits commute with countable limits in $\Spc$.

   For the second step, given $c$ in $\cC$ and $i$ in $I$, write $I(i,-) \otimes c$ for the left Kan extension
   \[\begin{tikzcd}
    	{*} & \cC. \\
    	I
    	\arrow["c", from=1-1, to=1-2]
    	\arrow["i"', from=1-1, to=2-1]
    	\arrow[""{name=0, anchor=center, inner sep=0}, "{\mathrm{Lan}_i c}"', from=2-1, to=1-2]
    	\arrow[shorten >=2pt, Rightarrow, from=1-1, to=0]
    \end{tikzcd}\]
    Since the mapping spaces of $I$ are countable, this left Kan extension exists and is preserved by the inclusion $\Fun(I,\cC) \hookrightarrow \Fun(I,\Ind_{\omega_1}(\cC))$. We see that
    \[\Map(I(i,-) \otimes c, Y) \simeq \Map(c,Y(i))\]
    for any $Y$ in $\Fun(I,\Ind_{\omega_1}(\cC))$. Since the functors $\{\Map(c,-)\}_{c \in \cC}$ jointly detect equivalences in $\Ind_{\omega_1}(\cC)$, the result follows.
\end{proof}

By \cref{prop: ind-omega1-differentiable}, composing with the Yoneda embedding and then left Kan extending along $\cC \hookrightarrow \Ind_{\omega_1}(\cC)$ provides a fully faithful functor
\begin{equation}\label{eq:Lkany}
\begin{tikzcd}
    \mathrm{LKan} \circ y_* \colon \Fun^{\ast, \mathrm{seq}}(\sC, \sD) \ar[r, hook] & \Fun^{\ast, \mathrm{\omega}}(\Ind_{\omega_1}(\sC), \Ind_{\omega_1}(\sD)).
\end{tikzcd}
\end{equation}
Note that this inclusion preserves countable colimits since $\sD \hookrightarrow \Ind_{\omega_1}(\sD)$ does.
Writing $\qdiff$ for the 2-category of (small) quasi-differentiable categories and reduced quasi-finitary functors, we can assemble these fully faithful functors into a 2-functor as follows.

\begin{corollary}
    The construction $\cC \mapsto \Ind_{\omega_1}(\cC)$ defines a locally fully faithful inclusion
    \[\qdiff \hookrightarrow \diff\]
    of 2-categories, which on mapping categories is given by the functor \cref{eq:Lkany}.
\end{corollary}

\begin{proof}
    Let $\Cat$ denote the category of small categories and $\widehat{\Cat}{}^{\omega_1\text{-}\mathrm{fil}}$ the category of (possibly large) categories that admit $\omega_1$-filtered colimits and functors between them that preserve these.
    Then the 2-functor \[\Ind_{\omega_1} \colon \Cat \to \widehat{\Cat}{}^{\omega_1\text{-}\mathrm{fil}}\] has the desired description on mapping categories.
    (This 2-functor can be constructed using \cref{corollary:left-adjoints-constructed-locally}, by an argument similar to \cite[Corollary 5.3.6.10]{HTT}.)
    It follows from \cref{prop: ind-omega1-differentiable} that this functor takes the locally full subcategory $\qdiff \subset \Cat$ to the locally full subcategory $\diff \subset \widehat{\Cat}{}^{\omega_1\text{-}\mathrm{fil}}$, giving the desired functor.
\end{proof}

To prove a chain rule in this context, we will need to define a composition product for sequences of functors of the form
\[\{A_n \colon \cC^{\times n} \to \cD\}_{n \geq 0} \in \prod_{n \geq 0} \SymFun_n^{\omega_1{\text -}\mathrm{sm}}(\sC,\sD) \eqqcolon \SymFun^{\omega_1{\text -}\mathrm{sm}}(\sC,\sD),\]
where $\cC$ and $\cD$ are categories that admit countable colimits.
By \cref{prop: ind-omega1-differentiable}, we may view $\SymFun^{\omega_1{\text -}\mathrm{sm}}(\sC,\sD)$ as the full subcategory of $\SymFunL(\Ind_{\omega_1}(\cC),\Ind_{\omega_1}(\cD))$ spanned by those symmetric sequences $\{F_n\}_{n \geq 0}$ such that for each $n$, the composite
\[\begin{tikzcd}
    \cC^{\times n} \ar[r] & \Ind_{\omega_1}(\cC)^{\otimes n}_{h\Sigma_n} \ar[r,"F_n"] & \Ind_{\omega_1}(\cD)
\end{tikzcd}\]
lands in $\cD$.
Moreover, these subcategories are closed under countable colimits since $\cD \hookrightarrow \Ind_{\omega_1}(\cD)$ preserves these, so the composition product restricts to these subcategories.
This shows that we can make the following definition.

\begin{definition}[The 2-category $\Cat^{\omega_1,\Sym}$]\label{def:omega1-symseq}
    The 2-category $\Cat^{\omega_1,\Sym}$ is the locally full subcategory (see \cref{construction:1-full-subcategory}) of $\pressym$ whose objects are the presentable categories of the form $\Ind_{\omega_1}(\cC)$ where $\cC$ admits countable colimits, and whose mapping objects are the full subcategories $\SymFun^{\omega_1{\text -}\mathrm{sm}}(\sC,\sD)$ of $\SymFunL(\Ind_{\omega_1}(\cC),\Ind_{\omega_1}(\cD))$.
    We will write $\cC$ instead of $\Ind_{\omega_1}(\cC)$ for an object in this category. We write $\Cat^{\omega_1,\Sym}_{\mathrm{St}}$ for the full subcategory whose objects are stable and $\Cat^{\omega_1,\Sym}_{\mathrm{St},\geq 1}$ for the locally full subcategory whose 1-morphisms are concentrated in arity $\geq 1$.
\end{definition}

We will also need that taking derivatives commutes with passing to $\Ind_{\omega_1}$ in the appropriate sense.
For this we need the following lemma.

\begin{lemma} \label{lem: ind-and-stabilization}
    Let $\cC$ be a small quasi-differentiable category. Then there is an equivalence $\Sp(\Ind_{\omega_1}(\sC)) \simeq \Ind_{\omega_1}(\Sp(\sC))$.
\end{lemma}

\begin{proof}
    Recall that $\Sp(\sD) = \Fun^{*,1{\text -}\mathrm{exc}}(\Spc^{\mathrm{fin}}_*, \sD)$, where $\Spc^{\mathrm{fin}}_*$ is the category of finite pointed spaces.
    Since $\Spc^{\mathrm{fin}}_*$ is countable, \cref{lem: commute-fun-ind-countable} provides an equivalence 
    \[
    \Fun(\Spc^{\mathrm{fin}}_*, \Ind_{\omega_1}(\sC)) \simeq \Ind_{\omega_1}(\Fun(\Spc^{\mathrm{fin}}_*, \sC)).
    \]
    We will show that this restricts to an equivalence $\Sp(\Ind_{\omega_1}(\sC)) \simeq \Ind_{\omega_1}(\Sp(\sC))$.
    This amounts to showing that $F \colon \Spc_*^{\mathrm{fin}} \to \Ind_{\omega_1}(\cC)$ is linear if and only if it can be written as an $\omega_1$-filtered colimit of linear functors landing in $\cC$.
    The ``if'' direction follows since $\omega_1$-filtered colimits commute with finite limits in $\Ind_{\omega_1}(\cC)$.
    For the converse, suppose $F \colon \Spc^{\mathrm{fin}} \to \Ind_{\omega_1}(\sC)$ is linear and write $F \simeq \colim_i F_i$ as an $\omega_1$-filtered colimit of functors landing in $\sC$.
    Then $F \simeq \colim_i P_1(\mathrm{red}(F_i))$ expresses $F$ as an $\omega_1$-filtered colimit of linear functors landing in $\cC$.
\end{proof}

We obtain \cref{thm:intro-quasi-diff} by combining all the previous results.

\begin{theorem}\label{thm:chainrule-quasi-diff}
    The lax 2-functor $\partial_* \colon \diff \to \pressymstc$ from \cref{thm: lax-derivatives} restricts to a lax 2-functor
    \[
    \begin{tikzcd}
    \partial_* \colon \qdiff \ar[r] & \Cat^{\omega_1,\Sym}_{\mathrm{St},\geq 1}
    \end{tikzcd}
    \]
    which takes a quasi-differentiable category $\cC$ to its stabilization $\Sp(\cC)$, and on mapping categories is given by the Goodwillie derivatives functor.
    
    Moreover, for any pair $\cC \xrightarrow{G} \cD \xrightarrow{F} \cE$ of reduced quasi-finitary functors between quasi-differentiable categories, the induced map
    \[\begin{tikzcd}
        \partial_* F \circ_{\partial_*{\id_\cD}} \partial_*G \ar[r] & \partial_*(FG)
    \end{tikzcd}\]
    is an equivalence.
\end{theorem}

\begin{proof}
    Given a functor $F \colon \cC \to \cD$, write $\overline{F} \colon \Ind_{\omega_1}(\cC) \to \Ind_{\omega_1}(\cD)$ for the functor obtained by composing with the Yoneda embedding and left Kan extending along $\cC \hookrightarrow \Ind_{\omega_1}(\cC)$.
    By \cref{lem: ind-and-stabilization}, the composite
    \begin{equation}\label{eq:qdiffpressymstc}\begin{tikzcd}
    \qdiff \ar[r, hook] & \diff \ar[r,"\partial_*"] & \pressymstc
    \end{tikzcd}\end{equation}
    takes $F \colon \cC \to \cD$ to the symmetric sequence $\partial_* \overline{F}$ in $\SymFunL(\Ind_{\omega_1}(\cC),\Ind_{\omega_1}(\cD))$. Since the full subcategory $\cD \subset \Ind_{\omega_1}(\cD)$ is closed under finite limits and countable colimits by \cref{prop: ind-omega1-differentiable}, it follows that for every $n$, the composite
    \[\begin{tikzcd}
        \cC^{\times n} \ar[r] & \Ind_{\omega_1}(\cC)^{\otimes n}_{h\Sigma_n} \ar[r,"\partial_n \overline{F}"] & \Ind_{\omega_1}(\cD)
    \end{tikzcd}\]
    lands in $\cD$.
    We therefore see that the composite \cref{eq:qdiffpressymstc} lands in the locally full subcategory $\Cat^{\omega_1,\Sym}_{\mathrm{St},\geq 1}$, producing the desired lax functor.
    Moreover, again since $\cD \subset \Ind_{\omega_1}(\cD)$ is closed under finite limits and countable colimits, this functor agrees with the Goodwillie derivatives functor from \cref{def:quasi-derivatives} on mapping categories.

    The fact that $\partial_* F \circ_{\partial_*{\id_\cD}} \partial_*G \to \partial_*(FG)$ is an equivalence follows directly from \cref{thm: unstable-chain-rule} and the fact that the inclusion
    \[\begin{tikzcd}
        \SymFun^{\omega_1\text{-}\mathrm{sm}}(\Sp(\sC), \Sp(\sE)) \ar[r, hook] & \SymFun^{\mathrm{L}}(\Ind_{\omega_1}(\Sp(\sC)), \Ind_{\omega_1}(\Sp(\sE)))
    \end{tikzcd}\]
    preserves countable colimits and hence geometric realizations.
\end{proof}

\begin{remark} \label{rem: enlarging-the-universe}
    The proof of \cref{thm:chainrule-quasi-diff} uses that $\Ind_{\omega_1}(\cC)$ is presentable, which requires $\cC$ to be small.
    To obtain a chain rule for functors between quasi-differentiable categories that are not necessarily small, one can enlarge the universe.
    Explicitly, suppose we are working in a Grothendieck universe $\sU$ that contains an infinite set.
    Let $\sV$ be a Grothendieck universe such that $\sU \in \sV$.
    Then any category that is large in $\sU$ becomes small in $\sV$, so by running the proof of \cref{thm:chainrule-quasi-diff} in $\sV$, it follows that the chain rule also holds for quasi-differentiable categories that are large with respect to the universe $\sU$.
\end{remark}

The same proof as \cref{theorem:Derivatives-as-strong-functor} shows that this chain rule can be reinterpreted as follows.

\begin{corollary}
    The lax functor $\partial_* \colon \qdiff \to \Cat^{\omega_1,\Sym}_{\mathrm{St},\geq 1}$ refines to a strong functor
    \[\begin{tikzcd}
        \partial_* \colon \qdiff \ar[r] & \MMor(\Cat^{\omega_1,\Sym}_{\mathrm{St},\geq 1})
    \end{tikzcd}\]
    into the Morita 2-category of $\Cat^{\omega_1,\Sym}_{\mathrm{St},\geq 1}$.
    On objects, this functor sends a quasi-differentiable category $\cC$ to the algebra $\partial_*{\id_\cC}$ in $\SymFun^{\omega_1{\text -}\mathrm{sm}}_{\geq 1}(\cC,\cC)$, and on 1-morphisms it is given by sending $F \colon \sC \to \sD$ to the $(\partial_*{\id_\sD},\partial_*{\id_\sC})$-bimodule $\partial_* F$.
\end{corollary}

\begin{remark}
    Since $\Cat^{\omega_1,\Sym}_{\mathrm{St},\geq 1} \hookrightarrow \pressymstc$ is a locally fully faithful inclusion that preserves geometric realizations on mapping categories, it follows that $\MMor(\Cat^{\omega_1,\Sym}_{\mathrm{St},\geq 1}) \hookrightarrow \MMor(\pressymstc)$ is a locally fully faithful inclusion.
    Essentially by construction, we have a commutative square
    \[\begin{tikzcd}
    	\qdiff & {\MMor(\Cat^{\omega_1,\Sym}_{\mathrm{St},\geq 1})} \\
    	\diff & {\MMor(\pressymstc)}
    	\arrow["{\partial_*}", from=1-1, to=1-2]
    	\arrow[hook, from=1-1, to=2-1]
    	\arrow[hook, from=1-2, to=2-2]
    	\arrow["{\partial_*}"', from=2-1, to=2-2]
    \end{tikzcd}\]
    of strong 2-functors.
\end{remark}

\newpage

\appendix

\section{Background on 2-category theory and lax functors}\label{appendix:Background-2-categories}

The goal of this appendix is to review the basics of $(\infty,2)$-category theory. As in the main text, we will say \emph{2-category} when we mean $(\infty,2)$-category and \emph{1-category} or \emph{category} when we mean $(\infty,1)$-category. The main function of this appendix is expository, and the results of this appendix are surely known to experts on 2-category theory. However, since we suspect that not everyone who is interested in Goodwillie calculus has a background in 2-category theory, we decided to include this appendix. In our treatment, we tried to be as self-contained as possible, developing the theory of 2-categories internal to the well-established language of 1-category theory.

We start in \cref{sec:2-categories-as-Segal-objects} with defining $2$-categories as certain simplicial objects in the category of $1$-categories. We do this in two different (yet equivalent) ways, both of which will have their merits throughout the rest of this appendix.

In \cref{sec:1/2-categorical-notions}, we then extend various fundamental notions from $1$-category theory to the setting of $2$-categories, such as functor categories, cocartesian fibrations, (co)limits, etc. Throughout these sections, a useful perspective to take is that $2$-categories are like $1$-categories, but where one has a mapping \emph{category} between any two objects as opposed to a mapping \emph{space}. In particular, many definitions here are obtained by considering their $1$-categorical analogue and replacing all mapping spaces with mapping categories.

In \cref{sec:Lax-functors-between-2-categories}, we focus on a feature of $2$-categories that is not so visible from the aforementioned perspective: that of commutativity up to a \emph{non-invertible} 2-morphism, aka \emph{lax commutativity}. For these sections, it can be beneficial to think of a $2$-category as a collection of monoidal categories $\{\sV_s\}_{s \in S}$ together with a collection of bimodule categories $\{\sM_{s,t}\}_{s,t \in S}$, parametrized by some ``space of objects'' $S$. In these sections we particularly focus on \emph{lax functors} between 2-categories, show that these behave roughly like lax monoidal functors between monoidal categories, and discuss several methods for constructing them.

\subsection{2-categories as certain Segal objects in \texorpdfstring{$\Cat$}{Cat}}\ \label{sec:2-categories-as-Segal-objects}

In this section we will define 2-categories. The definition that we will use is a variation on the usual definition of 1-categories as complete Segal spaces. Namely, by replacing all spaces in the definition of a complete Segal space with categories, one can add an extra categorical dimension. By imposing an extra condition to ensure that one has a  
\emph{space} of objects and not a category of objects, we arrive at the definition of a 2-category. It turns out that there are two natural ways to impose such an extra condition, leading to two possible definitions of the notion of a 2-category. The first of these is described in \cref{subsec:Twofold-Segal-spaces} and will simply be called a \emph{2-category}. The other definition, which turns out to be equivalent to the first one, will be called a \emph{flipped 2-category} and is studied in \cref{subsec:Flipped-2-categories}. Each of these two definitions has its merits: flipped 2-categories are very useful for translating facts from 1-category theory to the world of 2-categories, as will be demonstrated in \cref{sec:1/2-categorical-notions}. On the other hand, the definition of a 2-category given in \cref{subsec:Twofold-Segal-spaces} will be very convenient for the study of lax functors carried out in \cref{sec:Lax-functors-between-2-categories} below.

\subsubsection{Twofold Segal spaces}\label{subsec:Twofold-Segal-spaces}

Let us start by recalling the definition of a (complete) Segal object.

\begin{definition}[Complete Segal objects]\label{definition:Complete-Segal-objects}
	Let $\sC$ be a category that admits finite limits. A \emph{Segal object in $\sC$} is a simplicial object $X_\bullet \colon \Delta^{\op} \to \sC$ such that for every $n \geq 2$, the Segal map
	\[X_n \to \underbrace{X_1 \times_{X_0} \ldots \times_{X_0} X_1}_{n\text{-times}}\]
	is an equivalence. A Segal object $X$ in $\sC$ is \emph{complete} if moreover the map
	\[X_0 \to X_3 \times_{X_1 \times X_1} (X_0 \times X_0)\]
	is an equivalence. The category of (complete) Segal objects in $\sC$ is defined as the full subcategory of $\Fun(\Delta^{\op}, \sC)$ spanned by the (complete) Segal objects.
\end{definition}

\begin{remark}
    One can think of Segal objects $X_\bullet$ in $\cC$ as ``categories internal to $\cC$'', where $X_0$ is the ``object of objects'' and $X_1$ is the ``object of arrows''. The Segal condition then provides a composition map
    \[X_1 \times_{X_0} X_1 \simeq X_2 \xrightarrow{d_1} X_1.\]
    The pullback $X_3 \times_{X_1 \times X_1} (X_0 \times X_0)$ appearing in the completeness condition is taken over the cospan
    \[X_3 \xrightarrow{(d_1d_0, d_1d_3)} X_1 \times X_1 \xleftarrow{s_0 \times s_0} X_0 \times X_0.\]
    In particular, we can think of this pullback as the collection of composable triples of arrows $(f,g,h)$ such that $fg$ and $gh$ are the identity; alternatively, it is the ``object of invertible arrows'' of $X_\bullet$. The completeness condition ensures that this ``object of invertible arrows'' agrees with the ``object of objects'' $X_0$.
\end{remark}

\begin{example}[Complete Segal spaces]
	The category of complete Segal objects in $\Spc$ will be denoted $\CSS$. By the work of Joyal--Tierney \cite{JoyalTierney2007QuasicategoriesVsSegal}, the functor $\Catinfty \to \Fun(\Delta^{\op}, \Spc)$ that sends a category $\sC$ to its nerve 
	\[\segNerve(\sC) \colon \Delta^{\op} \to \Spc; \quad [n] \mapsto \Map([n],\sC)\]
	is fully faithful with essential image given by the full subcategory of complete Segal spaces. In particular, the nerve functor gives an equivalence $\segNerve \colon \Catinfty \eqarrow \CSS$ of categories.
\end{example}

\begin{example}[Monoidal categories]\label{example:monoidal-categories}
	The category $\MonCat$ of \emph{monoidal categories} and strong monoidal functors is equivalent to the category of Segal objects $X$ in $\Catinfty$ with the property that $X_0 \simeq *$.
\end{example}

\begin{example}[Double categories]\label{example:Double-categories}
    More generally, a Segal object in $\Catinfty$ is called a \emph{double category}.
\end{example}

\begin{example}[Coskeleta]\label{example:coskeleton}
    Let $\cC$ be a category that admits finite products. Then the evaluation functor $\ev_0 \colon \Fun(\Delta^\op, \cC) \to \cC$ admits a right adjoint $\cosk_0$ given by right Kan extension along $\{[0]\} \hookrightarrow \Delta^\op$. Informally, it can be described by the formula $(\cosk_0 X)_n = X^{\times n+1}$, where the face maps delete copies of $X$ and the degeneracy maps are given by diagonals. Observe that $\cosk_0 X$ is always a Segal object.
\end{example}

We now introduce (locally complete) twofold Segal spaces and their complete version. The latter will be the model for 2-categories used throughout this paper, so we will simply call them 2-categories. For brevity, we will use the term 2-precategory for what is generally called a locally complete twofold Segal space.

\begin{definition}[2-precategories]\label{definition:twofold-Segal-space}
    A \emph{2-precategory} is a Segal object $\sX$ in $\Catinfty$ such that $\sX_0$ is a space. The full subcategory of $\Fun(\Delta^\op,\Catinfty)$ spanned by the 2-precategories is denoted $\twofSeg$ and its morphisms are called \emph{functors (of 2-precategories)}.
\end{definition}

\begin{definition}[2-categories]\label{definition:2-category}
	A \emph{2-category} is a 2-precategory that is moreover complete; that is, a complete Segal object $\sX$ in $\Catinfty$ such that $\sX_0$ is a space. A \emph{functor (of 2-categories)} $\sX \to \sY$ is defined as a map from $\sX$ to $\sY$ in $\Fun(\Delta^\op,\Catinfty)$. The category $\Cattwo$ is the full subcategory of $\Fun(\Delta^\op,\Catinfty)$ spanned by the 2-categories.
\end{definition}

\begin{example}[Viewing 1-categories as 2-categories]\label{example:1-categories-as-2-categories}
    The fully faithful inclusion $\Spc \hookrightarrow \Cat$ induces a fully faithful inclusion $\Cat \simeq \CSS \hookrightarrow \Cattwo$. We will generally identify $\Cat$ with its image under this inclusion. To distinguish the objects of $\Cat$ from those of $\Cattwo$, we will call them \emph{1-categories}, or simply \emph{categories} when no confusion can arise.
\end{example}

\begin{example}\label{example:monoidal-category-as-2-category}
    By \cref{example:monoidal-categories}, any monoidal category $\sC$ is a 2-precategory in the sense of \cref{definition:twofold-Segal-space}. This 2-precategory is complete precisely if the space of invertible objects in $\sC$ is contractible; that is, the only invertible object in $\sC$ is the unit, and the space of automorphisms of the unit is contractible.
\end{example}

\begin{remark}\label{remark:Twofold-Segal-spaces-as-bisimplicial-spaces}
    The equivalence $\Catinfty \simeq \CSS$ allows us to identify the category $\Cattwo$ of complete twofold Segal spaces with the full subcategory of $\Fun(\Delta^{\op} \times \Delta^{\op}, \Spc)$ spanned by those simplicial spaces $X_{\bullet,\bullet}$ for which
    \begin{enumerate}[(i),noitemsep]
	   \item $X_{n,\bullet}$ and $X_{\bullet,m}$ are complete Segal spaces for any $n,m \geq 0$, and
	   \item $X_{0,\bullet}$ is a constant simplicial space.
    \end{enumerate}
    For the second item, we use that the nerve $\segNerve \colon \Catinfty \eqarrow \CSS$ takes spaces to constant simplicial spaces.
\end{remark}

\begin{remark}
	In \cite[\S 14]{BarwickSchommerPries2021UnicityTheoryHigher}, an alternative definition of complete twofold Segal spaces is given where only $X_{0,\bullet}$ and $X_{\bullet,m}$ are required to be complete. By \cite[Lemma 2.8]{Johnson-FreydScheimbauer2017OpLaxNatural} this implies that $X_{n,\bullet}$ is complete for every $n \geq 0$, hence this definition is equivalent to the one given above.
\end{remark}

We now give a few basic definitions and constructions related to 2-categories.
We will state them in the slightly more general setting of 2-precategories. The zeroth object $\sX_0$ of a 2-precategory $\sX$ will be referred to as its \emph{space of objects}.

\begin{definition}[Mapping categories]\label{def:mapping-infty-category}
	Let $\sX$ be a 2-precategory. For any pair of objects $x,y \in \sX_0$, the \emph{mapping category} $\sX(x,y)$ is defined by the pullback
	\[\begin{tikzcd}
		\sX(x,y) \ar[r] \ar[dr,"\lrcorner",very near start, phantom] \ar[d] & \sX_1 \ar[d,"{(d_1,d_0)}"]\\
		* \ar[r,"{(x,y)}"'] & \sX_0^{\times 2}.
	\end{tikzcd}\]
	An object in one of the categories $\sX(x,y)$ will be referred to as a \emph{1-morphism} of $\sX$, and a map in one of these categories as a \emph{2-morphism} of $\sX$.
\end{definition}

\begin{definition}
	Let $F \colon \sX \to \sY$ be a functor between 2-precategories. Then $F$ is called \emph{fully faithful} if for any $x,y \in \sX_0$, the map $F_{x,y} \colon \sX(x,y) \to \sY(x,y)$ is an equivalence. The functor $F$ is called \emph{essentially surjective} if the induced map $\pi_0(\sX_0) \to \pi_0(\sY_0)$ is surjective.
\end{definition}

We leave the proof of the following to the reader.

\begin{lemma}\label{lemma:fully-faithful-ess-surjective-is-equivalence}
	A functor $\sX \to \sY$ between 2-categories is an equivalence if and only if it is fully faithful and essentially surjective.
\end{lemma}

\begin{remark}[Completion]\label{remark:Completion-2-precategories}
	The inclusion $\PreCattwo \hookrightarrow \Cattwo$ admits a left adjoint called \emph{completion}.
    The existence of this left adjoint follows since $\Cattwo$ is a Bousfield localization of $\PreCattwo$ by e.g. \cite[Proposition 7.16]{Haugseng2015RectificationEnrichedCategories}.
	That is, any 2-precategory admits a universal functor to a 2-category.
    Moreover, this universal functor is essentially surjective and fully faithful, hence $\Cattwo$ is the localization of $\PreCattwo$ at the essentially surjective and fully faithful functors.
\end{remark}

\begin{definition}[Composition of 1-morphisms]\label{remark:composition-functor-2-category}
	Let $\sX$ be a 2-precategory and consider the commutative square
	\[\begin{tikzcd}
		\sX_2 \ar[r,"d_1"] \ar[d] & \sX_1 \ar[d]\\
		\sX_0^{\times 3} \ar[r,"\mathrm{pr}_{13}"] & \sX_0^{\times 2}.
	\end{tikzcd}\]
	Suppose that three objects $x,y,z \in \sX_0$ are given. Pulling back along the map $(x,y,z) \colon * \to \sX_0^{\times 3}$ and using the equivalence $\sX_2 \simeq \sX_1 \times_{\sX_0} \sX_1$, we obtain the \emph{composition functor}
    \[\circ \colon \sX(y,z) \times \sX(x,y) \to \sX(x,z).\]
    Moreover, for any object $x \in \sX_0$, the map $* \xrightarrow{x} \sX_0 \xrightarrow{s_0} \sX_1$ factors canonically through a map $* \to \sX(x,x)$, yielding the \emph{unit} or \emph{identity} $\unit_x$ of $\sX$ at $x$. We say that a 1-morphism is \emph{invertible} or an \emph{equivalence} if it admits both a left and a right inverse for this composition functor.
\end{definition}

\begin{construction}[Opposites of 2-categories]\label{construction:opposites-of-2-categories}
    There are 2 natural ways to take the opposite of a 2-category $\sX$: one can invert the direction of its 1-morphisms or of its 2-morphisms.
    For the first of these, recall that $\Delta$ has an automorphism $\mathit{rev}$ that reverses the order on a linearly ordered set.
    Precomposing $\sX \colon \Delta^\op \to \Cat$ with this functor $\mathit{rev}$ yields the 2-category $\sX^\op$ in which the direction of the 1-morphisms is inverted.
    To invert the 2-morphisms, one instead needs to postcompose $\sX$ with the automorphism $(-)^\op \colon \Cat \to \Cat$. The resulting 2-category will be denoted by $\sX^\co$.
    We will write $\sX^\coop$ for the 2-category $(\sX^\co)^\op \simeq (\sX^\op)^\co$.
\end{construction}

Note that the inclusion $\Spc \hookrightarrow \Cat$ has a right adjoint $(-)^\simeq$. Applying this right adjoint levelwise, we obtain a right adjoint $(-)^{\leq 1}$ to the inclusion $\Cat \simeq \CSS \hookrightarrow \Cattwo$ from \cref{example:1-categories-as-2-categories}.

\begin{definition}\label{def:underlying-1-category}
    Let $\sX$ be a 2-category. Then its \emph{underlying 1-category} is the 1-category $\sX^{\leq 1}$. Note that as a complete Segal space, $\sX^{\leq 1}$ is given by $[n] \mapsto (\sX_\bullet)^\simeq$.
\end{definition}

\begin{example}[The 2-category of categories]\label{example:Cat-as-a-2-category}
    The principal example of a 2-category is the 2-category of 1-categories. In \cite[\S 4.2]{Nuiten2023StraighteningSegalSpaces}, an explicit construction of this 2-category is given. One can verify that its mapping categories are given by functor categories and that its underlying 1-category is the usual 1-category of 1-categories (see \cite[Remark 4.23]{Nuiten2023StraighteningSegalSpaces}). For this reason, we will simply denote this 2-category by $\Cat$. When relevant, it will always be clear from the context whether we consider $\Cat$ as a 1-category or as a 2-category.
\end{example}

Let us conclude by describing a method to construct subobjects of 2-categories.

\begin{definition}[Locally fully faithful inclusion]
    Let $F \colon \sY \to \sX$ be a functor between 2-precategories.
    Then $F$ is called a \emph{locally fully faithful inclusion} if $F_0 \colon \sY_0 \to \sX_0$ is an inclusion of path components and for any $y,y' \in \sY_0$, the functor $F_{y,y'} \colon \sY(y,y') \to \sX(Fy,Fy')$ is fully faithful.
\end{definition}

Up to equivalence, any locally fully faithful inclusion into $\sX$ is constructed as follows.

\begin{construction}[Locally full subcategories]\label{construction:1-full-subcategory}
    Let $\sX$ be a 2-precategory and suppose that $\sY_0 \subset \sX_0$ is a subobject and that for every $y,y' \in \sY_0$, we are given full subcategories $\sY(y,y') \subset \sX(y,y')$ closed under composition and containing all 1-morphisms in the essential image of $\sY_0 \to \sX_1$. For $n \geq 1$, define $\sY_n \subset \sX_n \simeq \sX_1 \times_{\sX_0} \cdots \times_{\sX_0} \sX_1$ to be the full subcategory whose objects are those sequences $(f_1,\ldots,f_n) \in \sX_1 \times_{\sX_0} \cdots \times_{\sX_0} \sX_1$ such that for each $i$, there exist $y,y' \in \sY_0$ such that $f_i$ lies in $\sY(y,y')$. Then one can verify that the simplicial structure maps of $\sX$ make  $\sY = \sY_\bullet$ into a 2-precategory. In this situation, $\sY$ will be called a \emph{locally full sub(pre)category} of $\sX$. If for any pair $y,y' \in \sY_0$, the full subcategory $\sY(y,y') \subset \sX(y,y')$ is all of $\sX(y,y')$, then we will call $\sY$ a \emph{full sub(pre)category}. Note that if $\sX$ is a 2-category, then $\sY$ is also a 2-category.
\end{construction}

\subsubsection{Flipped 2-categories}\label{subsec:Flipped-2-categories}

Let us now describe a slightly different perspective on 2-categories. In \cref{remark:Twofold-Segal-spaces-as-bisimplicial-spaces}, we saw that $\Cattwo$ can be described as a certain full subcategory of $\Fun(\Delta^\op \times \Delta^\op, \Spc)$. By swapping the simplicial coordinates, we arrive at the following definition.

\begin{definition}[Flipped 2-categories]\label{definition:flipped-2-categories}
    A \emph{flipped 2-category} is a complete Segal object $\mathscr{Z}$ in $\Catinfty$ such that
    \begin{equation}\label{equation:Simplicial-groupoid-core}
		\Delta^{\op} \to \Catinfty; \quad [n] \mapsto (\mathscr{Z}_n)^\simeq
	\end{equation}
    is a constant simplicial space. The category $\Cattwofl$ is the full subcategory of $\Fun(\Delta^{\op}, \Cat)$ spanned by these flipped 2-categories.
\end{definition}

We let $(\cdot)^\fl \colon \Cattwo \to \Cattwofl$ denote the functor defined by viewing a 2-category as a bisimplicial space and swapping the simplicial coordinates. It follows immediately that this is an equivalence:

\begin{proposition}\label{proposition:Equivalence-with-flipped-2-categories}
	The functor $(\cdot)^\fl \colon \Cattwo \to \Cattwofl$ is an equivalence of categories.
\end{proposition}

\begin{proof}
	This follows since under the equivalence $N \colon \Catinfty \eqarrow \CSS$, the groupoid core $\sC^\simeq$ of a category $\sC$ corresponds to the zeroth space of $N(\sC)$.
\end{proof}

Given a 2-category $\sX$, let us briefly describe what the categories $\sX^\fl_n$ look like. The zeroth category $\sX^\fl_0$ of the simplicial object $\sX^\fl$ is the underlying 1-category of $\sX$, obtained by forgetting all non-invertible 2-morphisms. The category $\sX^\fl_n$ should more generally be thought of as the 1-category whose objects are those of $\sX$, and whose arrows are of the form
\[\begin{tikzcd}
	x \arrow[rr, "F_0"name=a0, looseness = 2, bend left = 70pt ]\arrow[rr, "F_1"'name=a1, bend left = 45pt]\arrow[rr, "F_{n-1}"name=am1, bend right = 45pt]\arrow[rr, "F_n"'name=am, looseness = 2, bend right = 70pt] & & y.	
	\arrow[from=a0,to=a1, Rightarrow, "\alpha_1", shorten <= 3pt, shorten >= 3pt]
	\arrow[from=a1,to=am1, phantom, "\vdots", yshift = 3pt]
	\arrow[from=am1,to=am, Rightarrow,"\alpha_n", shorten <= 3pt, shorten >= 3pt] 
\end{tikzcd}\]
Here the $F_i$ are 1-morphisms and the $\alpha_i$ are 2-morphisms in $\sX$.

\begin{remark}\label{remark:constant-flipped-2-categories}
   This informal description of $\sX^\fl_n$ suggests that 1-categories can be viewed as constant flipped 2-categories. Indeed, it is easily checked that the composition of the inclusion $\Cat \hookrightarrow \Cattwo$ from \cref{example:1-categories-as-2-categories} with the equivalence $\fl \colon \Cattwo \to \Cattwofl$ sends a category $\sC$ to the constant simplicial object in $\Fun(\Delta^\op,\Cat)$ with value $\sC$.
\end{remark}

\begin{remark}\label{remark:mapping-category-flipped-2-category}
    Let $\sX$ be a 2-category. It follows from \cref{remark:Twofold-Segal-spaces-as-bisimplicial-spaces,def:mapping-infty-category} that the mapping category $\sX(x,y)$, when viewed as a complete Segal space, is given by the pullback
    \[\begin{tikzcd}
		\sX(x,y)_\bullet \ar[r] \ar[dr,"\lrcorner",very near start, phantom] \ar[d] & \sX_{1,\bullet} \ar[d,"{(d_1,d_0)}"]\\
		* \ar[r,"{(x,y)}"'] & \sX_{0,\bullet}^{\times 2}.
	\end{tikzcd}\]
    Swapping the simplicial coordinates, it follows that $\sX(x,y)_\bullet \simeq (\sX^\fl_\bullet)(x,y)$.
\end{remark}

\begin{remark}\label{remark:opposite-flipped-2-category}
    Since the flipped 2-category $\sX^\fl$ is defined by swapping the two simplicial indexing categories of $\sX$, it follows from \cref{construction:opposites-of-2-categories} that $(\sX^\op)^\fl$ is obtained from $\sX^\fl$ by levelwise taking opposites.
    More precisely, $(\sX^\op)^\fl_\bullet \simeq (\sX^\fl_\bullet)^\op$.
\end{remark}

\subsection{Extending 1-categorical notions to 2-categories}\label{sec:1/2-categorical-notions}

We will now extend some basic notions from 1-category theory, such as functor categories, cocartesian fibrations, (co)limits, and so on, to the setting of 2-categories. While these extensions are known to the expert and have appeared in the literature before, we give a somewhat non-standard treatment of them. The general theme of this section is that when using flipped 2-categories, many 2-categorical notions agree with ``levelwise'' versions of their 1-categorical counterparts. This will allow us to deduce many facts about 2-categories from their 1-categorical analogues.

\subsubsection{The internal hom of \texorpdfstring{$\Cattwo$}{2-Cat}}

Given two bisimplicial objects $\sX$ and $\sY$ in $\Spc$, write $\sY^\sX$ for their internal mapping object in $\Fun(\Delta^{\op} \times \Delta^\op,\Spc)$. If $\sX$ and $\sY$ are complete Segal objects in all directions, or equivalently complete Segal objects in $\Cat$, then by \cite[Proposition 2.8]{Nuiten2023StraighteningSegalSpaces} this also holds for $\sY^\sX$. However, $\sY^\sX$ need not be a 2-category, even if $\sX$ and $\sY$ are. This is remedied by the following lemma. Let us write $\CSS(\Cat)$ for the category of complete Segal objects in $\Cat$.

\begin{lemma}[{\cite[Lemma 4.7]{Nuiten2023StraighteningSegalSpaces}}]
    The inclusion $\Cattwo \hookrightarrow \CSS(\Cat)$ admits a right adjoint $r$.
\end{lemma}

\begin{proof}
    We may equivalently show that the inclusion $\Cattwofl \hookrightarrow \CSS(\Cat)$ admits a right adjoint. Given an object in $\sY_\bullet$ in $\CSS(\Cat)$, let $(r\sY)_n$ denote the full subcategory of $\sY_n$ spanned by those objects that are in the essential image of the degeneracy map $\sY_0 \to \sY_n$.
    Then $r\sY$ defines a simplicial subobject of $\sY$.
    We leave it to the reader to verify that $r\sY$ is again a complete Segal object in $\Cat$.

    Now observe that the complete Segal condition on $\sY$ ensures that $n \mapsto (\sY_n)^\simeq$ is a complete Segal space. In particular, the map $(\sY_0)^\simeq \to (\sY_n)^\simeq$ is an inclusion of path components.
    Since $(r\sY)_n$ is defined as the full subcategory spanned by those path components in the essential image of $(\sY_0)^\simeq \to (\sY_n)^\simeq$, it follows that $(r\sY)_n^\simeq$ is a constant simplicial space.
    It therefore follows that $r\sY$ is a flipped 2-category. Moreover, if $\sX^\fl$ is a flipped 2-category, then any map $\sX^\fl \to \sY$ lands in $r\sY$ since $(\sX^\fl_n)^\simeq$ is constant. We conclude that $r$ is right adjoint to the inclusion $\Cattwo \simeq \Cattwofl \hookrightarrow \CSS(\Cat)$.
\end{proof}

Combined with the fact that for any pair of 2-categories $\sX$ and $\sY$, the object $\sY^\sX$ lies in $\CSS(\Cat)$, it follows that $r(\sY^\sX)$ defines an internal mapping object for $\Cattwo$.

\begin{proposition}\label{proposition:internal-hom-of-2Cat}
	The category $\Cattwo$ is cartesian closed. For a pair of 2-categories $\sX$ and $\sY$, the internal mapping object is given by $r(\sY^\sX)$.
\end{proposition}

\begin{definition}[Functor 2-categories]\label{definition:Functor-2-categories}
    Let $\sX$ and $\sY$ be 2-categories. Then the \emph{functor 2-category} $\Funtwo(\sX,\sY)$ is defined as the internal mapping object $r(\sY^\sX)$. The underlying 1-category of $\Funtwo(\sX,\sY)$ will be denoted by $\Fun(\sX,\sY)$.
\end{definition}

\begin{example}\label{example:flipped-internal-hom}
    Suppose $\sC$ is a 1-category, viewed as a (flipped) 2-category via \cref{example:1-categories-as-2-categories,remark:constant-flipped-2-categories}. From the explicit construction of the internal mapping object $\Funtwo(\sC,\sY)$ as $r(\sY^\sC)$, it follows that $\Fun(\sC,\sY)^\fl_n$ is equivalent to the full subcategory of $\Funtwo(\sC,\sY^\fl_n)$ spanned by those functors $\sC \to \sY^\fl_n$ that land in the essential image of $\sY^\fl_0 \to \sY^\fl_n$.
\end{example}

\begin{remark}
    The internal mapping object $\Fun$ of $\Cat$ also induces a $\Cat$-enrichment on $\Fun(\Delta^\op,\Cat)$. Since $\Cattwo$ and $\Cattwofl$ are both defined as full subcategories of $\Fun(\Delta^\op,\Cat)$, they both inherit a $\Cat$-enrichment. For the full subcategory $\Cattwofl$, this recovers the functor categories $\Fun(-,-)$ defined in \cref{definition:Functor-2-categories}. However, the $\Cat$-enrichment that $\Cattwo$ inherits from $\Fun(\Delta^\op,\Cat)$ does \textbf{not} agree with the functor categories $\Fun(-,-)$.
\end{remark}

\begin{example}[Slice categories]\label{example:arrow-2-categories}
    The \emph{arrow 2-category} $\Ar(\sX)$ of a 2-category $\sX$ is defined as $\Funtwo([1],\sX)$.
    The inclusions $\{0\},\{1\} \hookrightarrow [1]$ induce maps $\ev_0,\ev_1 \colon \Funtwo([1],\sX) \to \sX$.
    Given an object $x$ in $\sX$, we define the \emph{slice categories} $\sX_{x/}$ and $\sX_{/x}$ as the fibers $\ev_0^{-1}(x)$ and $\ev_1^{-1}(x)$, respectively.
    Using \cref{example:flipped-internal-hom} to unwind this construction, we see that $(\sX_{x/})^\fl_n$ is equivalent to the full subcategory of $(\sX_n^\fl)_{x/}$ spanned by those objects that are in the image of the degeneracy map $(\sX_0^\fl)_{x/} \to (\sX_n^\fl)_{x/}$, and similarly for $\sX_{/x}$.
	Using \cref{remark:mapping-category-flipped-2-category}, it follows that the mapping categories of $\sX_{x/}$ and $\sX_{/x}$ are given by
	\[\sX_{x/}(p \colon x \to y, q \colon x \to z) = \sX(y,z) \times_{\sX(x,z)} \{q\}\]
	and
	\[\sX_{/x}(p \colon y \to x, q \colon z \to x) = \sX(y,z) \times_{\sX(y,x)} \{p\},\]
	respectively.
\end{example}

\subsubsection{Cocartesian fibrations over 1-categories}

Let us now discuss the theory of cocartesian fibrations in the context of 2-categories, which will be used extensively when discussing symmetric monoidal 2-categories in \cref{sec:Algebras-in-2-categories}. In this paper we will only consider cocartesian fibrations where the base is a 1-category. Since this simplifies the theory significantly, we decided to focus only on this case.

\begin{definition}[Cocartesian fibrations]\label{definition:cocartesian-fibration}
    Let $\sC$ be a 1-category, $\sX$ a 2-category and $p \colon \sX \to \sC$ a functor. We say that a 1-morphism $f \colon x \to y$ in $\sX$ is a \emph{$p$-cocartesian lift} if for any object $z$ in $\sX$,
    \begin{equation}\label{diagram:pbsquare-cocartesian-lift}
        \begin{tikzcd}
            \sX(y,z) \ar[r,"f^*"] \ar[d] & \sX(x,z) \ar[d] \\
            \Map_\sC(p(y), p(z)) \ar[r,"p(f)^*"] & \Map_\sC(p(x),p(z))
        \end{tikzcd}
    \end{equation}
    is a pullback square in $\Cat$. The functor $p$ is a \emph{cocartesian fibration} if for any arrow $f' \colon c \to d$ in $\sC$ and any object $x$ in $\sX$ with $p(x) \simeq c$, there exists a $p$-cocartesian lift of $f'$ with domain $x$. We write $\Cocarttwo(\sC)$ for the category of cocartesian fibrations over $\sC$ and functors between them that preserve cocartesian lifts.
\end{definition}

\begin{remark}
    For the general definition of a cocartesian fibration between 2-categories $p \colon \sX \to \sY$, one needs to add an extra condition: namely that for any pair of objects $x,y$ in $\sX$, the map $\sX(x,y) \to \sY(p(x),p(y))$ is a cartesian fibration of 1-categories and that the composition maps of $\sX$ preserve cartesian lifts, see e.g. \cite[Definition 5.10]{Nuiten2023StraighteningSegalSpaces}. Observe that when $\sY$ is a 1-category, then these conditions are automatically satisfied.
\end{remark}

Using flipped 2-categories, one can deduce that $\Cocarttwo(\cC) \simeq \Fun(\cC,\Cattwo)$ from Lurie's straightening-unstraightening for 1-categories \cite[\S 3.2]{HTT}.

\begin{proposition}\label{proposition:2-cocartesian-straightening}
    Let $\sC$ be a 1-category. Then there is an equivalence $\Cocarttwo(\sC) \simeq \Fun(\sC,\Cattwo)$ that is natural in $\sC$.
\end{proposition}

\begin{proof}
    Using the equivalence $\Cattwo \simeq \Cattwofl \subset \Fun(\Delta^\op,\Cat)$ and swapping $\Delta^\op$ and $\sC$, we can view $\Fun(\sC,\Cattwo)$ as the full subcategory $\sA$ of $\Fun(\Delta^\op, \Fun(\sC,\Cat))$ spanned by those functors $F$ for which $F(-)(c)$ is always a flipped 2-category.
    Now recall the classical straightening equivalence $\Fun(\sC,\Cat) \simeq \Cocart(\sC)$, where $\Cocart(\sC)$ denotes the subcategory of $\Cat_{/\sC}$ spanned by cocartesian fibrations and functors over $\sC$ that preserve cocartesian lifts. Noting that 1-categories can be viewed as constant flipped 2-categories by \cref{remark:constant-flipped-2-categories}, it follows that the full subcategory $\sA$ can be identified with a (non-full) subcategory of $\Fun(\Delta^\op,\Cat)_{/\sC}$. An object $p \colon X_\bullet \to \sC$ lies in this subcategory precisely if
    \begin{enumerate}[(1)]
        \item\label{item1:cocartproof} for every $n \geq 0$, the functor $p_n \colon X_n \to \sC$ is a cocartesian fibration,
        \item\label{item2:cocartproof} the simplicial structure maps of $X_\bullet$ preserve cocartesian lifts, and
        \item\label{item3:cocartproof} for every $c$ in $\sC$, the fiber $X_\bullet \times_{\sC} \{c\}$ is a flipped 2-category,
    \end{enumerate}
    while a map $X_\bullet \to Y_\bullet$ over $\sC$ lies in this subcategory if and only if it levelwise preserves cocartesian lifts.
    
    Using that the functors $p_n$ are cocartesian fibrations and that the simplicial structure maps of $X_\bullet$ preserve cocartesian lifts, it follows that item \ref{item3:cocartproof} is equivalent to $X_\bullet$ being a flipped 2-category. It now follows from the description of the mapping categories of a flipped 2-category in \cref{remark:mapping-category-flipped-2-category} that items \ref{item1:cocartproof}-\ref{item2:cocartproof} are equivalent to the condition that $p \colon X_\bullet \to \sC$ is a cocartesian fibration in the sense of \cref{definition:cocartesian-fibration}. One similarly sees that a map $X_\bullet \to Y_\bullet$ levelwise preserves cocartesian lifts if and only if the corresponding map between 2-categories preserves cocartesian lifts, hence the equivalence $\Fun(\sC,\Cattwo) \simeq \Cocarttwo(\sC)$ follows. Observe that this equivalence is natural in $\sC$ since the classical straightening equivalence $\Fun(\sC,\Cat) \simeq \Cocart(\sC)$ is.
\end{proof}

\begin{remark}
    \Cref{definition:cocartesian-fibration,proposition:2-cocartesian-straightening} have obvious analogues for cartesian fibration.
\end{remark}

For recognizing when a functor is a cocartesian fibration, it will be convenient to also define locally cocartesian fibrations. This will be used in \cref{sec:The-2-category-PrL} for showing that the tensor product of presentable categories makes $\presl$ into a symmetric monoidal 2-category.

\begin{definition}[Locally cocartesian fibrations]
    Let $\cC$ be a 1-category, $\sX$ a 2-category, $p \colon \sX \to \sC$ a functor and $f \colon x \to y$ a 1-morphism in $\sX$. Consider the pullback
    \[\begin{tikzcd}
        {[1]} \times_{\sC} \sX \ar[r] \ar[d,"p'"] \ar[dr, phantom, very near start, "\lrcorner"] & \sX \ar[d,"p"] \\
        {[1]} \ar[r,"p(f)"] & \sC
    \end{tikzcd}\]
    in $\Cattwo$. We say that $f$ is a locally $p$-cocartesian lift if the corresponding 1-morphism in $[1] \times_{\sC} \sX$ is a $p'$-cocartesian lift. We say that $p$ is a locally cocartesian fibration if for $f' \colon c \to d$ in $\sC$ and any object $x$ in $\sX$ with $p(x) \simeq c$, there exists a locally $p$-cocartesian lift of $f'$ with domain $x$.
\end{definition}

\begin{lemma}\label{lemma:When-locally-cocartesian-fibration-is-cocartesian}
    Let $\sC$ be a 1-category, $\sX$ a 2-category and $p \colon \sX \to \sC$ a locally cocartesian fibration. Then $p$ is a cocartesian fibration if and only if for any pair of composable locally $p$-cocartesian lifts $f \colon x \to y$ and $g \colon y \to z$, the composite $gf \colon x \to z$ is again a locally $p$-cocartesian lift.
\end{lemma}

\begin{proof}
    Observe that by a similar argument as in the proof of \cref{proposition:2-cocartesian-straightening}, $p \colon \sX \to \sC$ is a (locally) cocartesian fibration if and only if for every $n \geq 0$, the functor $\sX^\fl_n \to \sC$ is a (locally) cocartesian fibration and the simplicial structure maps of $\sX^\fl$ preserve (locally) cocartesian lifts. Moreover, the locally $p$-cocartesian lifts of $\sX \to \sC$ are closed under composition precisely if for any $n \geq 0$, the locally $p_n$-cocartesian lifts of $p_n \colon \sX_n \to \sC$ are closed under composition. In particular, the result follows from \cite[Proposition 2.4.2.8]{HTT}.
\end{proof}

\subsubsection{Adjunctions between 2-categories}

The functor categories defined in \cref{definition:Functor-2-categories} allow us to define adjunctions between 2-categories in terms of units and/or counits. Let $F,G \colon \sX \to \sY$ be functors between 2-categories. Then a \emph{natural transformation $F \Rightarrow G$}, i.e. a 1-morphism from $F$ to $G$ in $\Fun(\sX, \sY)$, is simply a map $H \colon \sX \times [1] \to \sY$ such that $F = H|_{\sX \times \{0\}}$ and $G = H|_{\sX \times \{1\}}$. In particular, for any $x$ in $\sX$, we obtain a map $F(x) \to G(x)$.

\begin{definition}[Adjunctions between 2-categories]
    Let $\sX$ and $\sY$ be a pair of 2-categories. An \emph{adjunction between $\sX$ and $\sY$} is a pair of functors $F \colon \sX \to \sY$ and $G \colon \sY \to \sX$ together with a \emph{unit transformation} $\eta \colon \id_{\sX} \Rightarrow GF$ such that for any objects $x$ in $\sX$ and $y$ in $\sY$, the composite
    \[
    \sY(F(x),y) \to \sX(GF(x),G(y)) \xrightarrow{\eta^*} \sX(x,G(y))
    \]
    is an equivalence of categories. In this situation, we will call $F$ the \emph{left adjoint} and $G$ the \emph{right adjoint} and write $F \dashv G$.
\end{definition}

Note that $[1]$, when viewed as a flipped 2-category, is simply the constant simplicial object in $\Cat$ with value $[1]$. Now observe that under the straightening equivalence $\Fun(\Delta^\op,\Cat) \simeq \Cocart(\Delta^\op)$, the flipped 2-category $(\sX \times [1])^\fl$ is taken to $(\un \sX^\fl) \times [1]$, where $\un \sX^\fl$ denotes the cocartesian unstraightening of $\sX^\fl$. In particular, a natural transformation $F \Rightarrow G$ between functors of 2-categories $\sX \to \sY$ is equivalent to a natural transformation
\[\begin{tikzcd}[column sep = tiny]
	{\un \sX^\fl} && {\un \sY^\fl} \\
	& {\Delta^\op}
	\arrow[""{name=0, anchor=center, inner sep=0}, "F", curve={height=-9pt}, from=1-1, to=1-3]
	\arrow[""{name=1, anchor=center, inner sep=0}, "G"', curve={height=9pt}, from=1-1, to=1-3]
	\arrow["p"', from=1-1, to=2-2]
	\arrow["q", from=1-3, to=2-2]
	\arrow["\mu", shorten <=2pt, shorten >=2pt, Rightarrow, from=0, to=1]
\end{tikzcd}\]
such that $q(\mu)$ is invertible. From this we can deduce the following.

\begin{proposition}\label{proposition:2-adjunctions-vs-relative-adjunctions}
    Let $\sX$ and $\sY$ be 2-categories. Providing an adjunction between $\sX$ and $\sY$ is equivalent to providing a relative adjunction
    \[\begin{tikzcd}[column sep = tiny]
    	{\un \sX^\fl} && {\un \sY^\fl} \\
    	& {\Delta^\op}
    	\arrow[""{name=0, anchor=center, inner sep=0}, "F", curve={height=-6pt}, from=1-1, to=1-3]
    	\arrow[""{name=1, anchor=center, inner sep=0}, "G"', curve={height=6pt}, tail reversed, no head, from=1-1, to=1-3]
    	\arrow["p"', from=1-1, to=2-2]
    	\arrow["q", from=1-3, to=2-2]
    	\arrow["\dashv"{anchor=center, rotate=-90}, draw=none, from=0, to=1]
    \end{tikzcd}\]
    in the sense of \cite[Definition 7.3.2.2]{HA} such that $F$ and $G$ preserve all cocartesian lifts.
\end{proposition}

\begin{proof}
    First suppose we are given such a relative adjunction. Since $F$ and $G$ preserve cocartesian lifts, by straightening we obtain functors $F' : \sX \rightleftarrows \sY : G'$ and a natural transformation $\id_{\sX} \Rightarrow G'F'$. Since $F$ and $G$ form a relative adjunction, we see that for any $x$ in $\sX$ and any $y$ in $\sY$, the composite
    \[\sY^\fl_n(F'(x),y) \to \sX^\fl_n(G'F'(x),G'(y)) \rightarrow \sX^\fl_n(x,G'(y))\]
    is an equivalence. By \cref{remark:mapping-category-flipped-2-category}, it follows that
    \[\sY(F'(x),y) \to \sX(G'F'(x),G'(y)) \rightarrow \sX(x,G'(y))\]
    is an equivalence.

    Conversely, suppose we are given an adjunction $F' \dashv G'$ between 2-categories $\sX$ and $\sY$. Then $F$ and $G$ induce an adjunction between $\sX^\fl_n$ and $\sY^\fl_n$ for every $n$. The desired relative adjunction is now obtained from the dual of \cite[Proposition 7.3.2.6]{HA}.
\end{proof}

From this result we obtain a useful criterion for detecting when a functor $F \colon \sX \to \sY$ admits a right adjoint. Suppose
\begin{equation}\label{eq:Adjointable-square}
    \begin{tikzcd}
        \sC \ar[r,"f"] \ar[d,"h"'] & \sD \ar[d,"k"] \\
        \sC' \ar[r,"f'"'] & \sD'
    \end{tikzcd}
\end{equation}
is a commutative square in $\Cat$ where $f$ and $f'$ admit right adjoints $g$ and $g'$, respectively. Using the unit of $f' \dashv g'$ and the counit of $f \dashv g$, we obtain a lax commuting square
\[\begin{tikzcd}
    \sC \ar[d,"h"'] \ar[dr,Rightarrow, shorten <=7pt, shorten >=7pt,] & \sD \ar[d,"k"] \ar[l,"g"'] \\
    \sC' & \sD' \ar[l,"g'"]
\end{tikzcd}\]
A square of the form \cref{eq:Adjointable-square} is called \emph{right adjointable} if both $f$ and $f'$ have right adjoints and this lax square commutes strictly; that is, if the composite $hg \Rightarrow g'f'hg \simeq g'kfg \Rightarrow g'k$ is an invertible natural transformation.

\begin{proposition}\label{proposition:2-adjunctions-as-levelwise-adjunctions}
    Let $F \colon \sX \to \sY$ be a functor between 2-categories. Then $F$ admits a right adjoint if and only if for every $n \geq 0$, the functor $F_n \colon \sX^\fl_n \to \sY_n^\fl$ admits a right adjoint and for every $\alpha \colon [n] \to [m]$ in $\Delta$, the square
    \begin{equation}\label{eq:Adjointable-condition-2-adjunction}
        \begin{tikzcd}
            \sX^\fl_m \ar[r,"F_m"] \ar[d,"\alpha^*"] & \sY^\fl_m \ar[d,"\alpha^*"] \\
            \sX^\fl_n \ar[r,"F_n"] & \sY^\fl_n
        \end{tikzcd}
    \end{equation}
    is right adjointable.
\end{proposition}

\begin{proof}
    It follows from the dual of \cite[Proposition 7.3.2.6]{HTT} that $\un F \colon \un \sX^\fl \to \un \sY^\fl$ admits a relative right adjoint, so it suffices to show that this relative right adjoint preserves cocartesian lifts. We leave it to the reader to verify that this is the case precisely if all squares of the form \cref{eq:Adjointable-condition-2-adjunction} are right adjointable.
\end{proof}

\begin{remark}
    It follows dually that if $G \colon \sY \to \sX$ is a functor between 2-categories such that for every $n \geq 0$, the functor $G_n \colon \sY^\fl_n \to \sX_n^\fl$ admits a left adjoint and for every $\alpha \colon [n] \to [m]$ in $\Delta$, the square
    \begin{equation}
        \begin{tikzcd}
            \sY^\fl_m \ar[r,"G_m"] \ar[d,"\alpha^*"] & \sX^\fl_m \ar[d,"\alpha^*"] \\
            \sY^\fl_n \ar[r,"G_n"] & \sX^\fl_n
        \end{tikzcd}
    \end{equation}
    is left adjointable, then $G$ admits a left adjoint.
\end{remark}

\begin{remark}
    In fact, the assumptions of \cref{proposition:2-adjunctions-as-levelwise-adjunctions} can be weakened. Because of the Segal conditions on $\sX^\fl$ and $\sY^\fl$ and the fact that $(\sX^\fl_\bullet)^\simeq$ and $(\sY^\fl_\bullet)^\simeq$ are constant, one can show that it suffices for the square
    \[\begin{tikzcd}
            \sX^\fl_0 \ar[r,"F_0"] \ar[d,"s"] & \sY^\fl_0 \ar[d,"s"] \\
            \sX^\fl_1 \ar[r,"F_1"] & \sY^\fl_1
        \end{tikzcd}\]
    to be right adjointable.
\end{remark}

Recall that to specify a left adjoint to a functor $G \colon \sD \to \sC$ between 1-categories, it suffices to specify, for any $c$ in $\sC$, an object $F(c)$ in $\sD$ together with a map $c \to GF(c)$ in $\sC$ such that for any $d$ in $\sD$, the composite
\[\sD(F(c),d) \to \sC(GF(c),G(d)) \to \sC(c,G(d))\]
is an equivalence. Namely, one of the many equivalent ways to define adjunctions between 1-categories is as bicartesian fibrations over $\Delta^1$, see \cite[Definition 5.2.2.1]{HTT}. If we unstraighten $G \colon \sD \to \sC$ to a cartesian fibration $\cM \to \Delta^1$, then the condition stated above is precisely the condition that this functor $\cM \to \Delta^1$ is also a cocartesian fibration. It follows as in the proof of \cite[Proposition 5.2.2.8]{HTT} that the maps $c \to GF(c)$ can be refined to become the unit natural transformation $\id_{\sC} \Rightarrow GF$ of an adjunction $F \dashv G$. We will now prove the analogous result for 2-categories.

\begin{corollary}\label{corollary:left-adjoints-constructed-locally}
    Let $G \colon \sY \to \sX$ be a functor between 2-categories and suppose that for every $x$ in $\sX$, there exists an object $F(x)$ in $\sY$ and a map $x \to GF(x)$ such that for every $y$ in $\sY$, the composite
    \[\sY(F(x),y) \to \sY(GF(x),G(y)) \to \sY(x,G(y))\]
    is an equivalence. Then $F$ can be refined to a functor $\sX \to \sY$ that is left adjoint to $G$.
\end{corollary}

\begin{proof}
    The 1-categorical case shows that for every $n \geq 0$, the functor $G_n \colon \sY^\fl_n \to \sX^\fl_n$ admits a left adjoint $F_n$. We leave it to the reader to verify that the conditions of (the dual of) \cref{proposition:2-adjunctions-as-levelwise-adjunctions} are satisfied, hence that the functors $F_n$ assemble into a left adjoint $F \colon \sX \to \sY$ of $G$.
\end{proof}

\subsubsection{Colimits in 2-categories}

Having defined functor 2-categories and adjunctions between 2-categories, we can define (co)limits as follows.

\begin{definition}[(Co)limits in 2-categories]\label{definition:colimits-in-2-categories}
    Let $I$ and $\sX$ be 2-categories. Restriction along $I \to *$ induces a functor $\Delta \colon \sX \to \Funtwo(I,\sX)$. We say that \emph{$\sX$ admits colimits or limits indexed by $I$} if this functor admits a left or right adjoint, respectively. In this case we write $\colim_I$ for the left adjoint and $\lim_I$ for the right adjoint.
\end{definition}

For our study of free algebras in \cref{subsec:Free-algebras-in-2-category}, we will need to have a good understanding of colimits in 2-categories indexed by spaces. Note that we can view spaces as 2-categories through the inclusions $\Spc \hookrightarrow \Cat \hookrightarrow \Cattwo$.

\begin{definition}[Groupoidally (co)complete 2-categories]\label{definition:Groupoidally-cocomplete-category}
    A 2-category $\sX$ is said to be \emph{groupoidally (co)complete} if it admits all (co)limits indexed by small spaces.
\end{definition}

The following is a direct consequence of \cref{proposition:2-adjunctions-as-levelwise-adjunctions} and its dual.

\begin{proposition}\label{proposition:Characterization-groupoidally-cocomplete-categories}
    Let $\sX$ be a 2-category. Then $\sX$ is groupoidally (co)complete if and only if for every $n \geq 0$, the category $\sX^\fl_n$ admits (co)limits indexed by small spaces and for any $\alpha \colon [m] \to [n]$, the map $\alpha^* \colon \sX^\fl_n \to \sX^\fl_m$ preserves these.
\end{proposition}

\begin{proof}
    If $I$ is a space, then by \cref{example:flipped-internal-hom} it follows that $\Funtwo(I,\sX)^\fl_n = \Fun(I,\sX^\fl_n)$. Applying \cref{proposition:2-adjunctions-as-levelwise-adjunctions} (or its dual) now yields the result.
\end{proof}

\subsubsection{Tensors and cotensors in 2-categories}

In many cases, adjunctions and (co)limits can be detected on the level of the underlying 1-categories. The reason for this is the existence of certain (co)tensors. In what follows, we write $\Map(x,y)$ for the groupoid core $\sX(x,y)^\simeq$ of $\sX(x,y)$. Observe that the composition map of a 2-category $\sX$ induces a map $\Map(y,z) \times \sX(x,y) \to \sX(x,z)$ and hence a map $\Map(y,z) \to \Map(\sX(x,y),\sX(x,z))$. Similarly, we obtain a map $\Map(x,y) \to \Map(\sX(y,z),\sX(x,z))$.

\begin{definition}[(Co)tensors in 2-categories]\label{definition:Tensors-and-cotensors}
    Let $\sX$ be a 2-category and $\sC$ a 1-category. Then $\sX$ is said to \emph{admit (weak) tensors with $\sC$} if for any object $x$ in $\sX$, there exists an object $x \otimes \sC$ together with a map $\sC \to \sX(x, x \otimes \sC)$ such that for every $y$ in $\sX$, the composite
    \[\Map(x \otimes \sC, y) \to \Map(\sX(x, x\otimes\sC), \sX(x,y)) \to \Map(\sC, \sX(x,y))\]
    is an equivalence. Dually, $\sX$ \emph{admits (weak) cotensors with $\sC$} if for any object $y$ in $\sX$, there exists an object $y^\sC$ together with a map $\sC \to \sX(y^\sC, y)$ such that for every $x$ in $\sX$, the composite
    \[\Map(x, y^\sC) \to \Map(\sX(y^\sC,y), \sX(x,y)) \to \Map(\sC, \sX(x,y))\]
    is an equivalence.
\end{definition}

\begin{example}
    Any 2-category admits tensors and cotensors with $*$.
\end{example}

\begin{example}
    The 2-category $\Cat$ has all tensors and cotensors, given by the cartesian product and functor categories, respectively.
\end{example}

\begin{remark}
    The reason we call the above \emph{weak} (co)tensors is that in the usual definition of a (co)tensor, one asks for equivalences
    \[\sX(x \otimes \sC, y) \simeq \Fun(\sC,\sX(x,y)) \simeq \sX(x, y^\sC),\]
    and not just equivalences on the level of groupoid cores. However, for our purposes the weaker definition will suffice.
\end{remark}

As mentioned above, the main reason for introducing (co)tensors is the following result. Observe that if $G \colon \sY \to \sX$ is a functor between 2-categories that admit (weak) cotensors by $\sC$, then the composite $\sC \to \sY(y^{\sC}, y) \to \sX(G(y^{\sC}), G(y))$ induces a canonical map $G(y^{\sC}) \to G(y)^{\sC}$. We say that $G$ preserves cotensors by $\sC$ if for every $y$ in $\sY$, this map is an equivalence.

\begin{proposition}\label{proposition:2-adjunctions-detected-underlying}
    Let $\sX$ and $\sY$ be 2-categories that admit cotensors by $[1]$, and let $G \colon \sY \to \sX$ be a functor of 2-categories such that the underlying functor $G^{\leq 1} \colon \sY^{\leq 1} \to \sX^{\leq 1}$ of 1-categories admits a left adjoint. Then $G$ admits a left adjoint if and only if $G$ preserves cotensors by $[1]$.
\end{proposition}

\begin{proof}
    Let $F$ denote the left adjoint of $G^{\leq 1}$ and $\eta_x \colon x \to GFx$ the unit of this adjunction at $x$. Let us first assume that $G$ preserves cotensors by $[1]$. By \cref{corollary:left-adjoints-constructed-locally}, it suffices to show that the composite
    \[\sX(Fx,y) \xrightarrow{G} \sY(GFx, Gy) \xrightarrow{\eta_x^*} \sY(x,Gy)\]
    is an equivalence. By \cite[\href{https://kerodon.net/tag/01HG}{Theorem 01HG}]{kerodon}, it suffices to show that $\Map([1],\sX(Fx,y)) \to \Map([1],\sY(x,Gy))$ is an equivalence. A (somewhat tedious) diagram chase shows that
    \[\begin{tikzcd}
    	{\Map(Fx,y^{[1]})} & {\Map(x,G(y^{[1]}))} & {\Map(x,G(y)^{[1]})} \\
    	{\Map([1],\sX(Fx,y))} && {\Map([1],\sY(x,Gy))}
    	\arrow["\sim", from=1-1, to=1-2]
    	\arrow["\sim"', from=1-1, to=2-1]
    	\arrow[from=1-2, to=1-3]
    	\arrow["\sim", from=1-3, to=2-3]
    	\arrow[from=2-1, to=2-3]
    \end{tikzcd}\]
    commutes, where the top right horizontal arrow comes from the canonical map $G(y^{[1]}) \to G(y)^{[1]}$ constructed above. Since this map is an equivalence by assumption, the result follows. Conversely, if $G$ admits a left adjoint $F$, then the same diagram shows that for any $x \in \sX$, the map $\Map(x,G(y^{[1]})) \to \Map(x,G(y)^{[1]})$ is an equivalence. Since equivalences in a 2-category are detected in the underlying 1-category, the 1-categorical Yoneda lemma implies that $G(y^{[1]}) \to G(y)^{[1]}$ is an equivalence.
\end{proof}

\begin{remark}\label{remark:underlying-adjoint-is-adjoint}
    The proof of \cref{proposition:2-adjunctions-detected-underlying} moreover shows that the left adjoint of $G$ is a 2-categorical refinement of the left adjoint of $G^{\leq 1}$.
\end{remark}

\begin{corollary}\label{corollary:cotensors-detect-colimits-underlying}
    Let $\sX$ be a 2-category that admits cotensors by $[1]$ and let $I$ be a 1-category. Then $I$-indexed colimits in $\sX$ can be computed in the underlying category $\sX^{\leq 1}$.
\end{corollary}

\begin{proof}
    Combine \cref{definition:colimits-in-2-categories,proposition:2-adjunctions-detected-underlying,remark:underlying-adjoint-is-adjoint}.
\end{proof}

We also obtain dual versions of \cref{proposition:2-adjunctions-detected-underlying,corollary:cotensors-detect-colimits-underlying}.

\begin{proposition}\label{proposition:2-adjunctions-detected-underlying-tensor}
    Let $\sX$ and $\sY$ be 2-categories that admit tensors by $[1]$, and let $F \colon \sX \to \sY$ be a functor of 2-categories such that the underlying functor $F^{\leq 1} \colon \sX^{\leq 1} \to \sY^{\leq 1}$ of 1-categories admits a right adjoint. Then $F$ admits a right adjoint if and only if $F$ preserves tensors by $[1]$.
\end{proposition}

\begin{corollary}
\label{corollary:tensors-detect-limits-underlying}
    Let $\sX$ be a 2-category that admits tensors by $[1]$ and let $I$ be a 1-category. Then $I$-indexed limits in $\sX$ can be computed underlying in $\sX^{\leq 1}$.
\end{corollary}

\subsubsection{Representable functors}

In \cref{subsec:Lambda}, we define a certain 2-functor $\Lambda$ as a corepresentable functor.
There are several places in the literature where (co)representable functors of 2-categories are constructed, see for example \cite[Proposition 4.1.6]{Lurie2009InftyCategoriesGoodwillie}, \cite[\S 6.2.4]{Hinich2020YonedaLemmaEnriched}, or \cite[\S 4]{Nuiten2023StraighteningSegalSpaces}.
Here we state the basic properties of corepresentable functors that we will need, following \cite{Nuiten2023StraighteningSegalSpaces}.

\begin{proposition} \label{prop: corepresentable-functors-of-2-cats}
    Let $\sX$ be a 2-category and $x \in \sX_0$ an object.
    Then there exists a functor $\sX \to \Cat$ of 2-categories which is given on objects by $y \mapsto \sX(x,y)$, while on mapping categories it is given by the functor
    \[\begin{tikzcd}\sX(y,z) \ar[r] & \Fun(\sX(x,y),\sX(x,z))\end{tikzcd}\]
    adjoint to the composition $\circ \colon \sX(y,z) \times \sX(x,y) \to \sX(x,z)$.
\end{proposition}

\begin{proof}
    In \cite[Example 4.36]{Nuiten2023StraighteningSegalSpaces}, the functor
    \[\begin{tikzcd}
        \sX(-,-) \colon \sX^\op \times \sX \ar[r] &  \Cat;& (x,y) \ar[r,maps to]& \sX(x,y)
    \end{tikzcd}\]
    is constructed.
    Restriction along $\{x\} \times \sX \to \sX^\op \times \sX$ gives the desired functor.
    An unwinding of the proof of \cite[Corollary 4.33]{Nuiten2023StraighteningSegalSpaces} shows that it is of the stated form.
\end{proof}

\subsection{Lax functors between 2-categories}\label{sec:Lax-functors-between-2-categories}

In this section, we will study lax functors between 2-categories and describe several methods for constructing them. As mentioned at the start of this appendix, one possible perspective on 2-categories is to view them as a collection of monoidal categories $\{\sV_s\}_{s \in S}$ together with a collection of bimodule categories $\{\sM_{s,t}\}_{s,t \in S}$, parametrized by a space $S$. It will turn out that lax functors between 2-categories behave very much like lax monoidal functors between monoidal categories and bimodule categories, and many proofs will take a very similar form.

When taking such an ``algebraic'' perspective on 2-categories, the completeness condition is slightly awkward. It will turn out that when working with lax functors between 2-categories, it is indeed more convenient to work in the slightly more general setting of 2-precategories. For this reason, all the definitions and statements will be formulated in this more general context.

\subsubsection{Extracting bimodules from 2-categories}\label{subsubsection:extracting-bimodules}

As mentioned, one can informally think of a 2-(pre)category $\sX$ as a collection of monoidal categories $\sX(x,x)$ together with bimodule categories $\sX(x,y)$ and bilinear maps between them satisfying certain coherences, parametrized by the space $\sX_0$. We will now justify this perspective by showing that for any 2-precategory $\sX$, composition of 1-morphisms endows the categories $\sX(x,x)$ with a monoidal structure and the categories $\sX(x,y)$ with an $(\sX(y,y),\sX(x,x))$-bimodule structure.

For what follows, recall from \cref{example:coskeleton} that $\ev_0 \colon \Fun(\Delta^\op,\Cat) \to \Cat$ admits a right adjoint $\cosk_0$, the \emph{coskeleton functor}. Given an object $x$ in $\sX$, let $\overline{x} \colon * \to \cosk_0\sX_0$ denote the map corresponding to $x \colon * \to \sX_0$.

\begin{definition}[Monoidal endomorphism categories] \label{def:endomorphisms-infty-category}
	Let $\sX$ be a 2-precategory. For any object $x \in \sX_0$, the \emph{(monoidal) endomorphism category} $\sX(x,x)$ is defined by the pullback
	\[\begin{tikzcd}
		\sX(x,x) \ar[r,] \ar[d,] & \sX \ar[d,"\eta"]\\
		* \ar[r,"\overline x"] & \cosk_0 \sX_0
	\end{tikzcd}\]
	in $\Fun(\Delta^{\op},\Catinfty)$, where $\eta$ denotes the adjunction unit.
\end{definition}

Since $\cosk_0 \sX_0$ is a Segal object by \cref{example:coskeleton}, it follows that this pullback is indeed a monoidal category. The notation $\sX(x,x)$ for this monoidal category is slightly abusive, but justified by the fact that its underlying category is indeed the category $\sX(x,x)$ defined in \cref{def:mapping-infty-category}. We leave it to the reader to verify that the monoidal structure on $\sX(x,x)$ is indeed given by the composition functor from \cref{remark:composition-functor-2-category} and that its unit is $\unit_x$.

We now wish to endow $\sX(x,y)$ with the structure of a $(\sX(y,y),\sX(x,x))$-bimodule. Recall that the bimodule operad $\mathscr{BM}^\otimes$ is defined in \cite[Definition 4.3.1.6]{HA} as a certain $\infty$-operad with three colours $\mathfrak{a}_-$, $\mathfrak{m}$ and $\mathfrak{a}_+$. We will use a simplified version of the bimodule operad, namely the category $(\Delta_{/[1]})^\op$. In what follows, we write $[m|n]$ for the object of $\Delta_{/[1]}$ given by
\[\sigma \colon [m+n+1] \to [1]; \quad \sigma(i) =
\begin{cases}
    0 \quad &\text{if} \quad i \leq m\\
    1 \quad &\text{if} \quad i > m.
\end{cases}\]
Note that $m=-1$ and $n=-1$ are allowed, though not both at the same time.

\begin{proposition}\label{proposition:approximation-to-BM-operad}
    Let $\phi \colon (\Delta_{/[1]})^\op \to \mathscr{BM}^\otimes$ be the canonical functor given on objects by
    \[[m|n] \mapsto \begin{cases}
    (\underbrace{\mathfrak{a}_-,\ldots,\mathfrak{a}_-}_{m\text{-times}},\mathfrak{m},\underbrace{\mathfrak{a}_+,\ldots,\mathfrak{a}_+}_{n\text{-times}}) \quad &\text{if} \quad m,n \geq 0\\
    (\underbrace{\mathfrak{a}_-,\ldots,\mathfrak{a}_-}_{m\text{-times}}) \quad &\text{if} \quad n = -1\\
    (\underbrace{\mathfrak{a}_+,\ldots,\mathfrak{a}_+}_{n\text{-times}}) \quad &\text{if} \quad m = -1.
    \end{cases}\]
    Then $\phi$ is an approximation to $\mathscr{BM}^\otimes$ in the sense of \cite[Definition 2.3.3.6]{HA}.
\end{proposition}

\begin{proof}
    Combine Proposition 7.1.7 and Remark 7.1.9 of \cite{GepnerHaugseng2015EnrichedCategoriesNonsymmetric}.
\end{proof}

Recall that a bimodule category is defined as a $\mathscr{BM}^\otimes$-monoidal category in the sense of \cite[Definition 2.1.2.13]{HA}. By \cite[Example 2.4.2.4 \& Proposition 2.4.2.11]{HA}, there is an equivalence between $\mathscr{BM}^\otimes$-monoidal categories and Segal objects indexed by $(\Delta_{/[1]})^\op$, which is given by precomposition with $\phi \colon (\Delta_{/[1]})^\op \to \mathscr{BM}^\otimes$. We will therefore use Segal objects indexed by $(\Delta_{/[1]})^\op$ as our definition of a bimodule category. In order to state this definition precisely, observe that for a functor $F \colon (\Delta_{/[1]})^\op \to \Cat$, the maps $\{i-1,i\} \hookrightarrow [m+n+1]$ induce \emph{Segal maps}
\begin{equation}\label{eq:Segal-maps-bimodule-operad}
    F([m|n]) \to \begin{cases}
    F([1|-1])^{\times m} \times F([0|0]) \times F([-1|1])^{\times n} \quad &\text{if} \quad m,n \geq 0\\
    F([1|-1])^{\times m} \quad &\text{if} \quad n = -1\\
    F([-1|1])^{\times m} \quad &\text{if} \quad m = -1.
\end{cases}
\end{equation}
In the cases $[m|n] = [0|-1], [-1|0]$, we interpret the target of the Segal map as the terminal category $*$.
Also observe that there are canonical inclusions $\iota_{-}, \iota_{+} \colon \Delta^\op \hookrightarrow (\Delta_{/[1]})^\op$ sending $[n]$ to $[n|-1]$ and $[-1|n]$, respectively.

\begin{definition}[Bimodule categories]\label{definition:bimodule-categories}
    Let $\sC$ and $\sD$ be monoidal categories. Then a \emph{$(\sC,\sD)$-bimodule (category)} or \emph{$(\sC,\sD)$-bitensored category} is a functor $\sM \colon (\Delta_{/[1]})^\op \to \Cat$ together with equivalences $\sM \circ \iota_{-} \simeq \sC$ and $\sM \circ \iota_{+} \simeq \sD$ such that the Segal maps \cref{eq:Segal-maps-bimodule-operad} are equivalences. Given a $(\sC,\sD)$-bimodule category $\sM$, we will slightly abuse notation and write $\sM$ for the category $\sM([0|0])$ as well.
\end{definition}

\begin{example}\label{example:action-maps-of-bimodule}
    If $\sM \colon (\Delta_{/[1]})^\op \to \Cat$ is a $(\sC,\sD)$-bimodule object, then the \emph{action maps} $\sC \times \sM \to \sM$ and $\sM \times \sD \to \sM$ are given by the composites
    \[\begin{tikzcd}
        \sC \times \sM \simeq \sM([1|0]) \ar[r,"d_1"] & \sM([0|0]) \simeq \sM
    \end{tikzcd}\]
    and
    \[\begin{tikzcd}
        \sM \times \sD \simeq \sM([0|1]) \ar[r,"d_1"] & \sM([0|0]) \simeq \sM,
    \end{tikzcd}\]
    respectively.
\end{example}

\begin{remark}
    Let $\sN^\otimes \to \mathscr{BM}^\otimes$ be a $\mathscr{BM}^\otimes$-monoidal category, let $\sM \colon (\Delta_{/[1]})^\op \to \Cat$ be the corresponding bimodule category and let $q \colon \phi^* \sN^\otimes \simeq \un \sM \to (\Delta_{/[1]})^\op$ be its cocartesian unstraightening. It follows from \cite[Theorem 2.3.3.23]{HA} (and the universal property of the pullback $\phi^* \sN^\otimes$) that the category of bimodule objects in $\sN^\otimes$ is equivalent to the category of sections of $q$ that take inerts to cocartesian lifts. Here a map in $(\Delta_{/[1]})^\op$ is called inert if the underlying map in $\Delta^\op$ is inert.
\end{remark}

We will now describe how to obtain bimodule structures on the mapping categories $\sX(x,y)$ of a 2-precategory $\sX$. Let $p \colon (\Delta_{/[1]})^\op \to \Delta^\op$ denote the forgetful functor. Given a 2-precategory $\sX$, we obtain the functor $p^* \sX \colon (\Delta_{/[1]})^\op \to \Cat$. Write $\cosk_{0,0}$ for the right adjoint to the functor $\Fun((\Delta_{/[1]})^\op, \Catinfty) \to \Cat^{\times 2}$ that sends $F \colon (\Delta_{/[1]})^\op \to \Catinfty$ to the pair $(F([0|{-1}]),F([-1|0]))$. Let $\overline{(y,x)} \colon  * \to \cosk_{0,0}(\sX_0,\sX_0)$ denote the map corresponding to $(y,x) \colon (*,*) \to (\sX_0,\sX_0)$ under this adjunction and let $\eta \colon p^* \sX \to \cosk_{0,0}(\sX_0,\sX_0)$ be the map corresponding to the identity $(\sX_0,\sX_0) \to (\sX_0,\sX_0)$ under this adjunction.

\begin{proposition}[Bimodule mapping categories]\label{proposition:bimodule-mapping-categories}
    Let $\sX$ be a 2-precategory. For any $x,y \in \sX_0$, the pullback
    \[\begin{tikzcd}
    	P & {p^*\sX} \\
    	{*} & {\cosk_{0,0} (\sX_0,\sX_0)}
    	\arrow[from=1-1, to=2-1]
    	\arrow[from=2-1, to=2-2,"{\overline{(y,x)}}"]
    	\arrow[from=1-2, to=2-2,"\eta"]
    	\arrow[from=1-1, to=1-2]
    	\arrow["\lrcorner"{anchor=center, pos=0.125}, draw=none, from=1-1, to=2-2]
    \end{tikzcd}\]
    in $\Fun((\Delta_{/[1]})^\op,\Cat)$ defines an $(\sX(y,y),\sX(x,x))$-bimodule structure on $\sX(x,y)$. Moreover, the actions of $\sX(y,y)$ and $\sX(x,x)$ on $\sX(x,y)$ are given by the composition functor of \cref{remark:composition-functor-2-category}.
\end{proposition}

\begin{proof}
    It is easily verified that $P$ satisfies the Segal condition and hence that it is a bimodule category. Observe that for any $\sC$ and $\sD$ in $\Cat$, we have $\cosk_{0,0}(\sC,\sD)([m|n]) \simeq \sC^{\times m+1} \times \sD^{\times n+1}$. Applying this to $\sC = \sD = \sX_0$, we see that $P([0|0]) = \sX(x,y)$. Moreover, it follows that the restriction of $\cosk_{0,0}(\sX_0,\sX_0)$ along $\iota_+$ and $\iota_-$ recovers the simplicial object $\cosk_0(\sX_0)$. We therefore see that $P \circ \iota_- \simeq \sX(y,y)$ and that $P \circ \iota_+ \simeq \sX(x,x)$. Finally, it follows from \cref{example:action-maps-of-bimodule} that the actions of $\sX(y,y)$ and $\sX(x,x)$ on $\sX(x,y)$ are given by composition in $\sX$.
\end{proof}

\begin{example}
    Let $\sC$ be a monoidal category, viewed as a 2-precategory with space of objects given by $*$. Then \cref{proposition:bimodule-mapping-categories} exhibits $\sC$ as a bimodule over itself.
\end{example}

\begin{remark}\label{remark:simplified-left-right-module-operads}
    The category $(\Delta_{/[1]})^\op$ admits two notable full subcategories that can be used to describe left and right modules. Namely, one can consider the full subcategory $\LM$ spanned by objects of the form $[m|i]$ where $i = -1,0$ and $m \geq 0$, or the full subcategory $\RM$ spanned by objects of the form $[i|n]$ with $n \geq 0$ and $i = -1,0$. It is easy to see that both $\LM$ and $\RM$ are equivalent to $\Delta^\op \times [1]$, and it follows as in \cite[Remark 4.2.2.8]{HA} or \cite[Proposition 7.1.7 \& Remark 7.1.9]{GepnerHaugseng2015EnrichedCategoriesNonsymmetric} that the functor $\phi \colon (\Delta_{/[1]})^\op \to \mathscr{BM}^\otimes$ restricts to approximations $\LM \to \mathscr{LM}^\otimes$ and $\RM \to \mathscr{RM}^\otimes$ in the sense of \cite[Definition 2.3.3.6]{HA}. In particular, the categories $\LM$ and $\RM$ come with a notion of inert maps and left and right modules can be described as algebras over $\LM$ and $\RM$, respectively. Moreover, given a bimodule defined as an algebra over $(\Delta_{/[1]})^\op$, one obtains the corresponding left and right modules by restricting along the inclusions $\LM,\RM \to (\Delta_{/[1]})^\op$.
\end{remark}

\subsubsection{Lax functors}

We will now study lax functors between 2-(pre)categories and discuss their basic properties. Let $\un \sX \to \Delta^{\op}$ denote the cocartesian fibration corresponding to a 2-precategory $\sX \colon \Delta^\op \to \Cat$ under Lurie's straightening-unstraightening equivalence. Under this equivalence, a functor of 2-precategories corresponds to a functor between cocartesian fibrations that preserves cocartesian lifts. We define lax functors by weakening this condition. Recall that a map $\alpha \colon [n] \to [m]$ in $\Delta^{\op}$ is called \emph{inert} if it is the inclusion of a subinterval; i.e. if there exists a $0 \leq i \leq m - n$ such that $\alpha$ is the map $j \mapsto i+j$.

\begin{definition}[Lax functors]\label{definition:lax-functors}
	Let $\sX$ and $\sY$ be 2-precategories and denote their cocartesian unstraightenings by $\un \sX \to \Delta^{\op}$ and $\un \sY \to \Delta^{\op}$, respectively. A \emph{lax functor} $F \colon \sX \laxfun \sY$ is defined as a commutative diagram
	\[\begin{tikzcd}[column sep=small]
		\un \sX \ar[rr,"F"] \ar[dr] & & \un \sY \ar[dl] \\
		 & \Delta^{\op} &
	\end{tikzcd}\]
	such that $F$ preserves cocartesian lifts of inerts. Dually, an \emph{oplax functor} $F \colon \sX \oplaxfun \sY$ is defined as a functor between the cartesian unstraightenings of $\sX$ and $\sY$ that preserves cartesian lifts of inerts.
\end{definition}

To differentiate (op)lax functors from the functors defined in \cref{definition:twofold-Segal-space}, we will sometimes call the latter \emph{strong functors}. In what follows, we will restrict our attention to lax functors, but every statement admits an obvious dual for oplax functors.

\begin{remark}
	What we call a lax functor here is called a non-unital right-lax functor in \cite[\S 3.1]{GaitsgoryRozenblyum2017Volume1}.
\end{remark}

Part of the data of a lax functor $F \colon \sX \laxfun \sY$ is a collection of functors $F_n \colon \sX_n \to \sY_n$ for every $n$. These functors are natural in the inert maps, but for other maps they are only natural in a lax sense. More precisely, if $t \colon [n] \to [m]$ is not inert, then the square
\begin{equation}\label{diagram:Non-invertible-cells-lax-functor}
	\begin{tikzcd}
		{\sX_m} & {\sY_m} \\
		{\sX_n} & {\sY_n}
		\arrow["{t^*}"', from=1-1, to=2-1]
		\arrow[""{name=0, anchor=center, inner sep=0}, "{F_m}", from=1-1, to=1-2]
		\arrow["{t^*}", from=1-2, to=2-2]
		\arrow[""{name=1, anchor=center, inner sep=0}, "{F_n}"', from=2-1, to=2-2]
		\arrow[shorten <=8pt, shorten >=8pt, Rightarrow, "{\alpha_t}", from=1-2, to=2-1]
	\end{tikzcd}
\end{equation}
only needs to commute up to a non-invertible natural transformation $\alpha_t$. Suppose $x \to t^* x$ and $F_m(x) \to t^*F_m(x)$ are cocartesian lifts of $t$. Then the universal property of the latter cocartesian lift provides a map $\alpha_{t,x} \colon t^*F_n(x) \to F_n(t^*x)$, which is precisely the component of $\alpha_t$ at $x$. Note that it follows from this description that $F$ takes $x \to t^*x$ to a cocartesian lift if and only if $\alpha_{t,x}$ is an equivalence. In particular, $F$ preserves cocartesian lifts of $t$ precisely if $\alpha_t$ is invertible.

Given an object $x$ and a 1-cell $f$ of $\sX$, we will abusively write $F(x)$ and $F(f)$ for $F_0(x)$ and $F_1(f)$, respectively.

\begin{example}\label{example:Lax-functors-out-of-point}
	Let $*$ be the 2-precategory with exactly one object and only identity morphisms; that is, $*$ is the constant simplicial category with value the terminal category. Then a lax functor $* \laxfun \sX$ into another 2-precategory is the same as a map $\Delta^{\op} \to \un \sX$ over $\Delta^{\op}$ preserving cocartesian lifts of inerts; that is, a section of the cocartesian fibration $\un \sX \to \Delta^{\op}$ sending inerts to cocartesian lifts. Unwinding this, one sees that this corresponds precisely to the data of an object $x \in \sX_0$ together with an associative algebra in the monoidal category $\sX(x,x)$ from \cref{def:endomorphisms-infty-category}.
\end{example}

\begin{remark}
	Note that part of the data of a lax functor $F \colon \sX \laxfun \sY$ is the map $F_0 \colon \sX_0 \to \sY_0$. The completeness condition therefore ensures that $F$ takes equivalent objects of $\sX$ to equivalent objects of $\sY$. This is in contrast to the definition of a lax functor between ordinary $2$-categories, where this generally need not be the case. However, even though $F$ takes equivalent objects to equivalent objects, $F(f)$ does not need to be an equivalence when $f$ is. Namely, the square
	\[\begin{tikzcd}
		{\sX_0} & {\sY_0} \\
		{\sX_1} & {\sY_1}
		\arrow["{s_0}"', from=1-1, to=2-1]
		\arrow[""{name=0, anchor=center, inner sep=0}, "{F_0}", from=1-1, to=1-2]
		\arrow["{s_0}", from=1-2, to=2-2]
		\arrow[""{name=1, anchor=center, inner sep=0}, "{F_1}"', from=2-1, to=2-2]
		\arrow[shorten <=8pt, shorten >=8pt, Rightarrow, from=1-2, to=2-1]
	\end{tikzcd}\]
	only needs to commute up to a non-invertible natural transformation.
\end{remark}

A functor $\un \sX \to \un \sY$ induces a strong functor $\sX \to \sY$ precisely if it preserves cocartesian lifts. As explained below diagram \cref{diagram:Non-invertible-cells-lax-functor}, this is equivalent to the natural transformations $\alpha_t$ being invertible. There are two special cases of these natural transformations $\alpha_t$: when $t$ is $\delta^1 \colon [1] \to [2]$ and when $t$ is $\sigma^0 \colon [1] \to [0]$. For the first case, note that an object in $\sX_2$ consists of three elements $x,y,z \in \sX_0$ together with 1-cells $f \in \sX(x,y)$ and $g \in \sX(y,z)$. The component of $\alpha_t$ at this object then defines a map
\begin{equation}\label{map:composition-2-cell-of-lax-functor}
	Fg \circ Ff \to F(gf)
\end{equation}
in $\sY(F(x),F(z))$, witnessing that $F$ preserves composition up to a (possibly non-invertible) 2-cell. Similarly, the case where $t = \sigma^0$ gives us for every object $x \in \sX_0$ a map 
\begin{equation}\label{map:unit-2-cell-of-lax-functor}
	\unit_{F(x)} \to F(\unit_x)
\end{equation}
in $\sY(F(x),F(x))$.

\begin{lemma}\label{lemma:Criterion-when-lax-functors-are-strict}
	Let $F \colon \sX \laxfun \sY$ be a lax functor. Then $F$ is a strong functor if and only if for any pair of composable 1-cells $f,g$ in $\sX$ and any object $x \in \sX_0$, the 2-morphisms \cref{map:composition-2-cell-of-lax-functor} and \cref{map:unit-2-cell-of-lax-functor} are invertible in $\sY$.
\end{lemma}

\begin{proof}
	From the description of the components of the natural transformations $\alpha_t$ given just below diagram \cref{diagram:Non-invertible-cells-lax-functor}, it follows that $F$ preserves cocartesian edges precisely if for every $t \colon [m] \to [n]$, the transformation $\alpha_t$ is invertible. Since any map in $\Delta$ can be written as a composition of degeneracy maps $\sigma^i \colon [n+1] \to [n]$ and face maps $\delta^j \colon [m-1] \to [m]$, it suffices to show that $\alpha_t$ is invertible for these two cases. Note that $\delta^0,\delta^m \colon [m-1] \to [m]$ are inert, so $\alpha_t$ is invertible in those cases.
	
	Observing that the Segal maps $\sX_{k+l} \to \sX_k \times_{\sX_0} \sX_l$ are obtained from inert maps in $\Delta$, it follows that the square
	\[\begin{tikzcd}[column sep=huge]
		\sX_{k+l} \ar[r,"F_{k+l}"] \ar[d,"\sim" rot90] & \sY_{k+l} \ar[d,"\sim" rot90]\\
		\sX_k \times_{\sX_0} \sX_l \ar[r,"{F_k \times_{F_0} F_l}"] & \sY_k \times_{\sY_0} \sY_l
	\end{tikzcd}\]
	commutes. In particular, using that the Segal maps are equivalences, the problem reduces to showing that the natural transformations
	\[\begin{tikzcd}
		\sX_2 \ar[r,"F_2",""{name=0, anchor=center, inner sep=0}] \ar[d,"d_1"'] & \sY_2 \ar[d,"d_1"]\\
		\sX_1 \ar[r,"F_1"',""{name=1, anchor=center, inner sep=0}] & \sY_1
		\arrow[shorten <=8pt, shorten >=8pt, Rightarrow, from=1-2, to=2-1]
	\end{tikzcd}
	\quad \text{and} \quad
	\begin{tikzcd}
		\sX_0 \ar[r,"F_0",""{name=0, anchor=center, inner sep=0}] \ar[d,"s_0"'] & \sY_0 \ar[d,"s_0"]\\
		\sX_1 \ar[r,"F_1"',""{name=1, anchor=center, inner sep=0}] & \sY_1
		\arrow[shorten <=8pt, shorten >=8pt, Rightarrow, from=1-2, to=2-1]
	\end{tikzcd}\]
	are invertible. But this is exactly the assumption on $F$.
\end{proof}

If $\sC \curvearrowright \sM \curvearrowleft \sD$ and $\sC' \curvearrowright \sN \curvearrowleft \sD'$ are bimodule categories, viewed as cocartesian fibrations over $(\Delta_{/[1]})^\op$, then a lax bimodule functor $F$ from $\sM$ to $\sN$ is a functor between these cocartesian fibrations that preserves cocartesian lifts of inerts (but not necessarily cocartesian lifts of other maps).
Informally, such a lax bimodule functor consists of lax monoidal functors $\sC \to \sC'$ and $\sD \to \sD'$ together with a functor $\sM \to \sN$ that laxly commutes with the left and right actions.
Since a map in $(\Delta_{/[1]})^\op$ is inert precisely if its underlying map in $\Delta^\op$ is,
we obtain the following.

\begin{lemma}\label{lemma:lax-functor-gives-lax-bimodule-functor}
    Let $F \colon \sX \to \sY$ be a lax functor between 2-precategories and let $x$ and $y$ be objects in $\sX$. Then $F$ restricts to a lax bimodule functor
    \[\begin{tikzcd}[column sep = large]
        \sX(y,y) \curvearrowright \sX(x,y) \curvearrowleft \sX(x,x) \ar[r] & \sY(Fy,Fy) \curvearrowright \sY(Fx,Fy) \curvearrowleft \sY(Fx,Fx)
    \end{tikzcd}\]
\end{lemma}

\begin{proof}
    This follows since the pullbacks used to define the bimodule structures on $\sX(x,y)$ and $\sY(Fx,Fy)$ in \cref{proposition:bimodule-mapping-categories} can also be computed in the category of categories over $(\Delta_{/[1]})^\op$ that have cocartesian lifts of inerts.
\end{proof}

\subsubsection{Lax functors as local right adjoints}\label{ssec:local-adjunctions}

Let us now describe our main source of lax functors, namely as ``local right adjoints'' of (oplax) functors of 2-precategories.

\begin{definition}[Locally adjointable functor]\label{definition:locally-adjointable-functor}
    An (oplax) functor $F \colon \sX \oplaxfun \sY$ of 2-precategories is called \emph{locally right adjointable} if for any pair of objects $x,y \in \sX_0$, the induced functor $F_{x,y} \colon \sX(x,y) \to \sY(Fx,Fy)$ has a right adjoint.
\end{definition}

\begin{proposition}\label{prop: local-right-adjoint-is-lax}
    Let $F \colon \sX \oplaxfun \sY$ be a locally right adjointable oplax functor such that $F_0 \colon \sX_0 \to \sY_0$ is an equivalence. Then there exists a lax functor $R \colon \sY \laxfun \sX$ such that $R_0$ is inverse to $F_0$ and for any $x,y \in \sX_0$, the functor $R_{Fx,Fy} \colon \sY(Fx,Fy) \to \sX(x,y)$ is right adjoint to $F_{x,y} \colon \sX(x,y) \to \sY(Fx,Fy)$.
\end{proposition}

\begin{proof}
    Without loss of generality, we may assume that $\sX_0 = \sY_0$ and that $F_0$ is the identity. Consider the diagram
    \[\begin{tikzcd}[column sep=small]
    	\sX_1 \ar[rr,"F_1"] \ar[dr,"{(d_1,d_0)}"'] & & \sY_1 \ar[dl,"{(d_1,d_0)}"] \\
    		& \sX_0^{\times 2}, & 
    \end{tikzcd}\]
    which commutes since $d_0$ and $d_1$ are inert. Since $\sX_0$ is a space, this is a map between bicartesian fibrations. Since on fibers $F_1$ is exactly the functor $F_{x,y} \colon \sX(x,y) \to \sY(x,y)$, it follows from the dual of \cite[Proposition 7.3.2.6]{HA} that the right adjoints $R_{x,y}$ assemble into a relative right adjoint $R_1 \colon \sY_1 \to \sX_1$ of $F_1$. Now since $F_n$ is equivalent to the map
    \[\sX_n \simeq \sX_1 \times_{\sX_0} \cdots \times_{\sX_0} \sX_1 \xrightarrow{{(F_1,\ldots,F_1)}} \sY_1 \times_{\sX_0} \cdots \times_{\sX_0} \sY_1 \simeq \sY_n\]
    we see that it admits a right adjoint given by $(R_1,\ldots,R_1)$. Applying \cite[Theorem 3.1.11]{linskensLaxMonoidalAdjunctions} to the diagram
    \[\begin{tikzcd}[column sep=small]
    	\cartun \sX \ar[rr,"F"] \ar[dr] & & \cartun \sY \ar[dl] \\
    	& \Delta, &
    \end{tikzcd}\]
    where $\cartun$ denotes the cartesian unstraightening of a functor $\Delta^\op \to \Cat$, we obtain a functor $R \colon \un \sY \to \un \sX$ over $\Delta^\op$ which is fiberwise right adjoint to $F$. Using the Segal conditions of $\sX$ and $\sY$, an argument similar to the proof of \cite[Theorem 3.4.7]{linskensLaxMonoidalAdjunctions} yields the desired result.
\end{proof}

\begin{remark}\label{remark:lax-structure-maps-of-right-adjoint}
    Suppose $F \colon \sX \to \sY$ is a locally right adjointable oplax functor and let $R$ be constructed as in \cref{prop: local-right-adjoint-is-lax}. Write $\mu \colon F(g \circ f) \to F(g) \circ F(f)$ and $\nu \colon F(\unit_x) \to \unit_{F(x)}$ for the oplax structure maps of $F$ in the sense of \cref{map:composition-2-cell-of-lax-functor} and \cref{map:unit-2-cell-of-lax-functor}. By (the dual of) \cite[Proposition 3.2.7]{linskensLaxMonoidalAdjunctions}, it follows that the lax structure maps of $R$ are given by the composites
    \[Rg \circ Rf \xrightarrow{\eta} RF(Rg \circ Rf) \xrightarrow{R(\mu)} R(FRg \circ FRf) \xrightarrow{R(\varepsilon \circ \varepsilon)} R(g \circ f)\]
    and 
    \[\unit_{Rx} \xrightarrow{\eta} RF(\unit_{Rx}) \xrightarrow{R(\nu)} R(\unit_{FRx}) \simeq R(\unit_x),\]
    where $\eta$ and $\varepsilon$ denote the unit and counit of $F_{x,y} \dashv R_{x,y}$, respectively.
\end{remark}

\begin{remark}
Observe that the results of \cite{linskensLaxMonoidalAdjunctions} cited in the proof of \cref{prop: local-right-adjoint-is-lax} actually provide an equivalence between the category of right-adjointable oplax functors $\sX \oplaxfun \sY$ (that induce an equivalence on spaces of objects) and the opposite of the category of left-adjointable lax functors $\sY \laxfun \sX$ (that induce an equivalence on spaces of objects).
\end{remark}

In practice, not all locally right adjointable (oplax) functors $F \colon \sX \oplaxfun \sY$ induce an equivalence $F_0 \colon \sX_0 \to \sY_0$. To still construct its local right adjoint $R$, one needs to modify the space of objects of $\sY$. Recall the definition of $\cosk_0$ from \cref{example:coskeleton}.

\begin{construction}\label{construction:Modifying-space-of-objects}
    Let $\sX$ be a 2-precategory and let $\phi \colon S \to \sX_0$ be a map of spaces. Then $\phi^* \sX$ is the 2-precategory given by the pullback
    \[\begin{tikzcd}
		\phi^*\sX \ar[r] \ar[d] & \sX \ar[d,"\eta"]\\
		\cosk_0 S \ar[r,"\cosk_0 \phi"] & \cosk_0 \sX_0.
	\end{tikzcd}\]
    Note that given two points $s,t \in S$, the functor $\phi^* \sX \to \sX$ induces an equivalence $(\phi^*\sX)(s,t) \simeq \sX(\phi s, \phi t)$. In particular, the functor $\phi^*\sX \to \sX$ of 2-precategories is fully faithful.
\end{construction}

\begin{remark}
    Taking $S = *$, we recover \cref{def:endomorphisms-infty-category}.
\end{remark}

Given an oplax functor $F \colon \sX \oplaxfun \sY$, we apply \cref{construction:Modifying-space-of-objects} to the map of spaces $F_0 \colon \sX_0 \to \sY_0$ to obtain a 2-precategory $F_0^* \sY$ with the same space of objects as $\sX$.

\begin{corollary}
    Let $F \colon \sX \oplaxfun \sY$ be a locally right adjointable oplax functor. Then there exists a lax functor $R \colon F_0^* \sY \laxfun \sX$ that is the identity on spaces of objects and such that for any $x,y \in \sX_0$, the functor $R_{x,y} \colon (F_0^*\sY)(x,y) \simeq \sY(Fx, Fy) \to \sX(x,y)$ is right adjoint to $F_{x,y}$.
\end{corollary}

\begin{proof}
    Note that the cartesian unstraightening functor $\Fun(\Delta^\op, \Catinfty) \to \Catinfty_{/\Delta}$ preserves pullbacks. Since the diagram
    \[\begin{tikzcd}
		\cartun \sX \ar[r,"F"] \ar[d] & \cartun \sY \ar[d]\\
		\cartun \cosk_0 \sX_0 \ar[r,"\cosk_0 F_0"] & \cartun \cosk_0 \sY_0.
	\end{tikzcd}\]
    commutes and $F$ preserves cartesian lifts of inerts, we obtain a functor $F' \colon \cartun \sX \to \cartun F^*_0 \sY$ over $\Delta$ which preserves cartesian lifts of inerts. Applying \cref{prop: local-right-adjoint-is-lax} to $F'$, we obtain the desired result.
\end{proof}

\subsubsection{Localizations on mapping categories}

Given a monoidal category $\cC^\otimes$ and a reflective localization $L$ of the underlying category $\sC$, the full subcategory of local objects $L\sC$ is generally not closed under the monoidal structure of $\sC$. However, if $L$ is compatible with the monoidal structure in the sense of \cite[Definition 2.2.1.6]{HA}, then by \cite[Proposition 2.2.1.9]{HA} one does obtain a \emph{localized} monoidal structure on $L\sC$. Given a 2-category equipped with reflective localizations of its mapping categories, we will establish an analogous result.

\begin{definition}[Family of localizations]\label{definition:family-of-localizations}
    Let $\sX$ be a 2-precategory and suppose that we are given a full subcategory $\sX'_1 \subset \sX_1$ such that for any $x,y \in \sX_0$, the inclusion
    \[\sX'(x,y) \coloneqq \sX'_1 \times_{\sX_0^{\times 2}} \{(x,y)\} \hookrightarrow \sX(x,y)\]
    admits a left adjoint $L_{x,y}$. We will then call $\{L_{x,y}\}_{x,y \in \sX_0}$ a \emph{family of reflective localizations}.
\end{definition}

\begin{remark}
    In the previous definition, asking that for any $x,y \in \sX_0$, the inclusion $\sX'(x,y) \hookrightarrow \sX(x,y)$ admits a left adjoint is equivalent to asking that $\sX'_1 \hookrightarrow \sX_1$ admits a left adjoint. However, the definition given above will be more convenient in practice.
\end{remark}

\begin{definition}[Compatible family of localizations]\label{def: compatible family of localizations}
    Suppose $\{L_{x,y}\}_{x,y \in \sX_0}$ is a family of reflective localizations on a 2-precategory $\sX$. We will say that $\{L_{x,y}\}_{x,y \in \sX_0}$ is a \emph{compatible} family of reflective localizations if for any $x,y,z \in \sX$ and any $f \in \sX(x,y)$ and $g \in \sX(y,z)$, the maps
    \[L_{x,z}(g \circ f) \to L_{x,z}(L_{y,z} g \circ f) \quad \text{and} \quad L_{x,z}(g \circ f) \to L_{x,z}(g \circ L_{x,y}f) \]
    are equivalences in $\sX(x,z)$.
\end{definition}

Given a compatible family of reflective localizations $L = \{L_{x,y}\}_{x,y \in \sX_0}$ on $\sX$, consider the cocartesian unstraightening $\un \sX \to \Delta^\op$ of $\sX$. We can form the full subcategory $\sX'$ of $\un \sX$ whose fiber over $[0]$ is simply $\sX_0$ and whose fiber over $[n]$ is spanned by those objects $(f_1,\ldots,f_n)$ of $\sX_n \simeq \sX_1 \times_{\sX_0} \cdots \times_{\sX_0} \sX_1$ with the property that for every $i$, the 1-morphism $f_i$ is local.

\begin{proposition}[Localization of mapping categories]\label{proposition:local-localizations}
    Let $\sX$ be a 2-precategory and suppose that we are given a compatible family of reflective localizations $L = \{L_{x,y}\}_{x,y \in \sX_0}$. Then the full subcategory $\sX'$ of $\un \sX$ constructed above is a cocartesian fibration over $\Delta^\op$ whose straightening, denoted $L\sX$, is a 2-precategory. The composition in $L\sX$ is given by
    \[\circ^L \colon L\sX(y,z) \times L\sX(x,y) \to L\sX(x,z); \quad g \circ^L f = L_{x,z}(g \circ f)\]
    where $\circ$ denotes the composition of $\sX$, and the unit $\unit^L_x$ of $L \sX$ at $x$ is given by $L_{x,x}(\unit_x)$. Moreover, the inclusion $\sX' \hookrightarrow \sX$ defines a lax functor $L\sX \laxfun \sX$ and the localizations $L_{x,y}$ assemble into a strong functor $L \colon \sX \to L\sX$ that is locally left adjoint to $L\sX \laxfun \sX$.
\end{proposition}

This proposition can be proved in the same way as \cite[Proposition 2.2.1.9]{HA}. For the convenience of the reader, we include an alternative proof.

\begin{proof}
    Throughout this proof, we will write $L$ for $L_{x,y}$. Note that $\sX'$ has cocartesian lifts of inerts and that it satisfies the Segal condition by construction. In particular, to show that $\sX' \to \Delta^\op$ is a cocartesian fibration, it suffices to show that it has cocartesian lifts of those maps $[1] \to [n]$ that send $0$ to $0$ and $1$ to $n$. Let us first show that $\sX' \to \Delta^\op$ admits locally cocartesian lifts of such maps. For any sequence
    \[x_0 \xrightarrow{f_n} x_1 \xrightarrow{f_2} \cdots \xrightarrow{f_1} x_n\]
    of composable 1-morphisms in $\sX$ and any $g$ in $L\sX(x_0,x_n)$, we need to show that there exists an $f'$ in $L\sX(x_0,x_n)$ such that there is a natural equivalence $\Map_{\sX(x_0,x_n)}(f',g) \simeq \Map_{\sX(x_0,x_n)}(f_1 \circ \cdots \circ f_n, g)$. But such an $f'$ is given by $L(f_1 \circ \cdots \circ f_n)$. Note that this also includes the case where $n=0$ if we interpret $f_1 \circ \ldots \circ f_n$ as $\unit_{x_0}$ when $n=0$.

    To see that $\sX' \to \Delta^\op$ is a cocartesian fibration, we need to show by \cite[Proposition 2.4.2.8]{HTT} that the locally cocartesian lifts are closed under composition. For this, it suffices to show that for any sequence of composable morphisms $(f_1, \ldots, f_m)$ in $\sX$ and any sequence $0 = n_0 \leq n_1 \leq \cdots \leq n_m = n$, the map
    \[f_1 \circ \cdots \circ f_n \to L(f_{n_0 + 1} \circ \cdots \circ f_{n_1}) \circ \cdots \circ L(f_{n_{m-1}+1} \circ \cdots \circ f_{n_m})\]
    is an equivalence after applying $L$, where $f_{n_i+1} \circ \cdots \circ f_{n_{i+1}}$ is interpreted as $\unit_{x_{n_i}}$ if $n_i = n_{i+1}$. This follows since $\{L_{x,y}\}_{x,y \in \sX_0}$ is a compatible family of localizations.

    As remarked above, the straightening $L\sX$ of $\sX'$ satisfies the Segal condition by construction. Since $(L\sX)_0 = \sX_0$ is a space, we conclude that $L\sX$ is a 2-precategory. Since the cocartesian lifts of $\sigma^0 \colon [1] \to [0]$ and $\delta^1 \colon [1] \to [2]$ are given by maps of the form $x \to L (\unit_x)$ and $(f_1,f_2) \to L(f_1 \circ f_2)$, we see that the unit and composition of $L\sX$ are indeed as described in the proposition.

    Finally, the inclusion $\sX' \to \un \sX$ preserves cocartesian lifts of inerts by construction, hence we obtain a lax functor $L\sX \laxfun \sX$ which on mapping categories is given by the inclusion $L\sX(x,y) \hookrightarrow \sX(x,y)$ of the $L_{x,y}$-local objects. By the dual of \cref{prop: local-right-adjoint-is-lax}, we see that the localizations $L_{x,y}$ assemble into an oplax functor $\sX \oplaxfun L\sX$. By the duals of \cref{remark:lax-structure-maps-of-right-adjoint} and \cref{lemma:Criterion-when-lax-functors-are-strict}, this is a strong functor precisely if the maps $\unit_x \to L(\unit_x)$ and $f \circ g \to L(f) \circ L(g)$ are equivalences after applying $L$, which holds by assumption.
\end{proof}

\begin{remark}
    We will abusively write $L$ for the functor $\sX \to L\sX$ constructed in \cref{proposition:local-localizations}.
\end{remark}

\begin{remark}
    Even if we start with a 2-category $\sX$ and a compatible family of reflective localizations $L$, then $L\sX$ is generally not a 2-category. To obtain an actual 2-category, one needs to consider the completion $(L\sX)^\wedge$ of $L\sX$. However, while there is a lax functor $L\sX \laxfun \sX$, there need not be a lax functor $(L\sX)^\wedge \laxfun \sX$, so in practice it can be more convenient to work with the 2-precategory $L\sX$ as opposed to its completion.
\end{remark}

\begin{remark}[Weakly compatible family of localizations]\label{remark:weakly-compatible-family-of-localizations}
    Given a 2-precategory $\sX$ and a (not necessarily compatible) family of reflective localizations $L = \{ L_{x,y}\}_{x,y \in \sX_0}$, we can form a locally full subcategory of $\sZ$ of $\sX$ whose 1-morphisms are the units of $\sX$ and those 1-morphisms that can be obtained as finite compositions of $L$-local 1-morphisms in $\sX$. We will call $L$ a \emph{weakly compatible family of reflective localizations} if its restriction to $\sZ$ forms a compatible family of reflective localizations. In this case the construction of \cref{proposition:local-localizations} still yields a 2-precategory $L\sX \coloneqq L\sZ$, but instead of a strong functor we only obtain an oplax functor from $\sX$ to $L\sX$ (which is strong on the locally full subcategory $\sZ$).
\end{remark}

We conclude with the following recognition principle for 2-precategories of the form $L\sX$.

\begin{proposition}[Recognizing localizations]\label{proposition:Recognizing-localizations}
    Let $\phi \colon \sY \laxfun \sX$ be a lax functor between 2-precategories that is the identity on the spaces of objects $\sY_0 = \sX_0$.
    Suppose that for any $x,y \in \sY_0$, the functor $\sY(x,y) \to \sX(x, y)$ is fully faithful and admits a left adjoint $L = L_{x,y}$.
    Furthermore, suppose that for any $x$ in $\sY_0$ and any sequence of composable morphisms $(f_1,\ldots,f_n)$ in $\sY$, the lax structure maps
    \[\begin{tikzcd}
        \phi(f_n) \circ \cdots \circ \phi(f_1) \ar[r] & \phi(f_n \circ \cdots \circ f_1)
    \end{tikzcd}
    \quad \text{and} \quad
    \begin{tikzcd}
        \unit_x \ar[r] & \phi(\unit_x)
    \end{tikzcd}\]
    induce equivalences
    \[\begin{tikzcd}
        L(\phi(f_n) \circ \cdots \circ \phi(f_1)) \ar[r,"\sim"] & f_n \circ \cdots \circ f_1
    \end{tikzcd}
    \quad \text{and} \quad
    \begin{tikzcd}
        L(\unit_x) \ar[r,"\sim"] & \unit_x.
    \end{tikzcd}\]
    Then $L = \{L_{x,y}\}_{x,y \in \sX_0}$ is a weakly compatible family of localizations and $\sY$ is canonically equivalent to $L\sX$.
\end{proposition}

\begin{proof}
    Observe that $\phi$ induces a fully faithful map $\sY_1 \to \sX_1$. In particular, the localizations $L_{x,y}$ form a family of reflective localizations in the sense of \cref{definition:family-of-localizations}.
    It follows from the assumptions on the lax structure maps of $\phi$ that $L$ forms a weakly compatible family of localizations on $\sX$ in the sense of \cref{remark:weakly-compatible-family-of-localizations}.
    In particular, the category $L\sX$ exists.
    In the proof of \cref{proposition:local-localizations}, the cocartesian unstraightening $\un L\sX \to \Delta^\op$ of $L\sX$ is constructed as a certain full subcategory of $\un \sX$.
    It follows by construction that the objects of this full subcategory are precisely those in the essential image of $\phi \colon \un \sY \to \un \sX$.
    The result therefore follows if we can show that $\phi \colon \un \sY \to \un \sX$ is fully faithful.
    
    An inspection of the mapping spaces in $\un \sY$ and $\un \sX$ shows that this is the case precisely if for any sequence of composable morphisms
    \[\begin{tikzcd}
        x \ar[r,"f_1"] & \cdots \ar[r, "f_n"] & y
    \end{tikzcd}\]
    in $\sY$ (where we allow $n=0$) and any $g \colon x \to y$ in $\sY$, the lax structure maps of $\phi$ induce equivalences
    \[\Map(f_n \circ \cdots \circ f_1, g) \simeq \Map(\phi(f_n) \circ \cdots \circ \phi(f_1), \phi(g)).\]
    But this follows since $f_n \circ \cdots \circ f_1 \simeq L(\phi(f_n) \circ \cdots \circ \phi(f_1))$ by assumption.
\end{proof}

\newpage

\section{Orthofibrations, ends and mixed modules} \label{appendix:orthofibrations}

The purpose of this appendix is to prove two technical results that are used in our treatment of Koszul duality: we will give an explicit method for computing ends in the category $\Cat$ of categories and we will construct a natural equivalence $\LComod_Q(\RMod_A(\sM)) \simeq \RMod_A(\LComod_Q(\sM))$, where $\sC \curvearrowright \sM \curvearrowleft \sD$ is a bimodule object in $\Cat$.
Both proofs make use of \emph{orthofibrations}: these are functors $\sC \to \sA \times \sB$ which are cartesian over $\sA$ and cocartesian over $\sB$ in a compatible way.
Orthofibrations can be straightened to functors $\sA^\op \times \sB \to \Cat$.
They are useful in situations where one needs to mix two different directions, which is why we use them here.

We start by recalling the definition and some basic properties of orthofibrations in \cref{subsec:orthofibrations}.
We then prove our formula for computing ends in $\Cat$ in \cref{subsec:ends-using-orthofibrations}.
Finally, we prove the natural equivalence $\LComod_Q(\RMod_A(\sM)) \simeq \RMod_A(\LComod_Q(\sM))$ in \cref{subsec:Mixed-modules}, where we also deduce several other useful facts about (co)modules.

\subsection{Recollections on orthofibrations}\label{subsec:orthofibrations}

Let $\sA^\op \times \sB \to \Cat$ be a functor.
By unstraightening this functor in the second variable, we can equivalently view this as a functor $\sA^\op \to \Cocart(\sB)$.
Since $\Cocart(\sB)$ is a locally full subcategory of $\Cat_{/\sB} = \Funtwo([1],\Cat) \times_{\Cat} \{\sB\}$, we can view such a functor as an object of the 2-category
\[\Funtwo(\sA^\op, \Funtwo([1],\Cat) \times_{\Cat} \{\sB\} ) \simeq \Funtwo([1],\Funtwo(\sA^\op,\Cat)) \times_{\Funtwo(\sA^\op,\Cat)} \{\const_{\sB}\}. \]
If we now apply cartesian unstraightening to the variable $\sA$, we see that the right-hand side is equivalent to
\[\Funtwo([1],\Cart(\sA)) \times_{\Cart(\sA)} \{\sA \times \sB\} \simeq \Cart(\sA)_{/\sA \times \sB}.\]
In particular, we see that $\Funtwo(\sA^\op \times \sB, \Cat)$ is equivalent to a locally full subcategory of $\Cart(\sA)_{/\sA \times \sB} \subset \Cat_{/\sA \times \sB}$.
Unwinding the equivalences above, we see that a functor $p = (p_1,p_2) \colon \sC \to \sA \times \sB$ lies in this subcategory precisely if:
\begin{enumerate}[(a)]
    \item $p_1 \colon \sC \to \sA$ is cartesian, \label{ortho:a1}
    \item $p \colon \sC \to \sA \times \sB$ preserves cartesian lifts of maps in $\sA$, \label{ortho:a2}
    \item for every $a$ in $\sA$, the map on fibers $p_a \colon \sC_a \to (\sA \times \sB)_a = \sB$ is a cocartesian fibration, and \label{ortho:b1}
    \item for every $\alpha \colon a' \to a$ in $\sA$, the cartesian transport functor $\alpha^* \colon \sC_a \to \sC_{a'}$ preserves cocartesian lifts of maps in $\cB$. \label{ortho:b2}
\end{enumerate}

\begin{definition}\label{definition:orthofibration}
    A functor $\sC \to \sA \times \sB$ satisfying these properties is called an \emph{orthofibration}. Given orthofibrations $p \colon \sC \to \sA \times \sB$ and $q \colon \sD \to \sA \times \sB$, a \emph{map of orthofibrations} (from $p$ to $q$) is a functor $\sC \to \sD$ over $\sA \times \sB$ that preserves cartesian arrows over $\sA$ and cocartesian arrows over $\sB$. The 2-category $\Ortho(\sA,\sB)$ of orthofibrations over $\sA \times \sB$ is defined as the locally full subcategory of $\Cat_{/\sA \times \sB}$ whose objects are orthofibrations and whose 1-morphisms are maps of orthofibrations.
\end{definition}

\begin{remark}
    Observe that items \ref{ortho:a1}-\ref{ortho:a2} ensure that $p \colon \sC \to \sA \times \sB$ lies in $\Cart(\sA)_{/\sA \times \sB}$, while items \ref{ortho:b1}-\ref{ortho:b2} ensure that the straightened functor $\sA^\op \to \Cat_{/\cB}$ lands in $\Cocart(\cB)$.
\end{remark}

\begin{remark}
    What we call orthofibrations are called \emph{two-sided fibrations} in \cite{Stefanich2020HigherSheafTheory,RiehlVerity2022ElementsCategoryTheory} and \emph{mixed fibrations} in \cite{Torii2024Duoidal$infty$categories}. We will follow the terminology of \cite{linskensLaxMonoidalAdjunctions}.
\end{remark}

In the discussion preceding \cref{definition:orthofibration}, we showed the following.

\begin{theorem}[{\cite[Corollary 2.3.12]{linskensLaxMonoidalAdjunctions}, \cite[Proposition 2.1.16]{Stefanich2020HigherSheafTheory}}]\label{theorem:straightening-orthofibrations}
    Let $\sA$ and $\sB$ be categories. Then straightening in both variables yields an equivalence of 2-categories
    \[\Ortho(\sA,\sB) \simeq \Funtwo(\sA^\op \times \sB, \Cat).\]
\end{theorem}

Orthofibrations admit many equivalent characterizations, some of which are more symmetric in the variables $\sA$ and $\sB$ than the definition given above. For details, the reader is referred to Propositions 2.3.3 and 2.3.11 of \cite{linskensLaxMonoidalAdjunctions}.
We will not use these alternative characterizations.
However, let us just mention the following fact, which follows as in \cite[Propositions 2.3.3 \& 2.3.11]{linskensLaxMonoidalAdjunctions}.

\begin{lemma}\label{lemma:orthofibration-admits-both-cartesian-arrows}
    Let $p = (p_1,p_2) \colon \sC \to \sA \times \sB$ be an orthofibration.
    Then $\sC$ admits all $p$-cartesian lifts of maps in $\sA \times \sB^\simeq$ and all $p$-cocartesian lifts of maps in $\sA^\simeq \times \sB$.
    Moreover, cartesian transport over $p_1$ preserves $p_2$-cocartesian lifts, while cocartesian transport over $p_2$ preserves $p_1$-cartesian lifts.
\end{lemma}

In this situation, if we are given an object $c \in \sC$ and a map $\alpha \colon a \to p_1(c)$ in $\sA$, we will write $\alpha^* c \to c$ for a cartesian lift of $(\alpha,\id)$ at $c$. Similarly, if $\beta \colon p_2(c) \to b$ is a map in $\sB$, we shall write $c \to \beta_!c$ for a cocartesian lift of $(\id,\beta)$ at $c$.

\subsection{Computing ends as functors into orthofibrations}\label{subsec:ends-using-orthofibrations}

Let $H \colon \sA^\op \times \sA \to \sC$ be a functor. Then the \emph{end} of $H$ is defined (if it exists) as the limit of the composite
\[\begin{tikzcd}
    \TwL \sA \ar[r] & \sA^\op \times \sA \ar[r,"H"] & \sC
\end{tikzcd}\]
and is usually denoted by $\int_\sA H(a,a)$. For the basic properties of ends in the $\infty$-categorical setting, the reader is referred to \cite{Haugseng2022CoEndsCategories}.

One of the main reasons for considering ends is that they allow one to compute mapping spaces in functor categories. Namely, if $F,G \colon \cC \to \cD$ are functors, then by \cite[Proposition 5.1]{GepnerHaugsengea2017LaxColimitsFree} or \cite[Proposition 2.3]{Glasman2016SpectrumlevelHodgeFiltration} there is an equivalence
\[\Map_{\Fun(\cC,\cD)}(F,G) \simeq \int_\sC \Map(F(c),G(c)).\]
Ends can also be used to phrase the following variation on the Yoneda lemma, which is used in the proof of \cref{prop: comparison-bar-cobar-adjunctions}.
\begin{proposition}
\label{prop: parametrized-yoneda}
    For functors $H \colon \cA^\op \times \cB \to \Spc$ and $F \colon \cA \to \cB$, there is an equivalence
    \[
    \Map_{\Fun(\sA^\op \times \sB, \Spc)}(\Map_{\sB}(F(-),-), H) \simeq \int_{a \in \sA} H(a,Fa)
    \]
    of spaces.
\end{proposition}
\begin{proof}
     By combining the ``Fubini theorem'' for ends \cite[Proposition 2.21]{Haugseng2022CoEndsCategories} and the classical Yoneda lemma, it follows that
\begin{align*}
    &\Map_{\Fun(\sA^\op \times \sB, \Spc)}(\Map_{\sB}(F(-),-), H) \simeq \int_{(a,b) \in \sA^\op \times \sB} \Map_\Spc(\Map_\sB(Fa,b),H(a,b)) \\
    &\simeq \int_{a \in \sA^\op} \Map_{\Fun(\sB,\Spc)}(\Map_\sA(Fa,-),H(a,-)) \simeq \int_{a \in \sA} H(a,Fa). \qedhere
\end{align*}
\end{proof}
One can view this as a parametrized version of the Yoneda lemma.
In particular, if $\sA = *$ then this is the usual Yoneda lemma for presheaves on $\sB^\op$.

The goal of this section is to describe a general method for computing ends in the category $\Cat$ of categories and the category $\Spc$ of spaces using orthofibrations. This will require the following definition.

\begin{definition}[Orthocartesian arrows]\label{definition:orthocartesian-maps}
    Let $\sC \to \sA \times \sA$ be an orthofibration and suppose we are given a morphism $\gamma \colon c \to c'$ in $\sC$ lying over $(\alpha, \beta) \colon (a,b) \to (a',b')$ in $\sA \times \sA$.
    By \cref{lemma:orthofibration-admits-both-cartesian-arrows}, we may consider the (co)cartesian lifts $c \to \beta_!c$ and $\alpha^* c' \to c'$.
    Using the universal properties of these (co)cartesian lifts, we obtain the following essentially unique dotted arrows in $\sC$:
    \[\begin{tikzcd}[column sep = small]
    	& {\beta_!c} & \\
    	c && {c'} \\
    	& {\alpha^*c'} &
    	\arrow["{\exists!}"{description}, dashed, from=1-2, to=2-3]
    	\arrow["{\exists!}"'{pos=0.27}, dashed, from=1-2, to=3-2]
    	\arrow[from=2-1, to=1-2]
    	\arrow[from=2-1, to=2-3, crossing over,"\gamma"{pos=0.7}]
    	\arrow["{\exists !}"{description}, dashed, from=2-1, to=3-2]
    	\arrow[from=3-2, to=2-3]
    \end{tikzcd}\]
    We say that $\gamma$ is \emph{orthocartesian} if this morphism $\beta_!c \to \alpha^* c'$ is invertible.
\end{definition}

\begin{remark}
    What we call orthocartesian arrows are called bicartesian arrows in \cite[Definition 2.1.4]{Stefanich2020HigherSheafTheory}.
\end{remark}

\begin{remark}
    One can show that for orthofibrations $\sC \to \sA \times \sB$ and $\sD \to \sA \times \sB$, a map $\sC \to \sD$ over $\sA \times \sB$ is a map of orthofibrations if and only if it preserves orthocartesian morphisms.
\end{remark}

\begin{definition}
    Let $\sC \to \sA \times \sA$ be an orthofibration and suppose that we are given a commuting triangle of the form
    \[\begin{tikzcd}[column sep = tiny]
    	\sA && \sC \\
    	& {\sA \times \sA}
    	\arrow["F", from=1-1, to=1-3]
    	\arrow["\Delta"', from=1-1, to=2-2]
    	\arrow[from=1-3, to=2-2]
    \end{tikzcd}\]
    where $\Delta$ denotes the diagonal map.
    We say that $F$ is \emph{strict} if it takes any morphism in $\sA$ to an orthocartesian arrow in $\sC$.
    The full subcategory of $\Fun_{/\sA \times \sA}(\sA, \sC)$ spanned by the strict functors will be denoted $\Funstr{\sA}(\sA,\sC)$.
\end{definition}

The main result of this section is that the category of strict functors into an orthofibration agrees with the end of the corresponding functor.

\begin{proposition}\label{proposition:categories-ortho-section-end}
    Let $H \colon \sA^\op \times \sA \to \Cat$ be a functor and let $\sH \to \sA \times \sA$ be the corresponding orthofibration. Then there is an equivalence
    \[\int_{\sA} H \simeq \Funstr{\sA}(\sA,\sH)\]
    that is natural in $H$.
\end{proposition}

\begin{remark}
    In the proof of \cref{prop: comparison-bar-cobar-adjunctions}, we only require this proposition for functors $H$ that land in the full subcategory $\Spc \subset \Cat$ of spaces. This somewhat simplifies the situation (cf.\ \cref{cor: spaces-ortho-section-end} below), but we included the more general statement since we believe the result to be interesting in its own right.
\end{remark}

\begin{proof}[Proof of \cref{proposition:categories-ortho-section-end}]
    By \cite[Proposition 4.11]{GepnerHaugsengea2017LaxColimitsFree}, the inclusion of 2-categories $\Cart(\sA) \hookrightarrow \Cat_{/\sA}$ admits a (2-categorical) left adjoint given by $\cC \mapsto \sA^{[1]} \times_{\sA} \sC$. By slicing over the cartesian fibration $\pr_1 \colon \sA \times \sA \to \sA$, we obtain an induced adjunction of 2-categories
    \[
    \begin{tikzcd}[sep = large]
        \Cart(\sA)_{/ \pr_1} \ar[r, shift right, hook,""{name=0, anchor=center, inner sep=0}] &  \ar[l, bend right = 10, ""{name=1, anchor=center, inner sep=0}] (\Cat_{/\sA})_{/\pr_1} 
        \arrow["\dashv"{anchor=center, rotate=-90}, draw=none, from=0, to=1]
    \end{tikzcd}
    \simeq \Cat_{/\sA \times \sA}.
    \]
    The left adjoint takes a functor $(f_1,f_2) \colon \sB \to \sA \times \sA$ to the pullback
    \[\begin{tikzcd}
    	{\sA^{[1]} \times_{\sA} \sB} \ar[dr,phantom, very near start, "\lrcorner"] & \sB \\
    	{\sA^{[1]}} & \sA,
    	\arrow[from=1-1, to=1-2]
    	\arrow[from=1-1, to=2-1]
    	\arrow["{f_1}", from=1-2, to=2-2]
    	\arrow["{\ev_1}"', from=2-1, to=2-2]
    \end{tikzcd}\]
    which we view as a category over $\sA \times \sA$ via the map
    \[\begin{tikzcd}(\ev_0,f_2) \colon \sA^{[1]} \times_{\sA} \sB \ar[r] &\sA \times \sA.\end{tikzcd}\]
    Applying this left adjoint to the diagonal $\sA \to \sA \times \sA$, we obtain a natural equivalence
    \begin{align}\label{equivalence:1-free-orthofibration}
        \Fun_{/\sA \times \sA}^{1\text{-}\mathrm{cart}}(\sA^{[1]}&, \sH) \simeq \Fun_{/\sA \times \sA}(\sA,\sH) \\
        &\begin{tikzcd}[ampersand replacement = \&]
            F \ar[r, maps to] \& (a \mapsto F(\id_a))
        \end{tikzcd} \nonumber
    \end{align}
    where $\Fun_{/\sA \times \sA}^{1\text{-}\mathrm{cart}}$ denotes the mapping category of $\Cat(\sA)_{/\pr_1}$.
    (Observe that this mapping category consists of those functors over $\sA \times \sA$ that preserve cartesian lifts of arrows in the first variable.)
    
    Note that $\sA^{[1]} \to \sA \times \sA$ is in fact an orthofibration.
    In particular, we may consider the full subcategory $\Fun_{/\sA \times \sA}^{\mathrm{orth}}(\sA^{[1]},\sH)$ of $\Fun_{/\sA \times \sA}^{1\text{-}\mathrm{cart}}(\sA^{[1]}, \sH)$ spanned by the maps of orthofibrations; that is, the functors that furthermore preserve cocartesian lifts in the second variable.
    We will show that under the equivalence \cref{equivalence:1-free-orthofibration}, the category $\Fun_{\sA \times \sA}^\mathrm{orth}(\sA^{[1]},\sH)$ corresponds to the full subcategory $\Funstr{\sA}(\sA,\sH)$ of strict functors.
    Given a morphism $\beta \colon b \to c$ in $\sA$, one can construct cartesian and cocartesian lifts of $\beta$ in $\sA^{[1]}$ as
    \[\begin{tikzcd}
    	b & c \\
    	d & d
    	\arrow["\beta", from=1-1, to=1-2]
    	\arrow["{\gamma \circ \beta}"', from=1-1, to=2-1]
    	\arrow["\gamma", from=1-2, to=2-2]
    	\arrow[equal, from=2-1, to=2-2]
    \end{tikzcd}
    \quad \text{and} \quad
    \begin{tikzcd}
    	a & a \\
    	b & c,
    	\arrow[equal, from=1-1, to=1-2]
    	\arrow["\alpha"', from=1-1, to=2-1]
    	\arrow["{\beta \circ \alpha}", from=1-2, to=2-2]
    	\arrow["\beta"', from=2-1, to=2-2]
    \end{tikzcd}\]
    respectively.
    In particular, a functor $F$ in $\Fun_{/\sA \times \sA}^{1\text{-}\mathrm{cart}}(\sA^{[1]}, \sH)$ preserves cocartesian lifts in the second variable if and only if for every $\alpha \colon a \to b$ and $\beta \colon b \to c$ in $\sA$, the map $F(\alpha) \to F(\beta \circ \alpha)$ is cocartesian.
    Since $F$ preserves cartesian morphisms in the first variable, this map is equivalent to the map $\alpha^* F(\id_b) \to \alpha^* F(\beta)$, where $\alpha^*$ is the cartesian transport functor.
    Since $\alpha^*$ preserves cocartesian lifts by definition of an orthofibration, it follows that $F$ is a map of orthofibrations if and only if for any $\beta \colon b \to c$ in $\sA$, the map $F(\id_b) \to F(\beta)$ is cocartesian.
    This is the case if and only if $\beta_! F(\id_b) \to F(\beta) \simeq \beta^* F(\id_c)$ is an equivalence.
    But this is precisely the condition that the functor in $ \Fun_{/\sA \times \sA}(\sA,\sH)$ corresponding to $F$ is strict.

    To finish the proof, we note that by \cite[Corollary A.2.5]{HaugsengMelaniea2022ShiftedCoisotropicCorrespondences} and the dual of \cite[Proposition 5.2.1.11]{HA}, $\TwL(\sA)$ is the cocartesian unstraightening of the mapping space functor $\Map_{\sA}(-,-)$, while $\sA^{[1]}$ is the orthofibration corresponding to this functor.
    The equivalence of 2-categories $\Ortho(\sA,\sA) \simeq \Cocart(\sA^\op \times \sA)$ now yields that
    \[\Funstr{\sA}(\sA,\sH) \simeq \Fun_{/\sA \times \sA}^{\mathrm{orth}}(\sA^{[1]},\sH) \simeq \Fun_{/\sA^\op \times \sA}^{\mathrm{cocart}}(\TwL(\sA),\sH).\]
    By \cite[Corollary 3.3.3.2]{HTT}, we conclude that
    \[\Fun_{/\sA^\op \times \sA}^{\mathrm{cocart}}(\TwL(\sA),\sH) \simeq \lim_{\TwL(\sA)} H = \int_{a \in \sA} H(a,a). \qedhere\]
\end{proof}

If we are given a functor $H \colon \sA^\op \times \sA \to \Spc$ landing in spaces, then the fibers of the corresponding orthofibration $\sH \to \sA \times \sA$ are spaces. It follows that the morphism $\beta_!c \to \alpha^*c'$ constructed in \cref{definition:orthocartesian-maps} is always invertible, hence that any morphism in $\sH$ is orthocartesian. In particular, any functor $\sA \to \sH$ is strict and we conclude the following.

\begin{corollary}
\label{cor: spaces-ortho-section-end}
    Let $H \colon \sA^\op \times \sA \to \Spc$ be a functor landing in spaces and let $\sH \to \sA \times \sA$ be the corresponding orthofibration. Then there is an equivalence
    \[\int_{\sA} H(a,a) \simeq \Fun_{/\sA \times \sA}(\sA, \sH) \simeq \Map_{/\sA \times \sA}(\sA,\sH)\]
    that is natural in $H$.
\end{corollary}

\subsection{Interchanging module and comodule structures}\label{subsec:Mixed-modules}

Let $\sC \curvearrowright \cM \curvearrowleft \sD$ be a bimodule object in $\Cat$.
Lurie \cite[\S 4.3.2]{HA} showed that in this situation, $\LMod(\cM)$ inherits a right $\sD$-action and similarly $\RMod(\cM)$ inherits a left $\sC$-action.
He then shows \cite[Theorem 4.3.2.7]{HA} that there are natural equivalences
\[\LMod(\RMod(\cM)) \simeq \BMod(\cM) \simeq \RMod(\LMod(\cM)).\]
Since these equivalences are compatible with the projections to $\Alg(\sC) \times \Alg(\sD)$, this also gives natural equivalences $\LMod_A(\RMod_B(\cM)) \simeq \RMod_B(\LMod_A(\cM))$.

By dualizing this result, one obtains similar statements for comodules over coalgebras.
The goal of this section is to prove the following ``mixed'' version of this result.

\begin{theorem}\label{theorem:bicomod-interchange}
    Let $\cC \curvearrowright \cM \curvearrowleft \cD$ be a bimodule object in $\Cat$, let $Q$ be a coalgebra in $\cC$ and let $A$ be an algebra in $\cD$. Then there is an equivalence
    \[\LComod_Q(\RMod_A(\cM)) \simeq \RMod_A(\LComod_Q(\cM))\]
    that is natural in both $A$ and $Q$.
\end{theorem}

Since the proof of this result is very technical, we will first give a sketch. 

\subsubsection{A proof sketch}

By \cref{remark:simplified-left-right-module-operads}, we can define categories with a left action as functors $\LM \simeq \Delta^\op \times [1] \to \Cat$ satisfying a certain Segal condition---in particular, we may define them as certain maps $p \colon \sM_\bullet \to \sC_\bullet$ in $\Fun(\Delta^\op,\Cat)$.
Let us unwind the meaning of this Segal condition on $p$.
Write $i_n$ for the inclusion $\{n\} \hookrightarrow [n]$.
Then $p \colon \sM_\bullet \to \sC_\bullet$ defines a left action of $\sC \coloneqq \sC_\bullet$ on $\sM \coloneqq \sM_0$ if
\begin{itemize}
    \item $\sC$ is a monoidal category, and
    \item for every $n \geq 1$, the map $(p_n,i_n^*) \colon \sM_n \to \sC_n \times \sM_0$ is an equivalence.
\end{itemize}
We will generally write $\sC \curvearrowright \sM$ for such a left module category. By using the inclusions $\{0\} \hookrightarrow [n]$ instead of $\{n\} \hookrightarrow [n]$, one arrives at the definition of a category with a right action. Let us write $\LM(\Cat)$ (resp.\ $\RM(\Cat)$) for the category of categories with a left (resp.\ right) action.

Now suppose that $\cC \curvearrowright \cM \curvearrowleft \cD$ is a bimodule object in $\Cat$.
Then $\cC^{\times n} \times \cM$ admits a right $\cD$-action for every $n \geq 0$.
If we endow $\cC^{\times n}$ with the trivial right action of the category $*$, then we obtain a functor
\begin{align*}
    &\begin{tikzcd}[ampersand replacement = \&] \LM \simeq \Delta^\op \times [1] \ar[r] \& \RM(\Cat)\end{tikzcd} \\
    &\begin{tikzcd}[ampersand replacement = \&]
    ([n],i) \ar[r,mapsto] \& {\begin{cases}
        \cC^{\times n} \times \cM \curvearrowleft \cD \quad & \text{if} \quad i=0 \\
        \cC^{\times n} \curvearrowleft * \quad & \text{if} \quad i=1.
    \end{cases}} \end{tikzcd}
\end{align*}

We can now apply the functor $\RMod \colon \RM(\Cat) \to \Cat$ pointwise to obtain a diagram $\LM \to \Cat$ describing the left action of $\sC$ on $\RMod(\sM)$.
In particular, we can define $\LComod(\RMod(\sM))$. By switching the roles of $\LM$ and $\RM$, we instead obtain a functor $\RM \to \LM(\Cat)$ to which we can apply $\LComod$ pointwise.
This yields a right $\sD$-action on $\LComod(\sM)$, hence we can define $\RMod(\LComod(\sM))$.
The functors $\LM \to \RM(\Cat)$ and $\RM \to \LM(\Cat)$ both correspond to the same functor $F \colon \LM \times \RM \to \Cat$, namely the one given by
\begin{equation}\label{equation:informal-definition-F}
    F([m],i,[n],j) = \begin{cases}
        \sC^{\times m} \times \sM \times \sD^{\times n} \quad &\text{if} \quad i=j=0 \\
        \sC^{\times m} \quad &\text{if} \quad i=1, j=0 \\
        \sD^{\times n} \quad &\text{if} \quad i=0, j=1 \\
        * \quad &\text{if} \quad i=j=1
    \end{cases}
\end{equation}
By \cref{theorem:straightening-orthofibrations}, we obtain a corresponding orthofibration $\sF \to \LM^\op \times \RM$.
Let us call a section $\LM^\op \times \RM \to \sF$ that preserves cartesian lifts of inerts in the first variable and cocartesian lifts of inerts in the second variable a \emph{mixed module} and write $\MixMod(\sM)$ for their category.
One can then show that there are equivalences
\[\LComod(\RMod(\sM)) \simeq \MixMod(\sM) \simeq \RMod(\LComod(\sM)).\]
Moreover, these equivalences are compatible with the projections to $\coAlg(\sC) \times \Alg(\sD)$, producing the desired result.

\subsubsection{The actual proof}
We now make this proof sketch precise.
On the nose, our construction of the right action on $\LComod(\sM)$ and the left action on $\RMod(\sM)$ are different than those of Lurie \cite[\S 4.3.2]{HA}.
We will therefore compare our construction with that of Lurie in \cref{lemma:comparison-with-Lurie-bimod} below.

\begin{proof}[Proof of \cref{theorem:bicomod-interchange}]
    Let $\sC \curvearrowright \sM \curvearrowleft \sD$ be a bimodule object in $\Cat$.
    By \cref{proposition:approximation-to-BM-operad,definition:bimodule-categories}, we can view this as a functor $M \colon (\Delta_{/[1]})^\op \to \Cat$ satisfying the Segal condition. 
    
    \vspace{0.5em}
    \noindent\textbf{Step 1.} \textbf{Construction of the functor $F \colon \LM \times \RM \to \Cat$.}

    Write $\Delta^+$ for the augmented simplex category obtained by adjoining an initial object $[-1]$ to $\Delta$.
    Then $\Delta^+_{/[1]}$ is the category obtained by adjoining an initial object to $\Delta_{/[1]}$; let us denote this object by $[{-1}|{-1}]$.
    Moreover, concatenation of linearly ordered sets induces an equivalence
    \[\Delta^+ \times \Delta^+ \to \Delta^+_{/[1]}; \quad ([m],[n]) \mapsto [m|n].\]
    In particular, we can think of maps in $\Delta^+_{/[1]}$ as pairs $(\alpha,\beta)$ of maps in $\Delta^+$.
    Since $\Delta^+_{/[1]}$ is obtained by adjoining an initial object to $\Delta_{/[1]}$, there is a unique extension of $M$ to a functor
    \[\begin{tikzcd}
        M^+ \colon (\Delta^+_{/[1]})^\op \to \Cat
    \end{tikzcd}\]
    such that $M^+([{-1}|{-1}]) = *$.
    Define $\psi \colon \LM \times \RM \simeq \Delta^\op \times [1] \times \Delta^\op \times [1] \to (\Delta^+_{/[1]})^\op$ by
    \[\psi([m],i,[n],j) = \begin{cases}
        [m|n] \quad &\text{if} \quad i=j=0 \\
        [m|{-1}] \quad &\text{if} \quad i=1, j=0 \\
        [{-1}|n] \quad &\text{if} \quad i=0, j=1 \\
        [{-1}|{-1}] \quad &\text{if} \quad i=j=1
    \end{cases}\]
    on objects, and
    \[\psi(\alpha,i' \to i, \beta, j' \to j) = \begin{cases}
        (\alpha,\beta) \quad &\text{if} \quad i=j=0 \\
        (\alpha,*) \quad &\text{if} \quad i=1, j=0 \\
        (*,\beta) \quad &\text{if} \quad i=0, j=1 \\
        (*,*) \quad &\text{if} \quad i=j=1
    \end{cases}\]
    on morphisms. Here $*$ denotes, for any $[n]$ in $(\Delta^+)^\op$, the unique map to $[-1]$.
    The functor $F \colon \LM \times \RM \to \Cat$ is then defined as the composition $M^+ \circ \psi$. Unwinding this construction shows that $F$ is indeed of the form described in \cref{equation:informal-definition-F}. Write $\sF \to \LM^\op \times \RM$ for the corresponding orthofibration.
    
    \vspace{0.5em}
    \noindent\textbf{Step 2.} \textbf{Construction of $\sC \curvearrowright \RMod(\sM)$ and $\LComod(\sM) \curvearrowleft \sD$.}

    We start with constructing the left action $\sC \curvearrowright \RMod(\sM)$. Write $\pi_L \colon \LM^\op \times \RM \to \LM^\op$ for the projection, and given $([m],i)$ in $\LM$, write $\sF_{[m],i} \to \RM$ for the pullback $([m],i) \times_{\LM^\op} \sF$.
    Since $\pi_L$ is cocartesian, it is exponentiable by \cite[Lemma 3.2.1]{AyalaFrancis2020FibrationsInftyCategories} and hence $\pi_L^* \colon \Cat_{/\LM^\op} \to \Cat_{/\LM^\op \times \RM}$ admits a right adjoint $(\pi_L)_*$.
    Since $\pi_L^*$ preserves tensors by $[1]$, it follows from \cref{proposition:2-adjunctions-detected-underlying-tensor} that $(\pi_L)_*$ is a 2-categorical right adjoint.
    Observe that for any $([m],i)$ in $\LM$, we therefore have
    \begin{align*}
    ((\pi_L)_* \sF)_{([m],i)} &= \Fun_{/\LM^\op}(\{([m],i)\},(\pi_L)_* \sF) \simeq \Fun_{/\LM^\op \times \RM}(\{([m],i)\} \times \RM, \sF) \\
    &\simeq \Fun_{/\RM}(\RM,\sF_{[m],i}).
    \end{align*}
    In particular, an object of $(\pi_L)_* \sF$ consists of an object $([m],i)$ in $\LM$ together with a section $\RM \to \sF_{[m],i}$.
    We leave it to the reader to verify that $(\pi_L)_* \sF$ is a cartesian fibration over $\LM^\op$, where the cartesian morphisms are given by cartesian transport of sections $\RM \to \sF$.
    Let us write $(\pi_L^\inert)_* \sF$ for the full subcategory of $(\pi_L)_* \sF$ spanned by all such sections $\RM \to \sF_{[m],i}$ that preserve cocartesian lifts of inerts.
    Since cartesian transport preserves cocartesian lifts, it follows that $(\pi_L^\inert)_* \sF \to \LM^\op$ is again a cartesian fibration.
    Observe that $\sF_{[m],0}$ is the cocartesian fibration over $\RM$ describing the right action $\sC^{\times m} \times \sM \curvearrowleft \sD$.
    More precisely, this is the cartesian product in $\RM(\Cat)$ of the right actions $\sC^{\times m} \curvearrowleft *$ and $\sM \curvearrowleft \sD$.
    In particular, it follows that
    \[((\pi_L^\inert)_* \sF)_{([m],0)} = \RMod(\sC^{\times m} \curvearrowleft *) \times \RMod(\sM) \simeq \sC^{\times m} \times \RMod(\sM).\]
    We similarly see that $(\pi_L^\inert)_{([m],1)} = \sC^{\times m}$, hence the cartesian fibration $(\pi_L^\inert)_* \sF \to \LM^\op$ satisfies the Segal condition. In particular, it describes a left action of $\sC$ on $\RMod(\sM)$.

    To obtain the right action $\LComod(\sM) \curvearrowleft \sD$, note that $\LComod(\sM)$ is defined as $\LMod(\sC^\op \curvearrowright \sM^\op)^\op$, which is naturally equivalent to the category of sections over $\LM^\op$ that preserve cartesian lifts of inerts.
    We can therefore use an argument that is dual to the above: Let us write $\sF^{[n],j} \to \LM^\op$ for the pullback $\sF \times_{\RM} \{([n],j)\}$.
    As above, the pullback functor $(\pi_R)^* \colon \Cat_{/\RM} \to \Cat_{\LM^\op \times \RM}$ admits a right adjoint $(\pi_R)_*$.
    Define $(\pi_R^\inert)_* \sF$ as the full subcategory of $(\pi_R)_*\sF$ spanned by those sections $\LM^\op \to \sF^{[n],j}$ that preserve cartesian lifts of inert maps.
    It follows as above that $(\pi_R^\inert)_* \sF \to \RM$ is a cocartesian fibration describing a right action of $\sD$ on $\LComod(\sM)$.
    
    \vspace{0.5em}
    \noindent\textbf{Step 3.} \textbf{The equivalence $\LComod(\RMod(\sM)) \simeq \MixMod(\sM) \simeq \RMod(\LComod(\sM))$.}

    Given a category $\sB$, let us write $\pi_*$ for the functor $\Cat_{/\cB} \to \Cat$ sending a functor $\sA \to \sB$ to the category of sections $\Fun_{/\cB}(\sB,\sA)$.
    Note that this functor is right adjoint to pullback along $\pi \colon \sB \to *$.
    In case $\sA \to \sB$ is a (co)cartesian fibration and $\sB$ has a class of maps called ``inert'', we will write $\pi_*^\inert \sA$ for the full subcategory of sections that take all inert maps to (co)cartesian lifts.
    If $\sB = \LM^\op \times \RM$ and $\sA \to \LM^\op \times \RM$ is an orthofibration, then we will write $\pi_*^\inert \sA$ for the full subcategory of $\pi_* \sA$ spanned by those sections that preserve cartesian lifts of inerts in the first variable and cocartesian lifts of inerts in the second variable.

    It is now straightforward to show that the equivalences
    \[\pi_* (\pi_L)_* \sF \simeq \pi_* \sF \simeq \pi_* (\pi_R)_* \sF\]
    restrict to equivalences
    \[\pi_*^\inert (\pi_L^\inert)_* \sF \simeq \pi_*^\inert \sF \simeq \pi_*^\inert (\pi_R^\inert)_* \sF\]
    Defining $\MixMod(\sM)$ as $\pi_*^\inert\sF$, this gives the desired equivalences
    \[\LComod(\RMod(\sM)) \simeq \MixMod(\sM) \simeq \RMod(\LComod(\sM)).\]

    \vspace{0.5em}
    \noindent\textbf{Step 4.} \textbf{Compatibility with projections to $\coAlg(\sC) \times \Alg(\sD)$.}

    Note that the projections from $\RMod(\sM)$ and $\RMod(\LComod(\sM))$ to $\Alg(\sD)$ are obtained by restriction along the inclusion $i_R \colon \Delta^\op \hookrightarrow \RM$ given by $[n] \mapsto ([n],1)$.
    There is also an inclusion $j_R \colon \Delta^\op \hookrightarrow \LM^\op \times \RM$ given by $[n] \mapsto ([0],0,[n],1)$.
    Restriction along this inclusion induces a projection $\MixMod(\sM) \to \Alg(\sD)$.
    We will show that the diagram
    \[\begin{tikzcd}
    	{\LComod(\RMod(\sM))} & {\MixMod(\sM)} & {\RMod(\LComod(\sM))} \\
    	& {\Alg(\sD)}
    	\arrow[""{name=0, anchor=center, inner sep=0}, "\sim", from=1-1, to=1-2]
    	\arrow[""{name=1, anchor=center, inner sep=0}, from=1-1, to=2-2]
    	\arrow[from=1-2, to=2-2]
    	\arrow[""{name=2, anchor=center, inner sep=0}, "\sim"', from=1-3, to=1-2]
    	\arrow[""{name=3, anchor=center, inner sep=0}, from=1-3, to=2-2]
    	\arrow["{\text{\textcircled{a}}}"{pos=0.4}, draw=none, from=0, to=1]
    	\arrow["{\text{\textcircled{b}}}"'{pos=0.3}, draw=none, from=2, to=3]
    \end{tikzcd}\]
    commutes; an analogous argument can be used for the projection to $\coAlg(\sC)$.
    
    We first show that the triangle \textcircled{b} commutes.
    Consider the commutative diagram
    \[\begin{tikzcd}
    	{\LM^\op \times \RM} & \RM & {*}. \\
    	& {\Delta^\op}
    	\arrow["{\pi_R}", from=1-1, to=1-2]
    	\arrow["\pi", from=1-2, to=1-3]
    	\arrow["{j_R}", from=2-2, to=1-1]
    	\arrow["{i_R}", from=2-2, to=1-2]
    	\arrow["\pi"', from=2-2, to=1-3]
    \end{tikzcd}\]
    By pulling back and taking sections along these maps, we obtain a commutative diagram
    \[\begin{tikzcd}
    	{\pi^\inert_*(\pi^\inert_R)_*\sF} & {\pi^\inert_* i_R^* (\pi^\inert_R)_* \sF} \\
    	\pi_*^\inert \sF \ar[r]& {\pi^\inert_* j_R^* \sF}.
    	\arrow[from=1-1, to=1-2]
    	\arrow[from=1-1, to=2-1, "\sim" rot90]
    	\arrow[from=1-2, to=2-2]
    \end{tikzcd}\]
    It follows essentially by construction that $\pi_*^\inert i_R^*(\pi_R^\inert)_* \sF = \pi_*^\inert j_R^* \sF = \Alg(\sD)$ and that the top horizontal map is the projection $\RMod(\LComod(\sM)) \to \Alg(\sD)$, while the bottom horizontal map is the projection $\MixMod(\sM) \to \Alg(\sD)$. This shows that \textcircled{b} commutes.

    To see that \textcircled{a} commutes, let us first unwind where the projection $\LComod(\RMod(\sM)) \to \Alg(\sD)$ comes from:
    Observe that $\RMod(\sM) \to \Alg(\sD)$ can be viewed as a map in $\LM(\Cat)$, where we endow $\Alg(\sD)$ with the trivial left action $* \curvearrowright \Alg(\sD)$. Applying $\LComod$ to this then yields the projection
    \[\LComod(\RMod(\sM)) \to \LComod(* \curvearrowright \Alg(\sD)) = \Alg(\sD).\]
    Let $k \colon \Delta^\op \to \LM^\op \times \Delta^\op$ be the inclusion $[n] \mapsto ([0],0,[n])$, so that $i_R k = j_R$.
    Then the commutative diagram
    \[\begin{tikzcd}
	{\LM^\op \times \RM} & {\LM^\op} & {*} \\
	{\LM^\op \times \Delta^\op} \\
	{\Delta^\op}
	\arrow["{{\pi_L}}", from=1-1, to=1-2]
	\arrow["\pi", from=1-2, to=1-3]
	\arrow["{i_R}", from=2-1, to=1-1]
	\arrow["{\bar \pi_L}"', from=2-1, to=1-2]
	\arrow["{j_R}", bend left = 75, from=3-1, to=1-1]
	\arrow["\pi"', from=3-1, to=1-3]
	\arrow["k", hook, from=3-1, to=2-1]
\end{tikzcd}\]
    yields a commuting square
    \[\begin{tikzcd}
    	{\pi_*^\inert (\pi_L^\inert)_* \sF} & {\pi^\inert_* (\bar \pi_L^\inert)_* i_R^* \sF} \\
    	{\pi_*^\inert \sF} & {\pi_*^\inert j_R^* \sF}.
    	\arrow[from=1-1, to=1-2]
    	\arrow["\sim" rot90, from=1-1, to=2-1]
    	\arrow[from=1-2, to=2-2]
    	\arrow[from=2-1, to=2-2]
    \end{tikzcd}\]
    Since $(\bar\pi_L^\inert)_* i^*_R \sF$ is the category $\Alg(\sD)$ with the trivial left action $* \curvearrowright \sD$, it follows that the right vertical map is an equivalence $\pi_*^\inert(\pi_L^\inert)_* i^*_R \sF \simeq \pi_*^\inert j_R^* \sF = \Alg(\sD)$.
    Moreover, it follows by construction that the top horizontal map is the projection $\LComod(\RMod(\sM)) \to \Alg(\sD)$ and that the bottom horizontal map is $\MixMod(\sM) \to \Alg(\sD)$.
    We conclude that \textcircled{a} commutes.

    \vspace{0.5em}
    \noindent\textbf{Step 5.} \textbf{Naturality of the equivalence.}

    It follows by an argument similar to \cite[Corollary 4.3.3.2]{HA} that the projections from $\MixMod(\sM)$, $\LComod(\RMod(\sM))$ and $\RMod(\LComod(\sM))$ to $\Alg(\sD) \times \coAlg(\sC)$ are orthofibrations, where a map is an orthocartesian morphism in the sense of \cref{definition:orthocartesian-maps} precisely if the underlying map in $\sM$ is an equivalence.
    (Note that we switched $\Alg(\sD)$ and $\coAlg(\sC)$. This is because the projection is cartesian over $\Alg(\sD)$ and cocartesian over $\Alg(\sC)$.)
    We conclude that the equivalence $\LComod(\RMod(\sM)) \simeq \RMod(\LComod(\sM))$ is an equivalence of orthofibrations.
    By considering the fibers of the projections to $\Alg(\sD)$ and $\coAlg(\sC)$, we obtain that the equivalence
    \[\LComod_Q(\RMod_A(\sM)) \simeq \RMod_A(\LComod_Q(\sM))\]
    is natural with respect to restriction in the algebra direction and corestriction in the coalgebra direction.
\end{proof}

Note that in \cite[Construction 4.3.2.2]{HA}, Lurie gave a construction of a right $\sD$-action on $\LMod(\sM)$. Replacing left and right in this construction, one obtains a left $\sC$-action on $\RMod(\sM)$, allowing one to give an alternative construction of $\LComod(\RMod(\sM))$. We will now show that Lurie's construction agrees with ours. For this, it suffices to show the following.

\begin{lemma}\label{lemma:comparison-with-Lurie-bimod}
    Let $\sC \curvearrowright \sM \curvearrowleft \sD$ be a bimodule object in $\Cat$. Then the left $\sC$-action on $\RMod(\sM)$ constructed in \cite[Construction 4.3.2.2]{HA} agrees with the left $\sC$-action constructed on $\RMod(\sM)$ in the proof of \cref{theorem:bicomod-interchange} above.
\end{lemma}

\begin{proof}
    Let $\mathscr{LM}^\otimes$, $\mathscr{RM}^\otimes$ and $\mathscr{BM}^\otimes$ be the left/right/bimodule operads defined by Lurie in \cite[\S\S\ 4.2.1 \& 4.3.1]{HA}.
    Let $p_L \colon \LM \to \mathscr{LM}^\otimes$, $p_R \colon \RM \to \mathscr{RM}^\otimes$, and $q \colon \BM \to \mathscr{BM}^\otimes$ be the approximations of $\infty$-operads from \cref{proposition:approximation-to-BM-operad,remark:simplified-left-right-module-operads}.
    Let us write $\BM^+$ for the category $(\Delta^+_{/[1]})^\op$ obtained by adjoining a terminal object to $\BM$.
    Note that we can uniquely extend $q \colon \BM \to \mathscr{BM}^\otimes$ to a functor $q^+ \colon \BM^+ \to \mathscr{BM}^\otimes$ that sends the terminal object of $\BM^+$ to the terminal object of $\mathscr{BM}^\otimes$.
    The result follows by showing that the square
    \[\begin{tikzcd}
        \LM \times \RM \ar[r,"\psi"] \ar[d,"{(p_L,p_R)}"'] & \BM^+ \ar[d,"q^+"] \\
        \mathscr{LM}^\otimes \times \mathscr{RM}^\otimes \ar[r,"\mathbf{Pr}"] & \mathscr{BM}^\otimes
    \end{tikzcd}\]
    commutes, where $\mathbf{Pr}$ is the functor defined in \cite[Construction 4.3.2.1]{HA} and $\psi \colon \LM \times \RM \to \BM^+$ is the functor defined in the first step of the proof of \cref{theorem:bicomod-interchange}.
    We leave this verification to the reader.
\end{proof}

\subsubsection{Some related results}

We now state some auxiliary results that can be deduced from the proof above and that are helpful in \cref{sec: koszul-duality}.

\begin{lemma}\label{lem: Restriction-of-scalars-C-linear}
    Let $\sC \curvearrowright \sM \curvearrowleft \sD$ be a bimodule object in $\Cat$ and let $f \colon A \to B$ be a map in $\Alg(\sD)$. Then the restriction of scalars functor
    \[f^* \colon \RMod_B(\sM) \to \RMod_A(\sM)\]
    is left $\sC$-linear.
\end{lemma}

\begin{proof}
    Let us write $\MonCat$ for the category of monoidal categories. Note that a $\sC$-linear functor is the same as a map in
    \[\LMod_{\sC}(\Cat) = \LM(\Cat) \times_{\MonCat} \{\sC\}.\]
    We will show that the functor
    \begin{align*}
        &\qquad \begin{tikzcd}[ampersand replacement = \&] \Alg(\sD)^\op \ar[r] \& \Cat\end{tikzcd} \\
        &\begin{tikzcd}[ampersand replacement = \&]
        (f \colon A \to B) \ar[r,mapsto] \& (f^* \colon \RMod_B(\sM) \to \RMod_A(\sM)) \end{tikzcd}
    \end{align*}
    can be lifted to a functor that lands in $\LMod_{\sC}(\Cat)$.
    In the proof of \cref{theorem:bicomod-interchange}, we obtained the left action of $\sC$ on $\RMod(\sM)$ together with the projection $\RMod(\sM) \to \Alg(\sD)$ by constructing a functor of the form
    \begin{align*}
        &\begin{tikzcd}[ampersand replacement = \&] \LM \ar[r] \& \Fun([1],\Cat)\end{tikzcd} \\
        &\begin{tikzcd}[ampersand replacement = \&]
        ([n],i) \ar[r,mapsto] \& {\begin{cases}
            \RMod(\cC^{\times n} \times \cM) \to \Alg(\sD) \quad & \text{if} \quad i=0 \\
            \cC^{\times n} \to * \quad & \text{if} \quad i=1.
        \end{cases}} \end{tikzcd}
    \end{align*}
    Here the right $\sD$-action on $\sC^{\times n} \times \sM$ is the product of the trivial action $\sC^{\times n} \curvearrowleft *$ and the right action $\sM \curvearrowleft \sD$.
    Observe that by the right module version of \cite[Corollary 4.2.3.2]{HA}, the projections $\RMod(\sC^{\times n} \times \sM) \to \Alg(\sD)$ are cartesian fibrations, where a morphism in $\RMod(\sC^{\times n} \times \sM)$ is cartesian precisely if the underlying map in $\sC^{\times n} \times \sM$ is an equivalence.
    It follows that this functor $\LM \to \Fun([1],\Cat)$ lands in the (non-full) subcategory of cartesian fibrations. Applying cartesian straightening and swapping $\Alg(\sD)^\op$ with $\LM$, we obtain a diagram
    \begin{align*}
        &\begin{tikzcd}[ampersand replacement = \&] \Alg(\sD)^\op \ar[r] \& \LM(\Cat) \end{tikzcd} \\
        &\begin{tikzcd}[ampersand replacement = \&]
        A \ar[r,mapsto] \& {\begin{cases}
            \RMod_A(\cC^{\times n} \times \cM) \quad & \text{if} \quad i=0 \\
            \cC^{\times n} \quad & \text{if} \quad i=1.
        \end{cases}} \end{tikzcd}
    \end{align*}
    where the functoriality is given by restriction of scalars. Moreover, the composite
    \[\begin{tikzcd}
        \Alg(\sD)^\op \ar[r] & \LM(\Cat) \ar[r] & \MonCat
    \end{tikzcd}\]
    is by construction constant with value $\sC$. This gives the desired result.
\end{proof}

Let $\sC \curvearrowright \sM \curvearrowleft \sD$ and $\sC' \curvearrowright \sN \curvearrowleft \sD'$ be bimodule objects in $\Cat$, viewed as cocartesian fibrations over $(\Delta_{/[1]})^\op$ which we denote by $\overline \sM \to (\Delta_{/[1]})^\op$ and $\overline \sN \to (\Delta_{/[1]})^\op$, respectively. A \emph{lax bimodule functor} from $\sC \curvearrowright \sM \curvearrowleft \sD$ to $\sC' \curvearrowright \sN \curvearrowleft \sD'$ is then defined as a commuting triangle
\[\begin{tikzcd}[column sep = 0.5]
	\overline \sM && \overline \sN \\
	& {(\Delta_{/[1]})^\op}
	\arrow["F", from=1-1, to=1-3]
	\arrow[from=1-1, to=2-2]
	\arrow[from=1-3, to=2-2]
\end{tikzcd}\]
such that $F$ preserves cocartesian lifts of inert morphisms. We say that $F$ is \emph{strong on the left} if it furthermore preserves cocartesian lifts of any map of the form
\[\begin{tikzcd}
    (\alpha,\id_{[n]}) \colon [m|n] \to [m'|n].
\end{tikzcd}\]
This condition ensures that $F$ restricts to a strong monoidal functor $\sC \to \sC'$ and that it commutes strictly with left actions of $\sC$ and $\sC'$ on $\sM$ and $\sN$.
In particular, $F$ induces functors $\coAlg(\sC) \to \coAlg(\sC')$ and $\Alg(\sD) \to \Alg(\sD')$.
These functors lift to a functor $\MixMod(\sM) \to \MixMod(\sN)$:

\begin{lemma}
\label{lem: strong-lax-induced-map-mix-mod}
    Let $F$ be a lax bimodule functor from $\sC \curvearrowright \sM \curvearrowleft \sD$ to $\sC' \curvearrowright \sN \curvearrowleft \sD'$ and suppose that $F$ is strong on the left. Then $F$ induces a functor $\MixMod(\sM) \to \MixMod(\sN)$ such that the diagram
    \begin{equation}\label{diagram:Transporting-MixMod-along-lax-functor}
        \begin{tikzcd}
        	{\MixMod(\sM)} & {\MixMod(\sN)} \\
        	{\coAlg(\sC) \times \sM \times \Alg(\sD)} & {\coAlg(\sC') \times \sN \times \Alg(\sD')}
        	\arrow[from=1-1, to=1-2]
        	\arrow[from=1-1, to=2-1]
        	\arrow[from=1-2, to=2-2]
        	\arrow[from=2-1, to=2-2]
        \end{tikzcd}
    \end{equation}
    commutes.
\end{lemma}

\begin{proof}
    Write $M$ and $N$ for the functors $(\Delta_{/[1]})^\op \to \Cat$ corresponding to the bimodule objects $\sM$ and $\sN$.
    As in the proof of \cref{theorem:bicomod-interchange}, these extend uniquely to functors $M^+,N^+ \colon (\Delta_{/[1]}^+)^\op \to \Cat$ that preserve the terminal object.
    Write $\overline{\sM}$, $\overline{\sN}$, $\overline{\sM}{}^+$ and $\overline{\sN}{}^+$ for the cocartesian unstraightenings of these functors.
    The lax bimodule functor $F \colon \overline{\sM} \to \overline{\sN}$ extends uniquely to a functor $F^+ \colon \overline{\sM}{}^+ \to \overline{\sN}{}^+$ over $(\Delta_{/[1]}^+)^\op$.
    By pulling back along the map $\psi \colon \LM \times \RM \to (\Delta_{/[1]}^+)^\op$ constructed in the proof of \cref{theorem:bicomod-interchange}, we obtain the commutative triangle
    \[\begin{tikzcd}[column sep = tiny]
    	{\psi^*\overline \sM{}^+} && {\psi^*\overline \sN{}^+} \\
    	& {\LM \times \RM}
    	\arrow["{\psi^*F^+}", from=1-1, to=1-3]
    	\arrow[from=1-1, to=2-2]
    	\arrow[from=1-3, to=2-2]
    \end{tikzcd}\]
    where the diagonal arrows are cocartesian fibrations. The condition that $F$ is strong on the left implies that $\psi^*F$ preserves cocartesian lifts of all arrows of the form $(f,\id)$ in $\LM \times \RM$. By the same proof as \cite[Theorem 2.5.1]{linskensLaxMonoidalAdjunctions}, we obtain a commuting triangle
    \[\begin{tikzcd}
    	{\sO_\sM} && {\sO_\sN} \\
    	& {\LM^\op \times \RM}
    	\arrow["{F^\vee}", from=1-1, to=1-3]
    	\arrow[from=1-1, to=2-2]
    	\arrow[from=1-3, to=2-2]
    \end{tikzcd}\]
    where $\sO_\sM$ and $\sO_\sN$ are the orthofibrations corresponding to $\psi^* \overline\sM{}^+$ and $\psi^* \overline \sN{}^+$ and where $F^\vee$ preserves cartesian lifts of maps in $\LM$ and cocartesian lifts of inert maps in $\RM$. The functor $F^\vee$ induces the desired map
    \[\MixMod(\sM) \to \MixMod(\sN)\]
    by postcomposition, and the diagram \cref{diagram:Transporting-MixMod-along-lax-functor} commutes by construction.
\end{proof}
\newpage

\phantomsection
\printbibliography

\end{document}